\documentclass[twoside,
a4paper]{amsart}
\usepackage[T1,T2A]{fontenc}
\usepackage[koi8-r]{inputenc}
\usepackage[english,russian]{babel}
\providecommand{\cosp}{\loglike{cosp}}
\providecommand{\sinp}{\loglike{sinp}}

\providecommand{\alli}{\iota}
\providecommand{\ladder}[2][]{L_{#1}^{\!#2}}
\providecommand{\uir}[2][\relax]{\rho_{#2}^{#1}}
\usepackage{xr-hyper}
\usepackage{graphicx}

\PassOptionsToPackage{pdfauthor={Vladimir V. Kisil},%
    pdftitle={Erlangen Program at Large---2 1/2: Induced Representations and Hypercomplex Numbers},%
    pdfsubject={mathematics},%
    backref=page,%
  pdfkeywords={symmetries, dual numbers, double numbers, Moebius map,
    split-complex numbers}
}{hyperref}
\IfFileExists{eulervm.sty}{\usepackage[T1]{eulervm}
\renewcommand{\rmdefault}{ppl}
\newcommand{\Aprime}{A\!'}}
{\newcommand{\Aprime}{A'}}

\newcommand{\ppnum}[3]{\def\publname{\texttt{#2}, #3\\
    #1}}
\PII{\relax}
\copyrightinfo{}{\relax}

\newtheorem{thm}{Theorem}[section]
\newtheorem{prop}[thm]{Proposition}
\newtheorem{lem}[thm]{Lemma}
\newtheorem{principle}[thm]{Principle}

\newtheorem{thmr}{Теорема}[section]
\newtheorem{propr}[thm]{Предложение}
\newtheorem{lemr}[thm]{Лемма}
\newtheorem{principler}[thm]{Принцип}

\theoremstyle{definition}
\newtheorem{defn}[thm]{Definition}
\newtheorem{example}[thm]{Example}
\newtheorem{defnr}[thm]{Определение}
\newtheorem{exampler}[thm]{Пример}

\theoremstyle{remark}
\newtheorem{rem}[thm]{Remark}

\newtheorem{remr}[thm]{Замечание}

\newcommand{\AMSMSC}[2]{\subjclass[2000]{Primary #1; Secondary #2.}}

\usepackage[breaklinks=true,%
colorlinks=true,%
bookmarks=true,%
backref=page,%
pagebackref=true]{hyperref}

\usepackage[backrefs,msc-links
,nobysame]{myamsrefs}

\ppnum{
}{\href{http://arXiv.org/abs/0909.4464}{arXiv:0909.4464}}
{2009}

\input{mydef}
\usepackage[all]{xy}
 \IfFileExists{noweb.sty}{\usepackage{noweb}}{}
 \makeatletter
 \nwopt@smallcode
 \let\tab=&
 \def\idxexample#1{\nwix@id@uses#1}
 \makeatother
 
 \def\nwlbrace{\textbf{\texttt{\char123}}}
 \def\nwrbrace{\textbf{\texttt{\char125}}}
\newcommand{\CPP}{\texttt{C++}}
\newcommand{\NoWEB}{\texttt{noweb}}
\providecommand{\MetaPost}{\texttt{Meta}\-\texttt{Post}}
\providecommand{\GiNaC}{\textsf{GiNaC}}
\providecommand{\pyGiNaC}{\textsf{pyGiNaC}}

\providecommand{\clifford}[2][]{\ifcase #1 #2\or \tilde{#2} \or \breve{#2} \fi}

\providecommand{\rs}{\mathring{\sigma}}
\begin{document}
  \selectlanguage{english} 
\title
{
Induced Representations and Hypercomplex Numbers}

\author[Vladimir V. Kisil]%
{\href{http://www.maths.leeds.ac.uk/~kisilv/}{Vladimir V. Kisil}}
\thanks{On  leave from Odessa University.}

\address{%
School of Mathematics\\
University of Leeds\\
Leeds LS2\,9JT\\
UK
}

\email{\href{mailto:kisilv@maths.leeds.ac.uk}{kisilv@maths.leeds.ac.uk}}

\urladdr{\href{http://www.maths.leeds.ac.uk/~kisilv/}%
{http://www.maths.leeds.ac.uk/\~{}kisilv/}}

\begin{abstract}
  In the search for hypercomplex analytic functions on the half-plane,
  we review the construction of induced representations of the group
  \(G=\SL\). Firstly we note that \(G\)-action on the homogeneous
  space \(G/H\), where \(H\) is any one-dimensional subgroup of
  \(\SL\), is a linear-fractional transformation on hypercomplex
  numbers. Thus, we investigate various hypercomplex characters of
  subgroups \(H\). The correspondence between the structure of the
  group \(\SL\) and hypercomplex numbers can be illustrated in many
  other situations as well. We give examples of induced representations
  of \(\SL\) on spaces of hypercomplex valued functions, which are
  unitary in some sense. Raising/lowering operators for various
  subgroup prompt hypercomplex coefficients as well. 
\end{abstract}
\AMSMSC{22D30}{08A99, 13A50, 15A04, 20H05, 51M10} 

\maketitle

\tableofcontents

\section{Introduction}
\label{sec:introduction}

Analytic functions of a complex variable form a beautiful theory with
rich applications in many fields ranging from number theory to electrical
engineering. Thus, it is natural to look for its analogs and
generalisations in different directions. The most basic (or
fundamental?) situation appears if we replace the complex imaginary
unit \(\rmi^2=-1\) with either the hyperbolic one \(\rmh^2=+1\) or the
nilpotent \(\rmp^2=0\).  

Two-dimensional commutative associative algebra over reals generated
by \(1\) and \(\rmh\) consists of numbers \(x+\rmh y\), where \(x\),
\(y\in\Space{R}{}\).  They are known as split-complex, duplex,
hyperbolic or
\href{http://en.wikipedia.org/wiki/Split-complex_number}{\emph{double
    numbers}}~\cites{BocCatoniCannataNichZamp07,%
  KhrennikovSegre07a,Ulrych05a,Pilipchuk10a}. The algebra has zero
divisors \(\rmh_\pm=\frac{1}{2}(1\pm\rmh)\) with the properties
\(\rmh_\pm^2=\rmh_\pm\) and \(\rmh_+ \rmh_-=0\). Thus, double numbers
are isomorphic to \(\Space{R}{}\oplus \Space{R}{}\)---the direct sum
of two copies of the real line spanned by \(\rmh_+\) and \(\rmh_-\).
This explains the names ``split-complex'' and ``double''.

The analogous algebra associated to the nilpotent unit \(\rmp\)
consists of elements \(x+\rmp y\), which are called
\wiki{Dual_number}{\emph{dual
    numbers}}~\cites{Yaglom79,HerranzSantander02b,CatoniCannataNichelatti04}.
All zero divisors in the algebra are \(\rmp y\), \(y\in\Space{R}{}\).
Physical applications of hypercomplex numbers are scattered through
classical mechanics~\cites{Yaglom79,Kisil12c}, non-linear
dynamics~\cites{Pilipchuk11a,Pilipchuk10a},
relativity~\cites{BocCatoniCannataNichZamp07,Ulrych05a,Garasko09a},
cosmology~\cites{HerranzSantander02b,GromovKuratov05a} and quantum
mechanics~\cites{KhrennikovSegre07a,Kisil10a,Ulrych10a}.
\begin{rem}
  Unfortunately, there are no common notations for hypercomplex
  units. Moreover, it would be difficult simply to list the whole
  variety of symbols employed for this. Even the complex imaginary
  unit \(\rmi\) is oftenly written as \(j\) in engineering. The hyperbolic 
  unit is denoted by \(\rmh\) in many papers starting at least from the
  foundational article~\cite{VignauxDuranona35a}; although a different
  letter \(e\) is used in the remarkable book~\cite{Yaglom79}. The
  symbol \(\iota\) is used for the nilpotent unit
  in~\cites{GromovKuratov06a,Pimenov65a}, however we chose \(\rmp\)
  following Yaglom~\cite{Yaglom79}. The later notation is also
  suggestive in light of the following remark.
\end{rem}
\begin{rem}
  The parabolic unit \(\rmp\) is a close relative to the infinitesimal
  number \(\varepsilon\) from non-standard
  analysis~\cites{Devis77,Uspenskii88}. The former has the property
  that its square is \emph{exactly} zero, meanwhile the square of the
  latter is \emph{almost} zero at its own scale. In fact, there is a
  version of non-standard analysis~\cite{Bell08a} employing the
  nilpotent unit \(\rmp\) as an infinitesimal\footnote{I am grateful
    to the anonymous referee for pointing my attention to the book~\cite{Bell08a} by
    J.L.~Bell.}. Also, some non-standard proofs of the main
  calculus theorems are given in~\cite{CatoniCannataNichelatti04}. A
  similar property allows to obtain classical mechanics from the
  representations of the Heisenberg group~\cites{Kisil10a,Kisil12c}.
\end{rem}

\emph{What kind of ``analytic'' functions can be associated with dual and
double numbers}? Since this question is very natural it was addressed
over a prolonged period of time by various researchers. Many of them
were unaware of works of their predecessors, neither I can claim to
possess the complete knowledge.  Below is a brief summary of several
works known to me.

For double numbers, a systematic study was already accomplished
in~\cite{VignauxDuranona35a}, there are also numerous later
investigations and surveys,
see~\citelist{\cite{BocCatoniCannataNichZamp07}
  \cite{KhrennikovSegre07a} \cite{LavrentShabat77} \cite{MotterRosa98}
  \cite{Kravchenko09a}*{Part~IV} \cite{Pilipchuk10a}
  \cite{Sobczyk95a}  \cite{Libine07a}  \cite{Kisil97c}} and further references therein. The existing
consensus is based on the factorisation of the wave equation
\(\partial_x^2-\partial_y^2=(\partial_x-\rmh\partial_y)(\partial_x+\rmh\partial_y)\)
into a product of two linear differential operators. This is an analog
of the factorisation of the Laplacian
\(\partial_x^2+\partial_y^2=(\partial_x-\rmi\partial_y)(\partial_x+\rmi\partial_y)\)
into the product of the Cauchy--Riemann operator and its adjoint.
Thus, hyperbolic analytic functions are defined to be null solutions
of the operator \(\partial_x+\rmh\partial_y\). However, the split of
dual number in the basis \(\rmh_\pm\) reduces an ``analytic'' function
\(f(x,y)=\rmh_+f_+(x+y)- \rmh_-f_-(x-y)\) to the sum, where \(f_\pm\)
are two generic differentiable real-valued functions of a real
variable. This is related to the representation of the generic
solution of the wave equation on the infinite string as a sum of a
wave travelling to the left and another travelling to the right.

For the nilpotent unit \(\rmp\) the situation is even more trivial. The
above factorisation approach does not lead anywhere
useful\footnote{However, there is an interesting factorisation of a parabolic
  operator \(\partial_x^2+\partial_y\) which requires Clifford
  algebras of higher dimensions~\cite{CerejeirasKahlerSommen05a}.}, since
\((\partial_x-\rmp\partial_y)(\partial_x+\rmp\partial_y)=\partial_x^2\).
An attempt to define analytic functions in terms of power series of
\(x+\rmp y\) produces only functions of the form \(f(x)+\rmp y
f'(x)\) for a real-analytic function \(f(x)\) of a single
variable~\cite{CatoniCannataNichelatti04}.

In the series of previous
works~\cites{Kisil96c,Kisil97c,Kisil97a,Kisil02c} we proposed to
develop an analytic function theory in the spirit of the Erlangen
Programme of F.~Klein: from an appropriate group action. The present
paper is a further step towards that goal. We observe that classical
spaces of complex analytical functions---the Hardy and Bergman
spaces---are irreducible moduli under certain representations of the
group \(\SL\). Moreover, those representations are induced by 
complex-valued characters of the compact subgroup \(K\subset\SL\).
Thus, we hope to find \emph{other types of analytic function
  spaces among all irreducible (or primary) \(\SL\)-moduli under induced
  representations}.

The paper outline is as follows. We recall the structure of \(\SL\)
and list its tree non-isomorphic one-dimensional continuous subgroups
in Section~\ref{sec:intro-sl2}. This information is employed to
classify tree non-isomorphic two-dimensional \(\SL\)-homogeneous
spaces in Section~\ref{sec:induc-repr}. We discover that three kinds
of hypercomplex numbers are perfectly suited to describe
\(\SL\)-action on those homogeneous spaces. To induce representations
from one-dimensional subgroups we describe their characters (from the
algebraic point of view) in Section~\ref{sec:hyper-characters}. The
next Section~\ref{sec:second-attempt} provides an additional objects,
which can be also viewed as characters with geometrical spirit. An
induced representation can be linear if the inducing characters is
also linear. For geometric characters this is not straightforward and
we discuss various possibilities in
Sections~\ref{sec:rebu-algebr-struct}
and~\ref{sec:invar-line-algebra}. Finally, we are able to write down
various induced representations in Section~\ref{sec:concl-induc-repr}.
A convenient description of an irreducible \(\SL\)-module can be
produced in terms of a linear basis, consisting of eigenvectors of the
subgroup inducing the representation. The Lie algebra
\(\algebra{sl}_2\) of \(\SL\) acts transitively on this basis by means
of ladder operators, which are important in quantum mechanics. We
describe the respective structures in
Section~\ref{sec:correspondence}. The paper is concluded by an
(incomplete) list of interesting directions for further research.
\medskip

\textbf{Acknowledgement:}
I am grateful to the anonymous referee who made many useful
suggestions and comments, which helped to improve this paper.

\section{The Group $\SL$ and Its Subgroups}
\label{sec:intro-sl2}

Let \(\SL\) be \wiki{SL\%282\%2CR\%29}{the group of \(2\times 2\)
  matrices} with real entries and of determinant one~\cite{Lang85}.
This is the smallest semisimple Lie group. Any matrix in \(\SL\) admits a
(unique) decomposition of the form~\cite{HoweTan92}*{Exer.~I.14}:
\begin{equation}
  \label{eq:iwasawa-decomp}
  \begin{pmatrix}
    a&b \\c &d
  \end{pmatrix}= {\begin{pmatrix} \alpha & 0\\0&\alpha^{-1}
    \end{pmatrix}} {\begin{pmatrix} 1&\nu \\0&1
    \end{pmatrix}} {\begin{pmatrix}
      \cos\phi &  -\sin\phi\\
      \sin\phi & \cos\phi
    \end{pmatrix}},
\end{equation}
for the following real values of the parameters:
\begin{displaymath}
  \label{eq:iwasawa-param}
  \alpha=\sqrt{c^2+d^2}\in(0,\infty), \quad
  \nu=ac+bd\in(-\infty,\infty),\quad
  \phi = -\arctan \frac{c}{d}\in(-\pi,\pi].
\end{displaymath}
The formula~\eqref{eq:iwasawa-decomp}
rewritten in a way \(\SL=ANK\) is known as Iwasawa
decomposition~\cite{Lang85}*{\S~III.1} and can be generalised to any
semisimple Lie group.

Each out of the three types of matrices in the right-hand side
of~\eqref{eq:iwasawa-decomp} forms a one-parameter subgroup \(A\),
\(N\) and \(K\).  They are obtained by the exponentiation of
the respective zero-trace matrices:
\begin{eqnarray}
  \label{eq:a-subgroup}
  A&=&\left\{  
    \begin{pmatrix} e^t & 0\\0&e^{-t}
    \end{pmatrix}=\exp \begin{pmatrix} t & 0\\0&-t
    \end{pmatrix},\  t\in\Space{R}{}\right\},\\
  \label{eq:n-subgroup}
  N&=&\left\{   {\begin{pmatrix} 1&t \\0&1
      \end{pmatrix}=\exp \begin{pmatrix} 0 & t\\0&0
      \end{pmatrix},}\  t\in\Space{R}{}\right\},\\
  \label{eq:k-subgroup}
  K&=&\left\{ {\begin{pmatrix}
        \cos t &  \sin t\\
        -\sin t & \cos t
      \end{pmatrix}=   \exp \begin{pmatrix} 0& t\\-t&0
      \end{pmatrix}},\ t\in(-\pi,\pi]\right\}.
\end{eqnarray}
 The following simple result
have an instructive proof.
\begin{prop}
  Any continuous one-parameter subgroup of \(\SL\) is conjugate to one of
  subgroups \(A\), \(N\) or \(K\).
\end{prop}
\begin{proof}
  Any one-parameter subgroup is obtained through the exponentiation
  \begin{equation}
    \label{eq:exp-taylor}
    e^{tX}=\sum_{n=0}^\infty \frac{t^n}{n!}X^n
  \end{equation}
  of an element \(X\) of the Lie algebra \(\algebra{sl}_2\) of
  \(\SL\). Such \(X\) is a \(2\times 2\) matrix with the zero trace.
  The behaviour of the Taylor expansion~\eqref{eq:exp-taylor} depends
  from properties of powers \(X^n\). This can be classified by a
  straightforward calculation:
  \begin{lem}
    \label{le:square}
    The  square \(X^2\) of a traceless matrix \(X=
    \begin{pmatrix}
      a&b\\c&-a
    \end{pmatrix}\) is the identity matrix times \(a^2+bc=-\det X\).
    The factor can be negative, zero or positive, which corresponds to
    the three different types of the Taylor
    expansion~\eqref{eq:exp-taylor} of \(e^{tX}\).
  \end{lem}
  It is a simple exercise on characteristic polynomials to see that
  through the matrix similarity we can obtain from \(X\) a generator
  \begin{itemize}
  \item of the subgroup \(K\) if \((-\det X) <0\);
  \item of the subgroup \(N\) if \((-\det X) =0\);
  \item of the subgroup \(A\) if \((-\det X) >0\).
  \end{itemize}
  The determinant is invariant under the similarity, thus these
  cases are distinct.
\end{proof}
\begin{example}
  The following two subgroups are conjugated to \(A\) and \(N\) respectively:
\begin{eqnarray}
  \label{eq:ap-subgroup}
  \Aprime&=&\left\{  
    \begin{pmatrix} \cosh t & \sinh t\\ \sinh t& \cosh t
    \end{pmatrix}=\exp \begin{pmatrix} 0 & t\\t&0
    \end{pmatrix},\  t\in\Space{R}{}\right\},\\
  \label{eq:np-subgroup}
  N'&=&\left\{   {\begin{pmatrix} 1&0\\t&1
      \end{pmatrix}=\exp \begin{pmatrix} 0 & 0\\t&0
      \end{pmatrix},}\  t\in\Space{R}{}\right\}.
\end{eqnarray}
\end{example}

\section{Action of $\SL$ as a Source of Hypercomplex Numbers}
\label{sec:induc-repr}

Let \(H\) be a subgroup of a group \(G\).  Let \( \Omega=G / H\) be
the corresponding homogeneous space and \(s: \Omega \rightarrow G\) be
a smooth section~\cite{Kirillov76}*{\S~13.2}, which is a left inverse
to the natural projection \(p: G\rightarrow \Omega \). The choice of
\(s\) is inessential in the sense that by a smooth map
\(\Omega\rightarrow \Omega\) we can always reduce one to another.

Any \(g\in G\) has a unique decomposition of the form
\(g=s(\omega)h\), where \(\omega=p(g)\in \Omega\) and \(h\in H\).
Note that \(\Omega \) is a left homogeneous space with the
\(G\)-action defined in terms of \(p\) and \(s\) as follows:
\begin{equation}
  \label{eq:g-action}
  g: \omega  \mapsto g\cdot \omega=p(g* s(\omega)),
\end{equation}
where \(*\) is the multiplication on \(G\). This is also illustrated
by the following commutative diagram:
\begin{displaymath}
        \xymatrix{ 
        G \ar@<.5ex>[d]^p \ar[r]^{g*} & G \ar@<.5ex>[d]^p \\ 
        \Omega   \ar@<.5ex>[u]^s \ar[r]^{g\cdot}   & \Omega     \ar@<.5ex>[u]^s }
\end{displaymath}

For \(G=\SL\), as well as for other semisimple groups, it is common to
consider only the case of \(H\) being the maximal compact subgroup
\(K\). However in this paper we admit \(H\) to be any one-dimensional
subgroup. Due to the previous Proposition it is sufficient to take
\(H=K\), \(N'\) or \(\Aprime\). Then \(\Omega\) is a two-dimensional manifold
and for any choice of \(H\) we define~\cite{Kisil97c}*{Ex.~3.7(a)}:
\begin{equation}
  \label{eq:s-map}
  s: (u,v) \mapsto
  \frac{1}{\sqrt{v}}
  \begin{pmatrix}
    v & u \\ 0 & 1
  \end{pmatrix}, \qquad (u,v)\in\Space{R}{2},\  v>0.
\end{equation}
A direct (or computer algebra~\cite{Kisil07a}) calculation show that:
\begin{prop}
  \label{pr:sl2-act-brute}
  The \(\SL\) action~\eqref{eq:g-action} associated to the map
  \(s\)~\eqref{eq:s-map} is:
  \begin{equation}
    \label{eq:sl2-act-brute}
    \begin{pmatrix}
      a&b\\c&d
    \end{pmatrix}: (u,v)\mapsto
    \left(\frac{(au+b)(c u+d) -\sigma cav^2}{( c u+d)^2 -\sigma (cv)^2},
      \frac{v}{( c u+d)^2 -\sigma (cv)^2}\right),
  \end{equation}
  where \(\sigma=-1\), \(0\) and \(1\) for the subgroups \(K\), \(N'\)
  and \(\Aprime\) respectively.
\end{prop}
The expression in~\eqref{eq:sl2-act-brute} does not look very
appealing, however an introduction of hypercomplex numbers makes it
more attractive:
\begin{prop}
  \label{pr:sl2-act}
  Let a hypercomplex unit \(\alli\) be such that \(\alli^2=\sigma\), then the
  \(\SL\) action~\eqref{eq:sl2-act-brute} becomes:
  \begin{equation}
    \label{eq:sl2-act}
    \begin{pmatrix}
      a&b\\c&d
    \end{pmatrix}: w\mapsto \frac{aw+b}{cw+d}, \qquad
    \text{where } w=u+\alli v,
  \end{equation}
  for all three cases parametrised by \(\sigma\) as in Prop.~\ref{pr:sl2-act-brute}.
\end{prop}
\begin{rem}
  We wish to stress that the hypercomplex numbers were not introduced
  here by our intention, arbitrariness or ``generalising
  attitude''~\cite{Pontryagin86a}*{p.~4}. They were naturally created by
  the \(\SL\) action. 
\end{rem}

Notably the action~\eqref{eq:sl2-act} is a group homomorphism of the
group \(\SL\) into transformations of the ``upper half-plane'' on
hypercomplex numbers.  Although dual and double numbers are
algebraically trivial, the respective geometries in the spirit of
\wiki{Erlangen_program}{Erlangen programme} are refreshingly
inspiring~\cites{Kisil05a,Kisil12a,Kisil08a} and provide useful insights even in the
elliptic case~\cite{Kisil06a}.  In order to treat divisors of zero, we
need to consider M\"obius transformations~\eqref{eq:sl2-act} of
conformally completed plane~\cites{HerranzSantander02b,Kisil06b}.

Now we wish to linearise the action~\eqref{eq:g-action} through the
\wiki{Induced_representations}{induced
  representations}~\citelist{\cite{Kirillov76}*{\S~13.2}
  \cite{Kisil97c}*{\S~3.1}}.  We define a map \(r: G\rightarrow H\)
associated to \(p\) and \(s\) from the identities:
\begin{equation}
  \label{eq:r-map}
  r(g)={(s(\omega))}^{-1}g, \qquad \text{where }
  \omega=p(g)\in \Omega .
\end{equation}
Let \(\chi\) be an irreducible representation of \(H\) in a vector space \(V\), then
it induces a representation of \(G\) in the sense of
Mackey~\cite{Kirillov76}*{\S~13.2}.  This representation has the
realisation \(\uir{\chi}\) in the space 
of \(V\)-valued functions by the
formula~\cite{Kirillov76}*{\S~13.2.(7)--(9)}:
\begin{equation} 
  \label{eq:def-ind}
  [\uir{\chi}(g) f](\omega)= \chi(r(g^{-1} * s(\omega)))  f(g^{-1}\cdot \omega),
.
\end{equation}
where \(g\in G\), \(\omega\in\Omega\), \(h\in H\) and \(r: G
\rightarrow H\), \(s: \Omega \rightarrow G\) are maps defined
above; \(*\)~denotes multiplication on \(G\) and \(\cdot\) denotes the
action~\eqref{eq:g-action} of \(G\) on \(\Omega\).

In our consideration \(H\) is always one-dimensional. Traditionally,
an irreducible representation of such a subgroup is supposed to be a
complex valued character.  However, hypercomplex numbers naturally
appeared in the \(\SL\) action~\eqref{eq:sl2-act}, why shall we admit
only \(\rmi^2=-1\) to deliver a character then?

\section{Hypercomplex Characters---an Algebraic Approach}
\label{sec:hyper-characters}

\begin{figure}[htbp]
  \centering
  \includegraphics{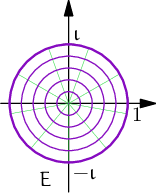}\qquad
  \includegraphics{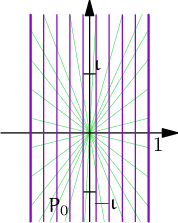}\qquad
  \includegraphics{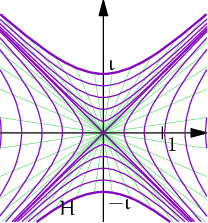}
  \caption[Rotations of wheels---Algebraic]{Rotations of algebraic
    wheels, i.e. the multiplication by \(e^{\alli t}\): elliptic
    (\(E\)), trivial parabolic (\(P_0\)) and hyperbolic (\(H\)). All
    blue orbits are defined by the identity
    \(x^2-\alli^2y^2=r^2\). Green ``spokes'' (straight lines from the
    origin to a point on the orbit) are ``rotated'' from the real axis.}
  \label{fig:rotations}
\end{figure}

As we already mentioned the typical discussion of induced
representations of \(\SL\) is centred around the case \(H=K\) and a
complex valued character of \(K\).  A linear transformation defined by
a matrix~\eqref{eq:k-subgroup} in \(K\) is a rotation of
\(\Space{R}{2}\) by the angle \(t\).  After identification
\(\Space{R}{2}=\Space{C}{}\) this action is given by the
multiplication \(e^{\rmi t}\), with \(\rmi^2=-1\).  The rotation
preserve the (elliptic) metric given by:
\begin{equation}
  \label{eq:ell-metric}
  x^2+y^2=(x+\rmi y)(x-\rmi y).
\end{equation}
Therefore the orbits of rotations are circles, any line passing the
origin (a ``spoke'') is rotated by the angle \(t\), see
Fig.~\ref{fig:rotations}(E).

Introduction of hypercomplex numbers produces the most straightforward
adaptation of this result.

\begin{prop}
  \label{pr:algebraic-char}
  The following table show correspondences between three types of
  algebraic characters:
  \begin{center}
    \begin{tabular}{||c|c|c||}
      \hline\hline
      Elliptic & Parabolic & Hyperbolic\\
      \hline\hline
      \(\strut^{\strut}\rmi^{2}=-1\)&    \(\rmp^2=0\)&    \(\rmh^2=1\)
      \\
      \hline
      \(w=x+\rmi y\) &\(w=x+\rmp y\) &\(w=x+\rmh y\) 
      \\
      \hline
      \(\bar{w}=x-\rmi y\) &\(\bar{w}=x-\rmp y\) &\(\bar{w}=x-\rmh y\) 
      \\
      \hline
      \(\strut^{\strut}e^{\rmi t} = \cos t +\rmi \sin t\)&\(e^{\rmp t} = 1 +\rmp t\)&\(e^{\rmh t} = \cosh t +\rmh \sinh t\)
      \\
      \hline
      \(\strut^{\strut}\modulus[e]{w}^{ 2}=w\bar{w}=x^2+y^2\) &\(\modulus[p]{w}^2=w\bar{w}=x^2\) &\(\modulus[h]{w}^2=w\bar{w}=x^2-y^2\) 
      \\
      \hline
      \(\arg w = \tan^{-1} \frac{y}{x} \frac{\strut}{\strut}\)& \(\arg w = \frac{y}{x}\)&\(\arg w = \tanh^{-1} \frac{y}{x}\)
      \\
      \hline
       unit circle \(\strut^{\strut} \modulus[e]{w}^{2}=1\) & ``unit'' strip \(x=\pm 1\) & unit hyperbola \(\modulus[h]{w}^2=1\)
      \\
      \hline\hline
    \end{tabular}
  \end{center}
  Geometrical action of multiplication by \(e^{\alli t}\) is drawn on
  Fig.~\ref{fig:rotations} for all three cases.
\end{prop}
Explicitly parabolic rotations associated with \(\rme^{\rmp t}\) acts on dual
numbers as follows:
\begin{equation}
  \label{eq:parab-rot-triv}
  \rme^{\rmp x}: a+\rmp b \mapsto a+\rmp (a x+b).
\end{equation}
This links the parabolic case with the Galilean
group~\cite{Yaglom79} of symmetries of the classic mechanics, with
the absolute time disconnected from space.

The obvious algebraic similarity and the connection to classical
kinematic is a wide spread justification for the following viewpoint on
the parabolic case, cf.~\cites{HerranzOrtegaSantander99a,Yaglom79}:
\begin{itemize}
\item the parabolic trigonometric functions are trivial:
  \begin{equation}
    \label{eq:par-trig-0}
    \cosp t =\pm 1, \qquad \sinp t=t;
  \end{equation}
\item the parabolic distance is independent from \(y\) if \(x\neq 0\):
  \begin{equation}
    \label{eq:par-metr-0}
    x^2=(x+\rmp y)(x-\rmp y);
  \end{equation}
\item the polar decomposition of a dual number is defined by~\cite{Yaglom79}*{App.~C(30')}:
  \begin{equation}
    \label{eq:p-polar-yaglom}
    u+\rmp v = u(1+\rmp \frac{v}{u}), \quad \text{ thus }
    \quad \modulus{u+\rmp v}=u, \quad \arg(u+\rmp v)=\frac{v}{u};
  \end{equation}
\item the parabolic wheel looks rectangular, see  Fig.~\ref{fig:rotations}(\(P_0\)).
\end{itemize}

Those algebraic analogies are quite explicit and widely accepted as an
ultimate source for parabolic
trigonometry~\cites{LavrentShabat77,HerranzOrtegaSantander99a,Yaglom79}.
Moreover, those three rotations are all non-isomorphic symplectic
linear transformations of the phase space, which makes them useful in
the context of classical and quantum mechanics~\cites{Kisil10a,Kisil11a}.

However we will see shortly that there exists an alternative with
geometric motivation and connection to equations of mathematical
physics.

\section{A Parabolic Wheel---a Geometrical Viewpoint}
\label{sec:second-attempt}

We make another attempt to describe parabolic rotations. If
multiplication (a linear transformation) is not sophisticated enough for
this we can advance to the next level of complexity:
linear-fractional.

Hypercomplex units do not need to be seen as abstract quantities. As
follows from Lem.~\ref{le:square} the generators of subgroup \(K\),
\(N\) and \(A\) represent units of complex, dual and double
numbers respectively. Their exponentiation to one-parameter
subgroups \(K\), \(N'\) and \(\Aprime\) of \(\SL\) produce matrix forms of
the Euler identities from the fifth row of the table in
Prop.~\ref{pr:algebraic-char}.

Thus we attempt to define characters of subgroups \(K\), \(N'\) and
\(\Aprime\) in term of geometric action of \(\SL\) by M\"obius
transformations. The action~\eqref{eq:sl2-act} is defined on the upper
half-plane and to relate it to unitary characters we wish to transfer
it to the unit disk. In the elliptic case this is done by the Cayley
transform, its action on the subgroup \(K\) is:
\begin{equation}
  \label{eq:ell-cayley}
  \frac{1}{2}
  \begin{pmatrix}
    1 & -\rmi \\ -\rmi &1
  \end{pmatrix}
  \begin{pmatrix}
    \cos t & -\sin t \\ \sin t & \cos t
  \end{pmatrix}
  \begin{pmatrix}
    1 & \rmi \\ \rmi &1
  \end{pmatrix} =
  \begin{pmatrix}
    \rme^{\rmi t} & 0 \\ 0 &  \rme^{-\rmi t}
  \end{pmatrix}.
\end{equation}
The diagonal matrix in the right hand side define the Moebius
transformation which reduces to multiplication by \(e^{2\rmi t}\),
i.e. the elliptic rotation.

A hyperbolic cousin of the Cayley transform is:
\begin{equation}
  \label{eq:hyp-cayley}
  \frac{1}{2}
  \begin{pmatrix}
    1 & \rmh \\ -\rmh &1
  \end{pmatrix}
  \begin{pmatrix}
    \cosh t & \sinh t \\ \sinh t & \cosh t
  \end{pmatrix}
  \begin{pmatrix}
    1 & -\rmh \\ \rmh &1
  \end{pmatrix} =
  \begin{pmatrix}
    \rme^{\rmh t} & 0 \\ 0 &  \rme^{-\rmh t}
  \end{pmatrix},
\end{equation}
similarly produces a Moebius transformation which is the
multiplication by \(e^{2\rmh t}\), which a unitary (Lorentz)
transformation of two-dimensional Minkowski space-time.  

In the parabolic case we use the similar pattern and define the Cayley
transform from the matrix:
\begin{displaymath}
  C_\rmp=
  \begin{pmatrix} 
    1 & -\rmp \\ -\rmp &1  
  \end{pmatrix}
\end{displaymath}
The Cayley transform of matrices~\eqref{eq:n-subgroup} from the
subgroup \(N\) is: 
\begin{equation}
  \label{eq:par-cayley}
  \begin{pmatrix}
    1 & -\rmp\\ -\rmp &1
  \end{pmatrix}
  \begin{pmatrix}
    1 & t \\ 0 & 1
  \end{pmatrix}
  \begin{pmatrix}
    1 & \rmp \\ \rmp &1
  \end{pmatrix} =
  \begin{pmatrix}
    1+\rmp t & t \\ 0 & 1-\rmp t
  \end{pmatrix}=
    \begin{pmatrix}
    \rme^{\rmp t} & t \\ 0 & \rme^{-\rmp t}
  \end{pmatrix}.
\end{equation}
This is not far from the previous identities~\eqref{eq:ell-cayley}
and~\eqref{eq:hyp-cayley}, however, the off-diagonal \((1,2)\)-term
destroys harmony. Nevertheless we will continue a unitary parabolic
rotation to be the M\"obius transformation with the
matrix~\eqref{eq:par-cayley}, which will not be a multiplication by a
scalar anymore.


\begin{figure}[htbp]
  \centering
  \includegraphics{parab-rot-k.eps}\hfill
  \includegraphics{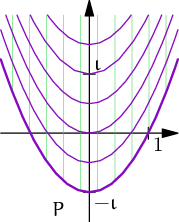}\hfill
  \includegraphics{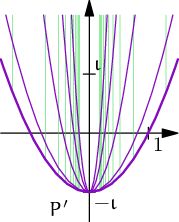}\hfill
  \includegraphics{parab-rot-a.eps}
  \caption[Rotation of wheels--- Geometric]{Rotation of geometric wheels:
    elliptic (\(E\)), two parabolic (\(P\) 
    and \(P'\)) and hyperbolic (\(H\)). Blue orbits are level lines for the
    respective moduli. Green straight lines join points with the same value of
    argument and are drawn with the constant ``angular
    step'' in each case.}
  \label{fig:p-rotations}
\end{figure}



\begin{example}
  \label{ex:upper-triang}
  The parabolic rotations
  with the upper-triangular matrices from the subgroup \(N\) becomes:
  \begin{equation}
    \label{eq:par-moebius}
    \begin{pmatrix}
      \rme^{\rmp t}&t\\0&\rme^{-\rmp t}
    \end{pmatrix}: -\rmp \mapsto t +\rmp(t^2-1).
  \end{equation}
  This coincides with the \emph{cyclic rotations} defined
  in~\cite{Yaglom79}*{\S~8}.  A comparison with the Euler formula
  seemingly confirms that \(\sinp t=t\), but suggests a new expression
  for \(\cosp t\):
  \begin{displaymath}
    \cosp t = 1-t^2, \qquad \sinp t= t.
  \end{displaymath}
  Therefore the parabolic Pythagoras' identity would be:
  \begin{equation}
    \label{eq:pyhagoras-p}
    \sinp^2 t +  \cosp t =1,
  \end{equation}
  which nicely fits in between the elliptic and hyperbolic versions:
  \begin{displaymath}
    \sin^2 t+\cos^2 t  =1, \qquad \sinh^2 t -  \cosh^2 t =-1.
  \end{displaymath}
  The identity~\eqref{eq:pyhagoras-p} is also less trivial than
  the version \( \cosp^2 t =1\)
  from~\cite{HerranzOrtegaSantander99a} (see
  also~\eqref{eq:par-trig-0}, \eqref{eq:par-metr-0}). 
\end{example}


\begin{example}
  \label{ex:lower-triang}
  There is the second option to define parabolic rotations for the
  lower-triangular matrices from the subgroup \(N'\). The important
  difference now is: the reference point cannot be \(-\rmp\) since it
  is a fixed point (as well as any point on the vertical axis).
  Instead we take \(\rmp^{-1}\), which is an ideal element (a point at
  infinity~\cite{Yaglom79}*{App.~C}) since \(\rmp\) is a divisor of
  zero. The proper compactifications by ideal elements for all three
  cases were discussed in~\cites{Kisil06b}.
  
  We get for the subgroup \(N'\):
  \begin{equation}
    \label{eq:par-moebius-prime}
    \begin{pmatrix}
      \rme^{-\rmp t}&0\\t&\rme^{\rmp t}
    \end{pmatrix}:\  \frac{1}{\rmp} \mapsto \frac{1}{t}+ \rmp \left(1-\frac{1}{t^2}\right). 
  \end{equation}
  
  A comparison with~\eqref{eq:par-moebius} shows that this form is
  obtained by the change \(t\mapsto t^{-1}\). The same transformation
  gives new expressions for parabolic trigonometric functions. The
  parabolic ``unit circle'' (or
  \emph{cycle}~\cites{Kisil05a,Kisil12a,Yaglom79}) is defined by the equation
  \(x^2-y=1\) in both cases, see Fig.~\ref{fig:p-rotations}(\(P\)) and
  (\(P'\)). However other orbits are different and we will give their
  description in the next Section.
\end{example}
Fig.~\ref{fig:p-rotations} illustrates M\"obius actions of
matrices~\eqref{eq:ell-cayley}, \eqref{eq:hyp-cayley}
and~\eqref{eq:par-cayley} on the respective ``unit disk'', which are
images of the upper half-planes under respective Cayley
transforms~\citelist{\cite{Kisil05a}*{\S~\ref{E-sec:unit-circles}} \cite{Kisil12a}*{Ch.~10}}.

\section{Rebuilding Algebraic Structures from Geometry}
\label{sec:rebu-algebr-struct}

We want induced representations to be linear, to this end the inducing
character shall be linear as well.  Rotations in elliptic and
hyperbolic cases are given by products of complex or double numbers
respectively and thus are linear. However non-trivial parabolic
rotations~\eqref{eq:par-moebius} and~\eqref{eq:par-moebius-prime}
(Fig.~\ref{fig:p-rotations}(\(P\)) and~(\(P'\))) are not linear.  Can
we find algebraic operations for dual numbers, which will linearise
those M\"obius transformations?

It is common in mathematics to ``revert a theorem into a definition''
and we will use this systematically to recover a compatible algebraic
structure.

\subsection{Modulus and Argument}
\label{sec:modulus-argument}

In the elliptic and hyperbolic cases orbits of rotations are points
with the constant norm (modulus): either \(x^2+y^2\) or
\(x^2-y^2\). In the parabolic case we employ this point of view
as well:
\begin{defn}
  \label{de:norm}
  Orbits of actions~\eqref{eq:par-moebius}
  and~\eqref{eq:par-moebius-prime} are contour lines for the following
  functions which we call respective moduli (norms):
  \begin{equation}
    \label{eq:parab-norm}
    \text{ for } N:\ \modulus{u+\rmp v}=u^2-v, \qquad \text{ for } N':\
    \modulus{u+\rmp v}'=\frac{u^2}{v+1}.
  \end{equation}
\end{defn}
\begin{rem}
  \begin{enumerate}
  \item The expression \( \modulus{(u,v)}=u^2-v\) represents a parabolic distance
    from \((0,\frac{1}{2})\) to \((u,v)\),
    see~\citelist{\cite{Kisil05a}*{Lem.~\ref{E-le:n-orbits-concentric}} \cite{Kisil12a}*{Ex.~10.6}},  which is in
    line with the ``parabolic Pythagoras' identity''~\eqref{eq:pyhagoras-p}.
  \item Modulus for \(N'\) expresses the parabolic focal length from \((0,-1)\)
    to  \((u,v)\) as described
    in~\citelist{\cite{Kisil05a}*{Lem.~\ref{E-le:np-orbits-p-confocal}} \cite{Kisil12a}*{Ex.~10.7}}.
  \end{enumerate}
\end{rem}

The only straight lines preserved by both the parabolic rotations
\(N\) and \(N'\) are vertical lines, thus we will treat them as
``spokes'' for parabolic ``wheels''. Elliptic spokes in mathematical
terms are ``points on the complex plane with the same argument'', thus
we again use this for the parabolic definition:
\begin{defn}
  \label{de:arg}
  Parabolic arguments are defined as follows:
  \begin{equation}
    \label{eq:parab-arg}
    \text{ for } N:\ \arg(u+\rmp v)=u, \qquad \text{ for } N':\
    \arg'(u+\rmp v)=\frac{1}{u}.
  \end{equation}
\end{defn}
Both Definitions~\ref{de:norm} and~\ref{de:arg} possess natural
properties with respect to parabolic rotations:
\begin{prop}
  \label{pr:rot-norm-arg}
  Let \(w_t\) be a parabolic rotation of \(w\) by an angle \(t\)
  in~\eqref{eq:par-moebius} or 
  in~\eqref{eq:par-moebius-prime}. Then:
  \begin{displaymath}
    \modulus{w_t}^{(\prime)}=\modulus{w}^{(\prime)},\qquad 
    \arg^{(\prime)} w_t=\arg^{(\prime)} w+t,
  \end{displaymath}
  where primed versions are used for subgroup \(N'\).
\end{prop}
All proofs in this and the following Sections were performed through
symbolic calculations on a computer. 
Details can be found in the earlier
(more heuristic) paper on this topic~\cite{Kisil07a}.

\begin{rem}
  Note that in the commonly accepted
  approach~\cite{Yaglom79}*{App.~C(30')} parabolic modulus and
  argument are given by expressions~\eqref{eq:p-polar-yaglom}, which
  are, in a sense, opposite to our agreements.
\end{rem}

\subsection{Rotation as Multiplication}
\label{sec:rotat-as-mult}

We revert again theorems into definitions to assign multiplication.
In fact, we consider parabolic rotations as multiplications by
unimodular numbers thus we define multiplication through an extension
of properties from Proposition~\ref{pr:rot-norm-arg}:
\begin{defn}
  \label{de:product}
  The product of vectors \(w_1\) and \(w_2\) is defined by the following
  two conditions:  
  \begin{enumerate}
  \item \(\arg^{(\prime)}(w_1 w_2)=\arg^{(\prime)} w_1 + \arg^{(\prime)} w_2\);
  \item \(\modulus{w_1 w_2}^{(\prime)} =\modulus{w_1}^{(\prime)}\cdot
    \modulus{w_2}^{(\prime)}\). 
  \end{enumerate}
\end{defn}
We also need a special form of parabolic conjugation, which coincides
with sign reversion of the argument.
\begin{defn}
  Parabolic conjugation is given by
  \begin{equation}
    \label{eq:conjugation}
    \overline{u+\rmp v}=-u+\rmp v.
  \end{equation}
\end{defn}
Obviously we have the properties:
\(\modulus{\overline{w}}^{(\prime)}=\modulus {w}^{(\prime)}\) and
\(\arg^{(\prime)}{\overline{w}}=-\arg^{(\prime)} {w}\). 
A combination of Definitions~\ref{de:norm}, \ref{de:arg}
and~\ref{de:product} uniquely determine expressions for products.
\begin{prop}
  The parabolic product of vectors is defined by formulae:
  \begin{align}
    \label{eq:parab-prod}
    \text{for } &N:& (u,v)*(u',v') & =
    (u+u',(u+u')^2-(u^2-v)(u'^2-v')); \\
    \label{eq:parab-prime-prod}
    \text{for } &N':&(u,v)*(u',v') & = \left(\frac{uu'}{u+u'},\frac{(v+1)(v'+1)}{(u+u')^2}-1\right).
  \end{align}
\end{prop}
Although both expressions look unusual they have many familiar properties:
\begin{prop}
  \label{pr:mul}
  Both products~\eqref{eq:parab-prod} and~\eqref{eq:parab-prime-prod}
  satisfy the following conditions:
  \begin{enumerate}
  \item\label{item:prod-comm-ass} They are commutative and
    associative;
  \item\label{it:rot-as-mult} The respective
    rotations~\eqref{eq:par-moebius} and~\eqref{eq:par-moebius-prime}
    are given by multiplications with a dual number with the unit norm.
  \item\label{item:prod-inv} The product \(w_1\bar{w}_2\) is invariant
    under respective rotations~\eqref{eq:par-moebius}
    and~\eqref{eq:par-moebius-prime}.
  \item\label{item:prod-norm-sq} 
    For any dual number \(w\) the following identity holds:
    \begin{displaymath}
      \modulus{w\bar{w}}= \modulus{w}^2.
    \end{displaymath}
  \end{enumerate}
\end{prop}
In particular, the property~(\ref{item:prod-inv}) will be crucial below
for an inner product~\eqref{eq:inner-product}, which makes induced
representations unitary.
\section{Invariant Linear Algebra}
\label{sec:invar-line-algebra}

Now we wish to define a linear structure on \(\Space{R}{2}\) which
would be invariant under point multiplication from the previous
Subsection (and thus under the parabolic rotations,
cf.~Prop.\ref{pr:mul}(\ref{it:rot-as-mult})). Multiplication by a real
scalar is straightforward (at least for a positive scalar): it should
preserve the argument and scale the norm of a vector. Thus we have
formulae for \(a>0\):
\begin{eqnarray}
  \label{eq:scalar-prod}
  a\cdot (u,v)&=&(u,a v+u^2(1-a))\quad\text{for } N,\\
  \label{eq:scalar-prod-prime}
  a\cdot (u,v)&=&\left(u,\frac{v+1}{a}-1\right)\quad\text{for } N'.
\end{eqnarray}

On the other hand, the addition of vectors can be done in several
different ways. One of them, related to tropical
mathematics~\cites{Litvinov05,Butkovic10a}, is outlined in
Example~\ref{ex:tropical}. Here, we present another alternative with
all due details.

Addition of vectors for both subgroups \(N\) and \(N'\) can be defined
by the common rules, where subtle differences are hidden within
corresponding Definitions~\ref{de:norm} (norms) and~\ref{de:arg}
(arguments).
\begin{defn}
  \label{de:p-add}
  Parabolic addition of vectors is defined by the following formulae:
  \begin{eqnarray}
    \label{eq:p-add-arg-exotic}
    \arg^{(\prime)}(w_1+w_2)&=&\frac{\arg^{(\prime)} w_1\cdot
      \modulus{w_1}^{(\prime)} 
      +\arg^{(\prime)} w_2\cdot\modulus{w_2}^{(\prime)}}{\modulus{w_1+w_2}^{(\prime)}},\\
    \label{eq:p-add-norm-exotic}
    \modulus{w_1+w_2}^{(\prime)}&=&\modulus{w_1}^{(\prime)}\pm\modulus{w_2}^{(\prime)},
  \end{eqnarray}
  where primed versions are used for the subgroup \(N'\).
\end{defn}

The rule for the norm of sum~\eqref{eq:p-add-norm-exotic} may look too
trivial at the first glance. We should say in its defence that it nicely
sits in between the elliptic \(\modulus{w+w'}\leq
\modulus{w}+\modulus{w'}\) and hyperbolic \(\modulus{w+w'}\geq
\modulus{w}+\modulus{w'}\) triangle inequalities for norms. 
  

The rule~\eqref{eq:p-add-arg-exotic} for argument of the sum is not
arbitrary as well.  From the Sine Theorem in the Euclidean geometry we
can deduce that:
  \begin{displaymath}
    \sin(\phi-\psi')=\frac{\modulus{w}\cdot\sin(\psi-\psi')}{\modulus{w+w'}},\qquad
    \sin(\psi'-\phi)=\frac{\modulus{w'}\cdot\sin(\psi-\psi')}{\modulus{w+w'}},
  \end{displaymath}
  where \(\psi^{(\prime)}=\arg w^{(\prime)}\) and \(\phi=\arg
  (w+w^{(\prime)})\). Using parabolic expression~\eqref{eq:par-trig-0}
  for the sine \(\sinp \theta=\theta\) we obtain the arguments addition
  formula~\eqref{eq:p-add-arg-exotic}.

A proper treatment of zeros in denominator
of~\eqref{eq:p-add-arg-exotic} can be achieved through a representation of
a dual number \(w=u+\rmp v\) as a pair of homogeneous polar coordinates
\([a,r]=[ \modulus{w}^{(\prime)} \cdot \arg^{(\prime)} w,
\modulus{w}^{(\prime)}]\) (dashed version for the subgroup
\(N'\)). Then the above addition is defined component-wise in the
homogeneous coordinates:
\begin{displaymath}
  w_1+w_2=[a_1+a_2, r_1+r_2], \qquad \text{ where } w_i=[a_i,r_i].
\end{displaymath}
The multiplication from Defn.~\ref{de:product} is given in the
homogeneous polar coordinates by:
\begin{displaymath}
  w_1\cdot w_2=[a_1r_2+a_2r_1, r_1 r_2], \qquad \text{ where } w_i=[a_i,r_i].
\end{displaymath}
Thus homogeneous coordinates linearise the
addition~\eqref{eq:p-add-arg-exotic}--\eqref{eq:p-add-norm-exotic} and
multiplication by a scalar~\eqref{eq:scalar-prod}. A transition to
other more transparent coordinates shall be treated withing birational
geometry framework~\cite{KollarMori08}.

Both
formulae~\eqref{eq:p-add-arg-exotic}--\eqref{eq:p-add-norm-exotic}
together uniquely define explicit expressions for addition of vectors.
However those expressions are rather cumbersome and not really much
needed. Instead we list properties of these operations:
\begin{prop}
  Vector additions for subgroups \(N\) and \(N'\) defined
  by~\eqref{eq:p-add-arg-exotic}--\eqref{eq:p-add-norm-exotic} satisfy
  the following conditions:
  \begin{enumerate}
  \item\label{item:add-is-comm-ass} They are commutative and
    associative.
  \item\label{item:distrib} They are distributive for
    multiplications~\eqref{eq:parab-prod} and~\eqref{eq:parab-prime-prod};
    consequently:
  \item\label{item:add-rot-inv} They are parabolic rotationally invariant;
  \item\label{item:distrib-scalar} They are distributive in both ways for the scalar
    multiplications~\eqref{eq:scalar-prod}
    and~\eqref{eq:scalar-prod-prime} respectively:
    \begin{displaymath}
      a\cdot(w_1+w_2)=a\cdot w_1+a\cdot w_2,\qquad
      (a+b)\cdot w=a\cdot w+b\cdot w.
    \end{displaymath}
  \end{enumerate}
\end{prop}
To complete the construction we need to define the zero vector and
the inverse. The inverse of \(w\) has the same argument as \(w\) and
the opposite norm.
\begin{prop}
  \begin{itemize}
  \item[(\(N\))] The zero vector is \((0,0)\) and consequently the inverse
    of \((u,v)\) is \((u,2u^2-v)\).
  \item[(\(N'\))] The zero vector is \((\infty,-1)\) and consequently the inverse
    of \((u,v)\) is \((u,-v-2)\).
  \end{itemize}
\end{prop}
Thereafter we can check that scalar
multiplications by negative reals are given by the same identities
~\eqref{eq:scalar-prod} and \eqref{eq:scalar-prod-prime} as for
positive ones.

\begin{rem}
  \label{re:conformality}
  The irrelevance of the standard linear structure for parabolic
  rotations manifests itself in many different ways, e.g. in an
  apparent ``non-conformality'' of lengths from parabolic foci, that
  is with the parameter \(\rs=0\)
  in~\citelist{\cite{Kisil05a}*{Prop.~\ref{E-it:conformity-length-foci}}
    \cite{Kisil12a}*{Ex.~7.14.iii}}.  An adjustment of notions to the
  proper framework restores the clear picture.

  The initial definition of
  conformality~\citelist{\cite{Kisil05a}*{Defn.~\ref{E-de:conformal}}
    \cite{Kisil12a}*{Defn.~7.13}} considered 
  the usual limit \(y'\rightarrow y\) along a straight line, i.e.
  ``spoke'' in terms of Fig.~\ref{fig:rotations}. This is justified in
  the elliptic and hyperbolic cases. However in the parabolic setting
  the proper ``spokes'' are vertical lines, see
  Fig.~\ref{fig:p-rotations}(\(P\)) and (\(P'\)), so the limit should
  be taken along
  them~\citelist{\cite{Kisil05a}*{Prop.~\ref{E-pr:parab-conf}}
    \cite{Kisil12a}*{\S~11.6.1}}. 
\end{rem}

\section{Induced Representations}
\label{sec:concl-induc-repr}

We discussed above various implementations of hypercomplex unitary
characters. Now we can return to consideration of induced
representations. We can notice that only the subgroup \(K\) requires a
complex valued character due to the fact of its compactness. For
subgroups \(N'\) and \(\Aprime\) we can consider characters of all three
types---elliptic, parabolic and hyperbolic. Moreover a parabolic
character can be taken either as algebraic~\eqref{eq:parab-rot-triv}
or any of two geometric~\eqref{eq:par-moebius}
and~\eqref{eq:par-moebius-prime}.  Therefore we have seven essentially
different induced representations, which multiply types to eleven (counting
flavours of parabolic characters).

\begin{example}
  Consider the subgroup \(H=K\), then we are limited to complex valued
  characters of \(K\) only. All of them are of the form \(\chi_k\):
  \begin{equation}
    \label{eq:k-character}
    \chi_k\begin{pmatrix}
      \cos t &  \sin t\\
      -\sin t & \cos t
    \end{pmatrix}=e^{-\rmi k t}, \qquad \text{ where }
    k\in\Space{Z}{}.
  \end{equation}
  Using the explicit form~\eqref{eq:s-map} of the map \(s\) we find 
  the map \(r\) given in~\eqref{eq:r-map} as follows:
  \begin{displaymath}
    r
    \begin{pmatrix}
      a&b\\c&d
    \end{pmatrix}
    =\frac{1}{\sqrt{c^2+d^2}}
    \begin{pmatrix}
      d&-c\\c&d
    \end{pmatrix}\in K.
  \end{displaymath}
  Therefore:
  \begin{displaymath}
    r(g^{-1} * s(u,v))  =  
    \frac{1}{\sqrt{(c u+d)^2 +(cv)^2}}
    \begin{pmatrix}
      cu+d&-cv\\cv&cu+d
    \end{pmatrix}, \text{ where } g^{-1}=    \begin{pmatrix}
      a&b\\c&d
    \end{pmatrix}.
  \end{displaymath}
  Substituting this into~\eqref{eq:k-character} and combining with the
  M\"obius transformation of the domain~\eqref{eq:sl2-act} we get the
  explicit realisation \(\uir{k}{}\) of the induced representation~\eqref{eq:def-ind}:
  \begin{equation}
    \label{eq:discrete}
    \uir{k}{}(g) f(w)=\frac{\modulus{cw+d}^k}{(cw+d)^k}f\left(\frac{aw+b}{cw+d}\right),
    \quad \text{ where } g^{-1}=\begin{pmatrix}a&b\\c&d
    \end{pmatrix}, \ w=u+\rmi v.
  \end{equation}
  This representation acts on complex valued functions in the upper half-plane
  \(\Space[+]{R}{2}=\SL/K\) and belongs to the discrete
  series~\cite{Lang85}*{\S~IX.2}. 
  \begin{prop}
    \label{pr:k-eigen}
    Let \(f_k(w)=\frac{\modulus{w-\rmi}^{k}}{(w-\rmi)^{k}}\) for \(k=2\), \(3\), \ldots, then
    \begin{enumerate}
    \item \(f_k\) is an eigenvector for any operator \(\uir{k}{}(h)\),
      where \(h\in K\), with the eigenvalue
      \(\chi_k(h)\)~\cite{Lang85}*{\S~IX.2}. 
    \item The function \(K(z,w)=\uir{k}{}(s(z)) f_k(w)\), where \(s(z)\)
      is defined in~\eqref{eq:s-map}, is the Bergman reproducing
      kernel up to the factor \(\modulus{\frac{z-\rmi}{w-\rmi}}^k\)
      in the upper half-plane~\cite{Kisil97c}*{\S~3.2}.
    \end{enumerate}
  \end{prop}
  Similarly we can get the Cauchy kernel for the limiting case \(k=1\)
  of the mock discrete series~\cite{Lang85}*{Ch.~IX}. There are many
  other important connections of representation~\eqref{eq:discrete} with
  complex analysis and operator theory. For example, M\"obius
  transformations of operators lead to Riesz-Dunford functional calculus
  and associated spectrum~\cite{Kisil02a}.
\end{example}

\begin{example}
  \label{ex:n-induced}
  In the case of the subgroup \(N\) there is a wider choice of
  possible characters.
  \begin{enumerate}
  \item Traditionally only complex valued characters of the subgroup
    \(N\) are considered, they are:
    \begin{equation}
      \label{eq:np-character}
      \chi_{\tau}^{\Space{C}{}}\begin{pmatrix}
        1 &  0\\
        t & 1
      \end{pmatrix}=e^{\rmi \tau t}, \qquad \text{ where }
      \tau\in\Space{R}{}.
    \end{equation}
    A direct calculation shows that:
    \begin{displaymath}
      r
      \begin{pmatrix}
        a&b\\c&d
      \end{pmatrix}
      =
      \begin{pmatrix}
        1&0\\\frac{c}{d}&1
      \end{pmatrix}\in N'.
    \end{displaymath}
    Thus:
    \begin{equation}
      \label{eq:np-char-part}
      r(g^{-1}*s(u,v))=
      \begin{pmatrix}
        1&0\\\frac{cv}{d+cu}&1
      \end{pmatrix}, \quad\text{ where } g^{-1}=    \begin{pmatrix}
        a&b\\c&d
      \end{pmatrix}.
    \end{equation}
    A substitution of this value into the
    character~\eqref{eq:np-character} together with the M\"obius
    transformation~\eqref{eq:sl2-act} we obtain the next realisation of~\eqref{eq:def-ind}:
    \begin{displaymath}
      \uir[{\Space{C}{}}]{\tau}(g) f(w)= \exp\left(\rmi\frac{\tau c v}{cu+d} \right)
      f\left(\frac{aw+b}{cw+d}\right), 
      \quad \text{where } w=u+\rmp v, \  g^{-1}=\begin{pmatrix}a&b\\c&d
      \end{pmatrix}.
    \end{displaymath}
    The representation acts on the space of \emph{complex} valued
    functions on the upper half-plane \(\Space[+]{R}{2}\), which is
    subset of \emph{dual} numbers as a homogeneous space \(\SL/N\).
    The mixture of complex and dual numbers in the same expression is
    confusing.
  \item The parabolic character \(\chi_{\tau}\) with the algebraic flavour
    is provided by multiplication~\eqref{eq:parab-rot-triv} with the
    dual number:
    \begin{displaymath}
      \chi_{\tau}\begin{pmatrix}
        1 &  0\\
        t & 1
      \end{pmatrix}=e^{\rmp \tau t}=1+\rmp \tau t, \qquad \text{ where }
      \tau\in\Space{R}{}.
    \end{displaymath}
    As before substitute the value~\eqref{eq:np-char-part} into this
    character we receive the representation:
    \begin{displaymath}
      \uir{\tau}{}(g) f(w)= \left(1+\rmp\frac{\tau c v}{cu+d} \right)
      f\left(\frac{aw+b}{cw+d}\right),
    \end{displaymath}
    where \(w\), \(\tau\) and \(g\) are as above.  The representation
    is defined on the space of dual numbers valued functions on the
    upper half-plane of dual numbers.  Thus expression contains only
    dual numbers with their usual algebraic operations. Thus it is
    linear with respect to them.
  \item The geometric character \(\chi_{\tau}^{g}\) is given by the
    action~\eqref{eq:par-moebius}. Then the corresponding
    representation acts again on the space of dual numbers valued
    functions on the upper half-plane of dual numbers as follows:
    \begin{displaymath}
      \uir[g]{\tau}(g) f(w)= \left(1+\rmp\frac{2\tau c v}{cu+d} \right)
      f\left(\frac{aw+b}{cw+d}\right)+\frac{\tau c v}{cu+d}+\rmp \frac{(\tau c v)^2}{(cu+d)^2}, 
    \end{displaymath}
    where \(w\), \(\tau\) and \(g\) are as above. This representation
    is linear with respect to operations~\eqref{eq:scalar-prod},
    \eqref{eq:p-add-arg-exotic} and~\eqref{eq:p-add-norm-exotic}.
  \end{enumerate}
\end{example}
All characters in the previous Example are unitary, the first two in a
conventional sense and the last one in the sense of
Prop.~\ref{pr:mul}. Then the general scheme of induced
representations~\cite{Kirillov76}*{\S~13.2} implies their unitarity in
proper senses.
\begin{thm}
  \label{th:unitarity}
  All three representations of \(\SL\) from Example~\ref{ex:n-induced}
  are unitary on the space of function on the upper half-plane
  \(\Space[+]{R}{2}\) of dual numbers with the inner product:
  \begin{equation}
    \label{eq:inner-product}
    \scalar{f_1}{f_2}=\int_{\Space[+]{R}{2}} f_1(w)
    \bar{f}_2(w)\,\frac{du\,dv}{v^2}, \qquad \text{ where } w=u+\rmp v,
  \end{equation}
  and we use
  \begin{enumerate}
  \item the conjugation and multiplication of functions'
    values in algebras of complex and dual numbers for
    representations \(\uir[{\Space{C}{}}]{\tau}\) and
    \(\uir{\tau}\) respectively;
  \item conjugation~\eqref{eq:conjugation} and
    multiplication~\eqref{eq:parab-prod} of functions' values for the
    representation \(\uir[g]{\tau}\).
  \end{enumerate}
\end{thm}
The inner product~\eqref{eq:inner-product} is positive defined for
the representation \(\uir[{\Space{C}{}}]{\tau}\) but is not for two
others. The respective spaces are parabolic cousins of the Krein
spaces~\cite{ArovDym08}, which are hyperbolic in our sense.

\section{Similarity and Correspondence: Ladder Operators}
\label{sec:correspondence}

From the above observation we can deduce the following empirical
principle, which has a heuristic value.

\begin{principle}[Similarity and correspondence]
  \label{pr:simil-corr-principle}
  \begin{enumerate}
  \item Subgroups \(K\), \(N\) and \(A\) play the similar r\^ole in a
    structure of the group \(\SL\) and its representations.
  \item The subgroups shall be swapped simultaneously with the
    respective replacement of hypercomplex unit \(\alli\).
  \end{enumerate}
\end{principle}
The first part of the Principle (similarity) does not look sound
alone. It is enough to mention that the subgroup \(K\) is compact (and
thus its spectrum is discrete) while two other subgroups are not.
However in a conjunction with the second part (correspondence) the
Principle have received the following confirmations so far:
\begin{itemize}
\item The action of \(\SL\) on the homogeneous space \(\SL/H\) for
  \(H=K\), \(N'\) or \(\Aprime\) is given by linear-fractional
  transformations of complex, dual or double numbers respectively
  (Prop.~\ref{pr:sl2-act}).  
\item Subgroups \(K\), \(N'\) or \(\Aprime\) are isomorphic to the groups of
  unitary rotations of respective unit cycles in complex, dual or
  double numbers (Prop.~\ref{pr:algebraic-char}).   
\item Representations induced from subgroups \(K\), \(N'\) or \(\Aprime\)
  are unitary if the inner product spaces of functions with values in
  complex, dual or double numbers (Thm.~\ref{th:unitarity}).
\end{itemize}
\begin{rem}
  The principle of similarity and correspondence resembles
  supersymmetry between bosons and fermions in particle physics, but
  we have similarity between three different types of entities in our case.
\end{rem}
Let us give another illustration to the Principle. Consider the Lie
algebra \(\algebra{sl}_2\) of the group \(\SL\). Pick up the following
basis in \(\algebra{sl}_2\)~\cite{MTaylor86}*{\S~8.1}:  
\begin{displaymath}
  A= \frac{1}{2}
  \begin{pmatrix}
    -1&0\\0&1
  \end{pmatrix},\quad 
  B= \frac{1}{2} \
  \begin{pmatrix}
    0&1\\1&0
  \end{pmatrix}, \quad 
  Z=
  \begin{pmatrix}
    0&1\\-1&0
  \end{pmatrix}.
\end{displaymath}
The commutation relation between those elements are:
\begin{equation}
  \label{eq:sl2-commutator}
  [Z,A]=2B, \qquad [Z,B]=-2A, \qquad [A,B]=- \frac{1}{2} Z.
\end{equation} 
Let \(\uir{}{}\) be a representation of the group \(\SL\) in a space
\(V\). Consider the derived representation \(d\uir{}{}\) of the Lie
algebra \(\algebra{sl}_2\)~\cite{Lang85}*{\S~VI.1} and denote
\(\tilde{X}=d\uir{}{}(X)\) for \(X\in\algebra{sl}_2\). To see the
structure of the representation \(\uir{}{}\) we can decompose the
space \(V\) into eigenspaces of the operator \(\tilde{X}\) for some
\(X\in \algebra{sl}_2\), cf.  Prop.~\ref{pr:k-eigen} or the Taylor
series in complex analysis.
\begin{example}
  \label{ex:ladder-op}
  It would not be surprising that we are going to consider three cases:
  \begin{enumerate}
  \item \label{it:ladder-ell}
    Let \(X=Z\) be a generator of the subgroup
    \(K\)~\eqref{eq:k-subgroup}. Since this is a compact subgroup the
    corresponding eigenspaces \(\tilde{Z} v_k=\rmi k v_k\) are
    parametrised by an integer \(k\in\Space{Z}{}\).  The
    raising/lowering (ladder) operators
    \(\ladder{\pm}\)~\citelist{\cite{Lang85}*{\S~VI.2}
      \cite{MTaylor86}*{\S~8.2}} are defined by the following
    commutation relations:
    \begin{equation}
      \label{eq:raising-lowering}
      [\tilde{Z},\ladder{\pm}]=\lambda_\pm \ladder{\pm}. 
    \end{equation}
    In other words \(\ladder{\pm}\) are eigenvectors for operators 
    \(\loglike{ad}Z\) of adjoint representation of \(\algebra{sl}_2\)~\cite{Lang85}*{\S~VI.2}.
    From the commutators~\eqref{eq:raising-lowering} we deduce that
    \(\ladder{+} v_k\) are eigenvectors of \(\tilde{Z}\) as well:
    \begin{eqnarray*}
      \tilde{Z}(\ladder{+} v_k)&=&(\ladder{+}\tilde{Z}+\lambda_+\ladder{+})v_k=\ladder{+}(\tilde{Z}v_k)+\lambda_+\ladder{+}v_k
      =\rmi k \ladder{+}v_k+\lambda_+\ladder{+}v_k\\
      &=&(\rmi k+\lambda_+)\ladder{+}v_k.
    \end{eqnarray*}
    Thus action of ladder operators on respective eigenspaces can be
    visualised by the diagram:
    \begin{displaymath}
      \xymatrix@1{
        \ldots\, \ar@<.4ex>[r]^{\ladder{+}} &
        \,V_{\rmi k-\lambda}\,  \ar@<.4ex>[l]^{\ladder{-}}\ar@<.4ex>[r]^{\ladder{+}} &
        \,V_{\rmi k}\, \ar@<.4ex>[l]^{\ladder{-}} \ar@<.4ex>[r]^{\ladder{+}} &
        \,V_{\rmi k+ \lambda}\,\ar@<.4ex>[l]^{\ladder{-}}  \ar@<.4ex>[r]^{\ladder{+}}
        &
        \,\ldots\ar@<.4ex>[l]^{\ladder{-}}}
    \end{displaymath}
    Assuming \(\ladder{+}=a\tilde{A}+b\tilde{B}+c\tilde{Z}\) from the
    relations~\eqref{eq:sl2-commutator} and defining
    condition~\eqref{eq:raising-lowering} we obtain linear equations
    with unknown \(a\), \(b\) and \(c\): 
    \begin{displaymath}
      c=0, \qquad 2a=\lambda_+ b, \qquad -2b=\lambda_+ a.
    \end{displaymath}
    The equations have a solution if and only if \(\lambda_+^2+4=0\),
    and the raising/lowering operators are \(\ladder{\pm}=\pm\rmi \tilde{A}+\tilde{B}\).
  \item \label{it:ladder-hyp}
    Consider the case \(X=2B\) of a generator of the subgroup
    \(\Aprime\)~\eqref{eq:ap-subgroup}. The subgroup is not compact
    and eigenvalues of the operator \(\tilde{B}\) can be arbitrary,
    however raising/lowering operators are still
    important~\citelist{\cite{HoweTan92}*{\S~II.1}
      \cite{Mazorchuk09a}*{\S~1.1}}.  We again seek a solution in the
    form \(\ladder[h]{+}=a\tilde{A}+b\tilde{B}+c\tilde{Z}\) for the commutator
    \([2\tilde{B},\ladder[h]{+}]=\lambda \ladder[h]{+}\). We will get the system:
    \begin{displaymath}
      4c=\lambda a,\qquad
      b=0,\qquad
      {a}=\lambda c.
    \end{displaymath}
    A solution exists if and only if \(\lambda^2=4\). There are
    obvious values \(\lambda=\pm 2\) with the ladder operators
    \(\ladder[h]{\pm}=\pm2\tilde{A}+\tilde{Z}\),
    see~\citelist{\cite{HoweTan92}*{\S~II.1}
      \cite{Mazorchuk09a}*{\S~1.1}}. Each indecomposable
    \(\algebra{sl}_2\)-module is formed by a one-dimensional chain of
    eigenvalues with a transitive action of ladder
    operators.
    
    Admitting double numbers we have an extra possibility to satisfy
    \(\lambda^2=4\) with values \(\lambda=\pm2\rmh\).  Then there is an
    additional pair of hyperbolic ladder operators
    \(\ladder[\rmh]{\pm}=\pm2\rmh\tilde{A}+\tilde{Z}\), which shift eigenvectors
    in the ``orthogonal'' direction to the standard operators \(\ladder[h]{\pm}\).
    Therefore an indecomposable \(\algebra{sl}_2\)-module can be
    parametrised by a two-dimensional lattice of eigenvalues on the
    double number plane:
\begin{displaymath}
  \xymatrix@R=2.5em@C=1.5em@M=.5em{
      & 
    \,\ldots\, \ar@<.4ex>[d]^{\ladder[\rmh]{+}} &  
    \,\ldots\, \ar@<.4ex>[d]^{\ladder[\rmh]{+}} & 
    \,\ldots\,  \ar@<.4ex>[d]^{\ladder[\rmh]{+}}  & 
    \\
    \ldots\, \ar@<.4ex>[r]^-{\ladder[h]{+}} & 
    \,V_{(n-2)+\rmh (k-2)}\,  \ar@<.4ex>[l]^-{\ladder[h]{-}}\ar@<.4ex>[r]^{\ladder[h]{+}}
     \ar@<.4ex>[u]^{\ladder[\rmh]{-}} \ar@<.4ex>[d]^{\ladder[\rmh]{+}} &  
    \,V_{n+\rmh (k-2)}\, \ar@<.4ex>[l]^{\ladder[h]{-}} \ar@<.4ex>[r]^{\ladder[h]{+}}
     \ar@<.4ex>[u]^{\ladder[\rmh]{-}} \ar@<.4ex>[d]^{\ladder[\rmh]{+}} & 
    \,V_{(n+2)+\rmh (k-2)}\,\ar@<.4ex>[l]^{\ladder[h]{-}}  \ar@<.4ex>[r]^-{\ladder[h]{+}}
     \ar@<.4ex>[u]^{\ladder[\rmh]{-}} \ar@<.4ex>[d]^{\ladder[\rmh]{+}}    & 
    \,\ldots\ar@<.4ex>[l]^-{\ladder[h]{-}}\\
    \ldots\, \ar@<.4ex>[r]^-{\ladder[h]{+}} & 
    \,V_{(n-2)+\rmh k}\,  \ar@<.4ex>[l]^-{\ladder[h]{-}}\ar@<.4ex>[r]^{\ladder[h]{+}}
     \ar@<.4ex>[u]^{\ladder[\rmh]{-}} \ar@<.4ex>[d]^{\ladder[\rmh]{+}} &  
    \,V_{n+\rmh k}\, \ar@<.4ex>[l]^{\ladder[h]{-}} \ar@<.4ex>[r]^{\ladder[h]{+}}
     \ar@<.4ex>[u]^{\ladder[\rmh]{-}} \ar@<.4ex>[d]^{\ladder[\rmh]{+}}& 
    \,V_{(n+2)+\rmh k}\,\ar@<.4ex>[l]^{\ladder[h]{-}}  \ar@<.4ex>[r]^-{\ladder[h]{+}}
     \ar@<.4ex>[u]^{\ladder[\rmh]{-}} \ar@<.4ex>[d]^{\ladder[\rmh]{+}}    & 
    \,\ldots\ar@<.4ex>[l]^-{\ladder[h]{-}}\\
    \ldots\, \ar@<.4ex>[r]^-{\ladder[h]{+}} & 
    \,V_{(n-2)+\rmh (k+2)}\,  \ar@<.4ex>[l]^-{\ladder[h]{-}}\ar@<.4ex>[r]^{\ladder[h]{+}}
     \ar@<.4ex>[u]^{\ladder[\rmh]{-}} \ar@<.4ex>[d]^{\ladder[\rmh]{+}} &  
    \,V_{n+\rmh (k+2)}\, \ar@<.4ex>[l]^{\ladder[h]{-}} \ar@<.4ex>[r]^{\ladder[h]{+}} 
     \ar@<.4ex>[u]^{\ladder[\rmh]{-}} \ar@<.4ex>[d]^{\ladder[\rmh]{+}}& 
    \,V_{(n+2)+\rmh (k+2)}\,\ar@<.4ex>[l]^{\ladder[h]{-}}  \ar@<.4ex>[r]^-{\ladder[h]{+}}
     \ar@<.4ex>[u]^{\ladder[\rmh]{-}} \ar@<.4ex>[d]^{\ladder[\rmh]{+}}    & 
    \,\ldots\ar@<.4ex>[l]^-{\ladder[h]{-}}\\
      & 
    \,\ldots\, \ar@<.4ex>[u]^{\ladder[\rmh]{-}} &  
    \,\ldots\, \ar@<.4ex>[u]^{\ladder[\rmh]{-}} & 
    \,\ldots\,  \ar@<.4ex>[u]^{\ladder[\rmh]{-}}  & }
\end{displaymath}
    
  \item \label{it:ladder-par}
    Finally consider the case of a generator \(X=-B+Z/2\) of the
    subgroup \(N'\)~\eqref{eq:np-subgroup}. According to the above
    procedure we get the equations:
    \begin{displaymath}
      b+2c=\lambda a,\qquad
      -a=\lambda b,\qquad
      \frac{a}{2}=\lambda c,
    \end{displaymath}
    which can be resolved if and only if \(\lambda^2=0\). If we
    restrict ourselves with the only real (complex) root
    \(\lambda=0\), then the corresponding operators
    \(\ladder[p]{\pm}=-\tilde{B}+\tilde{Z}/2\) will not affect eigenvalues and
    thus are useless in the above context.  However the dual number
    roots \(\lambda =\pm\rmp\) lead to the operators \(\ladder[\rmp]{\pm}=\pm
    \rmp\tilde{A}-\tilde{B}+\tilde{Z}/2\). These operators are
    suitable to build an \(\algebra{sl}_2\)-modules with a
    one-dimensional chain of eigenvalues.
  \end{enumerate}
\end{example}
\begin{rem}
  \label{re:hyper-number-necessity}
  It is noteworthy that:
  \begin{itemize}
  \item the introduction of complex numbers is a necessity for the
    \emph{existence} of ladder operators in the elliptic
    case;
  \item in the parabolic case we need dual numbers to make
    ladder operators \emph{useful};
  \item in the hyperbolic case double numbers are required for
    neither existence nor usability of ladder operators, but
    do provide an enhancement. 
  \end{itemize}
\end{rem}
We summarise the above consideration with a focus on the Principle of
similarity and correspondence:
\begin{prop}
  \label{pr:ladder-sim-eq}
  Let a vector \(X\in\algebra{sl}_2\) generates the subgroup \(K\),
  \(N'\) or \(\Aprime\), that is \(X=Z\), \(B-Z/2\), or
  \(B\) respectively. Let \(\alli\) be the respective hypercomplex unit.   

  Then raising/lowering operators \(\ladder{\pm}\) satisfying to the
  commutation relation:
  \begin{displaymath}
    [X,\ladder{\pm}]=\pm\alli \ladder{\pm},\qquad [\ladder{-},\ladder{+}]=2\alli X.
  \end{displaymath}
  are:
  \begin{displaymath}
    \ladder{\pm}=\pm\alli \tilde{A} +\tilde{Y}.
  \end{displaymath}
  Here \(Y\in\algebra{sl}_2\) is a linear combination of  \(B\) and
  \(Z\) with the properties:
  \begin{itemize}
  \item \(Y=[A,X]\).
  \item \(X=[A,Y]\).
  \item Killings form \(K(X,Y)\)~\cite{Kirillov76}*{\S~6.2} vanishes.
  \end{itemize}
  Any of the above properties defines the vector \(Y\in\loglike{span}\{B,Z\}\)
  up to a real constant factor.
\end{prop}

It is natural to expect that the usability of the Principle of similarity and
correspondence will not be limited to the considered examples
only. For example, similar types of ladder operators appeared in
relation to the Heisenberg group as well~\cite{Kisil11a}. 

\section{Open Problems}
\label{sec:open-problems}

We start from an illustration, that the invariant linear algebra
presented in Section~\ref{sec:invar-line-algebra} is not unique.

\begin{example}[Tropical form]
\label{ex:tropical}
Let us introduce the lexicographic order on \(\Space{R}{2}\):
\begin{displaymath}
  (u,v)\prec(u',v') \quad \text{if and only if} \quad
  \left\{\begin{array}{ll}
      \text{either}& u<u'; \\
      \text{or}& u=u',\  v<v'.
  \end{array}\right.
\end{displaymath}
One can define functions \(\min\) and \(\max\) of a pair of points on
\(\Space{R}{2}\) respectively. Then an addition of two vectors can be
defined either as their minimum or maximum. A similar definition is
used in \emph{tropical mathematics}, also known as Maslov
dequantisation or \(\Space[\min]{R}{}\) and \(\Space[\max]{R}{}\)
algebras, see~\cite{Litvinov05} for an energetic survey and
\cite{Butkovic10a} for a comprehensive coverage. It is easy to check
that such an addition is distributive with respect to scalar
multiplications~\eqref{eq:scalar-prod}---\eqref{eq:scalar-prod-prime}
and consequently is invariant under parabolic rotations. This approach
looks promising and definitely deserves a further careful investigation.
\end{example}

There are many other interesting questions to be investigated for the 
induced representations built in this paper. Here is a sketchy list of
some of them:
\begin{itemize}
\item All complex-valued irreducible unitary representations of
  \(\SL\) split into the three series: discrete, principal and
  complementary~\cites{Lang85,HoweTan92}. It is useful to find
  relations of new hypercomplex representations to these series.
\item A generic irreducible
  \(\algebra{sl}_2\)-module~\cites{HoweTan92,Mazorchuk09a} may not be
  unitarisable for any complex-valued inner product. How many of those
  non-unitarisable moduli become unitary for some hypercomplex-valued
  product? 
\item Most of the above irreducible hypercomplex-unitary
  \(\SL\)-moduli can be realised as spaces of functions either on the
  real line or the upper half-plane~\citelist{\cite{Kisil97c}
    \cite{Kisil11c}*{\S~5}}.  The concept of covariant
  transform~\citelist{\cite{Kisil10c} \cite{Kisil11c}*{\S~4}} provides
  the universal tool for a construction of associated integral formulae and
  reproducing kernels.
\item The described spaces of functions are null-solutions to certain
  differential operators, which possess \(\SL\)-symmetry. In
  particular, the covariant transform with a mother wavelet, which is
  annihilated by a ladder operator (or any its power), creates the
  space of null-solutions to an associated differential equation,
  see~\cite{Kisil11c}*{\S~5.3}. This shall be useful in analysis of
  solutions of partial differential
  equations~\cites{Konovenko09a,Pilipchuk10a}.
\item Analytical spaces of complex-valued functions work as models
  for functional calculus of normal operators, e.g Dunford--Riesz
  calculus. The newly build \(\SL\)-moduli support covariant
  functional calculi and associated spectra~\citelist{\cite{Kisil02a}
    \cite{Kisil95i} \cite{Kisil11c}*{\S~6}} for non-selfadjoint
  operators.
\end{itemize}

These directions can be viewed as parts of the Erlangen programme at
large~\cites{Kisil06a,Kisil11c}. The Principle of similarity and
correspondence may help to find most harmonious constructions.

{\small
\bibliography{abbrevmr,akisil,analyse,arare,aclifford,aphysics,algebra,ageometry}
}

\newpage
\renewcommand{\rmdefault}{cmr}
\selectlanguage{russian}
\renewcommand{\encodingdefault}{T2A}
\fontencoding{T2A}\selectfont
\gdef\proofname{Доказательство}
\setcounter{section}{0}
\setcounter{equation}{0}
\setcounter{page}{1}
\thispagestyle{plain}
\vspace*{30mm}
  \selectlanguage{russian}
\markboth{В. В. Кисиль}{Индуцированные представления и гиперкомплексные
числа}

\begin{center}
{\bfseries
 {\Large Индуцированных представлений группы $\SL$\\ и гиперкомплексные числа}\\[5mm]
 {\large
   \href{http://www.maths.leeds.ac.uk/~kisilv/}{В. В. Кисиль}}\\[3mm]
}
\end{center}

\date{\today}

\begin{abstract}
  В работе рассматривается конструкция индуцированных представлений
  для группы \(G=\SL\). Оказывается, что действие этой группы на
  однородном пространстве \(G/H\), где \(H\) является произвольной
  однопараметрической подгруппой \(\SL\), является дробно-линейным
  преобразованием двумерной алгебры гиперкомплексных чисел. Это
  наблюдение может быть распространено на дальнейшие соответствия
  между структурными компонентами \(\SL\) и гиперкомплексными
  системами.  Соответственно мы рассматриваем вопрос о
  гиперкомплексных характерах подгруппы \(H\). В частности мы
  приводим примеры индуцированных представлений группы \(\SL\) в
  пространствах функций с гиперкомплексными значениями, которые
  являются унитарными в определённом смысле.\\[2pt]
 \emph{Ключевые слова:} индуцированные представления, унитарные
   представления, полупростые группы Ли, комплексные числа, дуальные
   числа, двойные числа, преобразования Мёбиуса, дробно-линейные
   преобразования, операторы создания/уничтожения.
\end{abstract}
\vspace{1cm}
\par
\medskip

\section{Введение}
\label{sec:intro-ru}

Теория аналитических функций комплексного переменного---красивая
теория с многочисленными применениями в различных областях от теории
чисел до электродинамики. Естественно, происходит поиск различных
обобщений и аналогов такой успешной теории. Наиболее простая (или
основная?) попытка заключается в замене мнимой единицы \(\rmi^2=-1\)
или на гиперболическую \(\rmh^2=+1\) или нильпотентную \(\rmp^2=0\).

Двухмерная коммутативная ассоциативная алгебра порожденная \(1\) и
\(\rmh\) состоит из выражений вида \(x+\rmh y\) с действительными
параметрами \(x\) и \(y\). Как правило в русскоязычной литературе, они
называются \emph{двойными числами}~\cite{Yaglom79}*{Прил.~Б}, хотя
из-за влияния англоязычной литературы встречаются и названия вроде
\href{http://en.wikipedia.org/wiki/Split-complex_number}{паракомплексные
  числа}, расщепляемые комплексные числа, комплексные числа
гиперболического типа~\cites{BocCatoniCannataNichZamp07,%
  KhrennikovSegre07a,Ulrych05a,Pilipchuk10a}. Эта алгебра содержит
делители нуля \(\rmh_\pm=\frac{1}{2}(1\pm\rmh)\) и выполняется
тождества \(\rmh_\pm^2=\rmh_\pm\) и \(\rmh_+ \rmh_-=0\). Поэтому,
алгебраически двойные числа изоморфны  \(\Space{R}{}\oplus
\Space{R}{}\)---прямой сумме двух наборов действительных чисел
натянутых на \(\rmh_+\) и \(\rmh_-\). Отсюда и происходит название
<<двойные>>  или <<расщеплённые>>. 

Сходная алгебра порождённая \(1\) и \(\rmp\) состоит из элементов
\(x+\rmp y\), которые называются \wiki{Dual_number}{\emph{дуальными
    числами}}~\cites{Yaglom79,HerranzSantander02b,CatoniCannataNichelatti04}.
Среди них так же есть делители нуля вида \(\rmp y\),
\(y\in\Space{R}{}\). Физические приложения дуальных чисел обнаруживаются 
в классической механике~\cites{Yaglom79,Kisil12c}, нелинейной динамике~\cites{Pilipchuk11a,Pilipchuk10a},
теории относительности~\cites{BocCatoniCannataNichZamp07,Ulrych05a,Garasko09a},
космологии~\cites{HerranzSantander02b,GromovKuratov05a} и квантовой
механике~\cites{KhrennikovSegre07a,Kisil10a,Ulrych10a}.
\begin{remr}
  К сожалению, не существует устоявшихся обозначений для образующих
  гиперкомплесных алгебр. Более того, было бы затруднительно просто
  перечислить всё многообразие использованных символов и мы упомянем
  лишь некоторые из них. Стоит отметить, что даже наиболее
  традиционная мнимая единица \(\rmi\) обозначается как \(j\) во
  многих инжинерных текстах. Гиперболическая единица обозначается как
  \(\rmh\) во многочисленных работах начиная, по крайней мере, с
  основополагающего текста~\cite{VignauxDuranona35a}; однако, другое
  обозначение \(e\) употребляется в замечательной
  книге~\cite{Yaglom79}. Для нильпотентной единицы часто используется
  символ \(\iota\)~\cites{GromovKuratov06a,Pimenov65a}, но наш выбор
  \(\rmp\) позаимствован у Яглома~\cite{Yaglom79} и его
  последователей.  Нам это обозначение кажется удачным в силу
  следующего Замечания.
\end{remr}
\begin{remr}
  Нильпотентная единица \(\rmp\) близко связана с инфинитезимальным
  числом \(\varepsilon\) используемом в нестандартном
  анализе~\cites{Devis77,Uspenskii88}. Квадрат нилпотентно единицы в
  точности равен нулю в то время, как квадрат инфинитезимального числа
  есть <<практически>> ноль. Это сходство использовалось в
  работе~\cite{CatoniCannataNichelatti04} для вывода основных теорем
  дифференциального исчисления на основе нильпотентной единицы. Более
  того, существует изложение нестандартного анализа основанное на
  нильпотентной единице~\cite{Bell08a}\footnote{Я благодарен
    анонимному рецензенту за указание этой работы.}. То же самое
  свойство нильпотентности стоит за построением классической механики
  из представлений группы Гейзенберга~\cite{Kisil10a,Kisil12c}.
\end{remr}

Вполне естественен вопрос: \emph{какие <<аналитические>> функции связаны с
двойными и дуальными числами?} Не удивительно, что различные ответы на
этот вопрос уже длительное время предлагаются различными
исследователями. Так сложилось: что многие из них не были осведомлены
о работах предшественников. Не претендуя на полноту, мы укажем
известные автору работы, что бы повысить связанность в этой области. 

Для двойных чисел систематическое исследование было предпринято в
работе~\cite{VignauxDuranona35a}, среди последующих исследований и
обзоров укажем~\citelist{\cite{BocCatoniCannataNichZamp07}
  \cite{KhrennikovSegre07a} \cite{LavrentShabat77} \cite{MotterRosa98}
  \cite{Kravchenko09a}*{Part~IV} \cite{Pilipchuk10a} \cite{Sobczyk95a}
  \cite{Libine07a} \cite{Kisil97c}}, а так же дальнейшие ссылки в этих работах.
Доминирующий подход основан на факторизации волнового оператора
\(\partial_x^2-\partial_y^2=(\partial_x-\rmh\partial_y)(\partial_x+\rmh\partial_y)\)
в произведение двух линейных операторов. Это аналогично разложению
оператора Лапласа
\(\partial_x^2+\partial_y^2=(\partial_x-\rmi\partial_y)(\partial_x+\rmi\partial_y)\)
в произведение оператора Коши--Римана и его
сопряжённого. Соответственно, гиперболические аналитические функции
определяются как нулевые решения оператора
\(\partial_x+\rmh\partial_y\). Однако, разложение двойных чисел в базисе
\(\rmh_\pm\) представляет любую <<аналитическую>>  функцию 
\(f(x,y)=\rmh_+f_+(x+y)- \rmh_-f_-(x-y)\) в виде суммы двух
произвольных функций \(f_\pm\) действительной переменной. Это,
естественно, соответствует представлению общего решения волнового
уравнения в виде суммы двух волн приближающихся справа и слева.

Для нильпотентной единицы \(\rmp\) ситуация ещё более тривиальна.
Факторизация, сходная с предыдущими\footnote{Стоит, однако, отметить
  нетривиальную факторизацию волнового оператора
  \(\partial_x^2+\partial_y\) использующую высоко-размерные алгебры
  Клиффорда~\cite{CerejeirasKahlerSommen05a}.} вырождается
\((\partial_x-\rmp\partial_y)(\partial_x+\rmp\partial_y)=\partial_x^2\).
Попытка определить аналитические функции через разложения в степенной
ряд переменной \(x+\rmp y\) проводит к выражению \(f(x)+\rmp y
f'(x)\) порождённой произвольной дифференцируемой функцией \(f(x)\)
действительной переменной~\cite{CatoniCannataNichelatti04}.

В наших предыдущих работах~\cites{Kisil96c,Kisil97c,Kisil97a,Kisil02c}
мы предложили строить теорию аналитических функций в духе Эрлангенской
программы Феликса Клейна, то есть исходя из некоторого действия
группы. В данной работе мы продвигаемся в этом направлении. Заметим,
что классические пространства комплексно-аналитических
функций---пространства Харди и Бергмана---являются неприводимыми
модулями некоторых унитарных представлений группы \(\SL\). Далее, эти
представления индуцированы комплекснозначным характером компактной
подгруппы \(K\in\SL\). Поэтому, мы надеемся найти \emph{другие типы
  пространств аналитических функций среди всех неприводимых (или
  примарных) \(\SL\)-модулей для индуцированных представлений.}

Композиция работы такова. Мы напоминаем структуру группы \(\SL\) и
перечисляем все её неизоморфные однопараметрические непрерывные
подгруппы в разделе~\ref{r-sec:intro-sl2}. Этот результат позволяет описать три
неизоморфных типа \(\SL\)-однородных пространств в 
разделе~\ref{r-sec:induc-repr}. Как оказалось, три типа гиперкомплексных чисел наилучшим
образом описывают действие \(\SL\) на этих однородных посредством
дробно-линейных преобразований. Что бы индуцировать представления
\(\SL\) с одномерных подгрупп мы описываем их характеры в терминах
алгебр гиперкомплексных чисел в разделе~\ref{r-sec:introduction}. В следующем
разделе~\ref{r-sec:second-attempt} мы описываем объекты, которые могут трактоваться как
геометрические характеры однопараметрических подгрупп. Что бы
индуцированное представление было линейным, индуцированный характер
тоже должен быть линейным. Однако, линейность геометрических
характеров требует дополнительных определений, которые обсуждаются в
разделах~\ref{r-sec:rebu-algebr-struct}
и~\ref{r-sec:invar-line-algebra}. Далее, мы выписываем некоторые
примеры индуцированных представлений в
разделе~\ref{r-sec:concl-induc-repr}. Удобное описание неприводимого
\(\SL\)-модуля достигается в линейном базисе состоящем из собственных
векторов подгруппы индуцировавшей представление. Алгебра Ли
\(\algebra{sl}_2\) группы \(\SL\)  действует транзитивно на этом
базисе с помощью операторов повышения/понижения
(создания/уничтожения), которые широко известны в квантовой
механике. Мы описываем соответствующую структуру в
разделе~\ref{r-sec:correspondence}. В заключении мы даем некоторые
возможные направления дальнейших исследований.

\section{Группа $\SL$ и её подгруппы}
\label{r-sec:intro-sl2}

Группа \(\SL\) \wiki{SL\%282\%2CR\%29}{состоит из квадратных \(2\times
  2\) матриц} с действительными элементами и единичным
определителем. Групповая операция задаётся матричным умножением.
\(\SL\) является простейшей полупростой группой Ли. 
 Произвольный элемент \(\SL\) допускает единственное разложение
 вида~\cite{HoweTan92}*{Exer.~I.14}: 
\begin{equation}
  \label{r-eq:iwasawa-decomp}
  \begin{pmatrix}
    a&b \\c &d
  \end{pmatrix}= {\begin{pmatrix} \alpha & 0\\0&\alpha^{-1}
    \end{pmatrix}} {\begin{pmatrix} 1&\nu \\0&1
    \end{pmatrix}} {\begin{pmatrix}
      \cos\phi &  \sin\phi\\
      -\sin\phi & \cos\phi
    \end{pmatrix}},
\end{equation}
с некоторыми параметрами \(\alpha\in(0,\infty)\), \(\nu\in(-\infty,\infty)\)
и \(\phi\in(-\pi,\pi]\). Эта формула~\eqref{r-eq:iwasawa-decomp}
записанная в виде \(\SL=ANK\) известна как разложение
Ивасавы~\cite{Lang85}*{\S~III.1} и допускает обобщение на произвольную
полупростую группу Ли.

Семейства матриц каждого из трёх видов в правой
части~\eqref{r-eq:iwasawa-decomp} образуют однопараметрические подгруппы
обычно обозначаемые   \(A\),
\(N\) and \(K\).  Они получают экспоненциированием соответствующих
матриц с нулевым следом:
\begin{eqnarray}
  \label{r-eq:a-subgroup}
  A&=&\left\{  
    \begin{pmatrix} e^t & 0\\0&e^{-t}
    \end{pmatrix}=\exp \begin{pmatrix} t & 0\\0&-t
    \end{pmatrix},\  t\in\Space{R}{}\right\},\\
  \label{r-eq:n-subgroup}
  N&=&\left\{   {\begin{pmatrix} 1&t \\0&1
      \end{pmatrix}=\exp \begin{pmatrix} 0 & t\\0&0
      \end{pmatrix},}\  t\in\Space{R}{}\right\},\\
  \label{r-eq:k-subgroup}
  K&=&\left\{ {\begin{pmatrix}
        \cos t &  \sin t\\
        -\sin t & \cos t
      \end{pmatrix}=   \exp \begin{pmatrix} 0& t\\-t&0
      \end{pmatrix}},\ t\in(-\pi,\pi]\right\}.
\end{eqnarray}
Следующий простой результат имеет поучительное доказательство:
\begin{propr}
  \label{r-pr:three}
  Любая непрерывная однопараметрическая подгруппа \(\SL\) матрично
  сопряжена к одной из подгрупп \(A\), \(N\) или \(K\).
\end{propr}
\begin{proof}
  Любая непрерывная однопараметрическая подгруппа получается
  экспоненциированием 
  \begin{equation}
    \label{r-eq:exp-taylor}
    e^{tX}=\sum_{n=0}^\infty \frac{t^n}{n!}X^n
  \end{equation}
  элемента \(X\) алгебры Ли \(\algebra{sl}_2\) группы
  \(\SL\). Такой элемент представим  \(2\times 2\) матрицей с нулевым
  следом. Поведение ряда Тэйлора~\eqref{r-eq:exp-taylor} зависит от
  свойств степеней  \(X^n\). Последние легко классифицируются
  непосредственным вычислением:
  \begin{lemr}
    \label{r-le:square}
    Квадрат \(X^2\) матрицы \(X=
    \begin{pmatrix}
      a&b\\c&-a
    \end{pmatrix}\) с нулевым следом есть единичная матрица умноженная
    на \(a^2+bc=-\det X\). Этот множитель может быть отрицательным,
    нулевым или положительным, что соответствует трём различным типам
    ряда Тэйлора~\eqref{r-eq:exp-taylor} для \(e^{tX}\).
  \end{lemr}
  Несложно видеть, что матричное подобие приводит элемент \(X\) к
  генератору 
  \begin{itemize}
  \item подгруппы \(K\) если \((-\det X) <0\);
  \item подгруппы \(N\) если \((-\det X) =0\);
  \item подгруппы \(A\) если \((-\det X) >0\).
  \end{itemize}
  Так как матричное подобие не меняет определителя матрицы, то эти три
  случае не сводятся друг к другу.
\end{proof}
\begin{exampler}
  \label{r-ex:conjugate}
  Следующие две подгруппы сопряжены к \(A\) и \(N\) соответственно
  посредством матричного сопряжения с фиксированной матрицой \(
  \begin{pmatrix}
    0&1\\1&0
  \end{pmatrix}\):
  \begin{eqnarray}
    \label{r-eq:ap-subgroup}
    A'&=&\left\{  
      \begin{pmatrix} \cosh t & \sinh t\\ \sinh t& \cosh t
      \end{pmatrix}=\exp \begin{pmatrix} 0 & t\\t&0
      \end{pmatrix},\  t\in\Space{R}{}\right\},\\
    \label{r-eq:np-subgroup}
    N'&=&\left\{   {\begin{pmatrix} 1&0\\t&1
        \end{pmatrix}=\exp \begin{pmatrix} 0 & 0\\t&0
        \end{pmatrix},}\  t\in\Space{R}{}\right\}.
  \end{eqnarray}
  В дальнейшем мы будем, как правило, рассматривать подгруппы \(N'\) и
  \(A'\) вместо \(N\) и \(A\), так как это приведёт к более
  естественной геометрической картине.
\end{exampler}

\section{Действие  $\SL$ приводит к гиперкомплексным числам}
\label{r-sec:induc-repr}

Пусть \(H\) является подгруппой группы \(G\).  Для соответствующего
однородного пространства \( \Omega=G / H\) определим гладкое сечение
\(s: \Omega \rightarrow G\)~\cite{Kirillov76}*{\S~13.2}, которое
является левым обратным для естественной проекции \(p: G\rightarrow
\Omega \). Выбор такого сечения непринципиален в том смысле, что наши
дальнейшие построения определяются с точностью до гладкого отображения
\(\Omega\rightarrow \Omega\).

Каждый элемент \(g\in G\) может быть представлен единственным образом
в виде \(g=s(\omega)h\), где \(\omega=p(g)\in \Omega\) и \(h\in H\).
Соответственно, на \(\Omega \) определяется действие группы \(G\)
следующим образом:
\begin{equation}
  \label{r-eq:g-action}
  g: \omega  \mapsto g\cdot \omega=p(g* s(\omega)),
\end{equation}
где \(*\) обозначает групповое произведение. Это действие может быть
также проиллюстрировано коммутативной диаграммой:
\begin{displaymath}
        \xymatrix{ 
        G \ar@<.5ex>[d]^p \ar[r]^{g*} & G \ar@<.5ex>[d]^p \\ 
        \Omega   \ar@<.5ex>[u]^s \ar[r]^{g\cdot}   & \Omega     \ar@<.5ex>[u]^s }
\end{displaymath}

Для группы \(G=\SL\), как и для остальных полупростых групп, принято
рассматривать случай, когда \(H\) является максимальной компактной
подгруппой \(K\).  Однако в данной работе мы рассмотрим все три
возможности для однопараметрической подгруппы \(H\) описанных в
Предложении~\ref{r-pr:three} и Примере~\ref{r-ex:conjugate}, а именно
\(H=K\), \(N'\) или \(A'\). В этом случае многообразие \(\Omega\)
является двумерным и для всех трёх возможных подгрупп \(H\) мы
определим~\cite{Kisil97c}*{Ex.~3.7(a)}:
\begin{equation}
  \label{r-eq:s-map}
  s: (u,v) \mapsto
  \frac{1}{\sqrt{v}}
  \begin{pmatrix}
    v & u \\ 0 & 1
  \end{pmatrix}, \qquad (u,v)\in\Space{R}{2},\  v>0.
\end{equation}
Непосредственное вычисление или использование компьютерной
программы~\cite{Kisil07a} приводит к следующему результату:
\begin{propr}
  \label{r-pr:sl2-act-brute}
  Действие~\eqref{r-eq:g-action} группы \(\SL\) определённое сечением
  \(s\)~\eqref{r-eq:s-map} есть:
  \begin{equation}
    \label{r-eq:sl2-act-brute}
    \begin{pmatrix}
      a&b\\c&d
    \end{pmatrix}: (u,v)\mapsto
    \left(\frac{(au+b)(c u+d) -\sigma cav^2}{( c u+d)^2 -\sigma (cv)^2},
      \frac{v}{( c u+d)^2 -\sigma (cv)^2}\right),
  \end{equation}
  где \(\sigma=-1\), \(0\) или \(1\) для подгрупп \(K\), \(N'\)
  или \(A'\) соответственно.
\end{propr}
Выражение~\eqref{r-eq:sl2-act-brute} не выглядит привлекательно, однако,
введение гиперкоплексных чисел существенно его улучшает:
\begin{propr}
  \label{r-pr:sl2-act}
  Пусть единица \(\alli\) такова, что \(\alli^2=\sigma\), тогда
  действие группы \(\SL\)~\eqref{r-eq:sl2-act-brute} может быть записано
  как дробно-линейное (мёбиусово) преобразование:
  \begin{equation}
    \label{r-eq:sl2-act}
    \begin{pmatrix}
      a&b\\c&d
    \end{pmatrix}: w\mapsto \frac{aw+b}{cw+d}, \qquad
    \text{где } w=u+\alli v,
  \end{equation}
  для всех трёх возможных значений параметра \(\sigma\) из
  Предложения~\ref{r-pr:sl2-act-brute}. 
\end{propr}
\begin{remr}
  Отметим что гиперкомплексные единицы не вводились нами умышленно, по
  произволу или из <<чисто обобщательской
  попытки>>~\cite{Pontryagin86a}*{стр.~4}. Они естественно возникли из
  действия группы \(\SL\).
\end{remr}
Мнимая единица \(\rmi^2=-1\) порождает хорошо знакомые комплексные
числа. Двумерная алгебра выражений \(x+\rmh y\)
с образующей \(\rmh^2=1\) имеет много имён:
гиперболические, 
\href{r-http://en.wikipedia.org/wiki/Split-complex_number}{\emph{двойные 
    числа}}~\cites{BocCatoniCannataNichZamp07,KhrennikovSegre07a,Ulrych05a} и др.,
это один из простейших примеров
\href{r-http://en.wikipedia.org/wiki/Hypernumber}{гиперкомплексных
  чисел}. Параболические гиперкомплексные числа, часто так же называемые
\wiki{Dual_number}{\emph{двойными}} имеют вид \(x+\rmp y\)
с нильпотентной  единицей такой, что
\(\rmp^2=0\)~\cites{Yaglom79,HerranzSantander02b,CatoniCannataNichelatti04,Gromov90a}.
Мы используем символ \(\alli\) для обозначения любой из трёх 
единиц \(\rmi\), \(\rmp\) или \(\rmh\).

Действие~\eqref{r-eq:sl2-act} является групповым гомоморфизмом из
\(\SL\) в преобразование <<верхней полуплоскости>> гиперкомплексной
алгебры. Несмотря на то, что алгебраическая структура двойных и
дуальных чисел достаточно вырождена, они интересны как
однородные пространства для мёбиусовых преобразований. Построение
соответствующих геометрий в духе 
\wiki{Erlangen_program}{Эрлангенской программы} Ф.~Клейна ободряюще
необычно~\cite{Kisil05a} и приводит к новым результатам даже в хорошо
изученном случае геометрии Лобачевсеого~\cite{Kisil06a}.  

Из-за  делителей нуля корректное изучение мёбиусовых
преобразований~\eqref{r-eq:sl2-act} должно  рассматриваться  на
конформно пополненной плоскости~\cites{HerranzSantander02b,Kisil06b}.
Физические приложения гиперкомплексных чисел простираются от
классической механики~\cite{Yaglom79} и
специальной теории
относительности~\cites{BocCatoniCannataNichZamp07,Ulrych05a} до
космологии~\cites{HerranzSantander02b,GromovKuratov05a} и квантовой
механики~\cites{KhrennikovSegre07a,Kisil10a,GromovKuratov06a}. 

Стандартный способ линеаризации действия~\eqref{r-eq:g-action}
заключается в переходе к
\wiki{Induced_representations}{представлению индуцированному с
  подгруппы} \(H\)~\citelist{\cite{Kirillov76}*{\S~13.2}
  \cite{Kisil97c}*{\S~3.1}}.  Для этого мы определим отображение \(r: G\rightarrow H\)
порожденное отображениями \(p\) и \(s\) через соотношения:
\begin{equation}
  \label{r-eq:r-map}
  r(g)={(s(\omega))}^{-1}g, \qquad \text{где }
  \omega=p(g)\in \Omega .
\end{equation}
Пусть \(\chi\) является неприводимым представлением подгруппы \(H\)
в векторном пространстве \(V\), тогда оно  индуцирует
представление \(G\) в смысле Макки~\cite{Kirillov76}*{\S~13.2}. 
Существует несколько реализаций индуцированного представления, мы
будем рассматривать представление \(\uir{\chi}\) в пространстве
 \(V\)-значных функций по формуле~\cite{Kirillov76}*{\S~13.2.(7)--(9)}:
\begin{equation} 
  \label{r-eq:def-ind}
  [\uir{\chi}(g) f](\omega)= \chi(r(g^{-1} * s(\omega)))  f(g^{-1}\cdot \omega),
.
\end{equation}
где \(g\in G\), \(\omega\in\Omega\), \(h\in H\) и \(r: G
\rightarrow H\), \(s: \Omega \rightarrow G\) есть отображения
определённые выше; \(*\)~обозначает групповую операцию на \(G\) и
\(\cdot\) является действием~\eqref{r-eq:g-action} группы \(G\) на
\(\Omega\). 

Так как в нашем рассмотрении подгруппа \(H\) всегда одномерна, то её
неприводимые представления обычно предполагаются комплексными
характерами. Однако, если гиперкомплексные числа естественно
появляются в действии группы \(\SL\) на однородных
пространствах~\eqref{r-eq:sl2-act}, почему мы не должны также
рассмотреть гиперкомплексные характеры?

\section{Гиперкомплексные характеры---алгебраический подход}
\label{r-sec:introduction}

\begin{figure}[htbp]
  \centering
  \includegraphics{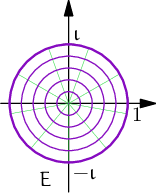}\qquad
  \includegraphics{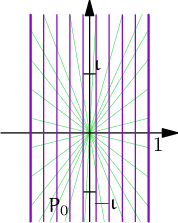}\qquad
  \includegraphics{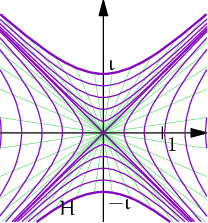}
  \caption[Унитарные вращения---алгебраическая версия]{
    Унитарные вращения, полученные алгебраическим способом,
    т.е. умножением на \(e^{\alli t}\): эллиптические
    (\(E\)), параболические (\(P_0\)) и гиперболические (\(H\)). Все
    фиолетовые орбиты определяются уравнением
    \(x^2-\alli^2y^2=r^2\). Зелёные  <<спицы>> (прямые линии из начала
    координат к точке на орбите) получаются вращением горизонтальной
    оси.} 
  \label{r-fig:rotations}
\end{figure}

Как уже отмечалось, традиционно основное внимание при изучении
представлений группы \(\SL\) уделяется случаю \(H=K\) и,
соответственно, комплекснозначным характерам \(K\). Линейное
преобразование определённое матрицами~\eqref{r-eq:k-subgroup} из \(K\)
представляют собой вращения \(\Space{R}{2}\) на угол \(t\).
Идентификация \(\Space{R}{2}=\Space{C}{}\) переводит это вращение в
умножение на \(e^{\rmi t}\), где \(\rmi^2=-1\).  Вращение является
унитарным для эллиптической метрики:
\begin{equation}
  \label{r-eq:ell-metric}
  x^2+y^2=(x+\rmi y)(x-\rmi y).
\end{equation}
Каждая орбита вращения являются окружностью (<<ободом колеса>>), и
любая прямая проходящая через начало координат (<<спица>>) переходит в
другую прямую повёрнутую на угол  \(t\), см.
Рис.~\ref{r-fig:rotations}(E).

Для гиперкомплексных чисел наиболее прямолинейное применение этой
модели для характеров выглядит так. 

\begin{propr}
  \label{r-pr:algebraic-char}
  Нижеследующая таблица показывает алгебраическое соответствие между
  гиперкомплексными характерами:
  \begin{center}
    \begin{tabular}{||c|c|c||}
      \hline\hline
      Эллиптический & Параболический & Гиперболический\\
      \hline\hline
      Subgroup \(K\)&Subgroup \(N'\)&Subgroup \(A'\)\\
      \hline
      \(\strut^{\strut}\rmi^{2}=-1\)&    \(\rmp^2=0\)&    \(\rmh^2=1\)
      \\
      \hline
      \(w=x+\rmi y\) &\(w=x+\rmp y\) &\(w=x+\rmh y\) 
      \\
      \hline
      \(\bar{w}=x-\rmi y\) &\(\bar{w}=x-\rmp y\) &\(\bar{w}=x-\rmh y\) 
      \\
      \hline
      \(\strut^{\strut}e^{\rmi t} = \cos t +\rmi \sin t\)&\(e^{\rmp t} = 1 +\rmp t\)&\(e^{\rmh t} = \cosh t +\rmh \sinh t\)
      \\
      \hline
      \(\strut^{\strut}\modulus[e]{w}^{ 2}=w\bar{w}=x^2+y^2\) &\(\modulus[p]{w}^2=w\bar{w}=x^2\) &\(\modulus[h]{w}^2=w\bar{w}=x^2-y^2\) 
      \\
      \hline
      \(\arg w = \arctan \frac{y}{x} \frac{\strut}{\strut}\)& \(\arg w = \frac{y}{x}\)&\(\arg w = \tanh^{-1} \frac{y}{x}\)
      \\
      \hline
       окружность \(\strut^{\strut} \modulus[e]{w}^{2}=1\) &
       две прямые \(x=\pm 1\) &  гипербола \(\modulus[h]{w}^2=1\)
      \\
      \hline\hline
    \end{tabular}
  \end{center}
  Геометрическое действие умножения на \(e^{\alli t}\) показано
  Рис.~\ref{r-fig:rotations} для всех трёх случаев. 
\end{propr}
Явным образом параболическое <<вращение>>
порождённое умножением на \(\rme^{\rmp t}\) действует на дуальные
числа так:
\begin{equation}
  \label{r-eq:parab-rot-triv}
  \rme^{\rmp x}: a+\rmp b \mapsto a+\rmp (a x+b).
\end{equation}
Это связывает параболический случай с галилеевыми симметриями
классической механики~\cite{Yaglom79}, в которой время абсолютно и
независимо от пространства. Такое кинематическое обоснование
предлагает следующие соответствия для параболического
случая~\cites{HerranzOrtegaSantander99a,Yaglom79}:
\begin{itemize}
\item Параболическая тригонометрия тривиальна:
  \begin{equation}
    \label{r-eq:par-trig-0}
    \cosp t =\pm 1, \qquad \sinp t=t;
  \end{equation}
\item параболическое расстояние (при ненулевом \(x\)) не зависит от \(y\):
  \begin{equation}
    \label{r-eq:par-metr-0}
    x^2=(x+\rmp y)(x-\rmp y);
  \end{equation}
\item полярное разложение дуального числа~\cite{Yaglom79}*{App.~C(30')}:
  \begin{equation}
    \label{r-eq:p-polar-yaglom}
    u+\rmp v = u(1+\rmp \frac{v}{u}), \quad \text{ поэтому }
    \quad \modulus{u+\rmp v}=u, \quad \arg(u+\rmp v)=\frac{v}{u};
  \end{equation}
\item параболическое <<колесо>> выглядит прямоугольным, см.
  Рис.~\ref{r-fig:rotations}(\(P_0\)). 
\end{itemize}

Алгебраическая аналогия и кинематическое обоснование кажутся достаточно
убедительными для единственно правильной параболической
тригонометрии~\cites{LavrentShabat77,HerranzOrtegaSantander99a,Yaglom79}.
Более того, данные три типа вращений представляют собой все возможные
линейные симплектические преобразования, что делает их и
соответствующие мнимые единицы полезными в классической и квантовой
механике~\cites{Kisil10a,Kisil11a}. 

Однако геометрический подход и симметрии параболических уравнений
математической физики предлагают менее вырожденную альтернативу
изложенную  ниже. 

\section{Параболическое колесо из геометрического действия}
\label{r-sec:second-attempt}

Попробуем ввести нетривиальные параболические вращения другим
способом. Если умножение в дуальных числах (линейные
преобразования) слишком тривиальны, мы перейдём к следующему уровню
сложности: дробно-линейным. 
linear-fractional.

Действительно, нет необходимости рассматривать гиперкомплексных
единицы как абстрактные сущности. Мы видели в Лемме~\ref{r-le:square},
что генераторы подгрупп \(K\), \(N\) и \(A\) по сути представляют
единицы для комплексных, дуальных и двойных чисел
соответственно.  Мы так же можем рассматривать выражения для
однопараметрических подгрупп \(K\), \(N'\) и \(A'\) как матричные
версии тождества Эйлера из пятой строки таблицы в
Предложении~\ref{r-pr:algebraic-char}.

Поэтому мы можем определить характеры подгрупп \(K\), \(N'\) и
\(A'\) через геометрическое действие \(\SL\) посредством
преобразований Мёбиуса. Однако действие~\eqref{r-eq:sl2-act} определено
в верхней полуплоскости, а было бы желательно получить его унитарную
форму, что возможно при конформном отображении на единичный круг. В
эллиптическом случае это достигается преобразованием Келли, которое
переводит подгруппу \(K\) в:
\begin{equation}
  \label{r-eq:ell-cayley}
  \frac{1}{2}
  \begin{pmatrix}
    1 & -\rmi \\ -\rmi &1
  \end{pmatrix}
  \begin{pmatrix}
    \cos t & -\sin t \\ \sin t & \cos t
  \end{pmatrix}
  \begin{pmatrix}
    1 & \rmi \\ \rmi &1
  \end{pmatrix} =
  \begin{pmatrix}
    \rme^{\rmi t} & 0 \\ 0 &  \rme^{-\rmi t}
  \end{pmatrix}.
\end{equation}
Диагональная матрица справа задаёт преобразование Мёбиуса, которое
сводится к умножению на\(e^{2\rmi t}\), то есть унитарному вращению
комплексной плоскости.

Гиперболическая версия преобразования Келли выглядит
так~\cite{Kisil05a}*{\S~\ref{E-sec:ellipt-hyperb-caley}}:
\begin{equation}
  \label{r-eq:hyp-cayley}
  \frac{1}{2}
  \begin{pmatrix}
    1 & \rmh \\ -\rmh &1
  \end{pmatrix}
  \begin{pmatrix}
    \cosh t & \sinh t \\ \sinh t & \cosh t
  \end{pmatrix}
  \begin{pmatrix}
    1 & -\rmh \\ \rmh &1
  \end{pmatrix} =
  \begin{pmatrix}
    \rme^{\rmh t} & 0 \\ 0 &  \rme^{-\rmh t}
  \end{pmatrix}.
\end{equation}
Преобразование Мёбиуса с полученной матрицей так же сводится к
умножению на \(e^{2\rmh t}\), что является унитарным преобразованием
Лоренца двумерного пространства-времени Минковского.  

В параболическом случае мы используем для преобразования Келли сходную
матрицу~\cite{Kisil05a}*{\S~\ref{E-sec:parab-cayl-transf}}: 
\begin{displaymath}
  C_\rmp=
  \begin{pmatrix} 
    1 & -\rmp \\ -\rmp &1  
  \end{pmatrix}
\end{displaymath}
Тогда образом матрицы~\eqref{r-eq:n-subgroup} из подгруппы \(N\) будет: 
\begin{equation}
  \label{r-eq:par-cayley}
  \begin{pmatrix}
    1 & -\rmp\\ -\rmp &1
  \end{pmatrix}
  \begin{pmatrix}
    1 & t \\ 0 & 1
  \end{pmatrix}
  \begin{pmatrix}
    1 & \rmp \\ \rmp &1
  \end{pmatrix} =
  \begin{pmatrix}
    1+\rmp t & t \\ 0 & 1-\rmp t
  \end{pmatrix}=
    \begin{pmatrix}
    \rme^{\rmp t} & t \\ 0 & \rme^{-\rmp t}
  \end{pmatrix}.
\end{equation}
Хотя этот результат не сильно отличается от двух предыдущих
случаев~\eqref{r-eq:ell-cayley} и~\eqref{r-eq:hyp-cayley}, но присутствие
ненулевого члена в позиции \((1,2)\) нарушает баланс. Тем не менее, мы
по прежнему определим параболическое вращение, как преобразование
Мёбиуса с матрицей~\eqref{r-eq:par-cayley}, которое уже не будет простым
умножением на скаляр.


\begin{figure}[htbp]
  \centering
  \includegraphics{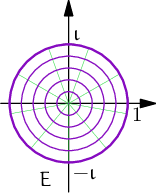}\hfill
  \includegraphics{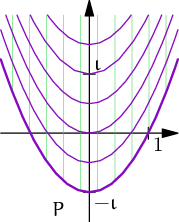}\hfill
  \includegraphics{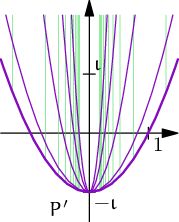}\hfill
  \includegraphics{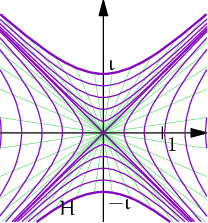}
  \caption[Унитарные вращение---геометрическая версия]{Геометрическая
    версия унитарных вращений:
    эллиптические (\(E\)), два варианта параболических (\(P\) 
    и \(P'\)), и гиперболические (\(H\)). Фиолетовые орбиты являются
    линиями уровня для соответствующих функций модуля. Зелёные прямые
    линии соединяют точки с одинаковым значением аргумента, они
    нарисованы с постоянным шагом приращения аргумента для каждого
    случая.} 
  \label{r-fig:p-rotations}
\end{figure}



\begin{exampler}
  \label{r-ex:upper-triang}
  Параболическое вращение с верхне-треугольной матрицей из подгруппы \(N\) есть:
  \begin{equation}
    \label{r-eq:par-moebius}
    \begin{pmatrix}
      \rme^{\rmp t}&t\\0&\rme^{-\rmp t}
    \end{pmatrix}: -\rmp \mapsto t +\rmp(t^2-1).
  \end{equation}
  Это совпадает с  \emph{циклическим вращением} определённым
  в~\cite{Yaglom79}*{\S~8}.  Сравнение с тождеством Эйлера по прежнему
  подтверждает определение параболического синуса как линейной функции
  \(\sinp t=t\), но предлагает новое значение для параболического
  косинуса (квадратичное вместо постоянного):
  \begin{displaymath}
    \cosp t = 1-t^2, \qquad \sinp t= t.
  \end{displaymath}
  Соответственно параболическая теорема Пифагора выглядит так:
  \begin{equation}
    \label{r-eq:pyhagoras-p}
    \sinp^2 t +  \cosp t =1,
  \end{equation}
  что неплохо соседствует с эллиптической с одной стороны и
  гиперболической---с другой версиями:
  \begin{displaymath}
    \sin^2 t+\cos^2 t  =1, \qquad \sinh^2 t -  \cosh^2 t =-1.
  \end{displaymath}
  Отметим также, что тождество~\eqref{r-eq:pyhagoras-p} менее вырождено,
  чем выражение вообще без синуса \( \cosp^2 t =1\)
  из~\cite{HerranzOrtegaSantander99a} (см. так
  же~\eqref{r-eq:par-trig-0}, \eqref{r-eq:par-metr-0}).
\end{exampler}


\begin{exampler}
  \label{r-ex:lower-triang}
  Существует вторая возможность определить параболические вращения
  используя нижне-треугольные матрицы из подгруппы
  \(N'\). Существенное отличие в том, что точка \(-\rmp\) теперь
  неподвижна (как и все точки вертикальной оси) и для получения орбиты
  мы будем использовать \(\rmp^{-1}\), который является идеальным
  элементом (точкой в бесконечности)~\cite{Yaglom79}*{App.~C}) потому,
  что \(\rmp\)---делитель нуля. Вопросы замыкания гиперкомплексных
  чисел идеальными элементами в физическом контексте рассматривались
  в~\cites{Kisil06b}. 
  
   Для подгруппы \(N'\) мы получаем:
  \begin{equation}
    \label{r-eq:par-moebius-prime}
    \begin{pmatrix}
      \rme^{-\rmp t}&0\\t&\rme^{\rmp t}
    \end{pmatrix}:\  \frac{1}{\rmp} \mapsto \frac{1}{t}+ \rmp \left(1-\frac{1}{t^2}\right). 
  \end{equation}
  
  Сравнение с~\eqref{r-eq:par-moebius} показывает, что обе версии
  связаны отображением \(t\mapsto t^{-1}\). Такое же соотношение
  связывает и параболические тригонометрические
  функции. Параболическая <<единичная окружность>>  (или
  \emph{цикл}~\cites{Kisil05a,Yaglom79}) определяется уравнением
  \(x^2-y=1\) в обоих случаях, см Рис.~\ref{r-fig:p-rotations}(\(P\)) и
  (\(P'\)). Однако остальные орбиты различаются и мы опишем их в
  следующей секции.
\end{exampler}
Рис.~\ref{r-fig:p-rotations} показывает преобразования Мёбиуса
определённых матрицами~\eqref{r-eq:ell-cayley}, \eqref{r-eq:hyp-cayley}
и~\eqref{r-eq:par-cayley} с соответствующими <<единичными кругами>>,
которые являются образами верхней полуплоскости под действием
соответствующих преобразований
Келли~\cite{Kisil05a}*{\S~\ref{E-sec:unit-circles}}.

\section{Восстановление алгебраической структуры из геометрии}
\label{r-sec:rebu-algebr-struct}

Мы хотим, что бы индуцированные представления были линейными, для
этого порождающий их характер должен также быть линейный. 
Вращения в эллиптическом и параболическом случае задаются
произведениями с комплексными и двойными числами соответственно и
потому---линейны.
Однако, нетривиальные параболические вращения~\eqref{r-eq:par-moebius}
и~\eqref{r-eq:par-moebius-prime} (Рис.~\ref{r-fig:p-rotations}(\(P\))
и~(\(P'\)))---нелинейны и мы введём новые алгебраические операции на
дуальных числах, которые линеаризуют преобразования Мёбиуса.  Для
этого мы будем систематически использовать хорошо известный приём
<<обращения теорем в определения>>.

\subsection{Модуль и аргумент}
\label{r-sec:modulus-argument}

В эллиптическом и гиперболическом случае орбиты вращений являются
точками с постоянной нормой (модулем): \(x^2+y^2\) или
\(x^2-y^2\) соответственно. В параболическом случае мы используем этот
же признак: 
\begin{defnr}
  \label{r-de:norm}
  Орбиты преобразований~\eqref{r-eq:par-moebius}
  и~\eqref{r-eq:par-moebius-prime} являются линиями уровня для следующих
  функций, которые мы будем называть соответствующими \emph{модулем}
  (\emph{нормой}): 
  \begin{equation}
    \label{r-eq:parab-norm}
    \text{ для } N:\ \modulus{u+\rmp v}=u^2-v, \qquad \text{ для } N':\
    \modulus{u+\rmp v}'=\frac{u^2}{v+1}.
  \end{equation}
\end{defnr}
\begin{remr}
  \begin{enumerate}
  \item Выражение \( \modulus{(u,v)}=u^2-v\) равно параболическому
    расстоянию из \((0,\frac{1}{2})\) to \((u,v)\) к точке \((u,v)\),
    см.~\cite{Kisil05a}*{Lem.~\ref{E-le:n-orbits-concentric}},  и
    так же соответствует  <<параболической теореме
    Пифагора>>~\eqref{r-eq:pyhagoras-p}. 
  \item Модуль для подгруппы \(N'\) выражает параболическое фокальное
    расстояние из точки \((0,-1)\)
    к \((u,v)\), см.~\cite{Kisil05a}*{Lem.~\ref{E-le:np-orbits-p-confocal}}.
  \end{enumerate}
\end{remr}

Единственные прямые, которые переходят в прямые под действием
параболических вращений с подгруппами
\(N\) и \(N'\), это вертикальные прямые. Поэтому мы их будем
использовать как множества точек с постоянным значением аргумента и
дадим следующее:
\begin{defnr}
  \label{r-de:arg}
  Параболический \emph{аргумент} определятся как:
  \begin{equation}
    \label{r-eq:parab-arg}
    \text{ для } N:\ \arg(u+\rmp v)=u, \qquad \text{ для } N':\
    \arg'(u+\rmp v)=\frac{1}{u}.
  \end{equation}
\end{defnr}
Оба определения~\ref{r-de:norm} и~\ref{r-de:arg} здают свойства хорошо
согласующиеся с параболическими вращениями:
\begin{propr}
  \label{r-pr:rot-norm-arg}
  Пусть точка \(w_t\) получается из \(w\) параболическим вращением на
  угол \(t\) в~\eqref{r-eq:par-moebius} или 
  в~\eqref{r-eq:par-moebius-prime}. Тогда:
  \begin{displaymath}
    \modulus{w_t}^{(\prime)}=\modulus{w}^{(\prime)},\qquad 
    \arg^{(\prime)} w_t=\arg^{(\prime)} w+t,
  \end{displaymath}
  где выражения со штрихом используются для подгруппы \(N'\).
\end{propr}
Все доказательства в этой и следующей секциях производятся 
вычислением либо непосредственно, либо с помощью пакета символических
расчётов~\cites{Kisil09c,Kisil07a}.
\begin{remr} Отметим, что в традиционном
  изложении~\cite{Yaglom79}*{App.~C(30')} параболический модуль и
  аргумент выражаются как в~\eqref{r-eq:p-polar-yaglom}, что в некотором
  смысле обратно к нашим соглашениям.
\end{remr}

\subsection{Параболический поворот как умножение}
\label{r-sec:rotat-as-mult}

Для определение параболического произведения мы снова будем
использовать известные свойства комплексного произведения в качестве
определения. Действительно, параболическое вращение можно
рассматривать как умножение на унимодулярное скаляр, таким образом мы
приходим к следующему 
\begin{defnr}
  \label{r-de:product}
  Произведение двух векторов \(w_1\) и \(w_2\) определяется такими
  двумя условиями:  
  \begin{enumerate}
  \item \(\arg^{(\prime)}(w_1 w_2)=\arg^{(\prime)} w_1 + \arg^{(\prime)} w_2\);
  \item \(\modulus{w_1 w_2}^{(\prime)} =\modulus{w_1}^{(\prime)}\cdot
    \modulus{w_2}^{(\prime)}\). 
  \end{enumerate}
\end{defnr}
Параболическое сопряжение так же требует специального определения,
которое основано на изменение знака аргумента. 
\begin{defnr}
  Сопряжённый вектор задаётся равенством:
  \begin{equation}
    \label{r-eq:conjugation}
    \overline{u+\rmp v}=-u+\rmp v.
  \end{equation}
\end{defnr}
Очевидно, мы имеем следующие свойства сопряжения:
\(\modulus{\overline{w}}^{(\prime)}=\modulus {w}^{(\prime)}\) и
\(\arg^{(\prime)}{\overline{w}}=-\arg^{(\prime)} {w}\). 
Определения~\ref{r-de:norm}, \ref{r-de:arg}
и~\ref{r-de:product} совместно задают значение для произведения двух
векторов единственным образом. 
\begin{propr}
  Параболическое произведение двух векторов имеет следующие явные
  выражения: 
  \begin{align}
    \label{r-eq:parab-prod}
    \text{для } &N:& (u,v)*(u',v') & =
    (u+u',(u+u')^2-(u^2-v)(u'^2-v')); \\
    \label{r-eq:parab-prime-prod}
    \text{для } &N':&(u,v)*(u',v') & = \left(\frac{uu'}{u+u'},\frac{(v+1)(v'+1)}{(u+u')^2}-1\right).
  \end{align}
\end{propr}
Не смотря на необычный вид оба произведения имеют хорошо знакомые
свойства:
\begin{propr}
  \label{r-pr:mul}
  Оба произведения~\eqref{r-eq:parab-prod} и~\eqref{r-eq:parab-prime-prod}
  удовлетворяют таким условиям:
  \begin{enumerate}
  \item\label{r-item:prod-comm-ass} Они коммутативны и ассоциативны;
  \item\label{r-it:rot-as-mult} Соответствующие параболические
    повороты~\eqref{r-eq:par-moebius} и~\eqref{r-eq:par-moebius-prime}
    совпадают с умножением на унимодулярный вектор.
  \item\label{r-item:prod-inv} Значение произведения \(w_1\bar{w}_2\) не
    изменится если оба вектора \(w_1\) и \(w_2\) будут повёрнуты на
    тот же угол соответствующим вращением~\eqref{r-eq:par-moebius}
    или~\eqref{r-eq:par-moebius-prime}.
  \item\label{r-item:prod-norm-sq} 
    Для любого вектора \(w\) выполняется тождество:
    \begin{displaymath}
      \modulus{w\bar{w}}= \modulus{w}^2.
    \end{displaymath}
  \end{enumerate}
\end{propr}
В частности, свойство~(\ref{r-item:prod-inv}) из этого Предложения будет
важно для скалярного произведения ~\eqref{r-eq:inner-product}, по
отношению к которому наши индуцированные представления будут
унитарны. 

\section{Инвариантная линейная алгебра}
\label{r-sec:invar-line-algebra}

Теперь мы хотим определить линейную структуру на \(\Space{R}{2}\),
которая будет совместима с параболическим умножением из предыдущей
секции и, потому, с параболическими вращениями,
см. Предложение~\ref{r-pr:mul}(\ref{r-it:rot-as-mult})). Умножение на
положительный 
действительный скаляр задаётся прямолинейно: оно не должно менять
аргумент и соответственно маштабировать норму вектора. Поэтому для
\(a>0\) мы имеем:
\begin{eqnarray}
  \label{r-eq:scalar-prod}
  a\cdot (u,v)&=&(u,a v+u^2(1-a))\quad\text{для } N,\\
  \label{r-eq:scalar-prod-prime}
  a\cdot (u,v)&=&\left(u,\frac{v+1}{a}-1\right)\quad\text{для } N'.
\end{eqnarray}

С другой стороны, сложение векторов может быть задано несколькими
существенно отличными свойствами. Мы приведём две возможности:
тропическую и экзотическую.

\subsection{Тропическая форма}
\label{r-sec:tropical-form}
\label{r-ex:tropical}
Введём лексиграфический порядок  на \(\Space{R}{2}\):
\begin{displaymath}
  (u,v)\prec(u',v') \quad \text{тогда и только тогда,} \quad
  \left\{\begin{array}{ll}
      \text{когда}& u<u'; \\
      \text{или}& u=u',\  v<v'.
  \end{array}\right.
\end{displaymath}
Соответственно мы можем задать функции \(\min\) и \(\max\) для любой
пары векторов из \(\Space{R}{2}\) соответственно. Тогда сумма двух
векторов определяется как их минимум (или максимум). Такое же
определение используется в \emph{тропической математике}, так же
известной как <<деквантизация>> Маслова или алгебры
\(\Space[\min]{R}{}\) и \(\Space[\max]{R}{}\), см. превосходный
обзор~\cite{Litvinov05}.  Легко проверить, что такое сложение
совместно с
умножением~\eqref{r-eq:scalar-prod}---\eqref{r-eq:scalar-prod-prime}
удовлетворяет дистрибутивному закону и потому будет инвариантно
относительно параболических вращений. Это направление выглядит
интересным и обещающим, однако мы больше не будем его рассматривать в
данной работе.

\subsection{Экзотическая форма}
\label{r-sec:exotic-form}
Сложение двух векторов может быть задано одной формулой для обоих
подгрупп  \(N\) и \(N'\) таким образом, что их отличие будет спрятано
в соответствующих Определениях~\ref{r-de:norm} (норма) и~\ref{r-de:arg}
(аргумент).
\begin{defnr}
  \label{r-de:p-add}
  Параболическое сложение векторов задаётся в (экзотической) форме: 
  \begin{eqnarray}
    \label{r-eq:p-add-arg-exotic}
    \arg^{(\prime)}(w_1+w_2)&=&\frac{\arg^{(\prime)} w_1\cdot
      \modulus{w_1}^{(\prime)} 
      +\arg^{(\prime)} w_2\cdot\modulus{w_2}^{(\prime)}}{\modulus{w_1+w_2}^{(\prime)}},\\
    \label{r-eq:p-add-norm-exotic}
    \modulus{w_1+w_2}^{(\prime)}&=&\modulus{w_1}^{(\prime)}\pm\modulus{w_2}^{(\prime)},
  \end{eqnarray}
  где штрихованные версии формул соответствуют подгруппе \(N'\).
\end{defnr}

На первый взгляд формула для нормы суммы~\eqref{r-eq:p-add-norm-exotic}
выглядит слишком тривиально. В её защиту мы можем отметить, что она
располагается как раз между эллиптическим \(\modulus{w+w'}\leq
\modulus{w}+\modulus{w'}\) и гиперболическим \(\modulus{w+w'}\geq
\modulus{w}+\modulus{w'}\) неравенствами треугольника для норм.
  

Правило для сложения аргументов~\eqref{r-eq:p-add-arg-exotic} тоже не
является произвольным. Из теоремы синусов евклидовой геометрии мы
можем получить 
\begin{displaymath}
  \sin(\phi-\psi')=\frac{\modulus{w}\cdot\sin(\psi-\psi')}{\modulus{w+w'}},\qquad
  \sin(\psi'-\phi)=\frac{\modulus{w'}\cdot\sin(\psi-\psi')}{\modulus{w+w'}},
\end{displaymath}
где \(\psi^{(\prime)}=\arg w^{(\prime)}\) и \(\phi=\arg
(w+w^{(\prime)})\). Используя значение~\eqref{r-eq:par-trig-0}
параболического синуса \(\sinp \theta=\theta\) 
мы получаем~\eqref{r-eq:p-add-arg-exotic}.

Аккуратность в обращении с нулями в знаменателе
выражения~\eqref{r-eq:p-add-arg-exotic} достигается при переходе от
вектора \(w=u+\rmp v\) к паре однородных полярных координат \([a,r]=[
\modulus{w}^{(\prime)} \cdot \arg^{(\prime)} w,
\modulus{w}^{(\prime)}]\) (штрихованная версия соответствует подгруппе
\(N'\)). Тогда сложение задаётся покомпонентно в однородных координатах:
\begin{displaymath}
  w_1+w_2=[a_1+a_2, r_1+r_2], \qquad \text{ где } w_i=[a_i,r_i].
\end{displaymath}
Умножение из Определения~\ref{r-de:product} в однородных координатах есть:
\begin{displaymath}
  w_1\cdot w_2=[a_1r_2+a_2r_1, r_1 r_2], \qquad \text{ где } w_i=[a_i,r_i].
\end{displaymath}
Таким образом, однородные координаты линеаризуют
сложение~\eqref{r-eq:p-add-arg-exotic}--\eqref{r-eq:p-add-norm-exotic} и
умножение на скаляр~\eqref{r-eq:scalar-prod}. Переход к более прозрачным
координатам может быть рассмотрен в контексте бирациональной
геометрии~\cite{KollarMori08}. 

Параболическое сложение определяется
условиями~\eqref{r-eq:p-add-arg-exotic}--\eqref{r-eq:p-add-norm-exotic}
однозначно, но явные формулы достаточно сложны. Во многих случаях без
них можно обойтись используя следующие свойства:
\begin{propr}
  Векторное сложение для подгрупп \(N\) и \(N'\) заданное
  выражением~\eqref{r-eq:p-add-arg-exotic}--\eqref{r-eq:p-add-norm-exotic}
  обладает следующими свойствами:
  \begin{enumerate}
  \item\label{r-item:add-is-comm-ass} Они коммутативны и ассоциативны.
  \item\label{r-item:distrib}Выполнятся дистрибутивный закон для
    умножений~\eqref{r-eq:parab-prod} и~\eqref{r-eq:parab-prime-prod}; 
    соответственно:
  \item\label{r-item:add-rot-inv} Сложения перестановочны с поворотами;
  \item\label{r-item:distrib-scalar} Выполняется дистрибутивные законы
    для умножения на скаляр~\eqref{r-eq:scalar-prod}
    и~\eqref{r-eq:scalar-prod-prime} соответственно:
    \begin{displaymath}
      a\cdot(w_1+w_2)=a\cdot w_1+a\cdot w_2,\qquad
      (a+b)\cdot w=a\cdot w+b\cdot w.
    \end{displaymath}
  \end{enumerate}
\end{propr}
Что бы завершить построение действительного векторного пространства
нам надо описать нулевой и обратный вектора. В частности, мы определим
обратный элемент как имеющий тот же аргумент и противоположную норму.
\begin{propr}
  \begin{itemize}
  \item[(\(N\))] Нулевой вектор есть \((0,0)\) и соответственно
    противоположный элемент для \((u,v)\) есть \((u,2u^2-v)\).
  \item[(\(N'\))] Нулевой вектор есть \((\infty,-1)\) и соответственно
    противоположный элемент для  \((u,v)\) есть \((u,-v-2)\).
  \end{itemize}
\end{propr}
Теперь используя равенство \(-w=(-1)\cdot w\) можно убедится, что та
же формула~\eqref{r-eq:scalar-prod} и \eqref{r-eq:scalar-prod-prime} может
быть использована для умножения на отрицательные действительные
скаляры.

\begin{remr}
  \label{r-re:conformality}
  Обычная линейная структура является неподходящей и в других
  вопросах, например, в кажущейся <<неконформности>> длин из
  параболического фокуса, то есть с параметром \(\rs=0\)
  в~\cite{Kisil05a}*{Prop.~\ref{E-it:conformity-length-foci}}. Однако
  подходящая выбор определений восстанавливает искомые свойства.

  Первоначальное определение
  конформности~\cite{Kisil05a}*{Defn.~\ref{E-de:conformal}}
  рассматривало традиционный предел \(y'\rightarrow y\) вдоль прямых
  линий, т.е. <<спиц>> изображённых на Рис.~\ref{r-fig:rotations}.
  Хорошо известно, что этот подход полностью работает в эллиптическом
  и гиперболическом случае. Однако, в параболическом случае
  вертикальные линии больше подходят на роль <<спиц>>, см.
  Рис.~\ref{r-fig:p-rotations}(\(P\)) и (\(P'\)), таким образом
  обоснованно рассматривать конформный предел вдоль
  них~\cite{Kisil05a}*{Prop.~\ref{E-pr:parab-conf}}. 
\end{remr}

\section{Индуцированные представления}
\label{r-sec:concl-induc-repr}

Теперь мы можем применить гиперкомплексные характеры построенные выше
к индуцированным представлениям заданным формулой~\eqref{r-eq:def-ind}.
Заметим, что только компактная подгруппа \(K\) требует
комплексно-значных характеров, так как для однозначности представления
образ матрицы из~\eqref{r-eq:k-subgroup} при \(t=2\pi\) должен быть
единицей. Для подгрупп \(N'\) и \(A'\) такого ограничения нет и мы
можем рассмотреть характеры всех трёх типов: эллиптические.
параболические и гиперболические.  
Более того,
параболические характеры могут быть любого типов:
алгебраические~\eqref{r-eq:parab-rot-triv} или одного из
геометрических~\eqref{r-eq:par-moebius} и~\eqref{r-eq:par-moebius-prime}.
Таким образом мы имеем семь существенно разных типов индуцированных
представлений, или даже одиннадцати (если отличать разновидности
параболических характеров). 

\begin{exampler}
  Рассмотрим подгруппу \(H=K\), как отмечено выше, в силу её
  компактности характер может быть только комплексно-значным. Тогда
  унитарный характер \(\chi_k\) имеет вид:
  \begin{equation}
    \label{r-eq:k-character}
    \chi_k\begin{pmatrix}
      \cos t &  \sin t\\
      -\sin t & \cos t
    \end{pmatrix}=e^{-\rmi k t}, \qquad \text{ где }
    k\in\Space{Z}{}.
  \end{equation}
  Используя явный вид~\eqref{r-eq:s-map} отображения \(s\) вычисляем
  функцию \(r\) заданную в~\eqref{r-eq:r-map}:
  \begin{displaymath}
    r
    \begin{pmatrix}
      a&b\\c&d
    \end{pmatrix}
    =\frac{1}{\sqrt{c^2+d^2}}
    \begin{pmatrix}
      d&-c\\c&d
    \end{pmatrix}\in K.
  \end{displaymath}
  Тогда:
  \begin{displaymath}
    r(g^{-1} * s(u,v))  =  
    \frac{1}{\sqrt{(c u+d)^2 +(cv)^2}}
    \begin{pmatrix}
      cu+d&-cv\\cv&cu+d
    \end{pmatrix}, \text{ где } g^{-1}=    \begin{pmatrix}
      a&b\\c&d
    \end{pmatrix}.
  \end{displaymath}
  Подставляя это выражение в~\eqref{r-eq:k-character} и комбинируя с
  мёбиусовым преобразованием аргумента~\eqref{r-eq:sl2-act}  получаем
  явную реализацию \(\uir{k}{}\) индуцированного представления~\eqref{r-eq:def-ind}:
  \begin{equation}
    \label{r-eq:discrete}
    \uir{k}{}(g) f(w)=\frac{\modulus{cw+d}^k}{(cw+d)^k}f\left(\frac{aw+b}{cw+d}\right),
    \quad \text{ где } g^{-1}=\begin{pmatrix}a&b\\c&d
    \end{pmatrix}, \ w=u+\rmi v.
  \end{equation}
  Это представление действует на комплексно-значных функциях в верхней
  полуплоскости \(\Space[+]{R}{2}=\SL/K\) и является унитарным из
  дискретной серии~\cite{Lang85}*{\S~IX.2}.
  \begin{propr}
    \label{r-pr:k-eigen}
    Пусть \(f_k(w)=\frac{\modulus{w-\rmi}^{k}}{(w-\rmi)^{k}}\) for
    \(k=2\), \(3\), \ldots, тогда 
    \begin{enumerate}
    \item \(f_k\) является собственным вектором операторов
      \(\uir{k}{}(h)\), для любого \(h\in K\), с собственным значением
      \(\chi_k(h)\)~\cite{Lang85}*{\S~IX.2}.
    \item Функция \(K(z,w)=\uir{k}{}(s(z)) f_k(w)\), где \(s(z)\)
      определяется в~\eqref{r-eq:s-map}, 
       есть воспроизводящее ядро Бергмана в
      верхней полуплоскости~\cite{Kisil97c}*{\S~3.2} умноженное на
      \(\modulus{\frac{z-\rmi}{w-\rmi}}^k\). 
    \end{enumerate}
  \end{propr}
  Сходным образом мы получаем ядро Коши если рассмотрим предельный
  случай \(k=1\) ложной дискретной серии~\cite{Lang85}*{Ch.~IX}.
  Существует много других связей представления~\eqref{r-eq:discrete} с
  комплексным анализом~\cite{Kisil97c} и теорией операторов. К
  примеру, преобразования Мёбиуса операторов определяют функциональное
  исчисление Данфорда--Рисса и соответствующий спектр~\cite{Kisil02a}.
\end{exampler}

\begin{exampler}
  \label{r-ex:n-induced}
  Для подгруппы \(N'\) возможен более широкий выбор характеров.
  \begin{enumerate}
  \item Традиционно рассматривается комплеснозначные характеры
    подгруппы \(N'\), т.е.:
    \begin{equation}
      \label{r-eq:np-character}
      \chi_{\tau}^{\Space{C}{}}\begin{pmatrix}
        1 &  0\\
        t & 1
      \end{pmatrix}=e^{\rmi \tau t}, \qquad \text{ где }
      \tau\in\Space{R}{}.
    \end{equation}
    Непосредственные вычисления в этом случае показывают, что:
    \begin{displaymath}
      r
      \begin{pmatrix}
        a&b\\c&d
      \end{pmatrix}
      =
      \begin{pmatrix}
        1&0\\\frac{c}{d}&1
      \end{pmatrix}\in N'.
    \end{displaymath}
    Поэтому:
    \begin{equation}
      \label{r-eq:np-char-part}
      r(g^{-1}*s(u,v))=
      \begin{pmatrix}
        1&0\\\frac{cv}{d+cu}&1
      \end{pmatrix}, \quad\text{ где } g^{-1}=    \begin{pmatrix}
        a&b\\c&d
      \end{pmatrix}.
    \end{equation}
    Подставляя это значение в характер~\eqref{r-eq:np-character} и
    соединяя с мёбиусовым преобразованием~\eqref{r-eq:sl2-act} получаем
    соответствующую реализацию формулы~\eqref{r-eq:def-ind}:
    \begin{displaymath}
      \uir[{\Space{C}{}}]{\tau}(g) f(w)= \exp\left(\rmi\frac{\tau c v}{cu+d} \right)
      f\left(\frac{aw+b}{cw+d}\right), 
      \quad \text{где } w=u+\rmp v, \  g^{-1}=\begin{pmatrix}a&b\\c&d
      \end{pmatrix}.
    \end{displaymath}
    Соответствующее индуцированное представление действует в
    пространстве комплекснозначных функций в верхней полуплоскости
    \(\Space[+]{R}{2}\), которая является подмножеством однородного
    пространства \(\SL/N\) состоящего из дуальных чисел.  Соседство
    комплексных и дуальных чисел в одном выражении является
    непривычным.
  \item Возьмём параболический характер \(\chi_{\tau}\)
    алгебраического типа заданный умножением дуальных
    чисел~\eqref{r-eq:parab-rot-triv}:
    \begin{displaymath}
      \chi_{\tau}\begin{pmatrix}
        1 &  0\\
        t & 1
      \end{pmatrix}=e^{\rmp \tau t}=1+\rmp \tau t, \qquad \text{ где }
      \tau\in\Space{R}{}.
    \end{displaymath}
    По прежнему подставляем в него значение~\eqref{r-eq:np-char-part} и
    получаем такое выражение:
    \begin{displaymath}
      \uir{\tau}{}(g) f(w)= \left(1+\rmp\frac{\tau c v}{cu+d} \right)
      f\left(\frac{aw+b}{cw+d}\right),
    \end{displaymath}
    где \(w\), \(\tau\) и \(g\) те же что и раньше. 

    Это индуцированное представление действует в пространстве функций
    определённых на верхней полуплоскости, как подмножестве дуальных
    чисел, и со значениями в дуальных числах. То есть, оно использует
    только дуальные числа и обычные алгебраические операции над ними.
    Естественным образом это представление является линейным.
  \item Рассмотрим параболический характер \(\chi_{\tau}^{g}\)
    заданный геометрическим действием~\eqref{r-eq:par-moebius}. Тогда
    представление действует в том же пространстве, что и в последнем
    случае, но формула будет другой:
    \begin{displaymath}
      \uir[g]{\tau}(g) f(w)= \left(1+\rmp\frac{2\tau c v}{cu+d} \right)
      f\left(\frac{aw+b}{cw+d}\right)+\frac{\tau c v}{cu+d}+\rmp \frac{(\tau c v)^2}{(cu+d)^2}, 
    \end{displaymath}
    где \(w\), \(\tau\) и \(g\) описаны выше. Это представление
    линейно по отношению к операциям~\eqref{r-eq:scalar-prod},
    \eqref{r-eq:p-add-arg-exotic} и~\eqref{r-eq:p-add-norm-exotic}.
  \end{enumerate}
\end{exampler}
Все характеры использованные в предыдущих примерах являются
унитарными%
: два первых в традиционном смысле, последний---как объяснено в
Предложении~\ref{r-pr:mul}. Соответственно общая конструкция
индуцированных представлений~\cite{Kirillov76}*{\S~13.2} гарантирует
унитарность полученных представлений в соответствующих смыслах.
\begin{thmr}
  \label{r-th:unitarity}
  Все три индуцированных представления \(\uir[{\Space{C}{}}]{\tau}\),
  \(\uir{\tau}\) и \(\uir[g]{\tau}\)
  группы \(\SL\) из
  Примера~\ref{r-ex:n-induced} являются унитарными в пространстве
  функций заданных на верхней полуплоскости \(\Space[+]{R}{2}\) как
  подмножестве дуальных чисел с внутренним произведением:
  \begin{equation}
    \label{r-eq:inner-product}
    \scalar{f_1}{f_2}=\int_{\Space[+]{R}{2}} f_1(w)
    \bar{f}_2(w)\,\frac{du\,dv}{v^2}, \qquad \text{ где } w=u+\rmp v,
  \end{equation}
  и мы используем
   \begin{enumerate}
  \item сопряжение и умножение значений функций в алгебре
    комплексных и дуальных чисел для представлений
    \(\uir[{\Space{C}{}}]{\tau}\) и \(\uir{\tau}\) соответственно;
  \item сопряжение~\eqref{r-eq:conjugation} и
    умножение~\eqref{r-eq:parab-prod} значений функций для представления
    \(\uir[g]{\tau}\).
  \end{enumerate}
\end{thmr}
Внутреннее произведение~\eqref{r-eq:inner-product} положительно
определено для представления \(\uir[{\Space{C}{}}]{\tau}\), но не для
двух других. 
Соответствующие пространства с вырожденным внутренним
произведением являются параболическими аналогами пространств
Крейна~\cite{ArovDym08}. Сами пространства Крейна в нашей
классификации являются гиперболическими.

\section{Принцип сходства и соответствия: операторы повышения}
\label{r-sec:correspondence}

Приведённые выше наблюдения позволяют сформулировать следующий
эмпирический принцип, который должен представлять эвристический
интерес:

\begin{principler}[Сходство и соответствие]
  \begin{enumerate}
  \item Подгруппы \(K\), \(N'\) и \(A'\) играют сходную роль в
    структуре группы \(\SL\) и её представлений;
  \item Перемена подгруппы должна производится с одновременной
    соответствующей заменой гиперкомплексных единиц \(\rmi\), \(\rmp\)
    или \(\rmh\).
  \end{enumerate}
\end{principler}
Безусловно, первая часть Принципа (сходство) способна вызвать удивление
у любого, кто знаком с группой \(\SL\). Однако, к настоящему моменту мы
видели, что совместно со второй частью (соответствие), принцип
уже проявился следующим образом:
\begin{itemize}
\item Действие \(\SL\) на однородном пространстве \(\SL/H\) для
  \(H=K\), \(N'\) или \(A'\) задаётся дробно-линейными
  преобразованиями комплексных, дуальных или двойных чисел
  соответственно (Предл.~\ref{r-pr:sl2-act}). 
\item Подгруппы \(A'\), \(N'\) или \(K\) изоморфны группам унитарных
  вращений соответствующих единичных <<окружностей>> на плоскостях
  двойных, дуальных или комплексных чисел
  (Предл.~\ref{r-pr:algebraic-char}).   
\item Представления индуцированные с подгруппы \(K\), \(N'\) или
  \(A\) являются унитарными если скалярное произведение функций
  определяются посредством произведения и сопряжения значений функций
  как комплексных, дуальных или двойных чисел (Теор.~\ref{r-th:unitarity}).
\end{itemize}
\begin{remr}
  Принцип сходства и соответствия близок к суперсимметрии между
  бозонами и фермионами в физике, но в нашем случае сходство
  устанавливается между тремя различными типами объектов.
\end{remr}
Приведём ещё одну иллюстрацию Принципа. Рассмотрим алгебру Ли
\(\algebra{sl}_2\)  группы \(\SL\). Один из возможных
базисов в \(\algebra{sl}_2\) таков~\cite{MTaylor86}*{\S~8.1}:  
\begin{displaymath}
  A= \frac{1}{2}
  \begin{pmatrix}
    -1&0\\0&1
  \end{pmatrix},\quad 
  B= \frac{1}{2} \
  \begin{pmatrix}
    0&1\\1&0
  \end{pmatrix}, \quad 
  Z=
  \begin{pmatrix}
    0&1\\-1&0
  \end{pmatrix}.
\end{displaymath}
Тогда коммутационные соотношения будут следующими:
\begin{equation}
  \label{r-eq:sl2-commutator}
  [Z,A]=2B, \qquad [Z,B]=-2A, \qquad [A,B]=- \frac{1}{2} Z.
\end{equation} 
Пусть у нас есть представление \(\uir{}{}\) группы \(\SL\) в
пространстве \(V\). Рассмотрим производное представление \(d\uir{}{}\)
алгебры \(\algebra{sl}_2\)~\cite{Lang85}*{\S~VI.1} и будем
использовать обозначениe \(\tilde{X}=d\uir{}{}(X)\) для
\(X\in\algebra{sl}_2\). Для анализа представления \(\uir{}{}\)
оказывается полезным разложить \(V\) по собственным векторам оператора
\(\tilde{X}\) для некоторого \(X\in \algebra{sl}_2\), см.
Предл.~\ref{r-pr:k-eigen} и ряд Тайлора в комплексном анализе.
\begin{exampler}
  Уже не должно быть удивительным, что мы собираемся рассмотреть три
  случая:
  \begin{enumerate}
  \item Пусть \(X=Z\) будет генератором подгруппы
    \(K\)~\eqref{r-eq:k-subgroup}. Так как эта подгруппа компактна, то
    собственные вектора \(\tilde{Z} v_k=\rmi k v_k\) параметризованы
    дискретным параметром \(k\in\Space{Z}{}\).  Особую роль играют
    операторы \(\ladder{\pm}\)
    повышения/понижения~\citelist{\cite{Lang85}*{\S~VI.2}
      \cite{MTaylor86}*{\S~8.2}}, определённые коммутационными
    соотношениями:
    \begin{equation}
      \label{r-eq:raising-lowering}
      [\tilde{Z},\ladder{\pm}]=\lambda_\pm \ladder{\pm}. 
    \end{equation}
    Таким образом \(\ladder{\pm}\) являются собственными для операторов
    \(\mathrm{ad}Z\) присоединённого
    представления \(\algebra{sl}_2\)~\cite{Lang85}*{\S~VI.2}.
    Важно, что вектор \(\ladder{+} v_k\) так же является собственным для
    \(\tilde{Z}\):
    \begin{eqnarray*}
      \tilde{Z}(\ladder{+} v_k)&=&(\ladder{+}\tilde{Z}+\lambda_+\ladder{+})v_k=\ladder{+}(\tilde{Z}v_k)+\lambda_+\ladder{+}v_k
      =\rmi k \ladder{+}v_k+\lambda_+\ladder{+}v_k\\
      &=&(\rmi k+\lambda_+)\ladder{+}v_k.
    \end{eqnarray*}

    Полагая \(\ladder{+}=a\tilde{A}+b\tilde{B}+c\tilde{Z}\) из коммутационного
    соотношений~\eqref{r-eq:sl2-commutator} и определяющего
    равенства~\eqref{r-eq:raising-lowering} получаем систему уравнений:
    \begin{displaymath}
      c=0, \qquad 2a=\lambda_+ b, \qquad -2b=\lambda_+ a.
    \end{displaymath}
    Эта система совместима тогда и только тогда, когда
    \(\lambda_+^2+4=0\). В этом случае операторы
    повышения/понижения\(\ladder{\pm}=\pm\rmi \tilde{A}+\tilde{B}\) действуют
    на одномерной цепочке собственных подпространств:
    \begin{displaymath}
             \xymatrix@1{
        \ldots\, \ar@<.4ex>[r]^{\ladder{+}} & 
        \,V_{k-2}\,  \ar@<.4ex>[l]^{\ladder{-}}\ar@<.4ex>[r]^{\ladder{+}} &  
        \,V_{k}\, \ar@<.4ex>[l]^{\ladder{-}} \ar@<.4ex>[r]^{\ladder{+}} & 
        \,V_{k+2}\,\ar@<.4ex>[l]^{\ladder{-}}  \ar@<.4ex>[r]^{\ladder{+}}
        & 
        \,\ldots\ar@<.4ex>[l]^{\ladder{-}}}
    \end{displaymath}
  \item Рассмотрим случай \(X=B\) генератора подгруппы
    \(A'\)~\eqref{r-eq:ap-subgroup}. Эта подгруппа некомпактна и
    собственные значения для \(\tilde{B}\) могут быть произвольным числом,
    однако операторы повышения/понижения по прежнему могут играть важную
    роль~\citelist{\cite{HoweTan92}*{\S~II.1}
      \cite{Mazorchuk09a}*{\S~1.1}}.  Будем опять искать решение в форме
    \(\ladder[h]{+}=a\tilde{A}+b\tilde{B}+c\tilde{Z}\) для
    коммутатора \([\tilde{B},\ladder[h]{+}]=\lambda \ladder{+}\). Получаем систему:
    \begin{displaymath}
      2c=\lambda a,\qquad
      b=0,\qquad
      \frac{a}{2}=\lambda c.
    \end{displaymath}
    Она совместима только если \(\lambda^2=1\). Очевидные значения
    \(\lambda=\pm 1\) приводят к операторам
    \(\ladder[h]{\pm}=\pm\tilde{A}+\tilde{Z}/2\), которые и
    используются в~\citelist{\cite{HoweTan92}*{\S~II.1}
      \cite{Mazorchuk09a}*{\S~1.1}}. В этом случае неприводимый
    \(\algebra{sl}_2\)-модуль представляется одномерной цепочкой
    собственных значений.
    
    Допуская двойные числа мы можем удовлетворить условие
    совместимости \(\lambda^2=1\) дополнительными значениями
    \(\lambda=\pm\rmh\).  Тогда возникает дополнительная пара
    операторов повышения/понижения
    \(\ladder[\rmh]{\pm}=\pm\rmh\tilde{A}+\tilde{Z}/2\), которая сдвигает
    собственные значения в <<ортогональном>> направлении к действию
    традиционных операторов \(\ladder[h]{\pm}\).  Следовательно неприводимый
    \(\algebra{sl}_2\)-модуль может параметризоваться уже двумерной
    решёткой собственных значений на плоскости двойных чисел:
\begin{displaymath}
  \xymatrix@R=2.5em@C=1.5em@M=.5em{
      & 
    \,\ldots\, \ar@<.4ex>[d]^{\ladder[\rmh]{+}} &  
    \,\ldots\, \ar@<.4ex>[d]^{\ladder[\rmh]{+}} & 
    \,\ldots\,  \ar@<.4ex>[d]^{\ladder[\rmh]{+}}  & 
    \\
    \ldots\, \ar@<.4ex>[r]^-{\ladder[h]{+}} & 
    \,V_{(n-2)+\rmh (k-2)}\,  \ar@<.4ex>[l]^-{\ladder[h]{-}}\ar@<.4ex>[r]^{\ladder[h]{+}}
     \ar@<.4ex>[u]^{\ladder[\rmh]{-}} \ar@<.4ex>[d]^{\ladder[\rmh]{+}} &  
    \,V_{n+\rmh (k-2)}\, \ar@<.4ex>[l]^{\ladder[h]{-}} \ar@<.4ex>[r]^{\ladder[h]{+}}
     \ar@<.4ex>[u]^{\ladder[\rmh]{-}} \ar@<.4ex>[d]^{\ladder[\rmh]{+}} & 
    \,V_{(n+2)+\rmh (k-2)}\,\ar@<.4ex>[l]^{\ladder[h]{-}}  \ar@<.4ex>[r]^-{\ladder[h]{+}}
     \ar@<.4ex>[u]^{\ladder[\rmh]{-}} \ar@<.4ex>[d]^{\ladder[\rmh]{+}}    & 
    \,\ldots\ar@<.4ex>[l]^-{\ladder[h]{-}}\\
    \ldots\, \ar@<.4ex>[r]^-{\ladder[h]{+}} & 
    \,V_{(n-2)+\rmh k}\,  \ar@<.4ex>[l]^-{\ladder[h]{-}}\ar@<.4ex>[r]^{\ladder[h]{+}}
     \ar@<.4ex>[u]^{\ladder[\rmh]{-}} \ar@<.4ex>[d]^{\ladder[\rmh]{+}} &  
    \,V_{n+\rmh k}\, \ar@<.4ex>[l]^{\ladder[h]{-}} \ar@<.4ex>[r]^{\ladder[h]{+}}
     \ar@<.4ex>[u]^{\ladder[\rmh]{-}} \ar@<.4ex>[d]^{\ladder[\rmh]{+}}& 
    \,V_{(n+2)+\rmh k}\,\ar@<.4ex>[l]^{\ladder[h]{-}}  \ar@<.4ex>[r]^-{\ladder[h]{+}}
     \ar@<.4ex>[u]^{\ladder[\rmh]{-}} \ar@<.4ex>[d]^{\ladder[\rmh]{+}}    & 
    \,\ldots\ar@<.4ex>[l]^-{\ladder[h]{-}}\\
    \ldots\, \ar@<.4ex>[r]^-{\ladder[h]{+}} & 
    \,V_{(n-2)+\rmh (k+2)}\,  \ar@<.4ex>[l]^-{\ladder[h]{-}}\ar@<.4ex>[r]^{\ladder[h]{+}}
     \ar@<.4ex>[u]^{\ladder[\rmh]{-}} \ar@<.4ex>[d]^{\ladder[\rmh]{+}} &  
    \,V_{n+\rmh (k+2)}\, \ar@<.4ex>[l]^{\ladder[h]{-}} \ar@<.4ex>[r]^{\ladder[h]{+}} 
     \ar@<.4ex>[u]^{\ladder[\rmh]{-}} \ar@<.4ex>[d]^{\ladder[\rmh]{+}}& 
    \,V_{(n+2)+\rmh (k+2)}\,\ar@<.4ex>[l]^{\ladder[h]{-}}  \ar@<.4ex>[r]^-{\ladder[h]{+}}
     \ar@<.4ex>[u]^{\ladder[\rmh]{-}} \ar@<.4ex>[d]^{\ladder[\rmh]{+}}    & 
    \,\ldots\ar@<.4ex>[l]^-{\ladder[h]{-}}\\
      & 
    \,\ldots\, \ar@<.4ex>[u]^{\ladder[\rmh]{-}} &  
    \,\ldots\, \ar@<.4ex>[u]^{\ladder[\rmh]{-}} & 
    \,\ldots\,  \ar@<.4ex>[u]^{\ladder[\rmh]{-}}  & }
\end{displaymath}
    
  \item Рассмотрим, наконец, случай \(X=-B+Z/2\) генератора подгруппы
    \(N'\)~\eqref{r-eq:np-subgroup}. По описанной ранее процедуре получаем
    систему:
    \begin{displaymath}
      b+2c=\lambda a,\qquad
      -a=\lambda b,\qquad
      \frac{a}{2}=\lambda c
    \end{displaymath}
    с условием совместимости \(\lambda^2=0\). Если мы ограничимся
    действительным (комплексным) корнем \(\lambda=0\), то
    соответствующие операторы \(\ladder[p]{\pm}=-\tilde{B}+\tilde{Z}/2\) не будут
    менять собственное значения и будут бесполезны в
    таком контексте. Однако дуальные значения \(\lambda =\pm\rmp\)
    позволяют использовать операторы \(\ladder[\rmp]{\pm}=\pm
    \rmp\tilde{A}-\tilde{B}+\tilde{Z}/2\) для построения
    \(\algebra{sl}_2\)-модулей с одномерной цепочкой
    собственных значений в дуальных числах.
  \end{enumerate}
\end{exampler}
\begin{remr}
  Стоит отметить, что 
  \begin{itemize}
  \item введение комплексных чисел необходимо для \emph{существования}
    операторов повышения/понижения в эллиптическом случае;
  \item в параболическом случае введение дуальных чисел необходимо для
    \emph{полезности} этих операторов;
  \item в гиперболическом случае двойственные числа выглядят
    необязательными. 
  \end{itemize}
\end{remr}
Подытожим рассмотренный пример подчеркнув в нём роль принципа
сходства и соответствия:  
\begin{propr}
  Пусть вектор \(X\in\algebra{sl}_2\) порождает подгруппу \(K\),
  \(N'\) или \(A'\), то есть векторами \(X=Z\), \(B-Z/2\), или
  \(B\) соответственно. Пусть \(\alli\) есть соответствующая
  гиперкомплексная единица.   

  Тогда операторы повышения/понижения \(\ladder{\pm}\) удовлетворяющие
  коммутационным соотношениям:
  \begin{displaymath}
    [X,\ladder{\pm}]=\pm\alli \ladder{\pm},\qquad [\ladder{-},\ladder{+}]=2\alli X.
  \end{displaymath}
  имеют вид:
  \begin{displaymath}
    \ladder{\pm}=\pm\alli \tilde{A} +\tilde{Y}.
  \end{displaymath}
  Здесь \(Y\in\algebra{sl}_2\) является линейной комбинацией \(B\) и
  \(Z\) а так же удовлетворяет  следующим условиям:
  \begin{itemize}
  \item \(Y=[A,X]\).
  \item \(X=[A,Y]\).
  \item форма Киллинга \(K(X,Y)\)~\cite{Kirillov76}*{\S~6.2} равна нулю.
  \end{itemize}
  Любое из приведённых выше условий совместно с \(Y\in\mathrm{span}\{B,Z\}\)
  определяет \(Y\) с точностью до действительного множителя.
\end{propr}

Кажется правдоподобным, что применение Принципа сходства и
соответствия не ограничивается приведёнными здесь случаями.
Например, сходным образом могут быть рассмотрены операторы
рождения/уничтожения для гармонического осциллятора и представлений
группы Гейзенберга~\cite{Kisil11a}.

\section{Заключение: открытые вопросы}
\label{r-sec:conclusion}

Введенные индуцированные представления заслуживают пристального
исследования. Среди важных вопросов можно отметить следующие: 
\begin{itemize}
\item Унитарные комплексно-значные представления \(\SL\) делятся на
  три серии: дискретную, непрерывную,
  дополнительную~\cites{Lang85,HoweTan92}.  Как эта классификация
  соотносится с введёнными гиперкомплексными представлениями
  группы \(\SL\).
\item Не каждый неприводимый \(\SL\)-модуль допускает
  введение комплексного скалярного произведения, в котором
  представление становится унитарным~\cites{HoweTan92,Mazorchuk09a}.
  Возможно, некоторые из таких модулей унитаризуются в
  гиперкомплексных числах.
\item Многие \(\SL\)-модули, допускающие гиперкомплексную унитаризацию,
  могут быть реализованы как пространства функций на действительной
  прямой или верхней полуплоскости~\citelist{\cite{Kisil97c}
    \cite{Kisil11c}*{\S~5}}. Концепция ковариантного преобразования~\citelist{\cite{Kisil10c}
    \cite{Kisil11c}*{\S~4}} позволяет строить в этих пространствах
  интегральные формулы и находить воспроизводящие ядра.
\item Упомянутые пространства функций состоят из нулевых векторов
  некоторых \(\SL\)-инвариантных дифференциальных операторов, которые
  являются аналогами операторов Лапласа и Коши--Римана. В частности,
  ковариантное преобразование с материнским всплеском, которое
  исчезает под действием оператора уничтожения (или любой его
  степени),  является решением соответствующего
  дифференциального уравнения~\cite{Kisil11c}*{\S~5.3}. Это может быть
  использовано для анализа решений уравнений в частных
  производных~\cites{Konovenko09a,Pilipchuk10a}.
\item Аналитические функции используются в качестве модели
  функционального исчисления операторов, в исчислении
  Гельфанада--Данфорда--Риса. Новые неприводимые \(\SL\)-модули
  являются подходящими моделями для исчисления и спектральной теории
  несамосопряжённых операторов~\citelist{\cite{Kisil02a}
    \cite{Kisil95i} \cite{Kisil11c}*{\S~6}}.
\end{itemize}
Эти направления являются частями Эрлангенской программы в широком
смысле~\cites{Kisil06a,Kisil11c} и должны быть продолжены в
последующих работах. Ожидается, что Принцип сходства и соответствия
окажет направляющую роль в поисках наиболее гармоничных конструкций.

\section*{Благодарности}
Автор благодарен проф.~Н.А.~Громову за полезное обсуждение данной работы и
многочисленные ценные замечания.

{\small
\providecommand{\noopsort}[1]{} \providecommand{\printfirst}[2]{#1}
  \providecommand{\singleletter}[1]{#1} \providecommand{\switchargs}[2]{#2#1}
  \providecommand{\irm}{\textup{I}} \providecommand{\iirm}{\textup{II}}
  \providecommand{\vrm}{\textup{V}} \providecommand{\cprime}{'}
  \providecommand{\eprint}[2]{\texttt{#2}}
  \providecommand{\myeprint}[2]{\texttt{#2}}
  \providecommand{\arXiv}[1]{\myeprint{http://arXiv.org/abs/#1}{arXiv:#1}}
  \providecommand{\doi}[1]{\href{http://dx.doi.org/#1}{doi:
  #1}}\providecommand{\CPP}{\texttt{C++}}
  \providecommand{\NoWEB}{\texttt{noweb}}
  \providecommand{\MetaPost}{\texttt{Meta}\-\texttt{Post}}
  \providecommand{\GiNaC}{\textsf{GiNaC}}
  \providecommand{\pyGiNaC}{\textsf{pyGiNaC}}
  \providecommand{\Asymptote}{\texttt{Asymptote}}

}

\newpage
  \selectlanguage{english} 

\appendix

\IfFileExists{parab-rotation.d}{
\section{Output of Symbolic Calculations}
\label{r-sec:outp-symb-calc} 
Here are the results of our symbolic calculations. The source code can
be obtained from this paper~\cite{Kisil09c} source at
\url{http://arXiv.org}. It uses Clifford algebra
facilities~\cite{Kisil05b} of the \GiNaC\ library~\cite{GiNaC}. The
source code is written in \NoWEB~\cite{NoWEB} literature programming
environment.

\noindent Calculations for subgroup \(N\) and straight spokes\\
Reference point: \(\left(\begin{array}{cc}0&-1\end{array}\right)\)\\
Reference point arg: \(0\)\\
Reference point norm: \(1\)\\
Cayley of the matrix X: \(\left(\begin{array}{cc}\clifford[0]{\mathbf{1}}+ \clifford[0]{e}^{{0} } x& \clifford[0]{\mathbf{1}} x\\0&\clifford[0]{\mathbf{1}}- \clifford[0]{e}^{{0} } x\end{array}\right)\)\\
Rotation by x: \(\left(\begin{array}{cc}u+x&2  u x+v+x^{2}\end{array}\right)\)\\
Rotation of \((u_0, v_0)\) by \(x\): \(\left(\begin{array}{cc}x&-1+x^{2}\end{array}\right)\)\\
Parabolic norm: \(u^{2}-v\)\\
Parabolic argument: \(u\)\\
Real number \(t_1\) as a dual number: \(\left(\begin{array}{cc}0&- t_1\end{array}\right)\)\\
Product: \(\left(\begin{array}{cc}u'+u&u'^{2}+2  u' u+u^{2}+ v' u^{2}- u'^{2} u^{2}- v' v+ u'^{2} v\end{array}\right)\)\\
Product by a scalar: \(\left(\begin{array}{cc}u&u^{2}+ a v- a u^{2}\end{array}\right)\)\\
Real part: \(\left(\begin{array}{cc}0&- u v-u^{2}+u^{3}+v\end{array}\right)\)\\
Imag part: \(\left(\begin{array}{cc}1&1+ u v-u^{3}\end{array}\right)\)\\
Zero vector: \(\left(\begin{array}{cc}0&0\end{array}\right)\)\\
Negative vector: \(\left(\begin{array}{cc}u&2 u^{2}-v\end{array}\right)\)\\
Dual number from argument \(0\) and norm \(n\): \(\left(\begin{array}{cc}0&- n\end{array}\right)\)\\
Dual number from argument \(a_1\) and norm \(n\): \(\left(\begin{array}{cc}a_1&a_1^{2}-n\end{array}\right)\)\\
Dual number from argument \(a_1\) and norm \(n\)---norm: \(n\)\\
Dual number from argument \(a_1\) and norm \(n\)---arg: \(a_1\)\\
Lin comb of two vectors \(a*(1, 0)+b*(-1, 0)\): \(\left(\begin{array}{cc}-\frac{b-a}{b+a}&\frac{b^{2}-b^{3}-a^{3}-3  b^{2} a+a^{2}-2  b a-3  b a^{2}}{{(b+a)}^{2}}\end{array}\right)\)\\
\(P\) is the sum \(\Re(P)\) and \(\Im(P)\): \textbf{true}\\
The real part of a real dual number is itself: \textbf{true}\\
norm is invariant under rotations: \textbf{true}\\
Product is invariant under rotations: \textbf{true}\\
Product \(w\bar{w}\) is norm squared: \textbf{true}\\
Product \((u, v)*(u_0, v_0)\) is \((u, v)\): \textbf{true}\\
Add is commutative: \textbf{true}\\
Add is associative: \textbf{true}\\
S-mult is commutative: \textbf{true}\\
S-mult is associative: \textbf{true}\\
S-mult is distributive 1: \textbf{true}\\
S-mult is distributive 2: \textbf{true}\\
Product is symmetric (commutative): \textbf{true}\\
Prod is associative: \textbf{true}\\
Product is distributive: \textbf{true}\\
\vspace{2mm}\hrule
Calculations for subgroup \(N'\) and straight spokes\\
Reference point: \(\left(\begin{array}{cc}\infty&-1+\infty^{2}\end{array}\right)\)\\
Reference point arg: \(0\)\\
Reference point norm: \(1\)\\
Cayley of the matrix X: \(\left(\begin{array}{cc}\clifford[0]{\mathbf{1}}- \clifford[0]{e}^{{0} } x&0\\ \clifford[0]{\mathbf{1}} x&\clifford[0]{\mathbf{1}}+ \clifford[0]{e}^{{0} } x\end{array}\right)\)\\
Rotation by x: \(\left(\begin{array}{cc}\frac{u}{1+ u x}&-\frac{2  u x+ u^{2} x^{2}-v}{1+2  u x+ u^{2} x^{2}}\end{array}\right)\)\\
Rotation of \((u_0, v_0)\) by \(x\): \(\left(\begin{array}{cc}\frac{1}{x}&-\frac{-1+x^{2}}{x^{2}}\end{array}\right)\)\\
Parabolic norm: \(\frac{u^{2}}{1+v}\)\\
Parabolic argument: \(\frac{1}{u}\)\\
Real number \(t_1\) as a dual number: \(\left(\begin{array}{cc}\infty&\frac{\infty^{2}-t_1}{t_1}\end{array}\right)\)\\
Product: \(\left(\begin{array}{cc}\frac{ u' u}{u'+u}&\frac{1+v'-u'^{2}-2  u' u-u^{2}+ v' v+v}{{(u'+u)}^{2}}\end{array}\right)\)\\
Product by a scalar: \(\left(\begin{array}{cc}u&-\frac{-1+a-v}{a}\end{array}\right)\)\\
Real part: \(\left(\begin{array}{cc}\infty&\frac{\infty^{2}-u^{2}+ \infty^{2} v+u}{ {(-1+u)} u}\end{array}\right)\)\\
Imag part: \(\left(\begin{array}{cc}1&-\frac{-1+u-v}{u}\end{array}\right)\)\\
Zero vector: \(\left(\begin{array}{cc}\infty&-1\end{array}\right)\)\\
Negative vector: \(\left(\begin{array}{cc}u&-2-v\end{array}\right)\)\\
Dual number from argument \(0\) and norm \(n\): \(\left(\begin{array}{cc}\infty&\frac{\infty^{2}-n}{n}\end{array}\right)\)\\
Dual number from argument \(a_1\) and norm \(n\): \(\left(\begin{array}{cc}\frac{1}{a_1}&-\frac{-1+ a_1^{2} n}{ a_1^{2} n}\end{array}\right)\)\\
Dual number from argument \(a_1\) and norm \(n\)---norm: \(n\)\\
Dual number from argument \(a_1\) and norm \(n\)---arg: \(a_1\)\\
Lin comb of two vectors \(a*(1, 0)+b*(-1, 0)\): \(\left(\begin{array}{cc}-\frac{b+a}{b-a}&\frac{b-b^{2}-a^{2}+2  b a+a}{{(b-a)}^{2}}\end{array}\right)\)\\
\(P\) is the sum \(\Re(P)\) and \(\Im(P)\): \textbf{true}\\
The real part of a real dual number is itself: \textbf{true}\\
norm is invariant under rotations: \textbf{true}\\
Product is invariant under rotations: \textbf{true}\\
Product \(w\bar{w}\) is norm squared: \textbf{true}\\
Product \((u, v)*(u_0, v_0)\) is \((u, v)\): \textbf{true}\\
Add is commutative: \textbf{true}\\
Add is associative: \textbf{true}\\
S-mult is commutative: \textbf{true}\\
S-mult is associative: \textbf{true}\\
S-mult is distributive 1: \textbf{true}\\
S-mult is distributive 2: \textbf{true}\\
Product is symmetric (commutative): \textbf{true}\\
Prod is associative: \textbf{true}\\
Product is distributive: \textbf{true}\\
\vspace{2mm}\hrule
Elliptic case of induced representations\\
map \(r(M)\): \(\left(\begin{array}{cc}\frac{d^{2}}{d^{2}+c^{2}}&-\frac{ c d}{d^{2}+c^{2}}\\\frac{ c d}{d^{2}+c^{2}}&\frac{d^{2}}{d^{2}+c^{2}}\end{array}\right)\)\\
map \(s^{-1}(M)\): \(\left(\begin{array}{cc}\frac{1}{d}&\frac{ b c^{2}+c+ b d^{2}}{d^{2}}\\0&\frac{d^{2}+c^{2}}{d}\end{array}\right)\)\\
character: \(\left(\begin{array}{cc}\frac{2  c d u+d^{2}+ c^{2} u^{2}}{2  c d u+d^{2}+ c^{2} u^{2}+ c^{2} v^{2}}&-\frac{ c {( c u+d)} v}{2  c d u+d^{2}+ c^{2} u^{2}+ c^{2} v^{2}}\\\frac{ c {( c u+d)} v}{2  c d u+d^{2}+ c^{2} u^{2}+ c^{2} v^{2}}&\frac{{( c u+d)}^{2}}{2  c d u+d^{2}+ c^{2} u^{2}+ c^{2} v^{2}}\end{array}\right)\)\\
Moebius map: \(\left(\begin{array}{cc}\frac{ c a v^{2}+ c a u^{2}+ b c u+ b d+ d a u}{2  c d u+d^{2}+ c^{2} u^{2}+ c^{2} v^{2}}&-\frac{ b c v- d a v}{2  c d u+d^{2}+ c^{2} u^{2}+ c^{2} v^{2}}\end{array}\right)\)\\
Moebius map is given by the imaginary unit: \textbf{true}\\
\vspace{2mm}\hrule
Parabolic (\(N^\prime\)) case of induced representations\\
map \(r(M)\): \(\left(\begin{array}{cc}1&0\\\frac{c}{d}&1\end{array}\right)\)\\
map \(s^{-1}(M)\): \(\left(\begin{array}{cc}\frac{1}{d}&b\\0&d\end{array}\right)\)\\
character: \(\left(\begin{array}{cc}1&0\\\frac{ c v}{ c u+d}&1\end{array}\right)\)\\
Moebius map: \(\left(\begin{array}{cc}\frac{b+ a u}{ c u+d}&-\frac{ b c v- d a v}{{( c u+d)}^{2}}\end{array}\right)\)\\
Moebius map is given by the imaginary unit: \textbf{true}\\
\vspace{2mm}\hrule
Hyperbolic case of induced representations\\
map \(r(M)\): \(\left(\begin{array}{cc}\frac{d^{2}}{d^{2}-c^{2}}&\frac{ c d}{d^{2}-c^{2}}\\\frac{ c d}{d^{2}-c^{2}}&\frac{d^{2}}{d^{2}-c^{2}}\end{array}\right)\)\\
map \(s^{-1}(M)\): \(\left(\begin{array}{cc}\frac{1}{d}&-\frac{ b c^{2}+c- b d^{2}}{d^{2}}\\0&\frac{d^{2}-c^{2}}{d}\end{array}\right)\)\\
character: \(\left(\begin{array}{cc}\frac{2  c d u+d^{2}+ c^{2} u^{2}}{2  c d u+d^{2}+ c^{2} u^{2}- c^{2} v^{2}}&\frac{ c {( c u+d)} v}{2  c d u+d^{2}+ c^{2} u^{2}- c^{2} v^{2}}\\\frac{ c {( c u+d)} v}{2  c d u+d^{2}+ c^{2} u^{2}- c^{2} v^{2}}&\frac{{( c u+d)}^{2}}{2  c d u+d^{2}+ c^{2} u^{2}- c^{2} v^{2}}\end{array}\right)\)\\
Moebius map: \(\left(\begin{array}{cc}-\frac{ c a v^{2}- c a u^{2}- b c u- b d- d a u}{2  c d u+d^{2}+ c^{2} u^{2}- c^{2} v^{2}}&-\frac{ b c v- d a v}{2  c d u+d^{2}+ c^{2} u^{2}- c^{2} v^{2}}\end{array}\right)\)\\
Moebius map is given by the imaginary unit: \textbf{true}\\
\vspace{2mm}\hrule

}{}

 \IfFileExists{dualnum.tex}{
 \section{Library for Dual Number Algebra} 
 \label{sec:lib-symb-calc}
 \label{r-sec:lib-symb-calc}
 This is the second version of the library which perform non-standard
 computations with dual number described in this paper. The first
 version was published in \url{arXiv.org} with the earlier more
 heuristic paper on this subject~\cite{Kisil07a}.  \LaTeX\ compilation
 of the file of this paper~\cite{Kisil09c} will produces the file
 \texttt{dualnum.nw} in the current directory.  This is a
 \NoWEB~\cite{NoWEB} code of the program. It uses Clifford algebra
 facilities~\cite{Kisil05b} of the \GiNaC\ library~\cite{GiNaC}.

\def\LA{\begingroup\maybehbox\bgroup\setupmodname\Rm$\langle$}\def\RA{$\rangle$\egroup\endgroup}\providecommand{\MM}{\kern.5pt\raisebox{.4ex}{\begin{math}\scriptscriptstyle-\kern-1pt-\end{math}}\kern.5pt}\providecommand{\PP}{\kern.5pt\raisebox{.4ex}{\begin{math}\scriptscriptstyle+\kern-1pt+\end{math}}\kern.5pt}\def\commopen{/\begin{math}\ast\,\end{math}}\def\commclose{\,\begin{math}\ast\end{math}\kern-.5pt/}\def\begcomm{\begingroup\maybehbox\bgroup\setupmodname}\def\endcomm{\egroup\endgroup}\nwfilename{dualnum.nw}\nwbegindocs{1}\nwdocspar

\nwenddocs{}\nwbegindocs{2}\nwdocspar
\nwenddocs{}\nwbegindocs{3}%
\nwenddocs{}\nwbegindocs{4}%
\nwenddocs{}\nwbegindocs{5}%
\nwenddocs{}\nwbegindocs{6}%
\nwenddocs{}\nwbegindocs{7}%
\nwenddocs{}\nwbegindocs{8}%
\nwenddocs{}\nwbegindocs{9}%
\nwenddocs{}\nwbegindocs{10}%
\nwenddocs{}\nwbegindocs{11}%
\nwenddocs{}\nwbegindocs{12}%
\nwenddocs{}\nwbegindocs{13}%
\nwenddocs{}\nwbegindocs{14}%
\nwenddocs{}\nwbegindocs{15}%
\nwenddocs{}\nwbegindocs{16}%
\nwenddocs{}\nwbegindocs{17}%
\nwenddocs{}\nwbegindocs{18}%
\nwenddocs{}\nwbegindocs{19}%
\nwdocspar
This piece of software is licensed under
\href{http://www.gnu.org/licenses/gpl.html}{GNU General Public
  License} (Version 3, 29 June 2007)~\cite{GNUGPL}.

\nwenddocs{}\nwbegindocs{20}The part related to geodesic spokes is still unfinished.

\nwenddocs{}\nwbegindocs{21}\nwdocspar
\subsection{Class {\Tt{}\Rm{}{\bf{}dual\_number}\nwendquote}}
\label{sec:class-dual_number}

\nwenddocs{}\nwbegindocs{22}\nwdocspar
\subsubsection{Public Methods}
\label{sec:public-methods}

\nwenddocs{}\nwbegindocs{23}\nwdocspar
A dual number can be created simply by listing its two components.
\nwenddocs{}\nwbegincode{24}\sublabel{NW4P20f9-1eKCCy-1}\nwmargintag{{\nwtagstyle{}\subpageref{NW4P20f9-1eKCCy-1}}}\moddef{Public methods~{\nwtagstyle{}\subpageref{NW4P20f9-1eKCCy-1}}}\endmoddef\Rm{}\nwstartdeflinemarkup\nwusesondefline{\\{NW4P20f9-4GrAIY-2}}\nwprevnextdefs{\relax}{NW4P20f9-1eKCCy-2}\nwenddeflinemarkup
        {\bf{}dual\_number}({\bf{}const} {\bf{}ex} & {\it{}a}, {\bf{}const} {\bf{}ex} & {\it{}b});

\nwalsodefined{\\{NW4P20f9-1eKCCy-2}\\{NW4P20f9-1eKCCy-3}\\{NW4P20f9-1eKCCy-4}\\{NW4P20f9-1eKCCy-5}\\{NW4P20f9-1eKCCy-6}\\{NW4P20f9-1eKCCy-7}\\{NW4P20f9-1eKCCy-8}\\{NW4P20f9-1eKCCy-9}\\{NW4P20f9-1eKCCy-A}\\{NW4P20f9-1eKCCy-B}\\{NW4P20f9-1eKCCy-C}}\nwused{\\{NW4P20f9-4GrAIY-2}}\nwidentuses{\\{{\nwixident{dual{\_}number}}{dual:unnumber}}}\nwindexuse{\nwixident{dual{\_}number}}{dual:unnumber}{NW4P20f9-1eKCCy-1}\nwendcode{}\nwbegindocs{25}Alternatively you can provide a \(1\times 2\) or \(2\times 1\)
matrix,  a list, another dual number {\Tt{}\Rm{}{\it{}P}\nwendquote} or a complex expression
with a non-zero imaginary part to give two components. 
If {\Tt{}\Rm{}{\it{}P}\nwendquote} does not have two components and is a real-valued
expression, it will be embedded into dual numbers with zero argument
and norm equal to {\Tt{}\Rm{}{\it{}P}\nwendquote}.
\nwenddocs{}\nwbegincode{26}\sublabel{NW4P20f9-1eKCCy-2}\nwmargintag{{\nwtagstyle{}\subpageref{NW4P20f9-1eKCCy-2}}}\moddef{Public methods~{\nwtagstyle{}\subpageref{NW4P20f9-1eKCCy-1}}}\plusendmoddef\Rm{}\nwstartdeflinemarkup\nwusesondefline{\\{NW4P20f9-4GrAIY-2}}\nwprevnextdefs{NW4P20f9-1eKCCy-1}{NW4P20f9-1eKCCy-3}\nwenddeflinemarkup
        {\bf{}dual\_number}({\bf{}const} {\bf{}ex} & {\it{}P});

\nwused{\\{NW4P20f9-4GrAIY-2}}\nwidentuses{\\{{\nwixident{dual{\_}number}}{dual:unnumber}}}\nwindexuse{\nwixident{dual{\_}number}}{dual:unnumber}{NW4P20f9-1eKCCy-2}\nwendcode{}\nwbegindocs{27} We can also obtain the module and argument of a {\Tt{}\Rm{}{\bf{}dual\_number}\nwendquote}.
\nwenddocs{}\nwbegincode{28}\sublabel{NW4P20f9-1eKCCy-3}\nwmargintag{{\nwtagstyle{}\subpageref{NW4P20f9-1eKCCy-3}}}\moddef{Public methods~{\nwtagstyle{}\subpageref{NW4P20f9-1eKCCy-1}}}\plusendmoddef\Rm{}\nwstartdeflinemarkup\nwusesondefline{\\{NW4P20f9-4GrAIY-2}}\nwprevnextdefs{NW4P20f9-1eKCCy-2}{NW4P20f9-1eKCCy-4}\nwenddeflinemarkup
        {\bf{}ex} {\it{}arg}() {\bf{}const};
        {\bf{}ex} {\it{}norm}() {\bf{}const};

\nwused{\\{NW4P20f9-4GrAIY-2}}\nwendcode{}\nwbegindocs{29}We define the conjugate of a {\Tt{}\Rm{}{\bf{}dual\_number}\nwendquote} as number
with the same modulus but negative value of the argument.
\nwenddocs{}\nwbegincode{30}\sublabel{NW4P20f9-1eKCCy-4}\nwmargintag{{\nwtagstyle{}\subpageref{NW4P20f9-1eKCCy-4}}}\moddef{Public methods~{\nwtagstyle{}\subpageref{NW4P20f9-1eKCCy-1}}}\plusendmoddef\Rm{}\nwstartdeflinemarkup\nwusesondefline{\\{NW4P20f9-4GrAIY-2}}\nwprevnextdefs{NW4P20f9-1eKCCy-3}{NW4P20f9-1eKCCy-5}\nwenddeflinemarkup
        {\bf{}ex} {\it{}conjugate}() {\bf{}const};

\nwused{\\{NW4P20f9-4GrAIY-2}}\nwendcode{}\nwbegindocs{31} Real part may be defined out of the formula
\(\Re(u,v)=\frac{1}{2}((u,v)+\overline{(u,v)})\), but it simply
reduces to the value of norm for the {\Tt{}\Rm{}{\bf{}dual\_number}\nwendquote}.
\nwenddocs{}\nwbegincode{32}\sublabel{NW4P20f9-1eKCCy-5}\nwmargintag{{\nwtagstyle{}\subpageref{NW4P20f9-1eKCCy-5}}}\moddef{Public methods~{\nwtagstyle{}\subpageref{NW4P20f9-1eKCCy-1}}}\plusendmoddef\Rm{}\nwstartdeflinemarkup\nwusesondefline{\\{NW4P20f9-4GrAIY-2}}\nwprevnextdefs{NW4P20f9-1eKCCy-4}{NW4P20f9-1eKCCy-6}\nwenddeflinemarkup
        {\bf{}ex} {\it{}real\_part}() {\bf{}const};
        {\bf{}ex} {\it{}imag\_part}() {\bf{}const}; 

\nwused{\\{NW4P20f9-4GrAIY-2}}\nwendcode{}\nwbegindocs{33}Negative of a {\Tt{}\Rm{}{\bf{}dual\_number}\nwendquote} and its power.
\nwenddocs{}\nwbegincode{34}\sublabel{NW4P20f9-1eKCCy-6}\nwmargintag{{\nwtagstyle{}\subpageref{NW4P20f9-1eKCCy-6}}}\moddef{Public methods~{\nwtagstyle{}\subpageref{NW4P20f9-1eKCCy-1}}}\plusendmoddef\Rm{}\nwstartdeflinemarkup\nwusesondefline{\\{NW4P20f9-4GrAIY-2}}\nwprevnextdefs{NW4P20f9-1eKCCy-5}{NW4P20f9-1eKCCy-7}\nwenddeflinemarkup
        {\bf{}dual\_number} {\it{}neg}() {\bf{}const};
        {\bf{}dual\_number} {\it{}power}({\bf{}const} {\bf{}ex} & {\it{}e}) {\bf{}const};

\nwused{\\{NW4P20f9-4GrAIY-2}}\nwidentuses{\\{{\nwixident{dual{\_}number}}{dual:unnumber}}}\nwindexuse{\nwixident{dual{\_}number}}{dual:unnumber}{NW4P20f9-1eKCCy-6}\nwendcode{}\nwbegindocs{35} We can also convert a {\Tt{}\Rm{}{\bf{}dual\_number}\nwendquote} to a {\Tt{}\Rm{}{\bf{}matrix}\nwendquote}.
\nwenddocs{}\nwbegincode{36}\sublabel{NW4P20f9-1eKCCy-7}\nwmargintag{{\nwtagstyle{}\subpageref{NW4P20f9-1eKCCy-7}}}\moddef{Public methods~{\nwtagstyle{}\subpageref{NW4P20f9-1eKCCy-1}}}\plusendmoddef\Rm{}\nwstartdeflinemarkup\nwusesondefline{\\{NW4P20f9-4GrAIY-2}}\nwprevnextdefs{NW4P20f9-1eKCCy-6}{NW4P20f9-1eKCCy-8}\nwenddeflinemarkup
        {\bf{}matrix} {\it{}to\_matrix}() {\bf{}const} {\nwlbrace} {\bf{}return} {\bf{}matrix} (1, 2, {\bf{}lst}({\it{}u\_comp}, {\it{}v\_comp})); {\nwrbrace}

\nwused{\\{NW4P20f9-4GrAIY-2}}\nwendcode{}\nwbegindocs{37}We define the rule for parabolic norm of a sum, see~\eqref{eq:p-add-norm-exotic}.
\nwenddocs{}\nwbegincode{38}\sublabel{NW4P20f9-1eKCCy-8}\nwmargintag{{\nwtagstyle{}\subpageref{NW4P20f9-1eKCCy-8}}}\moddef{Public methods~{\nwtagstyle{}\subpageref{NW4P20f9-1eKCCy-1}}}\plusendmoddef\Rm{}\nwstartdeflinemarkup\nwusesondefline{\\{NW4P20f9-4GrAIY-2}}\nwprevnextdefs{NW4P20f9-1eKCCy-7}{NW4P20f9-1eKCCy-9}\nwenddeflinemarkup
        {\bf{}ex} {\it{}add\_norms}({\bf{}const} {\bf{}dual\_number} & {\it{}P}) {\bf{}const};

\nwused{\\{NW4P20f9-4GrAIY-2}}\nwidentuses{\\{{\nwixident{dual{\_}number}}{dual:unnumber}}}\nwindexuse{\nwixident{dual{\_}number}}{dual:unnumber}{NW4P20f9-1eKCCy-8}\nwendcode{}\nwbegindocs{39}Upper half-plane action of the group \(\SL\) on the dual numbers.
\nwenddocs{}\nwbegincode{40}\sublabel{NW4P20f9-1eKCCy-9}\nwmargintag{{\nwtagstyle{}\subpageref{NW4P20f9-1eKCCy-9}}}\moddef{Public methods~{\nwtagstyle{}\subpageref{NW4P20f9-1eKCCy-1}}}\plusendmoddef\Rm{}\nwstartdeflinemarkup\nwusesondefline{\\{NW4P20f9-4GrAIY-2}}\nwprevnextdefs{NW4P20f9-1eKCCy-8}{NW4P20f9-1eKCCy-A}\nwenddeflinemarkup
        {\bf{}dual\_number} {\it{}sl2\_moebius}({\bf{}const} {\bf{}ex} & {\it{}a}, {\bf{}const} {\bf{}ex} & {\it{}b}, {\bf{}const} {\bf{}ex} & {\it{}c}, {\bf{}const} {\bf{}ex} & {\it{}d}) {\bf{}const};
        {\bf{}dual\_number} {\it{}sl2\_moebius}({\bf{}const} {\bf{}ex} & {\it{}M}) {\bf{}const};

\nwused{\\{NW4P20f9-4GrAIY-2}}\nwidentuses{\\{{\nwixident{dual{\_}number}}{dual:unnumber}}}\nwindexuse{\nwixident{dual{\_}number}}{dual:unnumber}{NW4P20f9-1eKCCy-9}\nwendcode{}\nwbegindocs{41}Similar action in the form of the group \(\mathrm{SU}\) on the unit disk.
\nwenddocs{}\nwbegincode{42}\sublabel{NW4P20f9-1eKCCy-A}\nwmargintag{{\nwtagstyle{}\subpageref{NW4P20f9-1eKCCy-A}}}\moddef{Public methods~{\nwtagstyle{}\subpageref{NW4P20f9-1eKCCy-1}}}\plusendmoddef\Rm{}\nwstartdeflinemarkup\nwusesondefline{\\{NW4P20f9-4GrAIY-2}}\nwprevnextdefs{NW4P20f9-1eKCCy-9}{NW4P20f9-1eKCCy-B}\nwenddeflinemarkup
        {\bf{}dual\_number} {\it{}su\_moebius}({\bf{}const} {\bf{}ex} & {\it{}a}, {\bf{}const} {\bf{}ex} & {\it{}b}, {\bf{}const} {\bf{}ex} & {\it{}c}, {\bf{}const} {\bf{}ex} & {\it{}d}) {\bf{}const};
        {\bf{}dual\_number} {\it{}su\_moebius}({\bf{}const} {\bf{}ex} & {\it{}M}) {\bf{}const};

\nwused{\\{NW4P20f9-4GrAIY-2}}\nwidentuses{\\{{\nwixident{dual{\_}number}}{dual:unnumber}}}\nwindexuse{\nwixident{dual{\_}number}}{dual:unnumber}{NW4P20f9-1eKCCy-A}\nwendcode{}\nwbegindocs{43}A particular, but important case of the above
\(\mathrm{SU}\)-action: unitary rotation by \(x\) for
subgroups \(N\) and \(N^\prime\).
\nwenddocs{}\nwbegincode{44}\sublabel{NW4P20f9-1eKCCy-B}\nwmargintag{{\nwtagstyle{}\subpageref{NW4P20f9-1eKCCy-B}}}\moddef{Public methods~{\nwtagstyle{}\subpageref{NW4P20f9-1eKCCy-1}}}\plusendmoddef\Rm{}\nwstartdeflinemarkup\nwusesondefline{\\{NW4P20f9-4GrAIY-2}}\nwprevnextdefs{NW4P20f9-1eKCCy-A}{NW4P20f9-1eKCCy-C}\nwenddeflinemarkup
        {\bf{}dual\_number} {\it{}rot}({\bf{}const} {\bf{}ex} & {\it{}x}) {\bf{}const};

\nwused{\\{NW4P20f9-4GrAIY-2}}\nwidentuses{\\{{\nwixident{dual{\_}number}}{dual:unnumber}}}\nwindexuse{\nwixident{dual{\_}number}}{dual:unnumber}{NW4P20f9-1eKCCy-B}\nwendcode{}\nwbegindocs{45}Algebraic operations are defined for {\Tt{}\Rm{}{\bf{}dual\_number}\nwendquote}s in a way
described in \S~\ref{sec:invar-line-algebra}. The standard \CPP\ operators
{\Tt{}\Rm{}+\nwendquote}, {\Tt{}\Rm{}-\nwendquote}, {\Tt{}\Rm{}\begin{math}\ast\end{math}\nwendquote}, {\Tt{}\Rm{}\begin{math}\div\end{math}\nwendquote} will be overloaded later in order to permit
natural expressions with {\Tt{}\Rm{}{\bf{}dual\_number}\nwendquote}s.
\nwenddocs{}\nwbegincode{46}\sublabel{NW4P20f9-1eKCCy-C}\nwmargintag{{\nwtagstyle{}\subpageref{NW4P20f9-1eKCCy-C}}}\moddef{Public methods~{\nwtagstyle{}\subpageref{NW4P20f9-1eKCCy-1}}}\plusendmoddef\Rm{}\nwstartdeflinemarkup\nwusesondefline{\\{NW4P20f9-4GrAIY-2}}\nwprevnextdefs{NW4P20f9-1eKCCy-B}{\relax}\nwenddeflinemarkup
        {\bf{}dual\_number} {\it{}add}({\bf{}const} {\bf{}dual\_number} & {\it{}P}) {\bf{}const};
        {\bf{}dual\_number} {\it{}sub}({\bf{}const} {\bf{}dual\_number} & {\it{}P}) {\bf{}const} {\nwlbrace} {\bf{}return} {\it{}add}({\it{}P}.{\it{}neg}()); {\nwrbrace}
        {\bf{}dual\_number} {\it{}mul}({\bf{}const} {\bf{}dual\_number} & {\it{}P}) {\bf{}const};
        {\bf{}dual\_number} {\it{}mul}({\bf{}const} {\bf{}ex} & {\it{}a}) {\bf{}const} {\nwlbrace} {\bf{}return} {\it{}mul}({\bf{}dual\_number}({\it{}a})); {\nwrbrace}
        {\bf{}dual\_number} {\it{}div}({\bf{}const} {\bf{}dual\_number} & {\it{}P}) {\bf{}const} {\nwlbrace}{\bf{}return} {\it{}mul}({\it{}P}.{\it{}power}(-1)); {\nwrbrace}
        {\bf{}dual\_number} {\it{}div}({\bf{}const} {\bf{}ex} & {\it{}a}) {\bf{}const} {\nwlbrace} {\bf{}return} {\it{}div}({\bf{}dual\_number}({\it{}a})); {\nwrbrace}

\nwused{\\{NW4P20f9-4GrAIY-2}}\nwidentuses{\\{{\nwixident{dual{\_}number}}{dual:unnumber}}}\nwindexuse{\nwixident{dual{\_}number}}{dual:unnumber}{NW4P20f9-1eKCCy-C}\nwendcode{}\nwbegindocs{47}We also define zero dual number.
\nwenddocs{}\nwbegincode{48}\sublabel{NW4P20f9-NQ5ML-1}\nwmargintag{{\nwtagstyle{}\subpageref{NW4P20f9-NQ5ML-1}}}\moddef{Additional routines declarations~{\nwtagstyle{}\subpageref{NW4P20f9-NQ5ML-1}}}\endmoddef\Rm{}\nwstartdeflinemarkup\nwusesondefline{\\{NW4P20f9-36Ytqo-1}}\nwprevnextdefs{\relax}{NW4P20f9-NQ5ML-2}\nwenddeflinemarkup
        {\bf{}dual\_number} {\it{}zero\_dual\_number}();
        {\bf{}dual\_number} {\it{}dn\_from\_arg\_mod}({\bf{}const} {\bf{}ex} & {\it{}a}, {\bf{}const} {\bf{}ex} & {\it{}n});

\nwalsodefined{\\{NW4P20f9-NQ5ML-2}\\{NW4P20f9-NQ5ML-3}\\{NW4P20f9-NQ5ML-4}\\{NW4P20f9-NQ5ML-5}}\nwused{\\{NW4P20f9-36Ytqo-1}}\nwidentuses{\\{{\nwixident{dual{\_}number}}{dual:unnumber}}}\nwindexuse{\nwixident{dual{\_}number}}{dual:unnumber}{NW4P20f9-NQ5ML-1}\nwendcode{}\nwbegindocs{49}\nwdocspar
\subsection{Algebraic Subroutines}
\label{sec:algebr-subr}
We need a couple of global variables which help to write uniformly
algebraic rules for both cases of subgroups \(N\) and \(N^\prime\).

\nwenddocs{}\nwbegindocs{50}Firstly, we need to consider separately cases of subgroup \(N\) and
\(N^\prime\), the following global variable keeps track on it.
\nwenddocs{}\nwbegincode{51}\sublabel{NW4P20f9-4EWpKB-1}\nwmargintag{{\nwtagstyle{}\subpageref{NW4P20f9-4EWpKB-1}}}\moddef{N-Nprime separation~{\nwtagstyle{}\subpageref{NW4P20f9-4EWpKB-1}}}\endmoddef\Rm{}\nwstartdeflinemarkup\nwusesondefline{\\{NW4P20f9-36Ytqo-1}}\nwprevnextdefs{\relax}{NW4P20f9-4EWpKB-2}\nwenddeflinemarkup
{\bf{}bool} {\it{}subgroup\_N}={\bf{}true};
{\bf{}bool} {\it{}straight\_spoke}={\bf{}true};

\nwalsodefined{\\{NW4P20f9-4EWpKB-2}\\{NW4P20f9-4EWpKB-3}}\nwused{\\{NW4P20f9-36Ytqo-1}}\nwendcode{}\nwbegindocs{52}In the case of the subgroup \(N^\prime\) the reference point lies at infinity, the
following {\Tt{}\Rm{}{\bf{}realsymbol}\nwendquote} variable represents it in the symbolic
calculations.
\nwenddocs{}\nwbegincode{53}\sublabel{NW4P20f9-4EWpKB-2}\nwmargintag{{\nwtagstyle{}\subpageref{NW4P20f9-4EWpKB-2}}}\moddef{N-Nprime separation~{\nwtagstyle{}\subpageref{NW4P20f9-4EWpKB-1}}}\plusendmoddef\Rm{}\nwstartdeflinemarkup\nwusesondefline{\\{NW4P20f9-36Ytqo-1}}\nwprevnextdefs{NW4P20f9-4EWpKB-1}{NW4P20f9-4EWpKB-3}\nwenddeflinemarkup
{\bf{}realsymbol} {\it{}Inf}({\tt{}"Inf"}, {\tt{}"{\char92}{\char92}infty"});

\nwused{\\{NW4P20f9-36Ytqo-1}}\nwendcode{}\nwbegindocs{54} We define the ``zero angle'', it is always zero.
\nwenddocs{}\nwbegincode{55}\sublabel{NW4P20f9-4EWpKB-3}\nwmargintag{{\nwtagstyle{}\subpageref{NW4P20f9-4EWpKB-3}}}\moddef{N-Nprime separation~{\nwtagstyle{}\subpageref{NW4P20f9-4EWpKB-1}}}\plusendmoddef\Rm{}\nwstartdeflinemarkup\nwusesondefline{\\{NW4P20f9-36Ytqo-1}}\nwprevnextdefs{NW4P20f9-4EWpKB-2}{\relax}\nwenddeflinemarkup
{\bf{}\char35{}define}{\tt{} Arg0 0}\nwindexdefn{\nwixident{Arg0}}{Arg0}{NW4P20f9-4EWpKB-3}

\nwused{\\{NW4P20f9-36Ytqo-1}}\nwidentdefs{\\{{\nwixident{Arg0}}{Arg0}}}\nwendcode{}\nwbegindocs{56}Here is the set of algebraic procedures representing definitions
made in this paper.

\nwenddocs{}\nwbegindocs{57}\nwdocspar
\subsubsection[Argument and u]{Argument and $u$}
\label{sec:argument-u}

\nwenddocs{}\nwbegindocs{58}We define the conjugate of a {\Tt{}\Rm{}{\bf{}dual\_number}\nwendquote} as number
with the same modulus but negative value of the argument.
\nwenddocs{}\nwbegincode{59}\sublabel{NW4P20f9-3jzoaG-1}\nwmargintag{{\nwtagstyle{}\subpageref{NW4P20f9-3jzoaG-1}}}\moddef{Algebraic procedures~{\nwtagstyle{}\subpageref{NW4P20f9-3jzoaG-1}}}\endmoddef\Rm{}\nwstartdeflinemarkup\nwusesondefline{\\{NW4P20f9-3HhVai-1}}\nwprevnextdefs{\relax}{NW4P20f9-3jzoaG-2}\nwenddeflinemarkup
{\bf{}ex} {\bf{}dual\_number}::{\it{}conjugate}() {\bf{}const} {\nwlbrace} {\bf{}return} {\it{}dn\_from\_arg\_mod}(-{\it{}arg}(), {\it{}norm}()); {\nwrbrace}

\nwalsodefined{\\{NW4P20f9-3jzoaG-2}\\{NW4P20f9-3jzoaG-3}\\{NW4P20f9-3jzoaG-4}\\{NW4P20f9-3jzoaG-5}\\{NW4P20f9-3jzoaG-6}}\nwused{\\{NW4P20f9-3HhVai-1}}\nwidentuses{\\{{\nwixident{dual{\_}number}}{dual:unnumber}}}\nwindexuse{\nwixident{dual{\_}number}}{dual:unnumber}{NW4P20f9-3jzoaG-1}\nwendcode{}\nwbegindocs{60}\nwdocspar
\subsubsection{Argument of a point}
\label{sec:argument-point}
The opposite task (finding argument of a point) is solved
similarly (a direct formula).
\nwenddocs{}\nwbegincode{61}\sublabel{NW4P20f9-4dxCds-1}\nwmargintag{{\nwtagstyle{}\subpageref{NW4P20f9-4dxCds-1}}}\moddef{Dual number class further implementation~{\nwtagstyle{}\subpageref{NW4P20f9-4dxCds-1}}}\endmoddef\Rm{}\nwstartdeflinemarkup\nwusesondefline{\\{NW4P20f9-3HhVai-1}}\nwprevnextdefs{\relax}{NW4P20f9-4dxCds-2}\nwenddeflinemarkup
{\bf{}ex} {\bf{}dual\_number}::{\it{}arg}() {\bf{}const} {\nwlbrace}
        {\bf{}if} ({\it{}subgroup\_N})
                {\bf{}return} {\it{}u\_comp};

\nwalsodefined{\\{NW4P20f9-4dxCds-2}\\{NW4P20f9-4dxCds-3}\\{NW4P20f9-4dxCds-4}\\{NW4P20f9-4dxCds-5}\\{NW4P20f9-4dxCds-6}\\{NW4P20f9-4dxCds-7}\\{NW4P20f9-4dxCds-8}\\{NW4P20f9-4dxCds-9}\\{NW4P20f9-4dxCds-A}\\{NW4P20f9-4dxCds-B}\\{NW4P20f9-4dxCds-C}\\{NW4P20f9-4dxCds-D}\\{NW4P20f9-4dxCds-E}\\{NW4P20f9-4dxCds-F}\\{NW4P20f9-4dxCds-G}\\{NW4P20f9-4dxCds-H}\\{NW4P20f9-4dxCds-I}\\{NW4P20f9-4dxCds-J}\\{NW4P20f9-4dxCds-K}\\{NW4P20f9-4dxCds-L}\\{NW4P20f9-4dxCds-M}}\nwused{\\{NW4P20f9-3HhVai-1}}\nwidentuses{\\{{\nwixident{dual{\_}number}}{dual:unnumber}}}\nwindexuse{\nwixident{dual{\_}number}}{dual:unnumber}{NW4P20f9-4dxCds-1}\nwendcode{}\nwbegindocs{62}Again in the case of \(N^\prime\) we need to consider cases of \(0\)\ldots
\nwenddocs{}\nwbegincode{63}\sublabel{NW4P20f9-4dxCds-2}\nwmargintag{{\nwtagstyle{}\subpageref{NW4P20f9-4dxCds-2}}}\moddef{Dual number class further implementation~{\nwtagstyle{}\subpageref{NW4P20f9-4dxCds-1}}}\plusendmoddef\Rm{}\nwstartdeflinemarkup\nwusesondefline{\\{NW4P20f9-3HhVai-1}}\nwprevnextdefs{NW4P20f9-4dxCds-1}{NW4P20f9-4dxCds-3}\nwenddeflinemarkup
        {\bf{}else} {\nwlbrace}
                {\bf{}if} ({\it{}straight\_spoke}) {\nwlbrace}
                        {\bf{}if} ({\it{}u\_comp}.{\it{}normal}().{\it{}is\_zero}())
                                {\bf{}return} {\it{}Inf};

\nwused{\\{NW4P20f9-3HhVai-1}}\nwendcode{}\nwbegindocs{64}\ldots and \(\infty\). We try to replace \(\frac{1}{\infty}\) by \(0\).
\nwenddocs{}\nwbegincode{65}\sublabel{NW4P20f9-4dxCds-3}\nwmargintag{{\nwtagstyle{}\subpageref{NW4P20f9-4dxCds-3}}}\moddef{Dual number class further implementation~{\nwtagstyle{}\subpageref{NW4P20f9-4dxCds-1}}}\plusendmoddef\Rm{}\nwstartdeflinemarkup\nwusesondefline{\\{NW4P20f9-3HhVai-1}}\nwprevnextdefs{NW4P20f9-4dxCds-2}{NW4P20f9-4dxCds-4}\nwenddeflinemarkup
                        {\bf{}else} 
                                {\bf{}try} {\nwlbrace}
                                        {\bf{}realsymbol} {\it{}t}({\tt{}"t"});
                                        {\bf{}return} {\it{}pow}({\it{}u\_comp}.{\it{}subs}({\it{}Inf} \begin{math}\equiv\end{math} {\it{}pow}({\it{}t}, -1)), -1).{\it{}normal}().{\it{}subs}({\it{}t} \begin{math}\equiv\end{math} 0).{\it{}normal}();
                                {\nwrbrace} {\bf{}catch} ({\it{}std}::{\it{}exception} &{\it{}p}) {\nwlbrace}
                                        {\bf{}return} {\it{}pow}({\it{}u\_comp}, -1);
                        {\nwrbrace}

\nwused{\\{NW4P20f9-3HhVai-1}}\nwendcode{}\nwbegindocs{66}We now calculate argument for the case geodesic spokes (a direct formula).
\nwenddocs{}\nwbegincode{67}\sublabel{NW4P20f9-4dxCds-4}\nwmargintag{{\nwtagstyle{}\subpageref{NW4P20f9-4dxCds-4}}}\moddef{Dual number class further implementation~{\nwtagstyle{}\subpageref{NW4P20f9-4dxCds-1}}}\plusendmoddef\Rm{}\nwstartdeflinemarkup\nwusesondefline{\\{NW4P20f9-3HhVai-1}}\nwprevnextdefs{NW4P20f9-4dxCds-3}{NW4P20f9-4dxCds-5}\nwenddeflinemarkup
                {\nwrbrace} {\bf{}else} 
                        {\bf{}return} (1-{\it{}pow}(2\begin{math}\ast\end{math}{\it{}v\_comp}+1, {\it{}half}))\begin{math}\div\end{math}{\it{}u\_comp};
        {\nwrbrace}
{\nwrbrace}

\nwused{\\{NW4P20f9-3HhVai-1}}\nwendcode{}\nwbegindocs{68}\nwdocspar

\subsubsection{Norm}
\label{sec:norm}

The corresponding value of the parabolic norm is calculated by the
formulae~\eqref{eq:parab-norm} (a direct formula).
\nwenddocs{}\nwbegincode{69}\sublabel{NW4P20f9-4dxCds-5}\nwmargintag{{\nwtagstyle{}\subpageref{NW4P20f9-4dxCds-5}}}\moddef{Dual number class further implementation~{\nwtagstyle{}\subpageref{NW4P20f9-4dxCds-1}}}\plusendmoddef\Rm{}\nwstartdeflinemarkup\nwusesondefline{\\{NW4P20f9-3HhVai-1}}\nwprevnextdefs{NW4P20f9-4dxCds-4}{NW4P20f9-4dxCds-6}\nwenddeflinemarkup
{\bf{}ex} {\bf{}dual\_number}::{\it{}norm}() {\bf{}const} {\nwlbrace}
        {\bf{}if} ({\it{}subgroup\_N})
                {\bf{}return} {\it{}pow}({\it{}u\_comp}, 2)-{\it{}v\_comp}; 

\nwused{\\{NW4P20f9-3HhVai-1}}\nwidentuses{\\{{\nwixident{dual{\_}number}}{dual:unnumber}}}\nwindexuse{\nwixident{dual{\_}number}}{dual:unnumber}{NW4P20f9-4dxCds-5}\nwendcode{}\nwbegindocs{70}The case of subgroup \(N^\prime\) require treatment of infinity. If
the denominator is zero than the treatment is the same for straight and geodesic spokes.
\nwenddocs{}\nwbegincode{71}\sublabel{NW4P20f9-4dxCds-6}\nwmargintag{{\nwtagstyle{}\subpageref{NW4P20f9-4dxCds-6}}}\moddef{Dual number class further implementation~{\nwtagstyle{}\subpageref{NW4P20f9-4dxCds-1}}}\plusendmoddef\Rm{}\nwstartdeflinemarkup\nwusesondefline{\\{NW4P20f9-3HhVai-1}}\nwprevnextdefs{NW4P20f9-4dxCds-5}{NW4P20f9-4dxCds-7}\nwenddeflinemarkup
        {\bf{}else} {\nwlbrace}
                {\bf{}ex} {\it{}denom}=(({\it{}straight\_spoke}?1:2)\begin{math}\ast\end{math}{\it{}v\_comp}+1).{\it{}normal}();
                {\bf{}if} ({\it{}denom}.{\it{}is\_zero}()) {\nwlbrace}
                        {\bf{}if} ({\it{}u\_comp}.{\it{}is\_zero}())
                                {\bf{}return} 1;
                        {\bf{}else} 
                                {\bf{}return} {\it{}Inf};

\nwused{\\{NW4P20f9-3HhVai-1}}\nwendcode{}\nwbegindocs{72}For a non-zero denominator we return the plain values.
\nwenddocs{}\nwbegincode{73}\sublabel{NW4P20f9-4dxCds-7}\nwmargintag{{\nwtagstyle{}\subpageref{NW4P20f9-4dxCds-7}}}\moddef{Dual number class further implementation~{\nwtagstyle{}\subpageref{NW4P20f9-4dxCds-1}}}\plusendmoddef\Rm{}\nwstartdeflinemarkup\nwusesondefline{\\{NW4P20f9-3HhVai-1}}\nwprevnextdefs{NW4P20f9-4dxCds-6}{NW4P20f9-4dxCds-8}\nwenddeflinemarkup
                {\nwrbrace} {\bf{}else} {\nwlbrace}
                        {\bf{}if} ({\it{}straight\_spoke})
                                {\bf{}return} ({\it{}pow}({\it{}u\_comp},2)\begin{math}\div\end{math}{\it{}denom}).{\it{}normal}();
                        {\bf{}else}
                                {\bf{}return} ({\it{}u\_comp}\begin{math}\div\end{math}{\it{}pow}({\it{}denom},{\it{}half})).{\it{}normal}();
                {\nwrbrace}
        {\nwrbrace}
{\nwrbrace}

\nwused{\\{NW4P20f9-3HhVai-1}}\nwendcode{}\nwbegindocs{74}\nwdocspar
\subsubsection[The value of v from the argument and norm]{The value of
  \(v\) from the argument and norm} 
\label{sec:value-v-from}
Often We need to find values of \(v\) such that for a given
value of the first component \(u\) point \((u, v)\) will have a given norm.
\nwenddocs{}\nwbegincode{75}\sublabel{NW4P20f9-3jzoaG-2}\nwmargintag{{\nwtagstyle{}\subpageref{NW4P20f9-3jzoaG-2}}}\moddef{Algebraic procedures~{\nwtagstyle{}\subpageref{NW4P20f9-3jzoaG-1}}}\plusendmoddef\Rm{}\nwstartdeflinemarkup\nwusesondefline{\\{NW4P20f9-3HhVai-1}}\nwprevnextdefs{NW4P20f9-3jzoaG-1}{NW4P20f9-3jzoaG-3}\nwenddeflinemarkup
{\bf{}ex} {\it{}v\_from\_u\_and\_norm}({\bf{}const} {\bf{}ex} & {\it{}u}, {\bf{}const} {\bf{}ex} & {\it{}n}) {\nwlbrace}
        {\bf{}realsymbol} {\it{}l}({\tt{}"l"});
        {\bf{}if} ({\it{}subgroup\_N})
                {\bf{}return} {\it{}lsolve}({\bf{}dual\_number}({\it{}u}, {\it{}l}).{\it{}norm}() \begin{math}\equiv\end{math} {\it{}n}, {\it{}l}).{\it{}normal}();

\nwused{\\{NW4P20f9-3HhVai-1}}\nwidentuses{\\{{\nwixident{dual{\_}number}}{dual:unnumber}}}\nwindexuse{\nwixident{dual{\_}number}}{dual:unnumber}{NW4P20f9-3jzoaG-2}\nwendcode{}\nwbegindocs{76}For subgroup \(N'\) we use powers to make equations linear.
\nwenddocs{}\nwbegincode{77}\sublabel{NW4P20f9-3jzoaG-3}\nwmargintag{{\nwtagstyle{}\subpageref{NW4P20f9-3jzoaG-3}}}\moddef{Algebraic procedures~{\nwtagstyle{}\subpageref{NW4P20f9-3jzoaG-1}}}\plusendmoddef\Rm{}\nwstartdeflinemarkup\nwusesondefline{\\{NW4P20f9-3HhVai-1}}\nwprevnextdefs{NW4P20f9-3jzoaG-2}{NW4P20f9-3jzoaG-4}\nwenddeflinemarkup
        {\bf{}else} {\nwlbrace}
                {\bf{}if} ({\it{}straight\_spoke}) 
                        {\bf{}return} {\it{}lsolve}({\it{}pow}({\bf{}dual\_number}({\it{}u}, {\it{}l}).{\it{}norm}(), -1) \begin{math}\equiv\end{math} {\it{}pow}({\it{}n}, -1), {\it{}l}).{\it{}normal}();
                {\bf{}else}
                        {\bf{}return} {\it{}lsolve}({\it{}pow}({\bf{}dual\_number}({\it{}u}, {\it{}l}).{\it{}norm}(), -2).{\it{}normal}() \begin{math}\equiv\end{math} {\it{}pow}({\it{}n}, -2), {\it{}l}).{\it{}normal}();
        {\nwrbrace}
{\nwrbrace}

\nwused{\\{NW4P20f9-3HhVai-1}}\nwidentuses{\\{{\nwixident{dual{\_}number}}{dual:unnumber}}}\nwindexuse{\nwixident{dual{\_}number}}{dual:unnumber}{NW4P20f9-3jzoaG-3}\nwendcode{}\nwbegindocs{78}\nwdocspar
\nwenddocs{}\nwbegincode{79}\sublabel{NW4P20f9-3jzoaG-4}\nwmargintag{{\nwtagstyle{}\subpageref{NW4P20f9-3jzoaG-4}}}\moddef{Algebraic procedures~{\nwtagstyle{}\subpageref{NW4P20f9-3jzoaG-1}}}\plusendmoddef\Rm{}\nwstartdeflinemarkup\nwusesondefline{\\{NW4P20f9-3HhVai-1}}\nwprevnextdefs{NW4P20f9-3jzoaG-3}{NW4P20f9-3jzoaG-5}\nwenddeflinemarkup
{\bf{}dual\_number} {\it{}zero\_dual\_number}() {\nwlbrace}
        {\bf{}return} ({\it{}subgroup\_N} ? {\bf{}dual\_number}(0, 0) : {\bf{}dual\_number}({\it{}Inf}, -1));
{\nwrbrace}

\nwused{\\{NW4P20f9-3HhVai-1}}\nwidentuses{\\{{\nwixident{dual{\_}number}}{dual:unnumber}}}\nwindexuse{\nwixident{dual{\_}number}}{dual:unnumber}{NW4P20f9-3jzoaG-4}\nwendcode{}\nwbegindocs{80}\nwdocspar
\nwenddocs{}\nwbegincode{81}\sublabel{NW4P20f9-4dxCds-8}\nwmargintag{{\nwtagstyle{}\subpageref{NW4P20f9-4dxCds-8}}}\moddef{Dual number class further implementation~{\nwtagstyle{}\subpageref{NW4P20f9-4dxCds-1}}}\plusendmoddef\Rm{}\nwstartdeflinemarkup\nwusesondefline{\\{NW4P20f9-3HhVai-1}}\nwprevnextdefs{NW4P20f9-4dxCds-7}{NW4P20f9-4dxCds-9}\nwenddeflinemarkup
{\bf{}bool} {\bf{}dual\_number}::{\it{}is\_zero}() {\bf{}const} {\nwlbrace}
        {\bf{}return} {\it{}is\_equal}({\it{}zero\_dual\_number}());
{\nwrbrace}

\nwused{\\{NW4P20f9-3HhVai-1}}\nwidentuses{\\{{\nwixident{dual{\_}number}}{dual:unnumber}}}\nwindexuse{\nwixident{dual{\_}number}}{dual:unnumber}{NW4P20f9-4dxCds-8}\nwendcode{}\nwbegindocs{82}\nwdocspar
\subsubsection{Real and Imaginary Parts}
\label{sec:real-imag-parts}

See \cite{Kisil07a} for a discussion of the real and
imaginary parts of dual numbers.
\nwenddocs{}\nwbegincode{83}\sublabel{NW4P20f9-3jzoaG-5}\nwmargintag{{\nwtagstyle{}\subpageref{NW4P20f9-3jzoaG-5}}}\moddef{Algebraic procedures~{\nwtagstyle{}\subpageref{NW4P20f9-3jzoaG-1}}}\plusendmoddef\Rm{}\nwstartdeflinemarkup\nwusesondefline{\\{NW4P20f9-3HhVai-1}}\nwprevnextdefs{NW4P20f9-3jzoaG-4}{NW4P20f9-3jzoaG-6}\nwenddeflinemarkup
{\bf{}ex} {\bf{}dual\_number}::{\it{}real\_part}() {\bf{}const} {\nwlbrace}
        {\bf{}return} {\it{}dn\_from\_arg\_mod}({\it{}Arg0}, (1-{\it{}arg}())\begin{math}\ast\end{math}{\it{}norm}());
{\nwrbrace}

{\bf{}ex} {\bf{}dual\_number}::{\it{}imag\_part}() {\bf{}const} {\nwlbrace}
        {\bf{}return} {\it{}dn\_from\_arg\_mod}(1, {\it{}arg}()\begin{math}\ast\end{math}{\it{}norm}());
{\nwrbrace}

\nwused{\\{NW4P20f9-3HhVai-1}}\nwidentuses{\\{{\nwixident{Arg0}}{Arg0}}\\{{\nwixident{dual{\_}number}}{dual:unnumber}}}\nwindexuse{\nwixident{Arg0}}{Arg0}{NW4P20f9-3jzoaG-5}\nwindexuse{\nwixident{dual{\_}number}}{dual:unnumber}{NW4P20f9-3jzoaG-5}\nwendcode{}\nwbegindocs{84}\nwdocspar
\subsubsection{Product of Two Points}
\label{sec:product-two-points}
We define now the product of two points according to the
Definition~\ref{de:product}. We also include a multiplication by a
scalar: if a factor is a scalar it is replaced by a vector with the
zero argument and norm equal to the scalar.
\nwenddocs{}\nwbegincode{85}\sublabel{NW4P20f9-4dxCds-9}\nwmargintag{{\nwtagstyle{}\subpageref{NW4P20f9-4dxCds-9}}}\moddef{Dual number class further implementation~{\nwtagstyle{}\subpageref{NW4P20f9-4dxCds-1}}}\plusendmoddef\Rm{}\nwstartdeflinemarkup\nwusesondefline{\\{NW4P20f9-3HhVai-1}}\nwprevnextdefs{NW4P20f9-4dxCds-8}{NW4P20f9-4dxCds-A}\nwenddeflinemarkup
{\bf{}dual\_number} {\bf{}dual\_number}::{\it{}mul}({\bf{}const} {\bf{}dual\_number} & {\it{}P}) {\bf{}const} {\nwlbrace}
        {\bf{}return} {\it{}dn\_from\_arg\_mod}({\it{}arg}()+{\it{}P}.{\it{}arg}(), {\it{}norm}()\begin{math}\ast\end{math}{\it{}P}.{\it{}norm}());
{\nwrbrace}

\nwused{\\{NW4P20f9-3HhVai-1}}\nwidentuses{\\{{\nwixident{dual{\_}number}}{dual:unnumber}}}\nwindexuse{\nwixident{dual{\_}number}}{dual:unnumber}{NW4P20f9-4dxCds-9}\nwendcode{}\nwbegindocs{86}\nwdocspar
\subsubsection{Vector Addition of Two Points}
\label{sec:vector-addition-two}

\nwenddocs{}\nwbegindocs{87}The rule for norm of a sum of two vectors
\nwenddocs{}\nwbegincode{88}\sublabel{NW4P20f9-4dxCds-A}\nwmargintag{{\nwtagstyle{}\subpageref{NW4P20f9-4dxCds-A}}}\moddef{Dual number class further implementation~{\nwtagstyle{}\subpageref{NW4P20f9-4dxCds-1}}}\plusendmoddef\Rm{}\nwstartdeflinemarkup\nwusesondefline{\\{NW4P20f9-3HhVai-1}}\nwprevnextdefs{NW4P20f9-4dxCds-9}{NW4P20f9-4dxCds-B}\nwenddeflinemarkup
{\bf{}ex} {\bf{}dual\_number}::{\it{}add\_norms}({\bf{}const} {\bf{}dual\_number} & {\it{}P}) {\bf{}const} {\nwlbrace} 
                {\bf{}return} ({\it{}norm}()+{\it{}P}.{\it{}norm}()).{\it{}normal}(); 
{\nwrbrace}

\nwused{\\{NW4P20f9-3HhVai-1}}\nwidentuses{\\{{\nwixident{dual{\_}number}}{dual:unnumber}}}\nwindexuse{\nwixident{dual{\_}number}}{dual:unnumber}{NW4P20f9-4dxCds-A}\nwendcode{}\nwbegindocs{89} The sum is calculated from the expression~\eqref{eq:p-add-arg-exotic}. 
\nwenddocs{}\nwbegincode{90}\sublabel{NW4P20f9-4dxCds-B}\nwmargintag{{\nwtagstyle{}\subpageref{NW4P20f9-4dxCds-B}}}\moddef{Dual number class further implementation~{\nwtagstyle{}\subpageref{NW4P20f9-4dxCds-1}}}\plusendmoddef\Rm{}\nwstartdeflinemarkup\nwusesondefline{\\{NW4P20f9-3HhVai-1}}\nwprevnextdefs{NW4P20f9-4dxCds-A}{NW4P20f9-4dxCds-C}\nwenddeflinemarkup
{\bf{}dual\_number} {\bf{}dual\_number}::{\it{}add}({\bf{}const} {\bf{}dual\_number} & {\it{}P}) {\bf{}const} {\nwlbrace}
        {\bf{}ex} {\it{}norms} = {\it{}add\_norms}({\it{}P});
        {\bf{}if} ({\it{}norms}.{\it{}normal}().{\it{}is\_zero}())
                {\bf{}return} {\it{}zero\_dual\_number}();
        {\bf{}else}
                {\bf{}return} {\it{}dn\_from\_arg\_mod}(({\it{}arg}()\begin{math}\ast\end{math}{\it{}norm}()+{\it{}P}.{\it{}arg}()\begin{math}\ast\end{math}{\it{}P}.{\it{}norm}())\begin{math}\div\end{math}{\it{}norms}, {\it{}norms});
{\nwrbrace}

\nwused{\\{NW4P20f9-3HhVai-1}}\nwidentuses{\\{{\nwixident{dual{\_}number}}{dual:unnumber}}}\nwindexuse{\nwixident{dual{\_}number}}{dual:unnumber}{NW4P20f9-4dxCds-B}\nwendcode{}\nwbegindocs{91}To construct a dual number from its argument and norm, we first
build a number with the given norm and then rotate it to obtain the
desired argument.
\nwenddocs{}\nwbegincode{92}\sublabel{NW4P20f9-3jzoaG-6}\nwmargintag{{\nwtagstyle{}\subpageref{NW4P20f9-3jzoaG-6}}}\moddef{Algebraic procedures~{\nwtagstyle{}\subpageref{NW4P20f9-3jzoaG-1}}}\plusendmoddef\Rm{}\nwstartdeflinemarkup\nwusesondefline{\\{NW4P20f9-3HhVai-1}}\nwprevnextdefs{NW4P20f9-3jzoaG-5}{\relax}\nwenddeflinemarkup
{\bf{}dual\_number} {\it{}dn\_from\_arg\_mod}({\bf{}const} {\bf{}ex} & {\it{}a}, {\bf{}const} {\bf{}ex} & {\it{}n}) {\nwlbrace}
        {\bf{}dual\_number} {\it{}P}={\bf{}dual\_number}(1, {\it{}v\_from\_u\_and\_norm}(1, {\it{}n}));
        {\bf{}ex} {\it{}angle}=({\it{}a}-{\it{}P}.{\it{}arg}()).{\it{}normal}();
        {\bf{}if} ({\it{}angle}.{\it{}is\_zero}())
                {\bf{}return} {\it{}P};
        {\bf{}else}
                {\bf{}return} {\it{}P}.{\it{}rot}({\it{}angle});
{\nwrbrace}

\nwused{\\{NW4P20f9-3HhVai-1}}\nwidentuses{\\{{\nwixident{dual{\_}number}}{dual:unnumber}}}\nwindexuse{\nwixident{dual{\_}number}}{dual:unnumber}{NW4P20f9-3jzoaG-6}\nwendcode{}\nwbegindocs{93} De Moivres Identity: 
\nwenddocs{}\nwbegincode{94}\sublabel{NW4P20f9-4dxCds-C}\nwmargintag{{\nwtagstyle{}\subpageref{NW4P20f9-4dxCds-C}}}\moddef{Dual number class further implementation~{\nwtagstyle{}\subpageref{NW4P20f9-4dxCds-1}}}\plusendmoddef\Rm{}\nwstartdeflinemarkup\nwusesondefline{\\{NW4P20f9-3HhVai-1}}\nwprevnextdefs{NW4P20f9-4dxCds-B}{NW4P20f9-4dxCds-D}\nwenddeflinemarkup
{\bf{}dual\_number} {\bf{}dual\_number}::{\it{}power}({\bf{}const} {\bf{}ex} & {\it{}e}) {\bf{}const} {\nwlbrace}
        {\bf{}return} {\it{}dn\_from\_arg\_mod}({\it{}arg}()\begin{math}\ast\end{math}{\it{}e}, {\it{}pow}({\it{}norm}(), {\it{}e}));
{\nwrbrace}

\nwused{\\{NW4P20f9-3HhVai-1}}\nwidentuses{\\{{\nwixident{dual{\_}number}}{dual:unnumber}}}\nwindexuse{\nwixident{dual{\_}number}}{dual:unnumber}{NW4P20f9-4dxCds-C}\nwendcode{}\nwbegindocs{95}Negative dual number for \((u,v)\) has the same argument and the
negative modulus.
\nwenddocs{}\nwbegincode{96}\sublabel{NW4P20f9-4dxCds-D}\nwmargintag{{\nwtagstyle{}\subpageref{NW4P20f9-4dxCds-D}}}\moddef{Dual number class further implementation~{\nwtagstyle{}\subpageref{NW4P20f9-4dxCds-1}}}\plusendmoddef\Rm{}\nwstartdeflinemarkup\nwusesondefline{\\{NW4P20f9-3HhVai-1}}\nwprevnextdefs{NW4P20f9-4dxCds-C}{NW4P20f9-4dxCds-E}\nwenddeflinemarkup
{\bf{}dual\_number} {\bf{}dual\_number}::{\it{}neg}() {\bf{}const} {\nwlbrace} 
        {\bf{}return} {\it{}dn\_from\_arg\_mod}({\it{}arg}(), -{\it{}norm}());
{\nwrbrace}

\nwused{\\{NW4P20f9-3HhVai-1}}\nwidentuses{\\{{\nwixident{dual{\_}number}}{dual:unnumber}}}\nwindexuse{\nwixident{dual{\_}number}}{dual:unnumber}{NW4P20f9-4dxCds-D}\nwendcode{}\nwbegindocs{97}Action of the group \(\SL\) on the upper halfplane of dual numbers 
(a direct formula).
\nwenddocs{}\nwbegincode{98}\sublabel{NW4P20f9-4dxCds-E}\nwmargintag{{\nwtagstyle{}\subpageref{NW4P20f9-4dxCds-E}}}\moddef{Dual number class further implementation~{\nwtagstyle{}\subpageref{NW4P20f9-4dxCds-1}}}\plusendmoddef\Rm{}\nwstartdeflinemarkup\nwusesondefline{\\{NW4P20f9-3HhVai-1}}\nwprevnextdefs{NW4P20f9-4dxCds-D}{NW4P20f9-4dxCds-F}\nwenddeflinemarkup
{\bf{}dual\_number} {\bf{}dual\_number}::{\it{}sl2\_moebius}({\bf{}const} {\bf{}ex} & {\it{}a}, {\bf{}const} {\bf{}ex} & {\it{}b}, {\bf{}const} {\bf{}ex} & {\it{}c}, {\bf{}const} {\bf{}ex} & {\it{}d}) {\bf{}const} {\nwlbrace}
        {\bf{}return} {\bf{}dual\_number}((({\it{}a}\begin{math}\ast\end{math}{\it{}u\_comp}+{\it{}b})\begin{math}\div\end{math}({\it{}c}\begin{math}\ast\end{math}{\it{}u\_comp}+{\it{}d})).{\it{}normal}(), 
                ({\it{}v\_comp}\begin{math}\ast\end{math}({\it{}a}\begin{math}\ast\end{math}{\it{}d}-{\it{}b}\begin{math}\ast\end{math}{\it{}c})\begin{math}\ast\end{math}{\it{}pow}({\it{}c}\begin{math}\ast\end{math}{\it{}u\_comp}+{\it{}d}, -2)).{\it{}normal}());
{\nwrbrace}

\nwused{\\{NW4P20f9-3HhVai-1}}\nwidentuses{\\{{\nwixident{dual{\_}number}}{dual:unnumber}}}\nwindexuse{\nwixident{dual{\_}number}}{dual:unnumber}{NW4P20f9-4dxCds-E}\nwendcode{}\nwbegindocs{99}The same transform, but given by a matrix
\nwenddocs{}\nwbegincode{100}\sublabel{NW4P20f9-4dxCds-F}\nwmargintag{{\nwtagstyle{}\subpageref{NW4P20f9-4dxCds-F}}}\moddef{Dual number class further implementation~{\nwtagstyle{}\subpageref{NW4P20f9-4dxCds-1}}}\plusendmoddef\Rm{}\nwstartdeflinemarkup\nwusesondefline{\\{NW4P20f9-3HhVai-1}}\nwprevnextdefs{NW4P20f9-4dxCds-E}{NW4P20f9-4dxCds-G}\nwenddeflinemarkup
{\bf{}dual\_number} {\bf{}dual\_number}::{\it{}sl2\_moebius}({\bf{}const} {\bf{}ex} & {\it{}M}) {\bf{}const} {\nwlbrace}
        {\bf{}if} ({\it{}is\_a}\begin{math}<\end{math}{\bf{}matrix}\begin{math}>\end{math}({\it{}M}) \begin{math}\wedge\end{math} {\it{}ex\_to}\begin{math}<\end{math}{\bf{}matrix}\begin{math}>\end{math}({\it{}M}).{\it{}rows}()\begin{math}\equiv\end{math}2  \begin{math}\wedge\end{math} {\it{}ex\_to}\begin{math}<\end{math}{\bf{}matrix}\begin{math}>\end{math}({\it{}M}).{\it{}cols}()\begin{math}\equiv\end{math}2) 
                {\bf{}return} {\it{}sl2\_moebius}({\it{}ex\_to}\begin{math}<\end{math}{\bf{}matrix}\begin{math}>\end{math}({\it{}M})(0,0), {\it{}ex\_to}\begin{math}<\end{math}{\bf{}matrix}\begin{math}>\end{math}({\it{}M})(0,1),
                                            {\it{}ex\_to}\begin{math}<\end{math}{\bf{}matrix}\begin{math}>\end{math}({\it{}M})(1,0), {\it{}ex\_to}\begin{math}<\end{math}{\bf{}matrix}\begin{math}>\end{math}({\it{}M})(1,1));
        {\bf{}else}
                {\bf{}throw}({\it{}std}::{\it{}invalid\_argument}({\tt{}"dual\_number.moebius(): parameter M"} 
                                                                        {\tt{}"should be a 2x2 matrix"}));
{\nwrbrace}

\nwused{\\{NW4P20f9-3HhVai-1}}\nwidentuses{\\{{\nwixident{dual{\_}number}}{dual:unnumber}}}\nwindexuse{\nwixident{dual{\_}number}}{dual:unnumber}{NW4P20f9-4dxCds-F}\nwendcode{}\nwbegindocs{101}Action of the group \(\SL\) conjugated with the Cayley transform on
the unit disk  of dual numbers (a direct formula).
\nwenddocs{}\nwbegincode{102}\sublabel{NW4P20f9-4dxCds-G}\nwmargintag{{\nwtagstyle{}\subpageref{NW4P20f9-4dxCds-G}}}\moddef{Dual number class further implementation~{\nwtagstyle{}\subpageref{NW4P20f9-4dxCds-1}}}\plusendmoddef\Rm{}\nwstartdeflinemarkup\nwusesondefline{\\{NW4P20f9-3HhVai-1}}\nwprevnextdefs{NW4P20f9-4dxCds-F}{NW4P20f9-4dxCds-H}\nwenddeflinemarkup
{\bf{}dual\_number} {\bf{}dual\_number}::{\it{}su\_moebius}({\bf{}const} {\bf{}ex} & {\it{}a}, {\bf{}const} {\bf{}ex} & {\it{}b}, {\bf{}const} {\bf{}ex} & {\it{}c}, {\bf{}const} {\bf{}ex} & {\it{}d}) {\bf{}const} {\nwlbrace}
        {\bf{}ex} {\it{}denom}=({\it{}c}\begin{math}\ast\end{math}{\it{}u\_comp}+{\it{}d}).{\it{}normal}(), {\it{}numer}, {\it{}summand}=({\it{}straight\_spoke}? {\bf{}ex}(1) : {\it{}half});

\nwused{\\{NW4P20f9-3HhVai-1}}\nwidentuses{\\{{\nwixident{dual{\_}number}}{dual:unnumber}}}\nwindexuse{\nwixident{dual{\_}number}}{dual:unnumber}{NW4P20f9-4dxCds-G}\nwendcode{}\nwbegindocs{103}We now define two different expressions for numerator depending from
the type of spokes we are using
\nwenddocs{}\nwbegincode{104}\sublabel{NW4P20f9-4dxCds-H}\nwmargintag{{\nwtagstyle{}\subpageref{NW4P20f9-4dxCds-H}}}\moddef{Dual number class further implementation~{\nwtagstyle{}\subpageref{NW4P20f9-4dxCds-1}}}\plusendmoddef\Rm{}\nwstartdeflinemarkup\nwusesondefline{\\{NW4P20f9-3HhVai-1}}\nwprevnextdefs{NW4P20f9-4dxCds-G}{NW4P20f9-4dxCds-I}\nwenddeflinemarkup
        {\bf{}if} ({\it{}straight\_spoke})
                {\it{}numer} = ({\it{}v\_comp}-{\it{}pow}({\it{}u\_comp}, 2)+1)\begin{math}\ast\end{math}({\it{}a}\begin{math}\ast\end{math}{\it{}d}-{\it{}b}\begin{math}\ast\end{math}{\it{}c})+{\it{}pow}({\it{}a}\begin{math}\ast\end{math}{\it{}u\_comp}+{\it{}b}, 2);
        {\bf{}else} 
                {\it{}numer} = {\it{}half}\begin{math}\ast\end{math}((2\begin{math}\ast\end{math}{\it{}v\_comp}-{\it{}pow}({\it{}u\_comp},2)+1)\begin{math}\ast\end{math}({\it{}a}\begin{math}\ast\end{math}{\it{}d}-{\it{}b}\begin{math}\ast\end{math}{\it{}c})+{\it{}pow}({\it{}a}\begin{math}\ast\end{math}{\it{}u\_comp}+{\it{}b}, 2));

\nwused{\\{NW4P20f9-3HhVai-1}}\nwendcode{}\nwbegindocs{105}For zero denominator we replace it by Inf.
\nwenddocs{}\nwbegincode{106}\sublabel{NW4P20f9-4dxCds-I}\nwmargintag{{\nwtagstyle{}\subpageref{NW4P20f9-4dxCds-I}}}\moddef{Dual number class further implementation~{\nwtagstyle{}\subpageref{NW4P20f9-4dxCds-1}}}\plusendmoddef\Rm{}\nwstartdeflinemarkup\nwusesondefline{\\{NW4P20f9-3HhVai-1}}\nwprevnextdefs{NW4P20f9-4dxCds-H}{NW4P20f9-4dxCds-J}\nwenddeflinemarkup
        {\bf{}if} ({\it{}denom}.{\it{}is\_zero}()) 
                {\bf{}return} {\bf{}dual\_number}(({\it{}a}\begin{math}\ast\end{math}{\it{}u\_comp}+{\it{}b})\begin{math}\ast\end{math}{\it{}Inf}, ({\it{}numer}\begin{math}\ast\end{math}{\it{}pow}({\it{}Inf},2)-{\it{}summand}).{\it{}normal}());
        {\bf{}else}

\nwused{\\{NW4P20f9-3HhVai-1}}\nwidentuses{\\{{\nwixident{dual{\_}number}}{dual:unnumber}}}\nwindexuse{\nwixident{dual{\_}number}}{dual:unnumber}{NW4P20f9-4dxCds-I}\nwendcode{}\nwbegindocs{107}A non-zero denominator may include infinity, which we treat separately.
\nwenddocs{}\nwbegincode{108}\sublabel{NW4P20f9-4dxCds-J}\nwmargintag{{\nwtagstyle{}\subpageref{NW4P20f9-4dxCds-J}}}\moddef{Dual number class further implementation~{\nwtagstyle{}\subpageref{NW4P20f9-4dxCds-1}}}\plusendmoddef\Rm{}\nwstartdeflinemarkup\nwusesondefline{\\{NW4P20f9-3HhVai-1}}\nwprevnextdefs{NW4P20f9-4dxCds-I}{NW4P20f9-4dxCds-K}\nwenddeflinemarkup
                {\bf{}if} ({\it{}u\_comp}.{\it{}has}({\it{}Inf})) {\nwlbrace}
                        {\bf{}realsymbol} {\it{}t}({\tt{}"t"});
                        {\bf{}return} {\bf{}dual\_number}((({\it{}a}\begin{math}\ast\end{math}{\it{}u\_comp}+{\it{}b})\begin{math}\div\end{math}{\it{}denom}).{\it{}subs}({\it{}Inf} \begin{math}\equiv\end{math} {\it{}pow}({\it{}t}, -1))
                                                           .{\it{}normal}().{\it{}subs}({\it{}t} \begin{math}\equiv\end{math} 0).{\it{}normal}(), 
                                                           ({\it{}numer}\begin{math}\ast\end{math}{\it{}pow}({\it{}denom}, -2)-{\it{}summand})
                                                           .{\it{}subs}({\it{}Inf} \begin{math}\equiv\end{math} {\it{}pow}({\it{}t}, -1)).{\it{}normal}().{\it{}subs}({\it{}t} \begin{math}\equiv\end{math} 0).{\it{}normal}());
                {\nwrbrace} {\bf{}else}

\nwused{\\{NW4P20f9-3HhVai-1}}\nwidentuses{\\{{\nwixident{dual{\_}number}}{dual:unnumber}}}\nwindexuse{\nwixident{dual{\_}number}}{dual:unnumber}{NW4P20f9-4dxCds-J}\nwendcode{}\nwbegindocs{109}For non-zero determinant we return the plain value.
\nwenddocs{}\nwbegincode{110}\sublabel{NW4P20f9-4dxCds-K}\nwmargintag{{\nwtagstyle{}\subpageref{NW4P20f9-4dxCds-K}}}\moddef{Dual number class further implementation~{\nwtagstyle{}\subpageref{NW4P20f9-4dxCds-1}}}\plusendmoddef\Rm{}\nwstartdeflinemarkup\nwusesondefline{\\{NW4P20f9-3HhVai-1}}\nwprevnextdefs{NW4P20f9-4dxCds-J}{NW4P20f9-4dxCds-L}\nwenddeflinemarkup
                        {\bf{}return} {\bf{}dual\_number}((({\it{}a}\begin{math}\ast\end{math}{\it{}u\_comp}+{\it{}b})\begin{math}\div\end{math}{\it{}denom}).{\it{}normal}(), 
                                                           ({\it{}numer}\begin{math}\ast\end{math}{\it{}pow}({\it{}denom}, -2)-{\it{}summand}).{\it{}normal}());
{\nwrbrace}

\nwused{\\{NW4P20f9-3HhVai-1}}\nwidentuses{\\{{\nwixident{dual{\_}number}}{dual:unnumber}}}\nwindexuse{\nwixident{dual{\_}number}}{dual:unnumber}{NW4P20f9-4dxCds-K}\nwendcode{}\nwbegindocs{111}A variant of the function if the matrix is given as a whole.
\nwenddocs{}\nwbegincode{112}\sublabel{NW4P20f9-4dxCds-L}\nwmargintag{{\nwtagstyle{}\subpageref{NW4P20f9-4dxCds-L}}}\moddef{Dual number class further implementation~{\nwtagstyle{}\subpageref{NW4P20f9-4dxCds-1}}}\plusendmoddef\Rm{}\nwstartdeflinemarkup\nwusesondefline{\\{NW4P20f9-3HhVai-1}}\nwprevnextdefs{NW4P20f9-4dxCds-K}{NW4P20f9-4dxCds-M}\nwenddeflinemarkup
{\bf{}dual\_number} {\bf{}dual\_number}::{\it{}su\_moebius}({\bf{}const} {\bf{}ex} & {\it{}M}) {\bf{}const} {\nwlbrace}
        {\bf{}if} ({\it{}is\_a}\begin{math}<\end{math}{\bf{}matrix}\begin{math}>\end{math}({\it{}M}) \begin{math}\wedge\end{math} {\it{}ex\_to}\begin{math}<\end{math}{\bf{}matrix}\begin{math}>\end{math}({\it{}M}).{\it{}rows}()\begin{math}\equiv\end{math}2  \begin{math}\wedge\end{math} {\it{}ex\_to}\begin{math}<\end{math}{\bf{}matrix}\begin{math}>\end{math}({\it{}M}).{\it{}cols}()\begin{math}\equiv\end{math}2) 
                {\bf{}return} {\it{}su\_moebius}({\it{}ex\_to}\begin{math}<\end{math}{\bf{}matrix}\begin{math}>\end{math}({\it{}M})(0,0), {\it{}ex\_to}\begin{math}<\end{math}{\bf{}matrix}\begin{math}>\end{math}({\it{}M})(0,1),
                                            {\it{}ex\_to}\begin{math}<\end{math}{\bf{}matrix}\begin{math}>\end{math}({\it{}M})(1,0), {\it{}ex\_to}\begin{math}<\end{math}{\bf{}matrix}\begin{math}>\end{math}({\it{}M})(1,1));
        {\bf{}else}
                {\bf{}throw}({\it{}std}::{\it{}invalid\_argument}({\tt{}"dual\_number.moebius(): parameter M should "}
                                                                        {\tt{}"be a 2x2 matrix"}));
{\nwrbrace}

\nwused{\\{NW4P20f9-3HhVai-1}}\nwidentuses{\\{{\nwixident{dual{\_}number}}{dual:unnumber}}}\nwindexuse{\nwixident{dual{\_}number}}{dual:unnumber}{NW4P20f9-4dxCds-L}\nwendcode{}\nwbegindocs{113}Parabolic rotations
\nwenddocs{}\nwbegincode{114}\sublabel{NW4P20f9-4dxCds-M}\nwmargintag{{\nwtagstyle{}\subpageref{NW4P20f9-4dxCds-M}}}\moddef{Dual number class further implementation~{\nwtagstyle{}\subpageref{NW4P20f9-4dxCds-1}}}\plusendmoddef\Rm{}\nwstartdeflinemarkup\nwusesondefline{\\{NW4P20f9-3HhVai-1}}\nwprevnextdefs{NW4P20f9-4dxCds-L}{\relax}\nwenddeflinemarkup
{\bf{}dual\_number} {\bf{}dual\_number}::{\it{}rot}({\bf{}const} {\bf{}ex} & {\it{}x}) {\bf{}const} {\nwlbrace}
        {\bf{}if} ({\it{}subgroup\_N})
                {\bf{}return} {\it{}su\_moebius}(1, {\it{}x}, 0, 1);
        {\bf{}else}
                {\bf{}return} {\it{}su\_moebius}(1, 0, {\it{}x}, 1);

{\nwrbrace}

\nwused{\\{NW4P20f9-3HhVai-1}}\nwidentuses{\\{{\nwixident{dual{\_}number}}{dual:unnumber}}}\nwindexuse{\nwixident{dual{\_}number}}{dual:unnumber}{NW4P20f9-4dxCds-M}\nwendcode{}\nwbegindocs{115}All algebraic routines are defined now.

\nwenddocs{}\nwbegindocs{116}\nwdocspar
\subsection{Library Outline}
\label{sec:library-outline}
Here is the outline how we use the above parts.

\nwenddocs{}\nwbegindocs{117}Routines for {\Tt{}\Rm{}{\bf{}dual\_number}\nwendquote} are collected in a separate library. We
start from the definition {\Tt{}\Rm{}{\bf{}dual\_number}\nwendquote} class in the header file.
\nwenddocs{}\nwbegincode{118}\sublabel{NW4P20f9-36Ytqo-1}\nwmargintag{{\nwtagstyle{}\subpageref{NW4P20f9-36Ytqo-1}}}\moddef{dualnum.h~{\nwtagstyle{}\subpageref{NW4P20f9-36Ytqo-1}}}\endmoddef\Rm{}\nwstartdeflinemarkup\nwenddeflinemarkup
 \LA{}Header initialisation~{\nwtagstyle{}\subpageref{NW4P20f9-1vbmXO-1}}\RA{}
 \LA{}N-Nprime separation~{\nwtagstyle{}\subpageref{NW4P20f9-4EWpKB-1}}\RA{}
 \LA{}Dual number class declaration~{\nwtagstyle{}\subpageref{NW4P20f9-4GrAIY-1}}\RA{}
 \LA{}Additional routines declarations~{\nwtagstyle{}\subpageref{NW4P20f9-NQ5ML-1}}\RA{}

\nwnotused{dualnum.h}\nwendcode{}\nwbegindocs{119}This is the initialisation part of the header file
\nwenddocs{}\nwbegincode{120}\sublabel{NW4P20f9-1vbmXO-1}\nwmargintag{{\nwtagstyle{}\subpageref{NW4P20f9-1vbmXO-1}}}\moddef{Header initialisation~{\nwtagstyle{}\subpageref{NW4P20f9-1vbmXO-1}}}\endmoddef\Rm{}\nwstartdeflinemarkup\nwusesondefline{\\{NW4P20f9-36Ytqo-1}}\nwenddeflinemarkup
{\bf{}\char35{}include}{\tt{} \begin{math}<\end{math}fstream\begin{math}>\end{math}}
 {\bf{}using} {\bf{}namespace} {\it{}std};

{\bf{}\char35{}include}{\tt{} \begin{math}<\end{math}ginac/ginac.h\begin{math}>\end{math}}
 {\bf{}using} {\bf{}namespace} {\it{}GiNaC};
{\bf{}\char35{}define}{\tt{} GINAC\_VERSION\_ATLEAST( major, minor) \begin{math}\backslash\end{math}}\nwindexdefn{\nwixident{GINAC{\_}VERSION{\_}ATLEAST}}{GINAC:unVERSION:unATLEAST}{NW4P20f9-1vbmXO-1}
        ({\it{}GINACLIB\_MAJOR\_VERSION} \begin{math}>\end{math} {\it{}major} \begin{math}\backslash\end{math}
        \begin{math}\vee\end{math} ({\it{}GINACLIB\_MAJOR\_VERSION} \begin{math}\equiv\end{math} {\it{}major} \begin{math}\wedge\end{math} {\it{}GINACLIB\_MINOR\_VERSION} \begin{math}\geq\end{math} {\it{}minor}))

\nwused{\\{NW4P20f9-36Ytqo-1}}\nwidentdefs{\\{{\nwixident{GINAC{\_}VERSION{\_}ATLEAST}}{GINAC:unVERSION:unATLEAST}}}\nwendcode{}\nwbegindocs{121}Here is the file with the implementation.
\nwenddocs{}\nwbegincode{122}\sublabel{NW4P20f9-3HhVai-1}\nwmargintag{{\nwtagstyle{}\subpageref{NW4P20f9-3HhVai-1}}}\moddef{dualnum.cpp~{\nwtagstyle{}\subpageref{NW4P20f9-3HhVai-1}}}\endmoddef\Rm{}\nwstartdeflinemarkup\nwenddeflinemarkup
{\bf{}\char35{}include}{\tt{} \begin{math}<\end{math}dualnum.h\begin{math}>\end{math}}
{\bf{}const} {\bf{}numeric} {\it{}half}={\bf{}numeric}(1,2);\nwindexdefn{\nwixident{numeric}}{numeric}{NW4P20f9-3HhVai-1}
 \LA{}Algebraic procedures~{\nwtagstyle{}\subpageref{NW4P20f9-3jzoaG-1}}\RA{}
 \LA{}Dual number class implementation~{\nwtagstyle{}\subpageref{NW4P20f9-3uuIiL-1}}\RA{}
 \LA{}Dual number class further implementation~{\nwtagstyle{}\subpageref{NW4P20f9-4dxCds-1}}\RA{}
 \LA{}Output routines~{\nwtagstyle{}\subpageref{NW4P20f9-4MrGQl-1}}\RA{}

\nwnotused{dualnum.cpp}\nwidentdefs{\\{{\nwixident{numeric}}{numeric}}}\nwendcode{}\nwbegindocs{123}\nwdocspar
\subsubsection{Output routines}
\label{sec:output-routines-lib}

We use standardised routines to output results of calculations.
\nwenddocs{}\nwbegincode{124}\sublabel{NW4P20f9-4MrGQl-1}\nwmargintag{{\nwtagstyle{}\subpageref{NW4P20f9-4MrGQl-1}}}\moddef{Output routines~{\nwtagstyle{}\subpageref{NW4P20f9-4MrGQl-1}}}\endmoddef\Rm{}\nwstartdeflinemarkup\nwusesondefline{\\{NW4P20f9-3HhVai-1}}\nwprevnextdefs{\relax}{NW4P20f9-4MrGQl-2}\nwenddeflinemarkup
{\bf{}void} {\it{}formula\_out}({\it{}string} {\it{}S}, {\bf{}const} {\bf{}ex} & {\it{}F}, {\bf{}bool} {\it{}lineend}) {\nwlbrace}\nwindexdefn{\nwixident{formula{\_}out}}{formula:unout}{NW4P20f9-4MrGQl-1}
        {\it{}cout} \begin{math}\ll\end{math} {\it{}S} \begin{math}\ll\end{math} ({\it{}latexout} ? {\tt{}"{\char92}{\char92}("} : {\tt{}""} ) \begin{math}\ll\end{math} {\it{}F} \begin{math}\ll\end{math} ({\it{}latexout} ? {\tt{}"{\char92}{\char92})"} : {\tt{}""} );
        {\bf{}if} ({\it{}lineend})
                {\it{}cout} \begin{math}\ll\end{math} ({\it{}latexout} ? {\tt{}"{\char92}{\char92}{\char92}{\char92}"} : {\tt{}""} ) \begin{math}\ll\end{math} {\it{}endl};
        {\bf{}else}
                {\it{}cout} \begin{math}\ll\end{math} {\tt{}"; "}; 
{\nwrbrace}

\nwalsodefined{\\{NW4P20f9-4MrGQl-2}}\nwused{\\{NW4P20f9-3HhVai-1}}\nwidentdefs{\\{{\nwixident{formula{\_}out}}{formula:unout}}}\nwendcode{}\nwbegindocs{125}This routine is used to check identities.
\nwenddocs{}\nwbegincode{126}\sublabel{NW4P20f9-4MrGQl-2}\nwmargintag{{\nwtagstyle{}\subpageref{NW4P20f9-4MrGQl-2}}}\moddef{Output routines~{\nwtagstyle{}\subpageref{NW4P20f9-4MrGQl-1}}}\plusendmoddef\Rm{}\nwstartdeflinemarkup\nwusesondefline{\\{NW4P20f9-3HhVai-1}}\nwprevnextdefs{NW4P20f9-4MrGQl-1}{\relax}\nwenddeflinemarkup
{\bf{}void} {\it{}test\_out}({\it{}string} {\it{}S}, {\bf{}const} {\bf{}ex} & {\it{}T}) {\nwlbrace}\nwindexdefn{\nwixident{test{\_}out}}{test:unout}{NW4P20f9-4MrGQl-2}
        {\it{}cout} \begin{math}\ll\end{math} {\it{}S} \begin{math}\ll\end{math} ({\it{}latexout} ? {\tt{}"{\char92}{\char92}textbf{\char123}"} : {\tt{}"*"} ) 
                 \begin{math}\ll\end{math} ({\it{}is\_a}\begin{math}<\end{math}{\bf{}dual\_number}\begin{math}>\end{math}({\it{}T}) ? {\it{}ex\_to}\begin{math}<\end{math}{\bf{}dual\_number}\begin{math}>\end{math}({\it{}T}).{\it{}normal}().{\it{}is\_zero}() :
                         {\it{}T}.{\it{}evalm}().{\it{}normal}().{\it{}is\_zero\_matrix}()) \begin{math}\ll\end{math} ({\it{}latexout} ? {\tt{}"{\char125}{\char92}{\char92}{\char92}{\char92}"} : {\tt{}"*"} )
                 \begin{math}\ll\end{math} {\it{}endl};
{\nwrbrace}

\nwused{\\{NW4P20f9-3HhVai-1}}\nwidentdefs{\\{{\nwixident{test{\_}out}}{test:unout}}}\nwidentuses{\\{{\nwixident{dual{\_}number}}{dual:unnumber}}}\nwindexuse{\nwixident{dual{\_}number}}{dual:unnumber}{NW4P20f9-4MrGQl-2}\nwendcode{}\nwbegindocs{127}Here is declarations of additional routines for the header file.
\nwenddocs{}\nwbegincode{128}\sublabel{NW4P20f9-NQ5ML-2}\nwmargintag{{\nwtagstyle{}\subpageref{NW4P20f9-NQ5ML-2}}}\moddef{Additional routines declarations~{\nwtagstyle{}\subpageref{NW4P20f9-NQ5ML-1}}}\plusendmoddef\Rm{}\nwstartdeflinemarkup\nwusesondefline{\\{NW4P20f9-36Ytqo-1}}\nwprevnextdefs{NW4P20f9-NQ5ML-1}{NW4P20f9-NQ5ML-3}\nwenddeflinemarkup
{\bf{}bool} {\it{}latexout};

\nwused{\\{NW4P20f9-36Ytqo-1}}\nwendcode{}\nwbegindocs{129}Function interface to global variables. Two work properly the
library shall fix which kind rotations, either \(N\) or \(N'\), is used.
\nwenddocs{}\nwbegincode{130}\sublabel{NW4P20f9-NQ5ML-3}\nwmargintag{{\nwtagstyle{}\subpageref{NW4P20f9-NQ5ML-3}}}\moddef{Additional routines declarations~{\nwtagstyle{}\subpageref{NW4P20f9-NQ5ML-1}}}\plusendmoddef\Rm{}\nwstartdeflinemarkup\nwusesondefline{\\{NW4P20f9-36Ytqo-1}}\nwprevnextdefs{NW4P20f9-NQ5ML-2}{NW4P20f9-NQ5ML-4}\nwenddeflinemarkup
{\bf{}ex} {\it{}arg0}() {\nwlbrace} {\bf{}return} {\it{}Arg0}; {\nwrbrace}
{\bf{}bool} {\it{}is\_subgroup\_N}() {\nwlbrace} {\bf{}return} {\it{}subgroup\_N}; {\nwrbrace}
{\bf{}void} {\it{}set\_subgroup}({\bf{}bool} {\it{}N});

\nwused{\\{NW4P20f9-36Ytqo-1}}\nwidentuses{\\{{\nwixident{Arg0}}{Arg0}}\\{{\nwixident{set{\_}subgroup}}{set:unsubgroup}}}\nwindexuse{\nwixident{Arg0}}{Arg0}{NW4P20f9-NQ5ML-3}\nwindexuse{\nwixident{set{\_}subgroup}}{set:unsubgroup}{NW4P20f9-NQ5ML-3}\nwendcode{}\nwbegindocs{131} We also define procedures to setup which spoke are used with the
subgroup \(N'\), either vertical lines or geodesics. This does not
affect calculations with the subgroup \(N\).
\nwenddocs{}\nwbegincode{132}\sublabel{NW4P20f9-NQ5ML-4}\nwmargintag{{\nwtagstyle{}\subpageref{NW4P20f9-NQ5ML-4}}}\moddef{Additional routines declarations~{\nwtagstyle{}\subpageref{NW4P20f9-NQ5ML-1}}}\plusendmoddef\Rm{}\nwstartdeflinemarkup\nwusesondefline{\\{NW4P20f9-36Ytqo-1}}\nwprevnextdefs{NW4P20f9-NQ5ML-3}{NW4P20f9-NQ5ML-5}\nwenddeflinemarkup
{\bf{}bool} {\it{}is\_spoke\_straight}() {\nwlbrace} {\bf{}return} {\it{}straight\_spoke}; {\nwrbrace}
{\bf{}void} {\it{}set\_straight\_spoke}({\bf{}bool} {\it{}S});

\nwused{\\{NW4P20f9-36Ytqo-1}}\nwidentuses{\\{{\nwixident{set{\_}straight{\_}spoke}}{set:unstraight:unspoke}}}\nwindexuse{\nwixident{set{\_}straight{\_}spoke}}{set:unstraight:unspoke}{NW4P20f9-NQ5ML-4}\nwendcode{}\nwbegindocs{133}Nicer output for testing routines.
\nwenddocs{}\nwbegincode{134}\sublabel{NW4P20f9-NQ5ML-5}\nwmargintag{{\nwtagstyle{}\subpageref{NW4P20f9-NQ5ML-5}}}\moddef{Additional routines declarations~{\nwtagstyle{}\subpageref{NW4P20f9-NQ5ML-1}}}\plusendmoddef\Rm{}\nwstartdeflinemarkup\nwusesondefline{\\{NW4P20f9-36Ytqo-1}}\nwprevnextdefs{NW4P20f9-NQ5ML-4}{\relax}\nwenddeflinemarkup
{\bf{}void} {\it{}test\_out}({\it{}string} {\it{}S}, {\bf{}const} {\bf{}ex} & {\it{}T});
{\bf{}void} {\it{}formula\_out}({\it{}string} {\it{}S}, {\bf{}const} {\bf{}ex} & {\it{}F}, {\bf{}bool} {\it{}lineend}={\bf{}true});

\nwused{\\{NW4P20f9-36Ytqo-1}}\nwidentuses{\\{{\nwixident{formula{\_}out}}{formula:unout}}\\{{\nwixident{test{\_}out}}{test:unout}}}\nwindexuse{\nwixident{formula{\_}out}}{formula:unout}{NW4P20f9-NQ5ML-5}\nwindexuse{\nwixident{test{\_}out}}{test:unout}{NW4P20f9-NQ5ML-5}\nwendcode{}\nwbegindocs{135}\nwdocspar
\subsection{Header and Implementation of the {\Tt{}\Rm{}{\bf{}dual\_number}\nwendquote} Class}
\label{sec:impl-dual_n-class}

\nwenddocs{}\nwbegindocs{136}\nwdocspar

\subsubsection{Header File for {\Tt{}\Rm{}{\bf{}dual\_number}\nwendquote}}
\label{sec:head-file-dual_n}

\nwenddocs{}\nwbegindocs{137}We use the standard \GiNaC\ machinery do define {\Tt{}\Rm{}{\bf{}dual\_number}\nwendquote}s as
derived of the class {\Tt{}\Rm{}{\bf{}basic}\nwendquote}.
\nwenddocs{}\nwbegincode{138}\sublabel{NW4P20f9-4GrAIY-1}\nwmargintag{{\nwtagstyle{}\subpageref{NW4P20f9-4GrAIY-1}}}\moddef{Dual number class declaration~{\nwtagstyle{}\subpageref{NW4P20f9-4GrAIY-1}}}\endmoddef\Rm{}\nwstartdeflinemarkup\nwusesondefline{\\{NW4P20f9-36Ytqo-1}}\nwprevnextdefs{\relax}{NW4P20f9-4GrAIY-2}\nwenddeflinemarkup
{\bf{}class} {\bf{}dual\_number} : {\bf{}public} {\bf{}basic}
{\nwlbrace}
        {\it{}GINAC\_DECLARE\_REGISTERED\_CLASS}({\bf{}dual\_number}, {\bf{}basic})
{\bf{}\char35{}if}{\tt{} GINAC\_VERSION\_ATLEAST(1,5)}
{\bf{}\char35{}else}{\tt{}}
        {\bf{}static} {\bf{}const} {\it{}tinfo\_static\_t} {\it{}return\_type\_tinfo\_static}[256];
{\bf{}\char35{}endif}{\tt{}}

\nwalsodefined{\\{NW4P20f9-4GrAIY-2}\\{NW4P20f9-4GrAIY-3}\\{NW4P20f9-4GrAIY-4}\\{NW4P20f9-4GrAIY-5}\\{NW4P20f9-4GrAIY-6}}\nwused{\\{NW4P20f9-36Ytqo-1}}\nwidentuses{\\{{\nwixident{dual{\_}number}}{dual:unnumber}}\\{{\nwixident{GINAC{\_}VERSION{\_}ATLEAST}}{GINAC:unVERSION:unATLEAST}}\\{{\nwixident{tinfo{\_}static{\_}t}}{tinfo:unstatic:unt}}}\nwindexuse{\nwixident{dual{\_}number}}{dual:unnumber}{NW4P20f9-4GrAIY-1}\nwindexuse{\nwixident{GINAC{\_}VERSION{\_}ATLEAST}}{GINAC:unVERSION:unATLEAST}{NW4P20f9-4GrAIY-1}\nwindexuse{\nwixident{tinfo{\_}static{\_}t}}{tinfo:unstatic:unt}{NW4P20f9-4GrAIY-1}\nwendcode{}\nwbegindocs{139}Public methods (constructors, algebraic, etc.)
\nwenddocs{}\nwbegincode{140}\sublabel{NW4P20f9-4GrAIY-2}\nwmargintag{{\nwtagstyle{}\subpageref{NW4P20f9-4GrAIY-2}}}\moddef{Dual number class declaration~{\nwtagstyle{}\subpageref{NW4P20f9-4GrAIY-1}}}\plusendmoddef\Rm{}\nwstartdeflinemarkup\nwusesondefline{\\{NW4P20f9-36Ytqo-1}}\nwprevnextdefs{NW4P20f9-4GrAIY-1}{NW4P20f9-4GrAIY-3}\nwenddeflinemarkup
{\bf{}public}:
        \LA{}Public methods~{\nwtagstyle{}\subpageref{NW4P20f9-1eKCCy-1}}\RA{}
        \LA{}Technical methods~{\nwtagstyle{}\subpageref{NW4P20f9-2XWNbs-1}}\RA{}
        \LA{}Dual number operators~{\nwtagstyle{}\subpageref{NW4P20f9-56t00-1}}\RA{}

\nwused{\\{NW4P20f9-36Ytqo-1}}\nwendcode{}\nwbegindocs{141}We redefine protected methods for printing only.
\nwenddocs{}\nwbegincode{142}\sublabel{NW4P20f9-4GrAIY-3}\nwmargintag{{\nwtagstyle{}\subpageref{NW4P20f9-4GrAIY-3}}}\moddef{Dual number class declaration~{\nwtagstyle{}\subpageref{NW4P20f9-4GrAIY-1}}}\plusendmoddef\Rm{}\nwstartdeflinemarkup\nwusesondefline{\\{NW4P20f9-36Ytqo-1}}\nwprevnextdefs{NW4P20f9-4GrAIY-2}{NW4P20f9-4GrAIY-4}\nwenddeflinemarkup
{\bf{}protected}:
        {\bf{}void} {\it{}do\_print}({\bf{}const} {\it{}print\_context} & {\it{}c}, {\bf{}unsigned} {\it{}level}) {\bf{}const};
        {\bf{}void} {\it{}do\_print\_latex}({\bf{}const} {\it{}print\_latex} & {\it{}c}, {\bf{}unsigned} {\it{}level}) {\bf{}const};

\nwused{\\{NW4P20f9-36Ytqo-1}}\nwendcode{}\nwbegindocs{143}Private members: two components of a {\Tt{}\Rm{}{\bf{}dual\_number}\nwendquote}.
\nwenddocs{}\nwbegincode{144}\sublabel{NW4P20f9-4GrAIY-4}\nwmargintag{{\nwtagstyle{}\subpageref{NW4P20f9-4GrAIY-4}}}\moddef{Dual number class declaration~{\nwtagstyle{}\subpageref{NW4P20f9-4GrAIY-1}}}\plusendmoddef\Rm{}\nwstartdeflinemarkup\nwusesondefline{\\{NW4P20f9-36Ytqo-1}}\nwprevnextdefs{NW4P20f9-4GrAIY-3}{NW4P20f9-4GrAIY-5}\nwenddeflinemarkup
{\bf{}protected}: 
        {\bf{}ex} {\it{}u\_comp};
        {\bf{}ex} {\it{}v\_comp};
{\nwrbrace};

\nwused{\\{NW4P20f9-36Ytqo-1}}\nwendcode{}\nwbegindocs{145} We overload standard algebraic operations for {\Tt{}\Rm{}{\bf{}dual\_number}\nwendquote}.
\nwenddocs{}\nwbegincode{146}\sublabel{NW4P20f9-56t00-1}\nwmargintag{{\nwtagstyle{}\subpageref{NW4P20f9-56t00-1}}}\moddef{Dual number operators~{\nwtagstyle{}\subpageref{NW4P20f9-56t00-1}}}\endmoddef\Rm{}\nwstartdeflinemarkup\nwusesondefline{\\{NW4P20f9-4GrAIY-2}}\nwenddeflinemarkup
//const dual\_number operator+(const dual\_number & rh);
//const dual\_number operator-(const dual\_number & rh);
//const dual\_number operator*(const dual\_number & rh);
//const dual\_number operator*(const ex & rh);
//const dual\_number operator/(const ex & rh);
//const dual\_number operator/(const dual\_number & rh);

\nwused{\\{NW4P20f9-4GrAIY-2}}\nwidentuses{\\{{\nwixident{dual{\_}number}}{dual:unnumber}}}\nwindexuse{\nwixident{dual{\_}number}}{dual:unnumber}{NW4P20f9-56t00-1}\nwendcode{}\nwbegindocs{147} We overload standard algebraic operations for {\Tt{}\Rm{}{\bf{}dual\_number}\nwendquote}.
\nwenddocs{}\nwbegincode{148}\sublabel{NW4P20f9-2jj2ci-1}\nwmargintag{{\nwtagstyle{}\subpageref{NW4P20f9-2jj2ci-1}}}\moddef{Dual number spare operators~{\nwtagstyle{}\subpageref{NW4P20f9-2jj2ci-1}}}\endmoddef\Rm{}\nwstartdeflinemarkup\nwenddeflinemarkup
{\bf{}const} {\bf{}dual\_number} {\bf{}operator}+({\bf{}const} {\bf{}dual\_number} & {\it{}lh}, {\bf{}const} {\bf{}dual\_number} & {\it{}rh});\nwindexdefn{\nwixident{dual{\_}number}}{dual:unnumber}{NW4P20f9-2jj2ci-1}
{\bf{}const} {\bf{}dual\_number} {\bf{}operator}-({\bf{}const} {\bf{}dual\_number} & {\it{}lh}, {\bf{}const} {\bf{}dual\_number} & {\it{}rh});\nwindexdefn{\nwixident{dual{\_}number}}{dual:unnumber}{NW4P20f9-2jj2ci-1}
{\bf{}const} {\bf{}dual\_number} {\bf{}operator}\begin{math}\ast\end{math}({\bf{}const} {\bf{}dual\_number} & {\it{}lh}, {\bf{}const} {\bf{}dual\_number} & {\it{}rh});\nwindexdefn{\nwixident{dual{\_}number}}{dual:unnumber}{NW4P20f9-2jj2ci-1}
{\bf{}const} {\bf{}dual\_number} {\bf{}operator}\begin{math}\ast\end{math}({\bf{}const} {\bf{}dual\_number} & {\it{}lh}, {\bf{}const} {\bf{}ex} & {\it{}rh});\nwindexdefn{\nwixident{dual{\_}number}}{dual:unnumber}{NW4P20f9-2jj2ci-1}
{\bf{}const} {\bf{}dual\_number} {\bf{}operator}\begin{math}\ast\end{math}({\bf{}const} {\bf{}ex} & {\it{}lh}, {\bf{}const} {\bf{}dual\_number} & {\it{}rh});\nwindexdefn{\nwixident{dual{\_}number}}{dual:unnumber}{NW4P20f9-2jj2ci-1}
{\bf{}const} {\bf{}dual\_number} {\bf{}operator}\begin{math}\div\end{math}({\bf{}const} {\bf{}dual\_number} & {\it{}lh}, {\bf{}const} {\bf{}ex} & {\it{}rh});\nwindexdefn{\nwixident{dual{\_}number}}{dual:unnumber}{NW4P20f9-2jj2ci-1}
{\bf{}const} {\bf{}dual\_number} {\bf{}operator}\begin{math}\div\end{math}({\bf{}const} {\bf{}ex} & {\it{}lh}, {\bf{}const} {\bf{}dual\_number} & {\it{}rh});\nwindexdefn{\nwixident{dual{\_}number}}{dual:unnumber}{NW4P20f9-2jj2ci-1}
{\bf{}const} {\bf{}dual\_number} {\bf{}operator}\begin{math}\div\end{math}({\bf{}const} {\bf{}dual\_number} & {\it{}lh}, {\bf{}const} {\bf{}ex} & {\it{}rh});\nwindexdefn{\nwixident{dual{\_}number}}{dual:unnumber}{NW4P20f9-2jj2ci-1}

\nwnotused{Dual number spare operators}\nwidentdefs{\\{{\nwixident{dual{\_}number}}{dual:unnumber}}}\nwendcode{}\nwbegindocs{149} We overload standard algebraic operations for {\Tt{}\Rm{}{\bf{}dual\_number}\nwendquote}.
\nwenddocs{}\nwbegincode{150}\sublabel{NW4P20f9-4GrAIY-5}\nwmargintag{{\nwtagstyle{}\subpageref{NW4P20f9-4GrAIY-5}}}\moddef{Dual number class declaration~{\nwtagstyle{}\subpageref{NW4P20f9-4GrAIY-1}}}\plusendmoddef\Rm{}\nwstartdeflinemarkup\nwusesondefline{\\{NW4P20f9-36Ytqo-1}}\nwprevnextdefs{NW4P20f9-4GrAIY-4}{NW4P20f9-4GrAIY-6}\nwenddeflinemarkup
//const dual\_number operator*(const ex & lh, const dual\_number & rh);
//const dual\_number operator/(const ex & lh, const dual\_number & rh);

\nwused{\\{NW4P20f9-36Ytqo-1}}\nwidentuses{\\{{\nwixident{dual{\_}number}}{dual:unnumber}}}\nwindexuse{\nwixident{dual{\_}number}}{dual:unnumber}{NW4P20f9-4GrAIY-5}\nwendcode{}\nwbegindocs{151}The following methods are needed for \GiNaC\ to work properly.
\nwenddocs{}\nwbegincode{152}\sublabel{NW4P20f9-2XWNbs-1}\nwmargintag{{\nwtagstyle{}\subpageref{NW4P20f9-2XWNbs-1}}}\moddef{Technical methods~{\nwtagstyle{}\subpageref{NW4P20f9-2XWNbs-1}}}\endmoddef\Rm{}\nwstartdeflinemarkup\nwusesondefline{\\{NW4P20f9-4GrAIY-2}}\nwenddeflinemarkup
        {\bf{}dual\_number} {\it{}normal}() {\bf{}const} {\nwlbrace} {\bf{}return} {\bf{}dual\_number}({\it{}u\_comp}.{\it{}normal}(), {\it{}v\_comp}.{\it{}normal}()); {\nwrbrace}
        {\bf{}dual\_number} {\it{}subs}({\bf{}const} {\bf{}ex} & {\it{}e}, {\bf{}unsigned} {\it{}options} = 0) {\bf{}const};
        {\bf{}bool} {\it{}is\_zero}() {\bf{}const};
        {\bf{}bool} {\it{}is\_equal}({\bf{}const} {\bf{}ex} & {\it{}other}) {\bf{}const};
        {\it{}size\_t} {\it{}nops}() {\bf{}const} {\nwlbrace} {\bf{}return} 2; {\nwrbrace}
        {\bf{}ex} {\it{}op}({\it{}size\_t} {\it{}i}) {\bf{}const};
        {\bf{}ex} & {\it{}let\_op}({\it{}size\_t} {\it{}i});
{\bf{}\char35{}if}{\tt{} GINAC\_VERSION\_ATLEAST(1,5)}
        {\bf{}dual\_number}({\bf{}const} {\it{}archive\_node} &{\it{}n}, {\bf{}lst} &{\it{}sym\_lst});
        {\bf{}void} {\it{}archive}({\it{}archive\_node} &{\it{}n}) {\bf{}const};
        {\bf{}ex} {\it{}unarchive}({\bf{}const} {\it{}archive\_node} &{\it{}n}, {\bf{}lst} &{\it{}sym\_lst});
        {\it{}return\_type\_t} {\it{}return\_type\_tinfo}() {\bf{}const} {\nwlbrace}{\bf{}return} {\it{}make\_return\_type\_t}\begin{math}<\end{math}{\bf{}dual\_number}\begin{math}>\end{math}();{\nwrbrace}
{\bf{}\char35{}endif}{\tt{}}

\nwused{\\{NW4P20f9-4GrAIY-2}}\nwidentuses{\\{{\nwixident{dual{\_}number}}{dual:unnumber}}\\{{\nwixident{GINAC{\_}VERSION{\_}ATLEAST}}{GINAC:unVERSION:unATLEAST}}}\nwindexuse{\nwixident{dual{\_}number}}{dual:unnumber}{NW4P20f9-2XWNbs-1}\nwindexuse{\nwixident{GINAC{\_}VERSION{\_}ATLEAST}}{GINAC:unVERSION:unATLEAST}{NW4P20f9-2XWNbs-1}\nwendcode{}\nwbegindocs{153}Additional routines.
\nwenddocs{}\nwbegincode{154}\sublabel{NW4P20f9-4GrAIY-6}\nwmargintag{{\nwtagstyle{}\subpageref{NW4P20f9-4GrAIY-6}}}\moddef{Dual number class declaration~{\nwtagstyle{}\subpageref{NW4P20f9-4GrAIY-1}}}\plusendmoddef\Rm{}\nwstartdeflinemarkup\nwusesondefline{\\{NW4P20f9-36Ytqo-1}}\nwprevnextdefs{NW4P20f9-4GrAIY-5}{\relax}\nwenddeflinemarkup
//dual\_number dn\_from\_arg\_mod(const ex & a, const ex & n);
{\bf{}ex} {\it{}v\_from\_u\_and\_norm}({\bf{}const} {\bf{}ex} & {\it{}u}, {\bf{}const} {\bf{}ex} & {\it{}n});

// End of "header"

\nwused{\\{NW4P20f9-36Ytqo-1}}\nwidentuses{\\{{\nwixident{dual{\_}number}}{dual:unnumber}}}\nwindexuse{\nwixident{dual{\_}number}}{dual:unnumber}{NW4P20f9-4GrAIY-6}\nwendcode{}\nwbegindocs{155}\nwdocspar
\subsubsection{Standard Implementation Part}
\label{sec:stand-impl-part}

\nwenddocs{}\nwbegindocs{156}The implementation uses standard \GiNaC\ technique.
\nwenddocs{}\nwbegincode{157}\sublabel{NW4P20f9-3uuIiL-1}\nwmargintag{{\nwtagstyle{}\subpageref{NW4P20f9-3uuIiL-1}}}\moddef{Dual number class implementation~{\nwtagstyle{}\subpageref{NW4P20f9-3uuIiL-1}}}\endmoddef\Rm{}\nwstartdeflinemarkup\nwusesondefline{\\{NW4P20f9-3HhVai-1}}\nwprevnextdefs{\relax}{NW4P20f9-3uuIiL-2}\nwenddeflinemarkup

{\it{}GINAC\_IMPLEMENT\_REGISTERED\_CLASS\_OPT}({\bf{}dual\_number}, {\bf{}basic}, 
                                                  {\it{}print\_func}\begin{math}<\end{math}{\it{}print\_context}\begin{math}>\end{math}(&{\bf{}dual\_number}::{\it{}do\_print}).
                                                  {\it{}print\_func}\begin{math}<\end{math}{\it{}print\_latex}\begin{math}>\end{math}(&{\bf{}dual\_number}::{\it{}do\_print\_latex}))

{\bf{}\char35{}if}{\tt{} GINAC\_VERSION\_ATLEAST(1,5)}
{\bf{}\char35{}else}{\tt{}}
{\bf{}const} {\it{}tinfo\_static\_t} {\bf{}dual\_number}::{\it{}return\_type\_tinfo\_static}[256] = {\nwlbrace}{\nwlbrace}{\nwrbrace}{\nwrbrace};\nwindexdefn{\nwixident{tinfo{\_}static{\_}t}}{tinfo:unstatic:unt}{NW4P20f9-3uuIiL-1}
{\bf{}\char35{}endif}{\tt{}}

\nwalsodefined{\\{NW4P20f9-3uuIiL-2}\\{NW4P20f9-3uuIiL-3}\\{NW4P20f9-3uuIiL-4}\\{NW4P20f9-3uuIiL-5}\\{NW4P20f9-3uuIiL-6}\\{NW4P20f9-3uuIiL-7}\\{NW4P20f9-3uuIiL-8}\\{NW4P20f9-3uuIiL-9}\\{NW4P20f9-3uuIiL-A}\\{NW4P20f9-3uuIiL-B}\\{NW4P20f9-3uuIiL-C}\\{NW4P20f9-3uuIiL-D}\\{NW4P20f9-3uuIiL-E}\\{NW4P20f9-3uuIiL-F}\\{NW4P20f9-3uuIiL-G}\\{NW4P20f9-3uuIiL-H}\\{NW4P20f9-3uuIiL-I}\\{NW4P20f9-3uuIiL-J}\\{NW4P20f9-3uuIiL-K}}\nwused{\\{NW4P20f9-3HhVai-1}}\nwidentdefs{\\{{\nwixident{tinfo{\_}static{\_}t}}{tinfo:unstatic:unt}}}\nwidentuses{\\{{\nwixident{dual{\_}number}}{dual:unnumber}}\\{{\nwixident{GINAC{\_}VERSION{\_}ATLEAST}}{GINAC:unVERSION:unATLEAST}}}\nwindexuse{\nwixident{dual{\_}number}}{dual:unnumber}{NW4P20f9-3uuIiL-1}\nwindexuse{\nwixident{GINAC{\_}VERSION{\_}ATLEAST}}{GINAC:unVERSION:unATLEAST}{NW4P20f9-3uuIiL-1}\nwendcode{}\nwbegindocs{158}\nwdocspar
\subsubsection{Implementation of Constructors}
\label{sec:impl-constr}
Default constructor.
\nwenddocs{}\nwbegincode{159}\sublabel{NW4P20f9-3uuIiL-2}\nwmargintag{{\nwtagstyle{}\subpageref{NW4P20f9-3uuIiL-2}}}\moddef{Dual number class implementation~{\nwtagstyle{}\subpageref{NW4P20f9-3uuIiL-1}}}\plusendmoddef\Rm{}\nwstartdeflinemarkup\nwusesondefline{\\{NW4P20f9-3HhVai-1}}\nwprevnextdefs{NW4P20f9-3uuIiL-1}{NW4P20f9-3uuIiL-3}\nwenddeflinemarkup
{\bf{}\char35{}if}{\tt{} GINAC\_VERSION\_ATLEAST(1,5)}
{\bf{}dual\_number}::{\bf{}dual\_number}() : {\it{}inherited}(), {\it{}u\_comp}(0), {\it{}v\_comp}(0)
{\bf{}\char35{}else}{\tt{}}
{\bf{}dual\_number}::{\bf{}dual\_number}() : {\it{}inherited}(&{\bf{}dual\_number}::{\it{}tinfo\_static}), {\it{}u\_comp}(0), {\it{}v\_comp}(0)
{\bf{}\char35{}endif}{\tt{}}
{\nwlbrace}
        {\it{}setflag}({\it{}status\_flags}::{\it{}not\_shareable});
{\nwrbrace}

\nwused{\\{NW4P20f9-3HhVai-1}}\nwidentuses{\\{{\nwixident{dual{\_}number}}{dual:unnumber}}\\{{\nwixident{GINAC{\_}VERSION{\_}ATLEAST}}{GINAC:unVERSION:unATLEAST}}}\nwindexuse{\nwixident{dual{\_}number}}{dual:unnumber}{NW4P20f9-3uuIiL-2}\nwindexuse{\nwixident{GINAC{\_}VERSION{\_}ATLEAST}}{GINAC:unVERSION:unATLEAST}{NW4P20f9-3uuIiL-2}\nwendcode{}\nwbegindocs{160}Constructor from two components.
\nwenddocs{}\nwbegincode{161}\sublabel{NW4P20f9-3uuIiL-3}\nwmargintag{{\nwtagstyle{}\subpageref{NW4P20f9-3uuIiL-3}}}\moddef{Dual number class implementation~{\nwtagstyle{}\subpageref{NW4P20f9-3uuIiL-1}}}\plusendmoddef\Rm{}\nwstartdeflinemarkup\nwusesondefline{\\{NW4P20f9-3HhVai-1}}\nwprevnextdefs{NW4P20f9-3uuIiL-2}{NW4P20f9-3uuIiL-4}\nwenddeflinemarkup
{\bf{}\char35{}if}{\tt{} GINAC\_VERSION\_ATLEAST(1,5)}
{\bf{}dual\_number}::{\bf{}dual\_number}({\bf{}const} {\bf{}ex} & {\it{}a}, {\bf{}const} {\bf{}ex} & {\it{}b}) : {\it{}inherited}(), 
{\it{}u\_comp}({\it{}a}), {\it{}v\_comp}({\it{}b})
{\bf{}\char35{}else}{\tt{}}
{\bf{}dual\_number}::{\bf{}dual\_number}({\bf{}const} {\bf{}ex} & {\it{}a}, {\bf{}const} {\bf{}ex} & {\it{}b}) : {\it{}inherited}(&{\bf{}dual\_number}::{\it{}tinfo\_static}), 
{\it{}u\_comp}({\it{}a}), {\it{}v\_comp}({\it{}b})
{\bf{}\char35{}endif}{\tt{}}
{\nwlbrace}
{\nwrbrace}

\nwused{\\{NW4P20f9-3HhVai-1}}\nwidentuses{\\{{\nwixident{dual{\_}number}}{dual:unnumber}}\\{{\nwixident{GINAC{\_}VERSION{\_}ATLEAST}}{GINAC:unVERSION:unATLEAST}}}\nwindexuse{\nwixident{dual{\_}number}}{dual:unnumber}{NW4P20f9-3uuIiL-3}\nwindexuse{\nwixident{GINAC{\_}VERSION{\_}ATLEAST}}{GINAC:unVERSION:unATLEAST}{NW4P20f9-3uuIiL-3}\nwendcode{}\nwbegindocs{162}Constructor from a single expression. It may contain two components\ldots
\nwenddocs{}\nwbegincode{163}\sublabel{NW4P20f9-3uuIiL-4}\nwmargintag{{\nwtagstyle{}\subpageref{NW4P20f9-3uuIiL-4}}}\moddef{Dual number class implementation~{\nwtagstyle{}\subpageref{NW4P20f9-3uuIiL-1}}}\plusendmoddef\Rm{}\nwstartdeflinemarkup\nwusesondefline{\\{NW4P20f9-3HhVai-1}}\nwprevnextdefs{NW4P20f9-3uuIiL-3}{NW4P20f9-3uuIiL-5}\nwenddeflinemarkup
{\bf{}\char35{}if}{\tt{} GINAC\_VERSION\_ATLEAST(1,5)}
{\bf{}dual\_number}::{\bf{}dual\_number}({\bf{}const} {\bf{}ex} & {\it{}P}) : {\it{}inherited}()
{\bf{}\char35{}else}{\tt{}}
{\bf{}dual\_number}::{\bf{}dual\_number}({\bf{}const} {\bf{}ex} & {\it{}P}) : {\it{}inherited}(&{\bf{}dual\_number}::{\it{}tinfo\_static})
{\bf{}\char35{}endif}{\tt{}}
{\nwlbrace}
        {\bf{}if} ({\it{}is\_a}\begin{math}<\end{math}{\bf{}lst}\begin{math}>\end{math}({\it{}P}) \begin{math}\vee\end{math} {\it{}is\_a}\begin{math}<\end{math}{\bf{}matrix}\begin{math}>\end{math}({\it{}P}) \begin{math}\vee\end{math} {\it{}is\_a}\begin{math}<\end{math}{\bf{}dual\_number}\begin{math}>\end{math}({\it{}P})) {\nwlbrace}
                {\it{}u\_comp} = {\it{}P}.{\it{}op}(0);
                {\it{}v\_comp} = {\it{}P}.{\it{}op}(1);

\nwused{\\{NW4P20f9-3HhVai-1}}\nwidentuses{\\{{\nwixident{dual{\_}number}}{dual:unnumber}}\\{{\nwixident{GINAC{\_}VERSION{\_}ATLEAST}}{GINAC:unVERSION:unATLEAST}}}\nwindexuse{\nwixident{dual{\_}number}}{dual:unnumber}{NW4P20f9-3uuIiL-4}\nwindexuse{\nwixident{GINAC{\_}VERSION{\_}ATLEAST}}{GINAC:unVERSION:unATLEAST}{NW4P20f9-3uuIiL-4}\nwendcode{}\nwbegindocs{164}\ldots if it is a real expression we embed it into {\Tt{}\Rm{}{\bf{}dual\_number}\nwendquote}\ldots
\nwenddocs{}\nwbegincode{165}\sublabel{NW4P20f9-3uuIiL-5}\nwmargintag{{\nwtagstyle{}\subpageref{NW4P20f9-3uuIiL-5}}}\moddef{Dual number class implementation~{\nwtagstyle{}\subpageref{NW4P20f9-3uuIiL-1}}}\plusendmoddef\Rm{}\nwstartdeflinemarkup\nwusesondefline{\\{NW4P20f9-3HhVai-1}}\nwprevnextdefs{NW4P20f9-3uuIiL-4}{NW4P20f9-3uuIiL-6}\nwenddeflinemarkup
        {\nwrbrace} {\bf{}else} {\bf{}if} ({\it{}P}.{\it{}imag\_part}().{\it{}normal}().{\it{}is\_zero}()) {\nwlbrace}
                {\bf{}dual\_number} {\it{}t}={\it{}dn\_from\_arg\_mod}({\it{}Arg0}, {\it{}P});
                {\it{}u\_comp} = {\it{}t}.{\it{}op}(0);
                {\it{}v\_comp} = {\it{}t}.{\it{}op}(1);

\nwused{\\{NW4P20f9-3HhVai-1}}\nwidentuses{\\{{\nwixident{Arg0}}{Arg0}}\\{{\nwixident{dual{\_}number}}{dual:unnumber}}}\nwindexuse{\nwixident{Arg0}}{Arg0}{NW4P20f9-3uuIiL-5}\nwindexuse{\nwixident{dual{\_}number}}{dual:unnumber}{NW4P20f9-3uuIiL-5}\nwendcode{}\nwbegindocs{166}\ldots or if its a complex expression we decompose it into the real and
imaginary parts.
\nwenddocs{}\nwbegincode{167}\sublabel{NW4P20f9-3uuIiL-6}\nwmargintag{{\nwtagstyle{}\subpageref{NW4P20f9-3uuIiL-6}}}\moddef{Dual number class implementation~{\nwtagstyle{}\subpageref{NW4P20f9-3uuIiL-1}}}\plusendmoddef\Rm{}\nwstartdeflinemarkup\nwusesondefline{\\{NW4P20f9-3HhVai-1}}\nwprevnextdefs{NW4P20f9-3uuIiL-5}{NW4P20f9-3uuIiL-7}\nwenddeflinemarkup
        {\nwrbrace} {\bf{}else} {\nwlbrace}
                {\it{}u\_comp} = {\it{}P}.{\it{}real\_part}();
                {\it{}v\_comp} = {\it{}P}.{\it{}imag\_part}();
        {\nwrbrace}
{\nwrbrace}

\nwused{\\{NW4P20f9-3HhVai-1}}\nwendcode{}\nwbegindocs{168}\nwdocspar
\subsubsection[Service methods for the GiNaC infrastructure]{Service methods for the {\GiNaC} infrastructure}
\label{sec:serv-meth-ginac}

Standard parts involving archiving, comparison and printing of the {\Tt{}\Rm{}{\bf{}dual\_number}\nwendquote} class
\nwenddocs{}\nwbegincode{169}\sublabel{NW4P20f9-3uuIiL-7}\nwmargintag{{\nwtagstyle{}\subpageref{NW4P20f9-3uuIiL-7}}}\moddef{Dual number class implementation~{\nwtagstyle{}\subpageref{NW4P20f9-3uuIiL-1}}}\plusendmoddef\Rm{}\nwstartdeflinemarkup\nwusesondefline{\\{NW4P20f9-3HhVai-1}}\nwprevnextdefs{NW4P20f9-3uuIiL-6}{NW4P20f9-3uuIiL-8}\nwenddeflinemarkup

{\bf{}dual\_number}::{\bf{}dual\_number}({\bf{}const} {\it{}archive\_node} &{\it{}n}, {\bf{}lst} &{\it{}sym\_lst}) : {\it{}inherited}()
{\nwlbrace}
        {\it{}n}.{\it{}find\_ex}({\tt{}"u\_comp"}, {\it{}u\_comp}, {\it{}sym\_lst});
        {\it{}n}.{\it{}find\_ex}({\tt{}"v\_comp"}, {\it{}v\_comp}, {\it{}sym\_lst});
{\nwrbrace}

\nwused{\\{NW4P20f9-3HhVai-1}}\nwidentuses{\\{{\nwixident{dual{\_}number}}{dual:unnumber}}}\nwindexuse{\nwixident{dual{\_}number}}{dual:unnumber}{NW4P20f9-3uuIiL-7}\nwendcode{}\nwbegindocs{170}Archiving routine.
\nwenddocs{}\nwbegincode{171}\sublabel{NW4P20f9-3uuIiL-8}\nwmargintag{{\nwtagstyle{}\subpageref{NW4P20f9-3uuIiL-8}}}\moddef{Dual number class implementation~{\nwtagstyle{}\subpageref{NW4P20f9-3uuIiL-1}}}\plusendmoddef\Rm{}\nwstartdeflinemarkup\nwusesondefline{\\{NW4P20f9-3HhVai-1}}\nwprevnextdefs{NW4P20f9-3uuIiL-7}{NW4P20f9-3uuIiL-9}\nwenddeflinemarkup
{\bf{}void} {\bf{}dual\_number}::{\it{}archive}({\it{}archive\_node} &{\it{}n}) {\bf{}const}\nwindexdefn{\nwixident{dual{\_}number}}{dual:unnumber}{NW4P20f9-3uuIiL-8}
{\nwlbrace}
        {\it{}inherited}::{\it{}archive}({\it{}n});
        {\it{}n}.{\it{}add\_ex}({\tt{}"u\_comp"}, {\it{}u\_comp});
        {\it{}n}.{\it{}add\_ex}({\tt{}"v\_comp"}, {\it{}v\_comp});
{\nwrbrace}

\nwused{\\{NW4P20f9-3HhVai-1}}\nwidentdefs{\\{{\nwixident{dual{\_}number}}{dual:unnumber}}}\nwendcode{}\nwbegindocs{172}Archiving routine.
\nwenddocs{}\nwbegincode{173}\sublabel{NW4P20f9-3uuIiL-9}\nwmargintag{{\nwtagstyle{}\subpageref{NW4P20f9-3uuIiL-9}}}\moddef{Dual number class implementation~{\nwtagstyle{}\subpageref{NW4P20f9-3uuIiL-1}}}\plusendmoddef\Rm{}\nwstartdeflinemarkup\nwusesondefline{\\{NW4P20f9-3HhVai-1}}\nwprevnextdefs{NW4P20f9-3uuIiL-8}{NW4P20f9-3uuIiL-A}\nwenddeflinemarkup
{\bf{}ex} {\bf{}dual\_number}::{\it{}unarchive}({\bf{}const} {\it{}archive\_node} &{\it{}n}, {\bf{}lst} &{\it{}sym\_lst})
{\nwlbrace}
        {\bf{}return} ({\bf{}new} {\bf{}dual\_number}({\it{}n}, {\it{}sym\_lst}))\begin{math}\rightarrow\end{math}{\it{}setflag}({\it{}status\_flags}::{\it{}dynallocated});
{\nwrbrace}

\nwused{\\{NW4P20f9-3HhVai-1}}\nwidentuses{\\{{\nwixident{dual{\_}number}}{dual:unnumber}}}\nwindexuse{\nwixident{dual{\_}number}}{dual:unnumber}{NW4P20f9-3uuIiL-9}\nwendcode{}\nwbegindocs{174}Comparison routine.
\nwenddocs{}\nwbegincode{175}\sublabel{NW4P20f9-3uuIiL-A}\nwmargintag{{\nwtagstyle{}\subpageref{NW4P20f9-3uuIiL-A}}}\moddef{Dual number class implementation~{\nwtagstyle{}\subpageref{NW4P20f9-3uuIiL-1}}}\plusendmoddef\Rm{}\nwstartdeflinemarkup\nwusesondefline{\\{NW4P20f9-3HhVai-1}}\nwprevnextdefs{NW4P20f9-3uuIiL-9}{NW4P20f9-3uuIiL-B}\nwenddeflinemarkup
{\bf{}int} {\bf{}dual\_number}::{\it{}compare\_same\_type}({\bf{}const} {\bf{}basic} & {\it{}other}) {\bf{}const}\nwindexdefn{\nwixident{dual{\_}number}}{dual:unnumber}{NW4P20f9-3uuIiL-A}
{\nwlbrace}
        {\it{}GINAC\_ASSERT}({\it{}is\_a}\begin{math}<\end{math}{\bf{}dual\_number}\begin{math}>\end{math}({\it{}other}));
        {\bf{}const} {\bf{}dual\_number} &{\it{}o} = {\bf{}static\_cast}\begin{math}<\end{math}{\bf{}const} {\bf{}dual\_number} &\begin{math}>\end{math}({\it{}other});
        
        {\bf{}int} {\it{}cmpval} = {\it{}u\_comp}.{\it{}compare}({\it{}o}.{\it{}op}(0));
        {\bf{}if} ({\it{}cmpval}\begin{math}\neq\end{math}0) {\bf{}return} {\it{}cmpval};

        {\bf{}return} {\it{}v\_comp}.{\it{}compare}({\it{}o}.{\it{}op}(1));
{\nwrbrace}

\nwused{\\{NW4P20f9-3HhVai-1}}\nwidentdefs{\\{{\nwixident{dual{\_}number}}{dual:unnumber}}}\nwendcode{}\nwbegindocs{176}Equality of two dual numbers.
\nwenddocs{}\nwbegincode{177}\sublabel{NW4P20f9-3uuIiL-B}\nwmargintag{{\nwtagstyle{}\subpageref{NW4P20f9-3uuIiL-B}}}\moddef{Dual number class implementation~{\nwtagstyle{}\subpageref{NW4P20f9-3uuIiL-1}}}\plusendmoddef\Rm{}\nwstartdeflinemarkup\nwusesondefline{\\{NW4P20f9-3HhVai-1}}\nwprevnextdefs{NW4P20f9-3uuIiL-A}{NW4P20f9-3uuIiL-C}\nwenddeflinemarkup
{\bf{}bool} {\bf{}dual\_number}::{\it{}is\_equal}({\bf{}const} {\bf{}ex} & {\it{}other}) {\bf{}const}
{\nwlbrace}
        {\it{}GINAC\_ASSERT}({\it{}is\_a}\begin{math}<\end{math}{\bf{}dual\_number}\begin{math}>\end{math}({\it{}other}));
        {\bf{}const} {\bf{}dual\_number} &{\it{}o} = {\bf{}static\_cast}\begin{math}<\end{math}{\bf{}const} {\bf{}dual\_number} &\begin{math}>\end{math}({\it{}other});
        
        {\bf{}return} ({\it{}u\_comp}-{\it{}o}.{\it{}op}(0)).{\it{}normal}().{\it{}is\_zero}() \begin{math}\wedge\end{math} ({\it{}v\_comp}-{\it{}o}.{\it{}op}(1)).{\it{}normal}().{\it{}is\_zero}();
{\nwrbrace}

\nwused{\\{NW4P20f9-3HhVai-1}}\nwidentuses{\\{{\nwixident{dual{\_}number}}{dual:unnumber}}}\nwindexuse{\nwixident{dual{\_}number}}{dual:unnumber}{NW4P20f9-3uuIiL-B}\nwendcode{}\nwbegindocs{178}\nwdocspar
\subsubsection{Printing}
\label{sec:printing}
Default printing.
\nwenddocs{}\nwbegincode{179}\sublabel{NW4P20f9-3uuIiL-C}\nwmargintag{{\nwtagstyle{}\subpageref{NW4P20f9-3uuIiL-C}}}\moddef{Dual number class implementation~{\nwtagstyle{}\subpageref{NW4P20f9-3uuIiL-1}}}\plusendmoddef\Rm{}\nwstartdeflinemarkup\nwusesondefline{\\{NW4P20f9-3HhVai-1}}\nwprevnextdefs{NW4P20f9-3uuIiL-B}{NW4P20f9-3uuIiL-D}\nwenddeflinemarkup
{\bf{}void} {\bf{}dual\_number}::{\it{}do\_print}({\bf{}const} {\it{}print\_context} & {\it{}c}, {\bf{}unsigned} {\it{}level}) {\bf{}const}\nwindexdefn{\nwixident{dual{\_}number}}{dual:unnumber}{NW4P20f9-3uuIiL-C}
{\nwlbrace}
        {\it{}c}.{\it{}s} \begin{math}\ll\end{math} {\tt{}"["};
        {\it{}u\_comp}.{\it{}print}({\it{}c});
        {\it{}c}.{\it{}s} \begin{math}\ll\end{math} {\tt{}","};
        {\it{}v\_comp}.{\it{}print}({\it{}c});
        {\it{}c}.{\it{}s} \begin{math}\ll\end{math} {\tt{}"]"};
{\nwrbrace}

\nwused{\\{NW4P20f9-3HhVai-1}}\nwidentdefs{\\{{\nwixident{dual{\_}number}}{dual:unnumber}}}\nwendcode{}\nwbegindocs{180}\LaTeX\ printing.
\nwenddocs{}\nwbegincode{181}\sublabel{NW4P20f9-3uuIiL-D}\nwmargintag{{\nwtagstyle{}\subpageref{NW4P20f9-3uuIiL-D}}}\moddef{Dual number class implementation~{\nwtagstyle{}\subpageref{NW4P20f9-3uuIiL-1}}}\plusendmoddef\Rm{}\nwstartdeflinemarkup\nwusesondefline{\\{NW4P20f9-3HhVai-1}}\nwprevnextdefs{NW4P20f9-3uuIiL-C}{NW4P20f9-3uuIiL-E}\nwenddeflinemarkup
{\bf{}void} {\bf{}dual\_number}::{\it{}do\_print\_latex}({\bf{}const} {\it{}print\_latex} & {\it{}c}, {\bf{}unsigned} {\it{}level}) {\bf{}const}\nwindexdefn{\nwixident{dual{\_}number}}{dual:unnumber}{NW4P20f9-3uuIiL-D}
{\nwlbrace}
        {\it{}c}.{\it{}s} \begin{math}\ll\end{math} {\tt{}"{\char92}{\char92}left({\char92}{\char92}begin{\char123}array{\char125}{\char123}cc{\char125}"};
        {\it{}u\_comp}.{\it{}print}({\it{}c});
        {\it{}c}.{\it{}s} \begin{math}\ll\end{math} {\tt{}"&"};
        {\it{}v\_comp}.{\it{}print}({\it{}c});
        {\it{}c}.{\it{}s} \begin{math}\ll\end{math} {\tt{}"{\char92}{\char92}end{\char123}array{\char125}{\char92}{\char92}right)"};
{\nwrbrace}

\nwused{\\{NW4P20f9-3HhVai-1}}\nwidentdefs{\\{{\nwixident{dual{\_}number}}{dual:unnumber}}}\nwendcode{}\nwbegindocs{182}\nwdocspar
\subsubsection{Overloading algebraic operations}
\label{sec:overl-algebr-oper}
Addition.
\nwenddocs{}\nwbegincode{183}\sublabel{NW4P20f9-3uuIiL-E}\nwmargintag{{\nwtagstyle{}\subpageref{NW4P20f9-3uuIiL-E}}}\moddef{Dual number class implementation~{\nwtagstyle{}\subpageref{NW4P20f9-3uuIiL-1}}}\plusendmoddef\Rm{}\nwstartdeflinemarkup\nwusesondefline{\\{NW4P20f9-3HhVai-1}}\nwprevnextdefs{NW4P20f9-3uuIiL-D}{NW4P20f9-3uuIiL-F}\nwenddeflinemarkup
{\bf{}const} {\bf{}dual\_number} {\bf{}operator}+({\bf{}const} {\bf{}dual\_number} & {\it{}lh}, {\bf{}const} {\bf{}dual\_number} & {\it{}rh}) \nwindexdefn{\nwixident{dual{\_}number}}{dual:unnumber}{NW4P20f9-3uuIiL-E}
{\nwlbrace}
        {\bf{}return} {\it{}lh}.{\it{}add}({\it{}rh});
{\nwrbrace}

{\bf{}const} {\bf{}dual\_number} {\bf{}operator}-({\bf{}const} {\bf{}dual\_number} & {\it{}lh}, {\bf{}const} {\bf{}dual\_number} & {\it{}rh})\nwindexdefn{\nwixident{dual{\_}number}}{dual:unnumber}{NW4P20f9-3uuIiL-E}
{\nwlbrace}
        {\bf{}return} {\it{}lh}.{\it{}sub}({\it{}rh});
{\nwrbrace}

\nwused{\\{NW4P20f9-3HhVai-1}}\nwidentdefs{\\{{\nwixident{dual{\_}number}}{dual:unnumber}}}\nwendcode{}\nwbegindocs{184}Multiplication.
\nwenddocs{}\nwbegincode{185}\sublabel{NW4P20f9-3uuIiL-F}\nwmargintag{{\nwtagstyle{}\subpageref{NW4P20f9-3uuIiL-F}}}\moddef{Dual number class implementation~{\nwtagstyle{}\subpageref{NW4P20f9-3uuIiL-1}}}\plusendmoddef\Rm{}\nwstartdeflinemarkup\nwusesondefline{\\{NW4P20f9-3HhVai-1}}\nwprevnextdefs{NW4P20f9-3uuIiL-E}{NW4P20f9-3uuIiL-G}\nwenddeflinemarkup
{\bf{}const} {\bf{}dual\_number} {\bf{}operator}\begin{math}\ast\end{math}({\bf{}const} {\bf{}dual\_number} & {\it{}lh}, {\bf{}const} {\bf{}ex} & {\it{}rh})\nwindexdefn{\nwixident{dual{\_}number}}{dual:unnumber}{NW4P20f9-3uuIiL-F}
{\nwlbrace}
        {\bf{}return} {\it{}lh}.{\it{}mul}({\it{}rh});
{\nwrbrace}

{\bf{}const} {\bf{}dual\_number} {\bf{}operator}\begin{math}\ast\end{math}({\bf{}const} {\bf{}ex} & {\it{}lh}, {\bf{}const} {\bf{}dual\_number} & {\it{}rh})\nwindexdefn{\nwixident{dual{\_}number}}{dual:unnumber}{NW4P20f9-3uuIiL-F}
{\nwlbrace}
        {\bf{}return} {\it{}rh}.{\it{}mul}({\it{}lh});
{\nwrbrace}

{\bf{}const} {\bf{}dual\_number} {\bf{}operator}\begin{math}\ast\end{math}({\bf{}const} {\bf{}dual\_number} & {\it{}lh}, {\bf{}const} {\bf{}dual\_number} & {\it{}rh})\nwindexdefn{\nwixident{dual{\_}number}}{dual:unnumber}{NW4P20f9-3uuIiL-F}
{\nwlbrace}
        {\bf{}return} {\it{}lh}.{\it{}mul}({\it{}rh});
{\nwrbrace}

\nwused{\\{NW4P20f9-3HhVai-1}}\nwidentdefs{\\{{\nwixident{dual{\_}number}}{dual:unnumber}}}\nwendcode{}\nwbegindocs{186}Division.
\nwenddocs{}\nwbegincode{187}\sublabel{NW4P20f9-3uuIiL-G}\nwmargintag{{\nwtagstyle{}\subpageref{NW4P20f9-3uuIiL-G}}}\moddef{Dual number class implementation~{\nwtagstyle{}\subpageref{NW4P20f9-3uuIiL-1}}}\plusendmoddef\Rm{}\nwstartdeflinemarkup\nwusesondefline{\\{NW4P20f9-3HhVai-1}}\nwprevnextdefs{NW4P20f9-3uuIiL-F}{NW4P20f9-3uuIiL-H}\nwenddeflinemarkup
{\bf{}const} {\bf{}dual\_number} {\bf{}operator}\begin{math}\div\end{math}({\bf{}const} {\bf{}dual\_number} & {\it{}lh}, {\bf{}const} {\bf{}dual\_number} & {\it{}rh})\nwindexdefn{\nwixident{dual{\_}number}}{dual:unnumber}{NW4P20f9-3uuIiL-G}
{\nwlbrace}
        {\bf{}return} {\it{}lh}.{\it{}mul}({\it{}rh}.{\it{}power}(-1));
{\nwrbrace}

{\bf{}const} {\bf{}dual\_number} {\bf{}operator}\begin{math}\div\end{math}({\bf{}const} {\bf{}ex} & {\it{}lh}, {\bf{}const} {\bf{}dual\_number} & {\it{}rh})\nwindexdefn{\nwixident{dual{\_}number}}{dual:unnumber}{NW4P20f9-3uuIiL-G}
{\nwlbrace}
        {\bf{}return} {\it{}rh}.{\it{}power}(-1)\begin{math}\ast\end{math}{\it{}lh};
{\nwrbrace}

{\bf{}const} {\bf{}dual\_number} {\bf{}operator}\begin{math}\div\end{math}({\bf{}const} {\bf{}dual\_number} & {\it{}lh}, {\bf{}const} {\bf{}ex} & {\it{}rh})\nwindexdefn{\nwixident{dual{\_}number}}{dual:unnumber}{NW4P20f9-3uuIiL-G}
{\nwlbrace}
        {\bf{}return} {\it{}lh}.{\it{}mul}({\it{}pow}({\it{}rh}, -1));
{\nwrbrace}

\nwused{\\{NW4P20f9-3HhVai-1}}\nwidentdefs{\\{{\nwixident{dual{\_}number}}{dual:unnumber}}}\nwendcode{}\nwbegindocs{188}\nwdocspar

\subsubsection{Component-related functions}
\label{sec:comp-relat-funct}

\nwenddocs{}\nwbegincode{189}\sublabel{NW4P20f9-3uuIiL-H}\nwmargintag{{\nwtagstyle{}\subpageref{NW4P20f9-3uuIiL-H}}}\moddef{Dual number class implementation~{\nwtagstyle{}\subpageref{NW4P20f9-3uuIiL-1}}}\plusendmoddef\Rm{}\nwstartdeflinemarkup\nwusesondefline{\\{NW4P20f9-3HhVai-1}}\nwprevnextdefs{NW4P20f9-3uuIiL-G}{NW4P20f9-3uuIiL-I}\nwenddeflinemarkup
{\bf{}ex} {\bf{}dual\_number}::{\it{}op}({\it{}size\_t} {\it{}i}) {\bf{}const}
{\nwlbrace}
        {\it{}GINAC\_ASSERT}({\it{}i}\begin{math}<\end{math}{\it{}nops}());
        
        {\bf{}switch} ({\it{}i}) {\nwlbrace}
        {\bf{}case} 0:
                {\bf{}return} {\it{}u\_comp};
        {\bf{}case} 1:
                {\bf{}return} {\it{}v\_comp};
        {\bf{}default}:
                {\bf{}throw}({\it{}std}::{\it{}invalid\_argument}({\tt{}"dual\_number::op(): requested"} 
                                                                        {\tt{}" operand out of the range (2)"}));
        {\nwrbrace}
{\nwrbrace}

\nwused{\\{NW4P20f9-3HhVai-1}}\nwidentuses{\\{{\nwixident{dual{\_}number}}{dual:unnumber}}}\nwindexuse{\nwixident{dual{\_}number}}{dual:unnumber}{NW4P20f9-3uuIiL-H}\nwendcode{}\nwbegindocs{190}\nwdocspar
\nwenddocs{}\nwbegincode{191}\sublabel{NW4P20f9-3uuIiL-I}\nwmargintag{{\nwtagstyle{}\subpageref{NW4P20f9-3uuIiL-I}}}\moddef{Dual number class implementation~{\nwtagstyle{}\subpageref{NW4P20f9-3uuIiL-1}}}\plusendmoddef\Rm{}\nwstartdeflinemarkup\nwusesondefline{\\{NW4P20f9-3HhVai-1}}\nwprevnextdefs{NW4P20f9-3uuIiL-H}{NW4P20f9-3uuIiL-J}\nwenddeflinemarkup
{\bf{}dual\_number} {\bf{}dual\_number}::{\it{}subs}({\bf{}const} {\bf{}ex} & {\it{}e}, {\bf{}unsigned} {\it{}options}) {\bf{}const}
{\nwlbrace}
        {\it{}exmap} {\it{}m};
        {\bf{}if} ({\it{}e}.{\it{}info}({\it{}info\_flags}::{\it{}list})) {\nwlbrace}
                {\bf{}lst} {\it{}l} = {\it{}ex\_to}\begin{math}<\end{math}{\bf{}lst}\begin{math}>\end{math}({\it{}e});
                {\bf{}for} ({\bf{}lst}::{\it{}const\_iterator} {\it{}i} = {\it{}l}.{\it{}begin}(); {\it{}i} \begin{math}\neq\end{math} {\it{}l}.{\it{}end}(); \protect\PP{\it{}i})
                        {\it{}m}.{\it{}insert}({\it{}std}::{\it{}make\_pair}({\it{}i}\begin{math}\rightarrow\end{math}{\it{}op}(0), {\it{}i}\begin{math}\rightarrow\end{math}{\it{}op}(1)));
        {\nwrbrace} {\bf{}else} {\bf{}if} ({\it{}is\_a}\begin{math}<\end{math}{\bf{}relational}\begin{math}>\end{math}({\it{}e})) {\nwlbrace}
                {\it{}m}.{\it{}insert}({\it{}std}::{\it{}make\_pair}({\it{}e}.{\it{}op}(0), {\it{}e}.{\it{}op}(1)));
        {\nwrbrace} {\bf{}else}
                {\bf{}throw}({\it{}std}::{\it{}invalid\_argument}({\tt{}"dual\_number::subs(): the parameter"}
                                                                        {\tt{}" should be a relational or a lst"}));
        
 {\bf{}return} {\it{}ex\_to}\begin{math}<\end{math}{\bf{}dual\_number}\begin{math}>\end{math}({\it{}inherited}::{\it{}subs}({\it{}m}, {\it{}options}));
{\nwrbrace}

\nwused{\\{NW4P20f9-3HhVai-1}}\nwidentuses{\\{{\nwixident{dual{\_}number}}{dual:unnumber}}}\nwindexuse{\nwixident{dual{\_}number}}{dual:unnumber}{NW4P20f9-3uuIiL-I}\nwendcode{}\nwbegindocs{192}\nwdocspar
\nwenddocs{}\nwbegincode{193}\sublabel{NW4P20f9-3uuIiL-J}\nwmargintag{{\nwtagstyle{}\subpageref{NW4P20f9-3uuIiL-J}}}\moddef{Dual number class implementation~{\nwtagstyle{}\subpageref{NW4P20f9-3uuIiL-1}}}\plusendmoddef\Rm{}\nwstartdeflinemarkup\nwusesondefline{\\{NW4P20f9-3HhVai-1}}\nwprevnextdefs{NW4P20f9-3uuIiL-I}{NW4P20f9-3uuIiL-K}\nwenddeflinemarkup
{\bf{}ex} & {\bf{}dual\_number}::{\it{}let\_op}({\it{}size\_t} {\it{}i})
{\nwlbrace}
        {\it{}GINAC\_ASSERT}({\it{}i}\begin{math}<\end{math}{\it{}nops}());
        
        {\it{}ensure\_if\_modifiable}();
        {\bf{}switch} ({\it{}i}) {\nwlbrace}
        {\bf{}case} 0:
                {\bf{}return} {\it{}u\_comp};
        {\bf{}case} 1:
                {\bf{}return} {\it{}v\_comp};
        {\bf{}default}:
                {\bf{}throw}({\it{}std}::{\it{}invalid\_argument}({\tt{}"dual\_number::let\_op(): requested operand"}
                                                                        {\tt{}" out of the range (2)"}));
        {\nwrbrace}
{\nwrbrace}

\nwused{\\{NW4P20f9-3HhVai-1}}\nwidentuses{\\{{\nwixident{dual{\_}number}}{dual:unnumber}}}\nwindexuse{\nwixident{dual{\_}number}}{dual:unnumber}{NW4P20f9-3uuIiL-J}\nwendcode{}\nwbegindocs{194}\nwdocspar
\nwenddocs{}\nwbegincode{195}\sublabel{NW4P20f9-3uuIiL-K}\nwmargintag{{\nwtagstyle{}\subpageref{NW4P20f9-3uuIiL-K}}}\moddef{Dual number class implementation~{\nwtagstyle{}\subpageref{NW4P20f9-3uuIiL-1}}}\plusendmoddef\Rm{}\nwstartdeflinemarkup\nwusesondefline{\\{NW4P20f9-3HhVai-1}}\nwprevnextdefs{NW4P20f9-3uuIiL-J}{\relax}\nwenddeflinemarkup
{\bf{}void} {\it{}set\_subgroup}({\bf{}bool} {\it{}N}) {\nwlbrace}\nwindexdefn{\nwixident{set{\_}subgroup}}{set:unsubgroup}{NW4P20f9-3uuIiL-K}
        {\it{}subgroup\_N} = {\it{}N};
{\nwrbrace}

{\bf{}void} {\it{}set\_straight\_spoke}({\bf{}bool} {\it{}S}) {\nwlbrace}\nwindexdefn{\nwixident{set{\_}straight{\_}spoke}}{set:unstraight:unspoke}{NW4P20f9-3uuIiL-K}
        {\it{}straight\_spoke} = {\it{}S}; 
        {\bf{}if} ({\it{}not} {\it{}S}) 
                {\it{}subgroup\_N}={\bf{}false};
{\nwrbrace}

\nwused{\\{NW4P20f9-3HhVai-1}}\nwidentdefs{\\{{\nwixident{set{\_}straight{\_}spoke}}{set:unstraight:unspoke}}\\{{\nwixident{set{\_}subgroup}}{set:unsubgroup}}}\nwendcode{}

\nwixlogsorted{c}{{Additional routines declarations}{NW4P20f9-NQ5ML-1}{\nwixd{NW4P20f9-NQ5ML-1}\nwixu{NW4P20f9-36Ytqo-1}\nwixd{NW4P20f9-NQ5ML-2}\nwixd{NW4P20f9-NQ5ML-3}\nwixd{NW4P20f9-NQ5ML-4}\nwixd{NW4P20f9-NQ5ML-5}}}%
\nwixlogsorted{c}{{Algebraic procedures}{NW4P20f9-3jzoaG-1}{\nwixd{NW4P20f9-3jzoaG-1}\nwixd{NW4P20f9-3jzoaG-2}\nwixd{NW4P20f9-3jzoaG-3}\nwixd{NW4P20f9-3jzoaG-4}\nwixd{NW4P20f9-3jzoaG-5}\nwixd{NW4P20f9-3jzoaG-6}\nwixu{NW4P20f9-3HhVai-1}}}%
\nwixlogsorted{c}{{Dual number class declaration}{NW4P20f9-4GrAIY-1}{\nwixu{NW4P20f9-36Ytqo-1}\nwixd{NW4P20f9-4GrAIY-1}\nwixd{NW4P20f9-4GrAIY-2}\nwixd{NW4P20f9-4GrAIY-3}\nwixd{NW4P20f9-4GrAIY-4}\nwixd{NW4P20f9-4GrAIY-5}\nwixd{NW4P20f9-4GrAIY-6}}}%
\nwixlogsorted{c}{{Dual number class further implementation}{NW4P20f9-4dxCds-1}{\nwixd{NW4P20f9-4dxCds-1}\nwixd{NW4P20f9-4dxCds-2}\nwixd{NW4P20f9-4dxCds-3}\nwixd{NW4P20f9-4dxCds-4}\nwixd{NW4P20f9-4dxCds-5}\nwixd{NW4P20f9-4dxCds-6}\nwixd{NW4P20f9-4dxCds-7}\nwixd{NW4P20f9-4dxCds-8}\nwixd{NW4P20f9-4dxCds-9}\nwixd{NW4P20f9-4dxCds-A}\nwixd{NW4P20f9-4dxCds-B}\nwixd{NW4P20f9-4dxCds-C}\nwixd{NW4P20f9-4dxCds-D}\nwixd{NW4P20f9-4dxCds-E}\nwixd{NW4P20f9-4dxCds-F}\nwixd{NW4P20f9-4dxCds-G}\nwixd{NW4P20f9-4dxCds-H}\nwixd{NW4P20f9-4dxCds-I}\nwixd{NW4P20f9-4dxCds-J}\nwixd{NW4P20f9-4dxCds-K}\nwixd{NW4P20f9-4dxCds-L}\nwixd{NW4P20f9-4dxCds-M}\nwixu{NW4P20f9-3HhVai-1}}}%
\nwixlogsorted{c}{{Dual number class implementation}{NW4P20f9-3uuIiL-1}{\nwixu{NW4P20f9-3HhVai-1}\nwixd{NW4P20f9-3uuIiL-1}\nwixd{NW4P20f9-3uuIiL-2}\nwixd{NW4P20f9-3uuIiL-3}\nwixd{NW4P20f9-3uuIiL-4}\nwixd{NW4P20f9-3uuIiL-5}\nwixd{NW4P20f9-3uuIiL-6}\nwixd{NW4P20f9-3uuIiL-7}\nwixd{NW4P20f9-3uuIiL-8}\nwixd{NW4P20f9-3uuIiL-9}\nwixd{NW4P20f9-3uuIiL-A}\nwixd{NW4P20f9-3uuIiL-B}\nwixd{NW4P20f9-3uuIiL-C}\nwixd{NW4P20f9-3uuIiL-D}\nwixd{NW4P20f9-3uuIiL-E}\nwixd{NW4P20f9-3uuIiL-F}\nwixd{NW4P20f9-3uuIiL-G}\nwixd{NW4P20f9-3uuIiL-H}\nwixd{NW4P20f9-3uuIiL-I}\nwixd{NW4P20f9-3uuIiL-J}\nwixd{NW4P20f9-3uuIiL-K}}}%
\nwixlogsorted{c}{{Dual number operators}{NW4P20f9-56t00-1}{\nwixu{NW4P20f9-4GrAIY-2}\nwixd{NW4P20f9-56t00-1}}}%
\nwixlogsorted{c}{{Dual number spare operators}{NW4P20f9-2jj2ci-1}{\nwixd{NW4P20f9-2jj2ci-1}}}%
\nwixlogsorted{c}{{dualnum.cpp}{NW4P20f9-3HhVai-1}{\nwixd{NW4P20f9-3HhVai-1}}}%
\nwixlogsorted{c}{{dualnum.h}{NW4P20f9-36Ytqo-1}{\nwixd{NW4P20f9-36Ytqo-1}}}%
\nwixlogsorted{c}{{Header initialisation}{NW4P20f9-1vbmXO-1}{\nwixu{NW4P20f9-36Ytqo-1}\nwixd{NW4P20f9-1vbmXO-1}}}%
\nwixlogsorted{c}{{N-Nprime separation}{NW4P20f9-4EWpKB-1}{\nwixd{NW4P20f9-4EWpKB-1}\nwixd{NW4P20f9-4EWpKB-2}\nwixd{NW4P20f9-4EWpKB-3}\nwixu{NW4P20f9-36Ytqo-1}}}%
\nwixlogsorted{c}{{Output routines}{NW4P20f9-4MrGQl-1}{\nwixu{NW4P20f9-3HhVai-1}\nwixd{NW4P20f9-4MrGQl-1}\nwixd{NW4P20f9-4MrGQl-2}}}%
\nwixlogsorted{c}{{Public methods}{NW4P20f9-1eKCCy-1}{\nwixd{NW4P20f9-1eKCCy-1}\nwixd{NW4P20f9-1eKCCy-2}\nwixd{NW4P20f9-1eKCCy-3}\nwixd{NW4P20f9-1eKCCy-4}\nwixd{NW4P20f9-1eKCCy-5}\nwixd{NW4P20f9-1eKCCy-6}\nwixd{NW4P20f9-1eKCCy-7}\nwixd{NW4P20f9-1eKCCy-8}\nwixd{NW4P20f9-1eKCCy-9}\nwixd{NW4P20f9-1eKCCy-A}\nwixd{NW4P20f9-1eKCCy-B}\nwixd{NW4P20f9-1eKCCy-C}\nwixu{NW4P20f9-4GrAIY-2}}}%
\nwixlogsorted{c}{{Technical methods}{NW4P20f9-2XWNbs-1}{\nwixu{NW4P20f9-4GrAIY-2}\nwixd{NW4P20f9-2XWNbs-1}}}%
\nwixlogsorted{i}{{\nwixident{Arg0}}{Arg0}}%
\nwixlogsorted{i}{{\nwixident{dual{\_}number}}{dual:unnumber}}%
\nwixlogsorted{i}{{\nwixident{formula{\_}out}}{formula:unout}}%
\nwixlogsorted{i}{{\nwixident{GINAC{\_}VERSION{\_}ATLEAST}}{GINAC:unVERSION:unATLEAST}}%
\nwixlogsorted{i}{{\nwixident{numeric}}{numeric}}%
\nwixlogsorted{i}{{\nwixident{set{\_}straight{\_}spoke}}{set:unstraight:unspoke}}%
\nwixlogsorted{i}{{\nwixident{set{\_}subgroup}}{set:unsubgroup}}%
\nwixlogsorted{i}{{\nwixident{test{\_}out}}{test:unout}}%
\nwixlogsorted{i}{{\nwixident{tinfo{\_}static{\_}t}}{tinfo:unstatic:unt}}%
\nwbegindocs{196}This is the end of {\Tt{}\Rm{}{\bf{}dual\_number}\nwendquote} implementation. 
\nwenddocs{}

}{}
 \IfFileExists{parab-rotation.tex}{
 \section{Program for Symbolic Calculations} 
 \label{r-sec:progr-symb-calc}
 This is a documentation for our symbolic calculations supporting this
 paper. You can obtain
 the program itself from the
 \href{http://arXiv.org/abs/math/yymmnnn}{source files} of this
 paper~\cite{Kisil09c} at \url{arXiv.org}; \LaTeX\ compilation of it
 will produces the file \texttt{parab-rotation.nw} in the current directory.
 This is a \NoWEB~\cite{NoWEB} code of the program. It uses Clifford
 algebra facilities~\cite{Kisil05b} of the \GiNaC\
 library~\cite{GiNaC}.


\def\LA{\begingroup\maybehbox\bgroup\setupmodname\Rm$\langle$}\def\RA{$\rangle$\egroup\endgroup}\providecommand{\MM}{\kern.5pt\raisebox{.4ex}{\begin{math}\scriptscriptstyle-\kern-1pt-\end{math}}\kern.5pt}\providecommand{\PP}{\kern.5pt\raisebox{.4ex}{\begin{math}\scriptscriptstyle+\kern-1pt+\end{math}}\kern.5pt}\def\commopen{/\begin{math}\ast\,\end{math}}\def\commclose{\,\begin{math}\ast\end{math}\kern-.5pt/}\def\begcomm{\begingroup\maybehbox\bgroup\setupmodname}\def\endcomm{\egroup\endgroup}\nwfilename{parab-rotation.nw}\nwbegindocs{1}\nwdocspar
\nwenddocs{}\nwbegindocs{2}%
\nwenddocs{}\nwbegindocs{3}%
\nwenddocs{}\nwbegindocs{4}%
\nwenddocs{}\nwbegindocs{5}%
\nwenddocs{}\nwbegindocs{6}%
\nwenddocs{}\nwbegindocs{7}%
\nwenddocs{}\nwbegindocs{8}%
\nwenddocs{}\nwbegindocs{9}%
\nwenddocs{}\nwbegindocs{10}%
\nwenddocs{}\nwbegindocs{11}%
\nwenddocs{}\nwbegindocs{12}%
\nwenddocs{}\nwbegindocs{13}%
\nwenddocs{}\nwbegindocs{14}%
\nwenddocs{}\nwbegindocs{15}%
\nwenddocs{}\nwbegindocs{16}%
\nwenddocs{}\nwbegindocs{17}%
\nwenddocs{}\nwbegindocs{18}%
\nwdocspar
This piece of software is licensed under
\href{http://www.gnu.org/licenses/gpl.html}{GNU General Public
  License} (Version 3, 29 June 2007)~\cite{GNUGPL}.

\nwenddocs{}\nwbegindocs{19}\nwdocspar
\subsection{Calculation and Tests}
This Subsection contains code for calculation of various
expression. See~\cite{Kisil04c} or \GiNaC info for usage of Clifford
algebra functions.

\nwenddocs{}\nwbegindocs{20}\subsubsection{Calculation of Expressions}
\label{sec:calc-expr}
Firstly, we output the expression of the Cayley transform for a generic
element from subgroups \(N\) and \(N^\prime\).
\nwenddocs{}\nwbegincode{21}\sublabel{NWgqRer-1Pc9Jw-1}\nwmargintag{{\nwtagstyle{}\subpageref{NWgqRer-1Pc9Jw-1}}}\moddef{Show expressions~{\nwtagstyle{}\subpageref{NWgqRer-1Pc9Jw-1}}}\endmoddef\Rm{}\nwstartdeflinemarkup\nwusesondefline{\\{NWgqRer-nXe8t-5}}\nwprevnextdefs{\relax}{NWgqRer-1Pc9Jw-2}\nwenddeflinemarkup
{\bf{}ex} {\it{}XC}={\it{}canonicalize\_clifford}(({\it{}TC}\begin{math}\ast\end{math}{\it{}X}\begin{math}\ast\end{math}{\it{}TCI}).{\it{}evalm}());
{\it{}formula\_out}({\tt{}"Cayley of the matrix X: "}, {\it{}XC}.{\it{}subs}({\it{}sign}\begin{math}\equiv\end{math}0).{\it{}normal}());

\nwalsodefined{\\{NWgqRer-1Pc9Jw-2}\\{NWgqRer-1Pc9Jw-3}\\{NWgqRer-1Pc9Jw-4}\\{NWgqRer-1Pc9Jw-5}\\{NWgqRer-1Pc9Jw-6}\\{NWgqRer-1Pc9Jw-7}\\{NWgqRer-1Pc9Jw-8}\\{NWgqRer-1Pc9Jw-9}\\{NWgqRer-1Pc9Jw-A}\\{NWgqRer-1Pc9Jw-B}\\{NWgqRer-1Pc9Jw-C}\\{NWgqRer-1Pc9Jw-D}}\nwused{\\{NWgqRer-nXe8t-5}}\nwidentuses{\\{{\nwixident{formula{\_}out}}{formula:unout}}}\nwindexuse{\nwixident{formula{\_}out}}{formula:unout}{NWgqRer-1Pc9Jw-1}\nwendcode{}\nwbegindocs{22}Then we calculate M\"obius action of those matrix on a point.
\nwenddocs{}\nwbegincode{23}\sublabel{NWgqRer-1Pc9Jw-2}\nwmargintag{{\nwtagstyle{}\subpageref{NWgqRer-1Pc9Jw-2}}}\moddef{Show expressions~{\nwtagstyle{}\subpageref{NWgqRer-1Pc9Jw-1}}}\plusendmoddef\Rm{}\nwstartdeflinemarkup\nwusesondefline{\\{NWgqRer-nXe8t-5}}\nwprevnextdefs{NWgqRer-1Pc9Jw-1}{NWgqRer-1Pc9Jw-3}\nwenddeflinemarkup
{\bf{}dual\_number} {\it{}W}({\it{}clifford\_moebius\_map}({\it{}XC}, {\it{}P}.{\it{}to\_matrix}(), {\it{}e}).{\it{}subs}({\it{}sign}\begin{math}\equiv\end{math}0).{\it{}normal}()),
        {\it{}W1} = {\it{}W}.{\it{}subs}({\bf{}lst}({\it{}u}\begin{math}\equiv\end{math}{\it{}u1}, {\it{}v}\begin{math}\equiv\end{math}{\it{}v1}));
{\it{}formula\_out}({\tt{}"Rotation by x: "}, {\it{}W});
{\bf{}if} ({\it{}not} {\it{}W}.{\it{}is\_equal}({\it{}P}.{\it{}rot}({\it{}x})))
        {\it{}cout} \begin{math}\ll\end{math} {\tt{}"*** dualnumber::rot() gives wrong answer *** {\char92}{\char92}("} \begin{math}\ll\end{math}
                {\it{}P}.{\it{}rot}({\it{}x}) \begin{math}\ll\end{math} {\tt{}"{\char92}{\char92})"} \begin{math}\ll\end{math} {\it{}endl} \begin{math}\ll\end{math} {\it{}endl};

\nwused{\\{NWgqRer-nXe8t-5}}\nwidentuses{\\{{\nwixident{formula{\_}out}}{formula:unout}}}\nwindexuse{\nwixident{formula{\_}out}}{formula:unout}{NWgqRer-1Pc9Jw-2}\nwendcode{}\nwbegindocs{24}Next we specialise the above result to the reference point.
\nwenddocs{}\nwbegincode{25}\sublabel{NWgqRer-1Pc9Jw-3}\nwmargintag{{\nwtagstyle{}\subpageref{NWgqRer-1Pc9Jw-3}}}\moddef{Show expressions~{\nwtagstyle{}\subpageref{NWgqRer-1Pc9Jw-1}}}\plusendmoddef\Rm{}\nwstartdeflinemarkup\nwusesondefline{\\{NWgqRer-nXe8t-5}}\nwprevnextdefs{NWgqRer-1Pc9Jw-2}{NWgqRer-1Pc9Jw-4}\nwenddeflinemarkup
{\it{}formula\_out}({\tt{}"Rotation of {\char92}{\char92}((u\_0, v\_0){\char92}{\char92}) by {\char92}{\char92}(x{\char92}{\char92}): "},\nwindexdefn{\nwixident{formula{\_}out}}{formula:unout}{NWgqRer-1Pc9Jw-3}
                        {\it{}W}.{\it{}subs}({\bf{}lst}({\it{}u} \begin{math}\equiv\end{math} {\it{}u0}, {\it{}v} \begin{math}\equiv\end{math} {\it{}v0})).{\it{}subs}({\it{}Inf} \begin{math}\equiv\end{math} {\it{}pow}({\it{}y}, -1)).{\it{}normal}().{\it{}subs}({\it{}y} \begin{math}\equiv\end{math} 0).{\it{}normal}());

\nwused{\\{NWgqRer-nXe8t-5}}\nwidentdefs{\\{{\nwixident{formula{\_}out}}{formula:unout}}}\nwendcode{}\nwbegindocs{26}The expression for the parabolic norm.
\nwenddocs{}\nwbegincode{27}\sublabel{NWgqRer-1Pc9Jw-4}\nwmargintag{{\nwtagstyle{}\subpageref{NWgqRer-1Pc9Jw-4}}}\moddef{Show expressions~{\nwtagstyle{}\subpageref{NWgqRer-1Pc9Jw-1}}}\plusendmoddef\Rm{}\nwstartdeflinemarkup\nwusesondefline{\\{NWgqRer-nXe8t-5}}\nwprevnextdefs{NWgqRer-1Pc9Jw-3}{NWgqRer-1Pc9Jw-5}\nwenddeflinemarkup
{\it{}formula\_out}({\tt{}"Parabolic norm: "}, {\it{}P}.{\it{}norm}());
{\it{}formula\_out}({\tt{}"Parabolic argument: "}, {\it{}P}.{\it{}arg}());

\nwused{\\{NWgqRer-nXe8t-5}}\nwidentuses{\\{{\nwixident{formula{\_}out}}{formula:unout}}}\nwindexuse{\nwixident{formula{\_}out}}{formula:unout}{NWgqRer-1Pc9Jw-4}\nwendcode{}\nwbegindocs{28}Embedding of reals into dual numbers.
\nwenddocs{}\nwbegincode{29}\sublabel{NWgqRer-1Pc9Jw-5}\nwmargintag{{\nwtagstyle{}\subpageref{NWgqRer-1Pc9Jw-5}}}\moddef{Show expressions~{\nwtagstyle{}\subpageref{NWgqRer-1Pc9Jw-1}}}\plusendmoddef\Rm{}\nwstartdeflinemarkup\nwusesondefline{\\{NWgqRer-nXe8t-5}}\nwprevnextdefs{NWgqRer-1Pc9Jw-4}{NWgqRer-1Pc9Jw-6}\nwenddeflinemarkup
{\it{}possymbol} {\it{}t1}({\tt{}"t1"},{\tt{}"t\_1"}), {\it{}a1}({\tt{}"a1"}, {\tt{}"a\_1"}), {\it{}n}({\tt{}"n"});
{\it{}formula\_out}({\tt{}"Real number {\char92}{\char92}(t\_1{\char92}{\char92}) as a dual number: "}, {\bf{}dual\_number}({\it{}t1}));

\nwused{\\{NWgqRer-nXe8t-5}}\nwidentuses{\\{{\nwixident{formula{\_}out}}{formula:unout}}}\nwindexuse{\nwixident{formula{\_}out}}{formula:unout}{NWgqRer-1Pc9Jw-5}\nwendcode{}\nwbegindocs{30}The expression for the product of two points.
\nwenddocs{}\nwbegincode{31}\sublabel{NWgqRer-1Pc9Jw-6}\nwmargintag{{\nwtagstyle{}\subpageref{NWgqRer-1Pc9Jw-6}}}\moddef{Show expressions~{\nwtagstyle{}\subpageref{NWgqRer-1Pc9Jw-1}}}\plusendmoddef\Rm{}\nwstartdeflinemarkup\nwusesondefline{\\{NWgqRer-nXe8t-5}}\nwprevnextdefs{NWgqRer-1Pc9Jw-5}{NWgqRer-1Pc9Jw-7}\nwenddeflinemarkup
{\it{}formula\_out}({\tt{}"Product: "}, {\it{}P}\begin{math}\ast\end{math}{\it{}P1});

\nwused{\\{NWgqRer-nXe8t-5}}\nwidentuses{\\{{\nwixident{formula{\_}out}}{formula:unout}}}\nwindexuse{\nwixident{formula{\_}out}}{formula:unout}{NWgqRer-1Pc9Jw-6}\nwendcode{}\nwbegindocs{32}The expression of the product of a point and a scalar.
\nwenddocs{}\nwbegincode{33}\sublabel{NWgqRer-1Pc9Jw-7}\nwmargintag{{\nwtagstyle{}\subpageref{NWgqRer-1Pc9Jw-7}}}\moddef{Show expressions~{\nwtagstyle{}\subpageref{NWgqRer-1Pc9Jw-1}}}\plusendmoddef\Rm{}\nwstartdeflinemarkup\nwusesondefline{\\{NWgqRer-nXe8t-5}}\nwprevnextdefs{NWgqRer-1Pc9Jw-6}{NWgqRer-1Pc9Jw-8}\nwenddeflinemarkup
{\it{}formula\_out}({\tt{}"Product by a scalar: "}, {\it{}a}\begin{math}\ast\end{math}{\it{}P});

\nwused{\\{NWgqRer-nXe8t-5}}\nwidentuses{\\{{\nwixident{formula{\_}out}}{formula:unout}}}\nwindexuse{\nwixident{formula{\_}out}}{formula:unout}{NWgqRer-1Pc9Jw-7}\nwendcode{}\nwbegindocs{34}Expressions for the real and imaginary parts.
\nwenddocs{}\nwbegincode{35}\sublabel{NWgqRer-1Pc9Jw-8}\nwmargintag{{\nwtagstyle{}\subpageref{NWgqRer-1Pc9Jw-8}}}\moddef{Show expressions~{\nwtagstyle{}\subpageref{NWgqRer-1Pc9Jw-1}}}\plusendmoddef\Rm{}\nwstartdeflinemarkup\nwusesondefline{\\{NWgqRer-nXe8t-5}}\nwprevnextdefs{NWgqRer-1Pc9Jw-7}{NWgqRer-1Pc9Jw-9}\nwenddeflinemarkup
{\it{}formula\_out}({\tt{}"Real part: "}, {\it{}P}.{\it{}real\_part}());
{\it{}formula\_out}({\tt{}"Imag part: "}, {\it{}P}.{\it{}imag\_part}());

\nwused{\\{NWgqRer-nXe8t-5}}\nwidentuses{\\{{\nwixident{formula{\_}out}}{formula:unout}}}\nwindexuse{\nwixident{formula{\_}out}}{formula:unout}{NWgqRer-1Pc9Jw-8}\nwendcode{}\nwbegindocs{36}The expression for a sum of two points is too cumbersome to be
printed.
\nwenddocs{}\nwbegincode{37}\sublabel{NWgqRer-1Pc9Jw-9}\nwmargintag{{\nwtagstyle{}\subpageref{NWgqRer-1Pc9Jw-9}}}\moddef{Show expressions~{\nwtagstyle{}\subpageref{NWgqRer-1Pc9Jw-1}}}\plusendmoddef\Rm{}\nwstartdeflinemarkup\nwusesondefline{\\{NWgqRer-nXe8t-5}}\nwprevnextdefs{NWgqRer-1Pc9Jw-8}{NWgqRer-1Pc9Jw-A}\nwenddeflinemarkup
//formula\_out("Add is: ", (P+P1).normal());

\nwused{\\{NWgqRer-nXe8t-5}}\nwidentuses{\\{{\nwixident{formula{\_}out}}{formula:unout}}}\nwindexuse{\nwixident{formula{\_}out}}{formula:unout}{NWgqRer-1Pc9Jw-9}\nwendcode{}\nwbegindocs{38}Expression for zero is
\nwenddocs{}\nwbegincode{39}\sublabel{NWgqRer-1Pc9Jw-A}\nwmargintag{{\nwtagstyle{}\subpageref{NWgqRer-1Pc9Jw-A}}}\moddef{Show expressions~{\nwtagstyle{}\subpageref{NWgqRer-1Pc9Jw-1}}}\plusendmoddef\Rm{}\nwstartdeflinemarkup\nwusesondefline{\\{NWgqRer-nXe8t-5}}\nwprevnextdefs{NWgqRer-1Pc9Jw-9}{NWgqRer-1Pc9Jw-B}\nwenddeflinemarkup
{\it{}formula\_out}({\tt{}"Zero vector: "}, {\it{}zero\_dual\_number}());

\nwused{\\{NWgqRer-nXe8t-5}}\nwidentuses{\\{{\nwixident{formula{\_}out}}{formula:unout}}}\nwindexuse{\nwixident{formula{\_}out}}{formula:unout}{NWgqRer-1Pc9Jw-A}\nwendcode{}\nwbegindocs{40}Expression for zero is
\nwenddocs{}\nwbegincode{41}\sublabel{NWgqRer-1Pc9Jw-B}\nwmargintag{{\nwtagstyle{}\subpageref{NWgqRer-1Pc9Jw-B}}}\moddef{Show expressions~{\nwtagstyle{}\subpageref{NWgqRer-1Pc9Jw-1}}}\plusendmoddef\Rm{}\nwstartdeflinemarkup\nwusesondefline{\\{NWgqRer-nXe8t-5}}\nwprevnextdefs{NWgqRer-1Pc9Jw-A}{NWgqRer-1Pc9Jw-C}\nwenddeflinemarkup
{\it{}formula\_out}({\tt{}"Negative vector: "}, {\it{}P}.{\it{}neg}());

\nwused{\\{NWgqRer-nXe8t-5}}\nwidentuses{\\{{\nwixident{formula{\_}out}}{formula:unout}}}\nwindexuse{\nwixident{formula{\_}out}}{formula:unout}{NWgqRer-1Pc9Jw-B}\nwendcode{}\nwbegindocs{42}Expression for zero is
\nwenddocs{}\nwbegincode{43}\sublabel{NWgqRer-1Pc9Jw-C}\nwmargintag{{\nwtagstyle{}\subpageref{NWgqRer-1Pc9Jw-C}}}\moddef{Show expressions~{\nwtagstyle{}\subpageref{NWgqRer-1Pc9Jw-1}}}\plusendmoddef\Rm{}\nwstartdeflinemarkup\nwusesondefline{\\{NWgqRer-nXe8t-5}}\nwprevnextdefs{NWgqRer-1Pc9Jw-B}{NWgqRer-1Pc9Jw-D}\nwenddeflinemarkup
{\it{}formula\_out}({\tt{}"Dual number from argument {\char92}{\char92}(0{\char92}{\char92}) and norm {\char92}{\char92}(n{\char92}{\char92}): "}, \nwindexdefn{\nwixident{formula{\_}out}}{formula:unout}{NWgqRer-1Pc9Jw-C}
                        {\it{}dn\_from\_arg\_mod}(0, {\it{}n}));
{\bf{}dual\_number} {\it{}PP}={\it{}dn\_from\_arg\_mod}({\it{}a1}, {\it{}n});
{\it{}formula\_out}({\tt{}"Dual number from argument {\char92}{\char92}(a\_1{\char92}{\char92}) and norm {\char92}{\char92}(n{\char92}{\char92}): "}, \nwindexdefn{\nwixident{formula{\_}out}}{formula:unout}{NWgqRer-1Pc9Jw-C}
                        {\it{}PP}.{\it{}normal}());
{\it{}formula\_out}({\tt{}"Dual number from argument {\char92}{\char92}(a\_1{\char92}{\char92}) and norm {\char92}{\char92}(n{\char92}{\char92})---norm: "}, \nwindexdefn{\nwixident{formula{\_}out}}{formula:unout}{NWgqRer-1Pc9Jw-C}
                        {\it{}PP}.{\it{}norm}().{\it{}normal}());
{\it{}formula\_out}({\tt{}"Dual number from argument {\char92}{\char92}(a\_1{\char92}{\char92}) and norm {\char92}{\char92}(n{\char92}{\char92})---arg: "}, \nwindexdefn{\nwixident{formula{\_}out}}{formula:unout}{NWgqRer-1Pc9Jw-C}
                        {\it{}PP}.{\it{}arg}().{\it{}normal}().{\it{}normal}());

\nwused{\\{NWgqRer-nXe8t-5}}\nwidentdefs{\\{{\nwixident{formula{\_}out}}{formula:unout}}}\nwendcode{}\nwbegindocs{44}Linear combination of points \((1,0)\) and \((-1,0)\) with
coefficients \(a\) and \(b\), for the linearisation presented in
\cite{Kisil07a}.
\nwenddocs{}\nwbegincode{45}\sublabel{NWgqRer-1Pc9Jw-D}\nwmargintag{{\nwtagstyle{}\subpageref{NWgqRer-1Pc9Jw-D}}}\moddef{Show expressions~{\nwtagstyle{}\subpageref{NWgqRer-1Pc9Jw-1}}}\plusendmoddef\Rm{}\nwstartdeflinemarkup\nwusesondefline{\\{NWgqRer-nXe8t-5}}\nwprevnextdefs{NWgqRer-1Pc9Jw-C}{\relax}\nwenddeflinemarkup
{\it{}formula\_out}({\tt{}"Lin comb of two vectors {\char92}{\char92}(a*(1, 0)+b*(-1, 0){\char92}{\char92}): "},\nwindexdefn{\nwixident{formula{\_}out}}{formula:unout}{NWgqRer-1Pc9Jw-D}
                        ({\it{}a}\begin{math}\ast\end{math}{\bf{}dual\_number}(1,0)+{\it{}b}\begin{math}\ast\end{math}{\bf{}dual\_number}(-1,0)).{\it{}normal}());

\nwused{\\{NWgqRer-nXe8t-5}}\nwidentdefs{\\{{\nwixident{formula{\_}out}}{formula:unout}}}\nwendcode{}\nwbegindocs{46}\nwdocspar
\subsubsection{Checking Algebraic Identities }
\label{sec:check-algebr-ident}
In this Subsection we verify basic algebraic properties of the defined operations.

\nwenddocs{}\nwbegindocs{47}A dual number is the sum of its real and imaginary parts.
\nwenddocs{}\nwbegincode{48}\sublabel{NWgqRer-22TNkn-1}\nwmargintag{{\nwtagstyle{}\subpageref{NWgqRer-22TNkn-1}}}\moddef{Check identities~{\nwtagstyle{}\subpageref{NWgqRer-22TNkn-1}}}\endmoddef\Rm{}\nwstartdeflinemarkup\nwusesondefline{\\{NWgqRer-nXe8t-5}}\nwprevnextdefs{\relax}{NWgqRer-22TNkn-2}\nwenddeflinemarkup
{\it{}test\_out}({\tt{}"{\char92}{\char92}(P{\char92}{\char92}) is the sum {\char92}{\char92}({\char92}{\char92}Re(P){\char92}{\char92}) and {\char92}{\char92}({\char92}{\char92}Im(P){\char92}{\char92}): "},\nwindexdefn{\nwixident{test{\_}out}}{test:unout}{NWgqRer-22TNkn-1}
 {\it{}P}-({\it{}ex\_to}\begin{math}<\end{math}{\bf{}dual\_number}\begin{math}>\end{math}({\it{}P}.{\it{}real\_part}())+{\it{}ex\_to}\begin{math}<\end{math}{\bf{}dual\_number}\begin{math}>\end{math}({\it{}P}.{\it{}imag\_part}())));

\nwalsodefined{\\{NWgqRer-22TNkn-2}\\{NWgqRer-22TNkn-3}\\{NWgqRer-22TNkn-4}\\{NWgqRer-22TNkn-5}\\{NWgqRer-22TNkn-6}\\{NWgqRer-22TNkn-7}\\{NWgqRer-22TNkn-8}\\{NWgqRer-22TNkn-9}\\{NWgqRer-22TNkn-A}\\{NWgqRer-22TNkn-B}\\{NWgqRer-22TNkn-C}\\{NWgqRer-22TNkn-D}\\{NWgqRer-22TNkn-E}\\{NWgqRer-22TNkn-F}}\nwused{\\{NWgqRer-nXe8t-5}}\nwidentdefs{\\{{\nwixident{test{\_}out}}{test:unout}}}\nwendcode{}\nwbegindocs{49}A dual number maid out of a real \(a\) has the norm of real part
equal to \(a\).
\nwenddocs{}\nwbegincode{50}\sublabel{NWgqRer-22TNkn-2}\nwmargintag{{\nwtagstyle{}\subpageref{NWgqRer-22TNkn-2}}}\moddef{Check identities~{\nwtagstyle{}\subpageref{NWgqRer-22TNkn-1}}}\plusendmoddef\Rm{}\nwstartdeflinemarkup\nwusesondefline{\\{NWgqRer-nXe8t-5}}\nwprevnextdefs{NWgqRer-22TNkn-1}{NWgqRer-22TNkn-3}\nwenddeflinemarkup
{\it{}test\_out}({\tt{}"The real part of a real dual number is itself: "},\nwindexdefn{\nwixident{test{\_}out}}{test:unout}{NWgqRer-22TNkn-2}
                 {\it{}ex\_to}\begin{math}<\end{math}{\bf{}dual\_number}\begin{math}>\end{math}({\bf{}dual\_number}({\it{}a}).{\it{}real\_part}()).{\it{}norm}()-{\it{}a});

\nwused{\\{NWgqRer-nXe8t-5}}\nwidentdefs{\\{{\nwixident{test{\_}out}}{test:unout}}}\nwendcode{}\nwbegindocs{51}The norm is invariant under parabolic rotations, i.e. they are in
agreement with Defn.~\ref{de:norm}.
\nwenddocs{}\nwbegincode{52}\sublabel{NWgqRer-22TNkn-3}\nwmargintag{{\nwtagstyle{}\subpageref{NWgqRer-22TNkn-3}}}\moddef{Check identities~{\nwtagstyle{}\subpageref{NWgqRer-22TNkn-1}}}\plusendmoddef\Rm{}\nwstartdeflinemarkup\nwusesondefline{\\{NWgqRer-nXe8t-5}}\nwprevnextdefs{NWgqRer-22TNkn-2}{NWgqRer-22TNkn-4}\nwenddeflinemarkup
{\it{}test\_out}({\tt{}"norm is invariant under rotations: "}, {\it{}P}.{\it{}norm}()-{\it{}W}.{\it{}norm}());

\nwused{\\{NWgqRer-nXe8t-5}}\nwidentuses{\\{{\nwixident{test{\_}out}}{test:unout}}}\nwindexuse{\nwixident{test{\_}out}}{test:unout}{NWgqRer-22TNkn-3}\nwendcode{}\nwbegindocs{53}The product \(w_1\bar{w}_2\) is invariant under rotations,
Prop.~\ref{item:prod-inv}.
\nwenddocs{}\nwbegincode{54}\sublabel{NWgqRer-22TNkn-4}\nwmargintag{{\nwtagstyle{}\subpageref{NWgqRer-22TNkn-4}}}\moddef{Check identities~{\nwtagstyle{}\subpageref{NWgqRer-22TNkn-1}}}\plusendmoddef\Rm{}\nwstartdeflinemarkup\nwusesondefline{\\{NWgqRer-nXe8t-5}}\nwprevnextdefs{NWgqRer-22TNkn-3}{NWgqRer-22TNkn-5}\nwenddeflinemarkup
{\it{}test\_out}({\tt{}"Product is invariant under rotations: "}, {\it{}P}\begin{math}\ast\end{math}{\it{}P1}.{\it{}conjugate}()-{\it{}W}\begin{math}\ast\end{math}{\it{}W1}.{\it{}conjugate}());

\nwused{\\{NWgqRer-nXe8t-5}}\nwidentuses{\\{{\nwixident{test{\_}out}}{test:unout}}}\nwindexuse{\nwixident{test{\_}out}}{test:unout}{NWgqRer-22TNkn-4}\nwendcode{}\nwbegindocs{55}Product \(w\bar{w}\) is \((0, \modulus{w}^2)\), Prop.~\ref{item:prod-norm-sq}.
\nwenddocs{}\nwbegincode{56}\sublabel{NWgqRer-22TNkn-5}\nwmargintag{{\nwtagstyle{}\subpageref{NWgqRer-22TNkn-5}}}\moddef{Check identities~{\nwtagstyle{}\subpageref{NWgqRer-22TNkn-1}}}\plusendmoddef\Rm{}\nwstartdeflinemarkup\nwusesondefline{\\{NWgqRer-nXe8t-5}}\nwprevnextdefs{NWgqRer-22TNkn-4}{NWgqRer-22TNkn-6}\nwenddeflinemarkup
{\it{}test\_out}({\tt{}"Product {\char92}{\char92}(w{\char92}{\char92}bar{\char123}w{\char125}{\char92}{\char92}) is norm squared: "},\nwindexdefn{\nwixident{test{\_}out}}{test:unout}{NWgqRer-22TNkn-5}
   ({\it{}P}\begin{math}\ast\end{math}{\it{}P}.{\it{}conjugate}()-{\it{}dn\_from\_arg\_mod}({\it{}Arg0}, {\it{}pow}({\it{}P}.{\it{}norm}(), 2))));

\nwused{\\{NWgqRer-nXe8t-5}}\nwidentdefs{\\{{\nwixident{test{\_}out}}{test:unout}}}\nwendcode{}\nwbegindocs{57}The reference point is unit under multiplication.
\nwenddocs{}\nwbegincode{58}\sublabel{NWgqRer-22TNkn-6}\nwmargintag{{\nwtagstyle{}\subpageref{NWgqRer-22TNkn-6}}}\moddef{Check identities~{\nwtagstyle{}\subpageref{NWgqRer-22TNkn-1}}}\plusendmoddef\Rm{}\nwstartdeflinemarkup\nwusesondefline{\\{NWgqRer-nXe8t-5}}\nwprevnextdefs{NWgqRer-22TNkn-5}{NWgqRer-22TNkn-7}\nwenddeflinemarkup
{\it{}test\_out}({\tt{}"Product {\char92}{\char92}((u, v)*(u\_0, v\_0){\char92}{\char92}) is {\char92}{\char92}((u, v){\char92}{\char92}): "}, {\it{}P}\begin{math}\ast\end{math}{\it{}P0}-{\it{}P});

\nwused{\\{NWgqRer-nXe8t-5}}\nwidentuses{\\{{\nwixident{test{\_}out}}{test:unout}}}\nwindexuse{\nwixident{test{\_}out}}{test:unout}{NWgqRer-22TNkn-6}\nwendcode{}\nwbegindocs{59}Addition is commutative, Prop.~\ref{item:add-is-comm-ass}.
\nwenddocs{}\nwbegincode{60}\sublabel{NWgqRer-22TNkn-7}\nwmargintag{{\nwtagstyle{}\subpageref{NWgqRer-22TNkn-7}}}\moddef{Check identities~{\nwtagstyle{}\subpageref{NWgqRer-22TNkn-1}}}\plusendmoddef\Rm{}\nwstartdeflinemarkup\nwusesondefline{\\{NWgqRer-nXe8t-5}}\nwprevnextdefs{NWgqRer-22TNkn-6}{NWgqRer-22TNkn-8}\nwenddeflinemarkup
{\it{}test\_out}({\tt{}"Add is commutative: "}, ({\it{}P}+{\it{}P1})-({\it{}P1}+{\it{}P}));

\nwused{\\{NWgqRer-nXe8t-5}}\nwidentuses{\\{{\nwixident{test{\_}out}}{test:unout}}}\nwindexuse{\nwixident{test{\_}out}}{test:unout}{NWgqRer-22TNkn-7}\nwendcode{}\nwbegindocs{61}Addition is associative, Prop.~\ref{item:add-is-comm-ass}.
\nwenddocs{}\nwbegincode{62}\sublabel{NWgqRer-22TNkn-8}\nwmargintag{{\nwtagstyle{}\subpageref{NWgqRer-22TNkn-8}}}\moddef{Check identities~{\nwtagstyle{}\subpageref{NWgqRer-22TNkn-1}}}\plusendmoddef\Rm{}\nwstartdeflinemarkup\nwusesondefline{\\{NWgqRer-nXe8t-5}}\nwprevnextdefs{NWgqRer-22TNkn-7}{NWgqRer-22TNkn-9}\nwenddeflinemarkup
{\it{}test\_out}({\tt{}"Add is associative: "}, (({\it{}P}+{\it{}P1})+ {\it{}P2})-({\it{}P}+({\it{}P1}+{\it{}P2})));

\nwused{\\{NWgqRer-nXe8t-5}}\nwidentuses{\\{{\nwixident{test{\_}out}}{test:unout}}}\nwindexuse{\nwixident{test{\_}out}}{test:unout}{NWgqRer-22TNkn-8}\nwendcode{}\nwbegindocs{63}Multiplication by a scalar is commutative.
\nwenddocs{}\nwbegincode{64}\sublabel{NWgqRer-22TNkn-9}\nwmargintag{{\nwtagstyle{}\subpageref{NWgqRer-22TNkn-9}}}\moddef{Check identities~{\nwtagstyle{}\subpageref{NWgqRer-22TNkn-1}}}\plusendmoddef\Rm{}\nwstartdeflinemarkup\nwusesondefline{\\{NWgqRer-nXe8t-5}}\nwprevnextdefs{NWgqRer-22TNkn-8}{NWgqRer-22TNkn-A}\nwenddeflinemarkup
{\it{}test\_out}({\tt{}"S-mult is commutative: "}, {\it{}P}\begin{math}\ast\end{math}{\it{}a}-{\it{}a}\begin{math}\ast\end{math}{\it{}P});

\nwused{\\{NWgqRer-nXe8t-5}}\nwidentuses{\\{{\nwixident{test{\_}out}}{test:unout}}}\nwindexuse{\nwixident{test{\_}out}}{test:unout}{NWgqRer-22TNkn-9}\nwendcode{}\nwbegindocs{65}Multiplication by a scalar is associative.
\nwenddocs{}\nwbegincode{66}\sublabel{NWgqRer-22TNkn-A}\nwmargintag{{\nwtagstyle{}\subpageref{NWgqRer-22TNkn-A}}}\moddef{Check identities~{\nwtagstyle{}\subpageref{NWgqRer-22TNkn-1}}}\plusendmoddef\Rm{}\nwstartdeflinemarkup\nwusesondefline{\\{NWgqRer-nXe8t-5}}\nwprevnextdefs{NWgqRer-22TNkn-9}{NWgqRer-22TNkn-B}\nwenddeflinemarkup
{\it{}test\_out}({\tt{}"S-mult is associative: "}, {\it{}b}\begin{math}\ast\end{math}{\it{}P}\begin{math}\ast\end{math}{\it{}a}-{\it{}a}\begin{math}\ast\end{math}{\it{}P}\begin{math}\ast\end{math}{\it{}b});

\nwused{\\{NWgqRer-nXe8t-5}}\nwidentuses{\\{{\nwixident{test{\_}out}}{test:unout}}}\nwindexuse{\nwixident{test{\_}out}}{test:unout}{NWgqRer-22TNkn-A}\nwendcode{}\nwbegindocs{67}Distributive law \(a(w_1+w_2)=aw_1+aw_2\), Prop.~\ref{item:distrib-scalar}.
\nwenddocs{}\nwbegincode{68}\sublabel{NWgqRer-22TNkn-B}\nwmargintag{{\nwtagstyle{}\subpageref{NWgqRer-22TNkn-B}}}\moddef{Check identities~{\nwtagstyle{}\subpageref{NWgqRer-22TNkn-1}}}\plusendmoddef\Rm{}\nwstartdeflinemarkup\nwusesondefline{\\{NWgqRer-nXe8t-5}}\nwprevnextdefs{NWgqRer-22TNkn-A}{NWgqRer-22TNkn-C}\nwenddeflinemarkup
{\it{}test\_out}({\tt{}"S-mult is distributive 1: "}, {\it{}a}\begin{math}\ast\end{math}({\it{}P}+{\it{}P1})-({\it{}a}\begin{math}\ast\end{math}{\it{}P} +{\it{}a}\begin{math}\ast\end{math}{\it{}P1}));

\nwused{\\{NWgqRer-nXe8t-5}}\nwidentuses{\\{{\nwixident{test{\_}out}}{test:unout}}}\nwindexuse{\nwixident{test{\_}out}}{test:unout}{NWgqRer-22TNkn-B}\nwendcode{}\nwbegindocs{69}Distributive law \((a+b)w=aw+bw\), Prop.~\ref{item:distrib-scalar}.
\nwenddocs{}\nwbegincode{70}\sublabel{NWgqRer-22TNkn-C}\nwmargintag{{\nwtagstyle{}\subpageref{NWgqRer-22TNkn-C}}}\moddef{Check identities~{\nwtagstyle{}\subpageref{NWgqRer-22TNkn-1}}}\plusendmoddef\Rm{}\nwstartdeflinemarkup\nwusesondefline{\\{NWgqRer-nXe8t-5}}\nwprevnextdefs{NWgqRer-22TNkn-B}{NWgqRer-22TNkn-D}\nwenddeflinemarkup
{\it{}test\_out}({\tt{}"S-mult is distributive 2: "}, {\it{}P}\begin{math}\ast\end{math}({\it{}a}+{\it{}b})-({\it{}P}\begin{math}\ast\end{math}{\it{}a} + {\it{}P}\begin{math}\ast\end{math}{\it{}b}));

\nwused{\\{NWgqRer-nXe8t-5}}\nwidentuses{\\{{\nwixident{test{\_}out}}{test:unout}}}\nwindexuse{\nwixident{test{\_}out}}{test:unout}{NWgqRer-22TNkn-C}\nwendcode{}\nwbegindocs{71}Product is commutative, Prop.~\ref{item:prod-comm-ass}.
\nwenddocs{}\nwbegincode{72}\sublabel{NWgqRer-22TNkn-D}\nwmargintag{{\nwtagstyle{}\subpageref{NWgqRer-22TNkn-D}}}\moddef{Check identities~{\nwtagstyle{}\subpageref{NWgqRer-22TNkn-1}}}\plusendmoddef\Rm{}\nwstartdeflinemarkup\nwusesondefline{\\{NWgqRer-nXe8t-5}}\nwprevnextdefs{NWgqRer-22TNkn-C}{NWgqRer-22TNkn-E}\nwenddeflinemarkup
{\it{}test\_out}({\tt{}"Product is symmetric (commutative): "}, {\it{}P}\begin{math}\ast\end{math}{\it{}P1}-{\it{}P1}\begin{math}\ast\end{math}{\it{}P});

\nwused{\\{NWgqRer-nXe8t-5}}\nwidentuses{\\{{\nwixident{test{\_}out}}{test:unout}}}\nwindexuse{\nwixident{test{\_}out}}{test:unout}{NWgqRer-22TNkn-D}\nwendcode{}\nwbegindocs{73}Product is associative, Prop.~\ref{item:prod-comm-ass}.
\nwenddocs{}\nwbegincode{74}\sublabel{NWgqRer-22TNkn-E}\nwmargintag{{\nwtagstyle{}\subpageref{NWgqRer-22TNkn-E}}}\moddef{Check identities~{\nwtagstyle{}\subpageref{NWgqRer-22TNkn-1}}}\plusendmoddef\Rm{}\nwstartdeflinemarkup\nwusesondefline{\\{NWgqRer-nXe8t-5}}\nwprevnextdefs{NWgqRer-22TNkn-D}{NWgqRer-22TNkn-F}\nwenddeflinemarkup
{\it{}test\_out}({\tt{}"Prod is associative: "}, ({\it{}P}\begin{math}\ast\end{math}{\it{}P1})\begin{math}\ast\end{math}{\it{}P2}-{\it{}P}\begin{math}\ast\end{math}({\it{}P1}\begin{math}\ast\end{math}{\it{}P2}));

\nwused{\\{NWgqRer-nXe8t-5}}\nwidentuses{\\{{\nwixident{test{\_}out}}{test:unout}}}\nwindexuse{\nwixident{test{\_}out}}{test:unout}{NWgqRer-22TNkn-E}\nwendcode{}\nwbegindocs{75}Product and addition are distributive, Prop.~\ref{item:distrib}.
\nwenddocs{}\nwbegincode{76}\sublabel{NWgqRer-22TNkn-F}\nwmargintag{{\nwtagstyle{}\subpageref{NWgqRer-22TNkn-F}}}\moddef{Check identities~{\nwtagstyle{}\subpageref{NWgqRer-22TNkn-1}}}\plusendmoddef\Rm{}\nwstartdeflinemarkup\nwusesondefline{\\{NWgqRer-nXe8t-5}}\nwprevnextdefs{NWgqRer-22TNkn-E}{\relax}\nwenddeflinemarkup
{\it{}test\_out}({\tt{}"Product is distributive: "}, ({\it{}P}+{\it{}P1})\begin{math}\ast\end{math}{\it{}P2}-({\it{}P}\begin{math}\ast\end{math}{\it{}P2}+{\it{}P1}\begin{math}\ast\end{math}{\it{}P2}));

\nwused{\\{NWgqRer-nXe8t-5}}\nwidentuses{\\{{\nwixident{test{\_}out}}{test:unout}}}\nwindexuse{\nwixident{test{\_}out}}{test:unout}{NWgqRer-22TNkn-F}\nwendcode{}\nwbegindocs{77}\nwdocspar
\subsection{Induced Representations}
\label{sec:induc-repr1}
 Here we calculate the basic formulae for Section~\ref{sec:induc-repr}.

\nwenddocs{}\nwbegindocs{78}\nwdocspar
\subsubsection{Encoded formulae}
\label{sec:encoded-formulae}

\nwenddocs{}\nwbegindocs{79}This routine encodes the map \(s: \Space{R}{2}\rightarrow
\SL\)~\eqref{eq:s-map}.
\nwenddocs{}\nwbegincode{80}\sublabel{NWgqRer-11aTCz-1}\nwmargintag{{\nwtagstyle{}\subpageref{NWgqRer-11aTCz-1}}}\moddef{Induced representations routines~{\nwtagstyle{}\subpageref{NWgqRer-11aTCz-1}}}\endmoddef\Rm{}\nwstartdeflinemarkup\nwusesondefline{\\{NWgqRer-1p0Y9w-2}}\nwprevnextdefs{\relax}{NWgqRer-11aTCz-2}\nwenddeflinemarkup
{\bf{}ex} {\it{}s\_map}({\bf{}const} {\bf{}ex} & {\it{}u}, {\bf{}const} {\bf{}ex} & {\it{}v}) {\nwlbrace}
        {\bf{}return} {\bf{}matrix}(2, 2, {\bf{}lst}({\it{}v},{\it{}u},0,1));
{\nwrbrace}

{\bf{}ex} {\it{}s\_map}({\bf{}const} {\bf{}ex} & {\it{}P}) {\nwlbrace}
        {\bf{}if} ({\it{}P}.{\it{}nops}() \begin{math}\equiv\end{math} 2)
                {\bf{}return} {\it{}s\_map}({\it{}P}.{\it{}op}(0), {\it{}P}.{\it{}op}(1));
        {\it{}cerr} \begin{math}\ll\end{math} {\tt{}"s\_map() error: parameter should have two operands"} \begin{math}\ll\end{math} {\it{}endl};
        {\bf{}return} {\it{}s\_map}({\it{}P},1);
{\nwrbrace}

\nwalsodefined{\\{NWgqRer-11aTCz-2}\\{NWgqRer-11aTCz-3}\\{NWgqRer-11aTCz-4}}\nwused{\\{NWgqRer-1p0Y9w-2}}\nwendcode{}\nwbegindocs{81}This routine encodes the map \(r: \SL \rightarrow
H\)~\eqref{eq:r-map}. The first parameter is an element of \(\SL\), the
second---is a generic element of subgroup \(H\). We look specific
matrix of the form \(K\) which makes the product \(MK\) belonging to
the image of {\Tt{}\Rm{}{\it{}s\_map}()\nwendquote}, i.e. its \((2,1)\) element should be zero.
\nwenddocs{}\nwbegincode{82}\sublabel{NWgqRer-11aTCz-2}\nwmargintag{{\nwtagstyle{}\subpageref{NWgqRer-11aTCz-2}}}\moddef{Induced representations routines~{\nwtagstyle{}\subpageref{NWgqRer-11aTCz-1}}}\plusendmoddef\Rm{}\nwstartdeflinemarkup\nwusesondefline{\\{NWgqRer-1p0Y9w-2}}\nwprevnextdefs{NWgqRer-11aTCz-1}{NWgqRer-11aTCz-3}\nwenddeflinemarkup
{\bf{}ex} {\it{}r\_map}({\bf{}const} {\bf{}ex} & {\it{}M}, {\bf{}const} {\bf{}ex} & {\it{}K}) {\nwlbrace}
        {\bf{}ex} {\it{}K1}={\it{}K}.{\it{}evalm}(), {\it{}K2};
        {\bf{}lst} {\it{}vars} = ({\it{}is\_a}\begin{math}<\end{math}{\bf{}symbol}\begin{math}>\end{math}({\it{}K1}.{\it{}op}(2)) ? {\bf{}lst}({\it{}K1}.{\it{}op}(2)) :  {\bf{}lst}({\it{}K1}.{\it{}op}(1)));
        {\bf{}if} ({\it{}is\_a}\begin{math}<\end{math}{\bf{}symbol}\begin{math}>\end{math}({\it{}K1}.{\it{}op}(3))) {\nwlbrace}
                {\it{}vars} = {\it{}vars}.{\it{}append}({\it{}K1}.{\it{}op}(3));
                {\it{}K2} = {\it{}K1}.{\it{}subs}({\it{}lsolve}({\bf{}lst}(({\it{}M}\begin{math}\ast\end{math}{\it{}K1}).{\it{}evalm}().{\it{}op}(2)\begin{math}\equiv\end{math}0), {\it{}vars})).{\it{}subs}({\it{}K1}.{\it{}op}(3)\begin{math}\equiv\end{math}1);
        {\nwrbrace} {\bf{}else}
                {\it{}K2} = {\it{}K1}.{\it{}subs}({\it{}lsolve}({\bf{}lst}(({\it{}M}\begin{math}\ast\end{math}{\it{}K1}).{\it{}evalm}().{\it{}op}(2)\begin{math}\equiv\end{math}0), {\it{}vars}));
        {\bf{}return} {\it{}pow}({\it{}K2}, -1).{\it{}evalm}();
{\nwrbrace}

\nwused{\\{NWgqRer-1p0Y9w-2}}\nwendcode{}\nwbegindocs{83}This is the inverse \(s^{-1}\) of the above map \(s\).
\nwenddocs{}\nwbegincode{84}\sublabel{NWgqRer-11aTCz-3}\nwmargintag{{\nwtagstyle{}\subpageref{NWgqRer-11aTCz-3}}}\moddef{Induced representations routines~{\nwtagstyle{}\subpageref{NWgqRer-11aTCz-1}}}\plusendmoddef\Rm{}\nwstartdeflinemarkup\nwusesondefline{\\{NWgqRer-1p0Y9w-2}}\nwprevnextdefs{NWgqRer-11aTCz-2}{NWgqRer-11aTCz-4}\nwenddeflinemarkup
{\bf{}ex} {\it{}p\_map}({\bf{}const} {\bf{}ex} & {\it{}M}, {\bf{}const} {\bf{}ex} & {\it{}K}) {\nwlbrace}
        {\bf{}ex} {\it{}x} = {\it{}K}.{\it{}evalm}().{\it{}op}(3);
        {\bf{}ex} {\it{}MK}=({\it{}M}\begin{math}\ast\end{math}{\it{}pow}({\it{}r\_map}({\it{}M},{\it{}K}),-1)).{\it{}evalm}();
        {\bf{}ex} {\it{}D}={\it{}MK}.{\it{}op}(3).{\it{}subs}({\it{}x}\begin{math}\equiv\end{math}1).{\it{}normal}();
        {\bf{}return} {\bf{}matrix}(1, 2, {\bf{}lst}(({\it{}MK}.{\it{}op}(1).{\it{}subs}({\it{}x}\begin{math}\equiv\end{math}1).{\it{}normal}()\begin{math}\div\end{math}{\it{}D}).{\it{}normal}(),
                                                        ({\it{}MK}.{\it{}op}(0).{\it{}subs}({\it{}x}\begin{math}\equiv\end{math}1).{\it{}normal}()\begin{math}\div\end{math}{\it{}D}).{\it{}normal}()));
{\nwrbrace}

\nwused{\\{NWgqRer-1p0Y9w-2}}\nwendcode{}\nwbegindocs{85}This is a matrix form of the above inverse map {\Tt{}\Rm{}{\it{}p\_map}()\nwendquote}.
\nwenddocs{}\nwbegincode{86}\sublabel{NWgqRer-11aTCz-4}\nwmargintag{{\nwtagstyle{}\subpageref{NWgqRer-11aTCz-4}}}\moddef{Induced representations routines~{\nwtagstyle{}\subpageref{NWgqRer-11aTCz-1}}}\plusendmoddef\Rm{}\nwstartdeflinemarkup\nwusesondefline{\\{NWgqRer-1p0Y9w-2}}\nwprevnextdefs{NWgqRer-11aTCz-3}{\relax}\nwenddeflinemarkup
{\bf{}ex} {\it{}p\_map\_m}({\bf{}const} {\bf{}ex} & {\it{}M}, {\bf{}const} {\bf{}ex} & {\it{}K}) {\nwlbrace}
 {\bf{}return} ({\it{}M}\begin{math}\ast\end{math}{\it{}pow}({\it{}r\_map}({\it{}M},{\it{}K}),-1)).{\it{}evalm}();
{\nwrbrace}

\nwused{\\{NWgqRer-1p0Y9w-2}}\nwendcode{}\nwbegindocs{87}\nwdocspar
\subsubsection{Caculation of induced representation formulae}
\label{sec:cacul-induc-repr}

\nwenddocs{}\nwbegindocs{88} Firstly we define a generic element {\Tt{}\Rm{}{\it{}M}\nwendquote} of \(\SL\).
\nwenddocs{}\nwbegincode{89}\sublabel{NWgqRer-1h89LW-1}\nwmargintag{{\nwtagstyle{}\subpageref{NWgqRer-1h89LW-1}}}\moddef{Induced representations~{\nwtagstyle{}\subpageref{NWgqRer-1h89LW-1}}}\endmoddef\Rm{}\nwstartdeflinemarkup\nwusesondefline{\\{NWgqRer-3sxui-4}}\nwprevnextdefs{\relax}{NWgqRer-1h89LW-2}\nwenddeflinemarkup
{\bf{}ex} {\it{}M}={\bf{}matrix}(2,2, {\bf{}lst}({\it{}a},{\it{}b},{\it{}c},{\it{}d})), {\it{}H};

\nwalsodefined{\\{NWgqRer-1h89LW-2}\\{NWgqRer-1h89LW-3}\\{NWgqRer-1h89LW-4}\\{NWgqRer-1h89LW-5}\\{NWgqRer-1h89LW-6}\\{NWgqRer-1h89LW-7}\\{NWgqRer-1h89LW-8}}\nwused{\\{NWgqRer-3sxui-4}}\nwendcode{}\nwbegindocs{90}We consider the three cases.
\nwenddocs{}\nwbegincode{91}\sublabel{NWgqRer-1h89LW-2}\nwmargintag{{\nwtagstyle{}\subpageref{NWgqRer-1h89LW-2}}}\moddef{Induced representations~{\nwtagstyle{}\subpageref{NWgqRer-1h89LW-1}}}\plusendmoddef\Rm{}\nwstartdeflinemarkup\nwusesondefline{\\{NWgqRer-3sxui-4}}\nwprevnextdefs{NWgqRer-1h89LW-1}{NWgqRer-1h89LW-3}\nwenddeflinemarkup
{\it{}string} {\it{}cases}[]={\nwlbrace}{\tt{}"Elliptic"}, {\tt{}"Parabolic ({\char92}{\char92}(N^{\char92}{\char92}prime{\char92}{\char92}))"}, {\tt{}"Hyperbolic"}{\nwrbrace};

\nwused{\\{NWgqRer-3sxui-4}}\nwendcode{}\nwbegindocs{92}In the those cases {\Tt{}\Rm{}{\it{}subgroups}\nwendquote} holds a generic element of a
subgroup \(H\), see~\eqref{eq:k-subgroup}, \eqref{eq:ap-subgroup}
and~\eqref{eq:np-subgroup}.
\nwenddocs{}\nwbegincode{93}\sublabel{NWgqRer-1h89LW-3}\nwmargintag{{\nwtagstyle{}\subpageref{NWgqRer-1h89LW-3}}}\moddef{Induced representations~{\nwtagstyle{}\subpageref{NWgqRer-1h89LW-1}}}\plusendmoddef\Rm{}\nwstartdeflinemarkup\nwusesondefline{\\{NWgqRer-3sxui-4}}\nwprevnextdefs{NWgqRer-1h89LW-2}{NWgqRer-1h89LW-4}\nwenddeflinemarkup
{\bf{}ex} {\it{}subgroups}={\bf{}lst}({\bf{}matrix}(2, 2, {\bf{}lst}({\it{}x},-{\it{}y},{\it{}y},{\it{}x})),
                                 {\bf{}matrix}(2, 2, {\bf{}lst}(1,0,{\it{}y},1)),
                                 {\bf{}matrix}(2, 2, {\bf{}lst}({\it{}x},{\it{}y},{\it{}y},{\it{}x})));

\nwused{\\{NWgqRer-3sxui-4}}\nwendcode{}\nwbegindocs{94}Now we run a cycle over the three cases\ldots
\nwenddocs{}\nwbegincode{95}\sublabel{NWgqRer-1h89LW-4}\nwmargintag{{\nwtagstyle{}\subpageref{NWgqRer-1h89LW-4}}}\moddef{Induced representations~{\nwtagstyle{}\subpageref{NWgqRer-1h89LW-1}}}\plusendmoddef\Rm{}\nwstartdeflinemarkup\nwusesondefline{\\{NWgqRer-3sxui-4}}\nwprevnextdefs{NWgqRer-1h89LW-3}{NWgqRer-1h89LW-5}\nwenddeflinemarkup
{\bf{}for}({\bf{}int} {\it{}i}=0; {\it{}i}\begin{math}<\end{math}3; {\it{}i}\protect\PP) {\nwlbrace} 
 {\it{}H}={\it{}subgroups}[{\it{}i}];
 {\it{}cout} \begin{math}\ll\end{math} {\it{}cases}[{\it{}i}] \begin{math}\ll\end{math} {\tt{}" case of induced representations{\char92}{\char92}{\char92}{\char92}"} \begin{math}\ll\end{math} {\it{}endl};
 //formula\_out("M*H: ", (M*H).evalm());

\nwused{\\{NWgqRer-3sxui-4}}\nwidentuses{\\{{\nwixident{formula{\_}out}}{formula:unout}}}\nwindexuse{\nwixident{formula{\_}out}}{formula:unout}{NWgqRer-1h89LW-4}\nwendcode{}\nwbegindocs{96}\ldots and output expression of \(r\)~\eqref{eq:r-map},\ldots
\nwenddocs{}\nwbegincode{97}\sublabel{NWgqRer-1h89LW-5}\nwmargintag{{\nwtagstyle{}\subpageref{NWgqRer-1h89LW-5}}}\moddef{Induced representations~{\nwtagstyle{}\subpageref{NWgqRer-1h89LW-1}}}\plusendmoddef\Rm{}\nwstartdeflinemarkup\nwusesondefline{\\{NWgqRer-3sxui-4}}\nwprevnextdefs{NWgqRer-1h89LW-4}{NWgqRer-1h89LW-6}\nwenddeflinemarkup
 {\it{}formula\_out}({\tt{}"map {\char92}{\char92}(r(M){\char92}{\char92}): "}, {\it{}r\_map}({\it{}M},{\it{}H}));

\nwused{\\{NWgqRer-3sxui-4}}\nwidentuses{\\{{\nwixident{formula{\_}out}}{formula:unout}}}\nwindexuse{\nwixident{formula{\_}out}}{formula:unout}{NWgqRer-1h89LW-5}\nwendcode{}\nwbegindocs{98}\ldots matrix form of the inverse \(s^{-1}\)~\eqref{eq:r-map},\ldots
\nwenddocs{}\nwbegincode{99}\sublabel{NWgqRer-1h89LW-6}\nwmargintag{{\nwtagstyle{}\subpageref{NWgqRer-1h89LW-6}}}\moddef{Induced representations~{\nwtagstyle{}\subpageref{NWgqRer-1h89LW-1}}}\plusendmoddef\Rm{}\nwstartdeflinemarkup\nwusesondefline{\\{NWgqRer-3sxui-4}}\nwprevnextdefs{NWgqRer-1h89LW-5}{NWgqRer-1h89LW-7}\nwenddeflinemarkup
 {\it{}formula\_out}({\tt{}"map {\char92}{\char92}(s^{\char123}-1{\char125}(M){\char92}{\char92}): "}, {\it{}p\_map\_m}({\it{}M},{\it{}H}).{\it{}subs}({\it{}a}\begin{math}\equiv\end{math}(1+{\it{}b}\begin{math}\ast\end{math}{\it{}c})\begin{math}\div\end{math}{\it{}d}).{\it{}normal}());

\nwused{\\{NWgqRer-3sxui-4}}\nwidentuses{\\{{\nwixident{formula{\_}out}}{formula:unout}}}\nwindexuse{\nwixident{formula{\_}out}}{formula:unout}{NWgqRer-1h89LW-6}\nwendcode{}\nwbegindocs{100}\ldots expression for the argument of the character in~\eqref{eq:def-ind},\ldots
\nwenddocs{}\nwbegincode{101}\sublabel{NWgqRer-1h89LW-7}\nwmargintag{{\nwtagstyle{}\subpageref{NWgqRer-1h89LW-7}}}\moddef{Induced representations~{\nwtagstyle{}\subpageref{NWgqRer-1h89LW-1}}}\plusendmoddef\Rm{}\nwstartdeflinemarkup\nwusesondefline{\\{NWgqRer-3sxui-4}}\nwprevnextdefs{NWgqRer-1h89LW-6}{NWgqRer-1h89LW-8}\nwenddeflinemarkup
 {\it{}formula\_out}({\tt{}"character: "}, {\it{}r\_map}({\it{}M}\begin{math}\ast\end{math}{\it{}s\_map}({\it{}P}),{\it{}H}));

\nwused{\\{NWgqRer-3sxui-4}}\nwidentuses{\\{{\nwixident{formula{\_}out}}{formula:unout}}}\nwindexuse{\nwixident{formula{\_}out}}{formula:unout}{NWgqRer-1h89LW-7}\nwendcode{}\nwbegindocs{102}\ldots and finally the action~\eqref{eq:g-action} of \(\SL\) on the
homogeneous space.
\nwenddocs{}\nwbegincode{103}\sublabel{NWgqRer-1h89LW-8}\nwmargintag{{\nwtagstyle{}\subpageref{NWgqRer-1h89LW-8}}}\moddef{Induced representations~{\nwtagstyle{}\subpageref{NWgqRer-1h89LW-1}}}\plusendmoddef\Rm{}\nwstartdeflinemarkup\nwusesondefline{\\{NWgqRer-3sxui-4}}\nwprevnextdefs{NWgqRer-1h89LW-7}{\relax}\nwenddeflinemarkup
 {\it{}formula\_out}({\tt{}"Moebius map: "}, {\it{}p\_map}({\it{}M}\begin{math}\ast\end{math}{\it{}s\_map}({\it{}P}.{\it{}to\_matrix}()),{\it{}H}));
 {\it{}test\_out}({\tt{}"Moebius map is given by the imaginary unit: "}, {\it{}p\_map}({\it{}M}\begin{math}\ast\end{math}{\it{}s\_map}({\it{}P}),{\it{}H}) -
                  {\it{}clifford\_moebius\_map}({\it{}a}\begin{math}\ast\end{math}{\it{}one}, {\it{}b}\begin{math}\ast\end{math}{\it{}one}, {\it{}c}\begin{math}\ast\end{math}{\it{}one}, {\it{}d}\begin{math}\ast\end{math}{\it{}one},{\it{}P}.{\it{}to\_matrix}(),{\it{}e}).{\it{}subs}({\it{}sign}\begin{math}\equiv\end{math}{\it{}i}-1));
 {\it{}cout} \begin{math}\ll\end{math}  ({\it{}latexout} ? {\tt{}"{\char92}{\char92}vspace{\char123}2mm{\char125}{\char92}{\char92}hrule"} :
                   {\tt{}"----------------------------------------"} ) \begin{math}\ll\end{math} {\it{}endl};
{\nwrbrace}

\nwused{\\{NWgqRer-3sxui-4}}\nwidentuses{\\{{\nwixident{formula{\_}out}}{formula:unout}}\\{{\nwixident{test{\_}out}}{test:unout}}}\nwindexuse{\nwixident{formula{\_}out}}{formula:unout}{NWgqRer-1h89LW-8}\nwindexuse{\nwixident{test{\_}out}}{test:unout}{NWgqRer-1h89LW-8}\nwendcode{}\nwbegindocs{104}\nwdocspar
\subsection{Program Outline}
\label{sec:program-outline}
Here is the outline how we use the above parts.

\nwenddocs{}\nwbegindocs{105}\nwdocspar
\subsubsection{Test program outline}
\label{sec:test-program-outline}

Firstly we load {\Tt{}\Rm{}{\bf{}dual\_number}\nwendquote} support.
\nwenddocs{}\nwbegincode{106}\sublabel{NWgqRer-1p0Y9w-1}\nwmargintag{{\nwtagstyle{}\subpageref{NWgqRer-1p0Y9w-1}}}\moddef{*~{\nwtagstyle{}\subpageref{NWgqRer-1p0Y9w-1}}}\endmoddef\Rm{}\nwstartdeflinemarkup\nwprevnextdefs{\relax}{NWgqRer-1p0Y9w-2}\nwenddeflinemarkup
{\bf{}\char35{}include}{\tt{} \begin{math}<\end{math}cycle.h\begin{math}>\end{math}}
{\bf{}\char35{}include}{\tt{} \begin{math}<\end{math}dualnum.h\begin{math}>\end{math}}

\nwalsodefined{\\{NWgqRer-1p0Y9w-2}}\nwnotused{*}\nwendcode{}\nwbegindocs{107}The rest of the program makes all checks.
\nwenddocs{}\nwbegincode{108}\sublabel{NWgqRer-1p0Y9w-2}\nwmargintag{{\nwtagstyle{}\subpageref{NWgqRer-1p0Y9w-2}}}\moddef{*~{\nwtagstyle{}\subpageref{NWgqRer-1p0Y9w-1}}}\plusendmoddef\Rm{}\nwstartdeflinemarkup\nwprevnextdefs{NWgqRer-1p0Y9w-1}{\relax}\nwenddeflinemarkup
 \LA{}Definition of variables~{\nwtagstyle{}\subpageref{NWgqRer-2iImO3-1}}\RA{}
 \LA{}Test routine~{\nwtagstyle{}\subpageref{NWgqRer-nXe8t-1}}\RA{}
 \LA{}Induced representations routines~{\nwtagstyle{}\subpageref{NWgqRer-11aTCz-1}}\RA{}
 \LA{}Main procedure~{\nwtagstyle{}\subpageref{NWgqRer-3sxui-1}}\RA{}

\nwendcode{}\nwbegindocs{109}\nwdocspar
\subsubsection{Variables}
\label{sec:variables}

\nwenddocs{}\nwbegindocs{110}These {\Tt{}\Rm{}{\bf{}realsymbol}\nwendquote}s are used in our calculations.
\nwenddocs{}\nwbegincode{111}\sublabel{NWgqRer-2iImO3-1}\nwmargintag{{\nwtagstyle{}\subpageref{NWgqRer-2iImO3-1}}}\moddef{Definition of variables~{\nwtagstyle{}\subpageref{NWgqRer-2iImO3-1}}}\endmoddef\Rm{}\nwstartdeflinemarkup\nwusesondefline{\\{NWgqRer-1p0Y9w-2}}\nwprevnextdefs{\relax}{NWgqRer-2iImO3-2}\nwenddeflinemarkup
{\bf{}const} {\bf{}numeric} {\it{}half}(1,2);\nwindexdefn{\nwixident{numeric}}{numeric}{NWgqRer-2iImO3-1}

\nwalsodefined{\\{NWgqRer-2iImO3-2}\\{NWgqRer-2iImO3-3}\\{NWgqRer-2iImO3-4}\\{NWgqRer-2iImO3-5}\\{NWgqRer-2iImO3-6}\\{NWgqRer-2iImO3-7}}\nwused{\\{NWgqRer-1p0Y9w-2}}\nwidentdefs{\\{{\nwixident{numeric}}{numeric}}}\nwendcode{}\nwbegindocs{112}Variables \(v\)s oftenly appear under square roots of the form
\(\sqrt{1+2v}\). To facilitate the simplifications of the type
\((\sqrt{1+2v})^2 = 1+2v\) we (falsely) define them to be positive
symbols.  
\nwenddocs{}\nwbegincode{113}\sublabel{NWgqRer-2iImO3-2}\nwmargintag{{\nwtagstyle{}\subpageref{NWgqRer-2iImO3-2}}}\moddef{Definition of variables~{\nwtagstyle{}\subpageref{NWgqRer-2iImO3-1}}}\plusendmoddef\Rm{}\nwstartdeflinemarkup\nwusesondefline{\\{NWgqRer-1p0Y9w-2}}\nwprevnextdefs{NWgqRer-2iImO3-1}{NWgqRer-2iImO3-3}\nwenddeflinemarkup
{\it{}possymbol} {\it{}v}({\tt{}"v"}), {\it{}v1}({\tt{}"v'"}), {\it{}v2}({\tt{}"v''"});

\nwused{\\{NWgqRer-1p0Y9w-2}}\nwendcode{}\nwbegindocs{114} Other real variables.
\nwenddocs{}\nwbegincode{115}\sublabel{NWgqRer-2iImO3-3}\nwmargintag{{\nwtagstyle{}\subpageref{NWgqRer-2iImO3-3}}}\moddef{Definition of variables~{\nwtagstyle{}\subpageref{NWgqRer-2iImO3-1}}}\plusendmoddef\Rm{}\nwstartdeflinemarkup\nwusesondefline{\\{NWgqRer-1p0Y9w-2}}\nwprevnextdefs{NWgqRer-2iImO3-2}{NWgqRer-2iImO3-4}\nwenddeflinemarkup
{\bf{}realsymbol} {\it{}u}({\tt{}"u"}), {\it{}u1}({\tt{}"u'"}), {\it{}u2}({\tt{}"u''"}), 
 {\it{}a}({\tt{}"a"}), {\it{}b}({\tt{}"b"}), {\it{}c}({\tt{}"c"}), {\it{}d}({\tt{}"d"}), {\it{}x}({\tt{}"x"}), {\it{}y}({\tt{}"y"}),

\nwused{\\{NWgqRer-1p0Y9w-2}}\nwendcode{}\nwbegindocs{116}Finally this variable keeps the signature of the metric space.
\nwenddocs{}\nwbegincode{117}\sublabel{NWgqRer-2iImO3-4}\nwmargintag{{\nwtagstyle{}\subpageref{NWgqRer-2iImO3-4}}}\moddef{Definition of variables~{\nwtagstyle{}\subpageref{NWgqRer-2iImO3-1}}}\plusendmoddef\Rm{}\nwstartdeflinemarkup\nwusesondefline{\\{NWgqRer-1p0Y9w-2}}\nwprevnextdefs{NWgqRer-2iImO3-3}{NWgqRer-2iImO3-5}\nwenddeflinemarkup
 {\it{}sign}({\tt{}"s"}, {\tt{}"{\char92}{\char92}sigma"});

\nwused{\\{NWgqRer-1p0Y9w-2}}\nwendcode{}\nwbegindocs{118}This an index used for the definition of Clifford units.
\nwenddocs{}\nwbegincode{119}\sublabel{NWgqRer-2iImO3-5}\nwmargintag{{\nwtagstyle{}\subpageref{NWgqRer-2iImO3-5}}}\moddef{Definition of variables~{\nwtagstyle{}\subpageref{NWgqRer-2iImO3-1}}}\plusendmoddef\Rm{}\nwstartdeflinemarkup\nwusesondefline{\\{NWgqRer-1p0Y9w-2}}\nwprevnextdefs{NWgqRer-2iImO3-4}{NWgqRer-2iImO3-6}\nwenddeflinemarkup
{\bf{}varidx} {\it{}mu}({\bf{}symbol}({\tt{}"mu"}, {\tt{}"{\char92}{\char92}mu"}), 1), {\it{}nu}({\bf{}symbol}({\tt{}"nu"}, {\tt{}"{\char92}{\char92}nu"}), 2);

\nwused{\\{NWgqRer-1p0Y9w-2}}\nwendcode{}\nwbegindocs{120}Three generic points which are used in calculations.
\nwenddocs{}\nwbegincode{121}\sublabel{NWgqRer-2iImO3-6}\nwmargintag{{\nwtagstyle{}\subpageref{NWgqRer-2iImO3-6}}}\moddef{Definition of variables~{\nwtagstyle{}\subpageref{NWgqRer-2iImO3-1}}}\plusendmoddef\Rm{}\nwstartdeflinemarkup\nwusesondefline{\\{NWgqRer-1p0Y9w-2}}\nwprevnextdefs{NWgqRer-2iImO3-5}{NWgqRer-2iImO3-7}\nwenddeflinemarkup
{\bf{}dual\_number} {\it{}P}({\it{}u}, {\it{}v}), {\it{}P1}({\it{}u1}, {\it{}v1}), {\it{}P2}({\it{}u2}, {\it{}v2});

\nwused{\\{NWgqRer-1p0Y9w-2}}\nwendcode{}\nwbegindocs{122}Here we define a parabolic Clifford units {\Tt{}\Rm{}{\it{}one}\nwendquote}, {\Tt{}\Rm{}{\it{}e0}\nwendquote}.
\nwenddocs{}\nwbegincode{123}\sublabel{NWgqRer-2iImO3-7}\nwmargintag{{\nwtagstyle{}\subpageref{NWgqRer-2iImO3-7}}}\moddef{Definition of variables~{\nwtagstyle{}\subpageref{NWgqRer-2iImO3-1}}}\plusendmoddef\Rm{}\nwstartdeflinemarkup\nwusesondefline{\\{NWgqRer-1p0Y9w-2}}\nwprevnextdefs{NWgqRer-2iImO3-6}{\relax}\nwenddeflinemarkup
{\bf{}ex} {\it{}e} = {\it{}clifford\_unit}({\it{}mu}, {\it{}diag\_matrix}({\bf{}lst}({\it{}sign}))),
 {\it{}e0} = {\it{}e}.{\it{}subs}({\it{}mu}\begin{math}\equiv\end{math}0),
 {\it{}one} = {\it{}dirac\_ONE}(),
 {\it{}e2} = {\it{}clifford\_unit}({\it{}nu}, {\it{}diag\_matrix}({\bf{}lst}(-1,{\it{}sign})));

\nwused{\\{NWgqRer-1p0Y9w-2}}\nwendcode{}\nwbegindocs{124}\nwdocspar
\subsubsection{Test routine}
\label{sec:test-routine}
This routine make the same sequence of checks for both cases of
subgroups \(N\) and \(N^\prime\).

\nwenddocs{}\nwbegindocs{125}First we define the reference point {\Tt{}\Rm{}({\it{}u0},{\it{}v0})\nwendquote}.
\nwenddocs{}\nwbegincode{126}\sublabel{NWgqRer-nXe8t-1}\nwmargintag{{\nwtagstyle{}\subpageref{NWgqRer-nXe8t-1}}}\moddef{Test routine~{\nwtagstyle{}\subpageref{NWgqRer-nXe8t-1}}}\endmoddef\Rm{}\nwstartdeflinemarkup\nwusesondefline{\\{NWgqRer-1p0Y9w-2}}\nwprevnextdefs{\relax}{NWgqRer-nXe8t-2}\nwenddeflinemarkup
{\bf{}void} {\it{}parab\_rot\_sub}({\bf{}bool} {\it{}N}, {\bf{}bool} {\it{}S}) {\nwlbrace}\nwindexdefn{\nwixident{parab{\_}rot{\_}sub}}{parab:unrot:unsub}{NWgqRer-nXe8t-1}
        {\it{}cout} \begin{math}\ll\end{math} {\tt{}"Calculations for subgroup {\char92}{\char92}(N"}
                 \begin{math}\ll\end{math} ({\it{}N} ? {\tt{}""} : {\tt{}"'"}) \begin{math}\ll\end{math} {\tt{}"{\char92}{\char92}) and "}
                 \begin{math}\ll\end{math} ({\it{}S} ? {\tt{}"straight"} : {\tt{}"geodesics"}) \begin{math}\ll\end{math} {\tt{}" spokes{\char92}{\char92}{\char92}{\char92}"}\begin{math}\ll\end{math}{\it{}endl};
        {\it{}set\_subgroup}({\it{}N});
        {\it{}set\_straight\_spoke}({\it{}S});
        {\bf{}ex} {\it{}X},
                {\it{}W0}={\it{}dn\_from\_arg\_mod}({\it{}Arg0}, 1),
                {\it{}u0}={\it{}W0}.{\it{}op}(0),
                {\it{}v0}={\it{}W0}.{\it{}op}(1),
                {\it{}P0}={\bf{}matrix}(1, 2, {\bf{}lst}({\it{}u0}, {\it{}v0}));
        
\nwalsodefined{\\{NWgqRer-nXe8t-2}\\{NWgqRer-nXe8t-3}\\{NWgqRer-nXe8t-4}\\{NWgqRer-nXe8t-5}}\nwused{\\{NWgqRer-1p0Y9w-2}}\nwidentdefs{\\{{\nwixident{parab{\_}rot{\_}sub}}{parab:unrot:unsub}}}\nwendcode{}\nwbegindocs{127}These two matrices define the Cayley transform and its inverse.
\nwenddocs{}\nwbegincode{128}\sublabel{NWgqRer-nXe8t-2}\nwmargintag{{\nwtagstyle{}\subpageref{NWgqRer-nXe8t-2}}}\moddef{Test routine~{\nwtagstyle{}\subpageref{NWgqRer-nXe8t-1}}}\plusendmoddef\Rm{}\nwstartdeflinemarkup\nwusesondefline{\\{NWgqRer-1p0Y9w-2}}\nwprevnextdefs{NWgqRer-nXe8t-1}{NWgqRer-nXe8t-3}\nwenddeflinemarkup
{\it{}cout} \begin{math}\ll\end{math} {\tt{}"Reference point: {\char92}{\char92}("} \begin{math}\ll\end{math} {\it{}P0} \begin{math}\ll\end{math} {\tt{}"{\char92}{\char92}){\char92}{\char92}{\char92}{\char92}"} \begin{math}\ll\end{math} {\it{}endl};
{\it{}cout} \begin{math}\ll\end{math} {\tt{}"Reference point arg: {\char92}{\char92}("} \begin{math}\ll\end{math} {\it{}ex\_to}\begin{math}<\end{math}{\bf{}dual\_number}\begin{math}>\end{math}({\it{}W0}).{\it{}arg}() \begin{math}\ll\end{math} {\tt{}"{\char92}{\char92}){\char92}{\char92}{\char92}{\char92}"} \begin{math}\ll\end{math} {\it{}endl};
{\it{}cout} \begin{math}\ll\end{math} {\tt{}"Reference point norm: {\char92}{\char92}("} \begin{math}\ll\end{math} {\it{}ex\_to}\begin{math}<\end{math}{\bf{}dual\_number}\begin{math}>\end{math}({\it{}W0}).{\it{}norm}() \begin{math}\ll\end{math} {\tt{}"{\char92}{\char92}){\char92}{\char92}{\char92}{\char92}"} \begin{math}\ll\end{math} {\it{}endl};

\nwused{\\{NWgqRer-1p0Y9w-2}}\nwendcode{}\nwbegindocs{129}For the subgroup \(N\) we consider upper-triangular matrices, for
\(N^\prime\)---lower-triangular.
\nwenddocs{}\nwbegincode{130}\sublabel{NWgqRer-nXe8t-3}\nwmargintag{{\nwtagstyle{}\subpageref{NWgqRer-nXe8t-3}}}\moddef{Test routine~{\nwtagstyle{}\subpageref{NWgqRer-nXe8t-1}}}\plusendmoddef\Rm{}\nwstartdeflinemarkup\nwusesondefline{\\{NWgqRer-1p0Y9w-2}}\nwprevnextdefs{NWgqRer-nXe8t-2}{NWgqRer-nXe8t-4}\nwenddeflinemarkup
{\bf{}if} ({\it{}subgroup\_N})
        {\it{}X}={\bf{}matrix}(2, 2, {\bf{}lst}({\it{}one}, {\it{}one}\begin{math}\ast\end{math}{\it{}x}, 0, {\it{}one}));
{\bf{}else}
        {\it{}X}={\bf{}matrix}(2, 2, {\bf{}lst}({\it{}one}, 0, {\it{}one}\begin{math}\ast\end{math}{\it{}x}, {\it{}one}));

\nwused{\\{NWgqRer-1p0Y9w-2}}\nwendcode{}\nwbegindocs{131}Two different types of Cayley transforms.
\nwenddocs{}\nwbegincode{132}\sublabel{NWgqRer-nXe8t-4}\nwmargintag{{\nwtagstyle{}\subpageref{NWgqRer-nXe8t-4}}}\moddef{Test routine~{\nwtagstyle{}\subpageref{NWgqRer-nXe8t-1}}}\plusendmoddef\Rm{}\nwstartdeflinemarkup\nwusesondefline{\\{NWgqRer-1p0Y9w-2}}\nwprevnextdefs{NWgqRer-nXe8t-3}{NWgqRer-nXe8t-5}\nwenddeflinemarkup
{\bf{}ex}   {\it{}TC}, {\it{}TCI};
{\bf{}if} ({\it{}S}) {\nwlbrace}
        {\it{}TC}={\bf{}matrix}(2, 2, {\bf{}lst}({\it{}one}, -{\it{}e0}, -{\it{}e0}, {\it{}one}));
        {\it{}TCI}={\bf{}matrix}(2, 2, {\bf{}lst}({\it{}one}, {\it{}e0}, {\it{}e0}, {\it{}one}));
{\nwrbrace} {\bf{}else}{\nwlbrace}
        {\it{}TC}={\bf{}matrix}(2, 2, {\bf{}lst}({\it{}one}, -{\it{}e0}\begin{math}\ast\end{math}{\it{}half}, -{\it{}e0}\begin{math}\ast\end{math}{\it{}half}, {\it{}one}));
        {\it{}TCI}={\bf{}matrix}(2, 2, {\bf{}lst}({\it{}one}, {\it{}e0}\begin{math}\ast\end{math}{\it{}half}, {\it{}e0}\begin{math}\ast\end{math}{\it{}half}, {\it{}one}));
{\nwrbrace}

\nwused{\\{NWgqRer-1p0Y9w-2}}\nwendcode{}\nwbegindocs{133}Common part of test routine.
\nwenddocs{}\nwbegincode{134}\sublabel{NWgqRer-nXe8t-5}\nwmargintag{{\nwtagstyle{}\subpageref{NWgqRer-nXe8t-5}}}\moddef{Test routine~{\nwtagstyle{}\subpageref{NWgqRer-nXe8t-1}}}\plusendmoddef\Rm{}\nwstartdeflinemarkup\nwusesondefline{\\{NWgqRer-1p0Y9w-2}}\nwprevnextdefs{NWgqRer-nXe8t-4}{\relax}\nwenddeflinemarkup
 \LA{}Show expressions~{\nwtagstyle{}\subpageref{NWgqRer-1Pc9Jw-1}}\RA{}
 \LA{}Check identities~{\nwtagstyle{}\subpageref{NWgqRer-22TNkn-1}}\RA{}
 {\it{}cout} \begin{math}\ll\end{math}  ({\it{}latexout} ? {\tt{}"{\char92}{\char92}vspace{\char123}2mm{\char125}{\char92}{\char92}hrule"} :
                   {\tt{}"----------------------------------------"} ) \begin{math}\ll\end{math} {\it{}endl};
{\nwrbrace}

\nwused{\\{NWgqRer-1p0Y9w-2}}\nwendcode{}\nwbegindocs{135}\nwdocspar
\subsubsection{Main procedure}
\label{sec:main-procedure}
It just calls the test routine, calculates the induced representation and
draws a few pictures.

\nwenddocs{}\nwbegindocs{136} We output formulae in \LaTeX\ mode.
\nwenddocs{}\nwbegincode{137}\sublabel{NWgqRer-3sxui-1}\nwmargintag{{\nwtagstyle{}\subpageref{NWgqRer-3sxui-1}}}\moddef{Main procedure~{\nwtagstyle{}\subpageref{NWgqRer-3sxui-1}}}\endmoddef\Rm{}\nwstartdeflinemarkup\nwusesondefline{\\{NWgqRer-1p0Y9w-2}}\nwprevnextdefs{\relax}{NWgqRer-3sxui-2}\nwenddeflinemarkup
{\bf{}int} {\it{}main}(){\nwlbrace}\nwindexdefn{\nwixident{main}}{main}{NWgqRer-3sxui-1}
 {\it{}latexout}={\bf{}true};

\nwalsodefined{\\{NWgqRer-3sxui-2}\\{NWgqRer-3sxui-3}\\{NWgqRer-3sxui-4}\\{NWgqRer-3sxui-5}}\nwused{\\{NWgqRer-1p0Y9w-2}}\nwidentdefs{\\{{\nwixident{main}}{main}}}\nwendcode{}\nwbegindocs{138}Preparation of output stream.
\nwenddocs{}\nwbegincode{139}\sublabel{NWgqRer-3sxui-2}\nwmargintag{{\nwtagstyle{}\subpageref{NWgqRer-3sxui-2}}}\moddef{Main procedure~{\nwtagstyle{}\subpageref{NWgqRer-3sxui-1}}}\plusendmoddef\Rm{}\nwstartdeflinemarkup\nwusesondefline{\\{NWgqRer-1p0Y9w-2}}\nwprevnextdefs{NWgqRer-3sxui-1}{NWgqRer-3sxui-3}\nwenddeflinemarkup
 {\it{}cout} \begin{math}\ll\end{math} {\it{}boolalpha};
 {\bf{}if} ({\it{}latexout})
         {\it{}cout} \begin{math}\ll\end{math} {\it{}latex};
 \LA{}Drawing pictures~{\nwtagstyle{}\subpageref{NWgqRer-Ghke2-1}}\RA{}

\nwused{\\{NWgqRer-1p0Y9w-2}}\nwendcode{}\nwbegindocs{140}Now we call the test routine for both \(N\) and \(N^\prime\) subgroups.
\nwenddocs{}\nwbegincode{141}\sublabel{NWgqRer-3sxui-3}\nwmargintag{{\nwtagstyle{}\subpageref{NWgqRer-3sxui-3}}}\moddef{Main procedure~{\nwtagstyle{}\subpageref{NWgqRer-3sxui-1}}}\plusendmoddef\Rm{}\nwstartdeflinemarkup\nwusesondefline{\\{NWgqRer-1p0Y9w-2}}\nwprevnextdefs{NWgqRer-3sxui-2}{NWgqRer-3sxui-4}\nwenddeflinemarkup
 {\it{}parab\_rot\_sub}({\bf{}true}, {\bf{}true});
 {\it{}parab\_rot\_sub}({\bf{}false}, {\bf{}true});
// parab\_rot\_sub(false, false); To work with geodesic spokes

\nwused{\\{NWgqRer-1p0Y9w-2}}\nwidentuses{\\{{\nwixident{parab{\_}rot{\_}sub}}{parab:unrot:unsub}}}\nwindexuse{\nwixident{parab{\_}rot{\_}sub}}{parab:unrot:unsub}{NWgqRer-3sxui-3}\nwendcode{}\nwbegindocs{142} Calculation of induced representations formulae.
\nwenddocs{}\nwbegincode{143}\sublabel{NWgqRer-3sxui-4}\nwmargintag{{\nwtagstyle{}\subpageref{NWgqRer-3sxui-4}}}\moddef{Main procedure~{\nwtagstyle{}\subpageref{NWgqRer-3sxui-1}}}\plusendmoddef\Rm{}\nwstartdeflinemarkup\nwusesondefline{\\{NWgqRer-1p0Y9w-2}}\nwprevnextdefs{NWgqRer-3sxui-3}{NWgqRer-3sxui-5}\nwenddeflinemarkup
        \LA{}Induced representations~{\nwtagstyle{}\subpageref{NWgqRer-1h89LW-1}}\RA{}

\nwused{\\{NWgqRer-1p0Y9w-2}}\nwendcode{}\nwbegindocs{144}And we finishing by drawing several pictures for
Figs.~\ref{fig:rotations} and~\ref{fig:p-rotations}.
\nwenddocs{}\nwbegincode{145}\sublabel{NWgqRer-3sxui-5}\nwmargintag{{\nwtagstyle{}\subpageref{NWgqRer-3sxui-5}}}\moddef{Main procedure~{\nwtagstyle{}\subpageref{NWgqRer-3sxui-1}}}\plusendmoddef\Rm{}\nwstartdeflinemarkup\nwusesondefline{\\{NWgqRer-1p0Y9w-2}}\nwprevnextdefs{NWgqRer-3sxui-4}{\relax}\nwenddeflinemarkup
{\nwrbrace}
\nwused{\\{NWgqRer-1p0Y9w-2}}\nwendcode{}\nwbegindocs{146}\nwdocspar
\subsection{Drawing Orbits}
\label{sec:drawing-orbits}
To draw cycles we use {\Tt{}\Rm{}{\bf{}cycle}\nwendquote} library~\cite{Kisil05b}.

\nwenddocs{}\nwbegindocs{147}Elliptic orbits (circles).
\nwenddocs{}\nwbegincode{148}\sublabel{NWgqRer-Ghke2-1}\nwmargintag{{\nwtagstyle{}\subpageref{NWgqRer-Ghke2-1}}}\moddef{Drawing pictures~{\nwtagstyle{}\subpageref{NWgqRer-Ghke2-1}}}\endmoddef\Rm{}\nwstartdeflinemarkup\nwusesondefline{\\{NWgqRer-3sxui-2}}\nwprevnextdefs{\relax}{NWgqRer-Ghke2-2}\nwenddeflinemarkup
{\it{}ofstream} {\it{}asymptote}({\tt{}"parab-rot-data.asy"});
{\it{}asymptote} \begin{math}\ll\end{math} {\tt{}"//This file is automatically created from the code in parab-rotation.nw"} 
                \begin{math}\ll\end{math} {\it{}endl} \begin{math}\ll\end{math} {\it{}endl}
                \begin{math}\ll\end{math} {\tt{}"path[] K="};
{\bf{}for}({\bf{}int} {\it{}i}=0; {\it{}i}\begin{math}<\end{math}6; {\it{}i}\protect\PP)
        {\bf{}cycle2D}({\bf{}lst}(0,0),{\it{}e2}.{\it{}subs}({\it{}sign}\begin{math}\equiv\end{math}-1),{\it{}i}\begin{math}\ast\end{math}{\it{}i}\begin{math}\ast\end{math}.04)
                .{\it{}asy\_path}({\it{}asymptote}, -1.5, 1.5, -2, 2, 0, ({\it{}i}\begin{math}>\end{math}0));
{\it{}asymptote} \begin{math}\ll\end{math} {\tt{}";"} \begin{math}\ll\end{math} {\it{}endl};
{\it{}asymptote} \begin{math}\ll\end{math} {\tt{}"path[] Kb="};
{\bf{}cycle2D}({\bf{}lst}(0,0),{\it{}e2}.{\it{}subs}({\it{}sign}\begin{math}\equiv\end{math}-1),1)
        .{\it{}asy\_path}({\it{}asymptote}, -1.5, 1.5, -2, 2, 0, {\bf{}false});
{\it{}asymptote} \begin{math}\ll\end{math} {\tt{}";"} \begin{math}\ll\end{math} {\it{}endl};

\nwalsodefined{\\{NWgqRer-Ghke2-2}\\{NWgqRer-Ghke2-3}\\{NWgqRer-Ghke2-4}\\{NWgqRer-Ghke2-5}\\{NWgqRer-Ghke2-6}\\{NWgqRer-Ghke2-7}\\{NWgqRer-Ghke2-8}}\nwused{\\{NWgqRer-3sxui-2}}\nwendcode{}\nwbegindocs{149}Hyperbolic orbits.
\nwenddocs{}\nwbegincode{150}\sublabel{NWgqRer-Ghke2-2}\nwmargintag{{\nwtagstyle{}\subpageref{NWgqRer-Ghke2-2}}}\moddef{Drawing pictures~{\nwtagstyle{}\subpageref{NWgqRer-Ghke2-1}}}\plusendmoddef\Rm{}\nwstartdeflinemarkup\nwusesondefline{\\{NWgqRer-3sxui-2}}\nwprevnextdefs{NWgqRer-Ghke2-1}{NWgqRer-Ghke2-3}\nwenddeflinemarkup
{\it{}asymptote} \begin{math}\ll\end{math} {\tt{}"path[] A="};
{\bf{}for}({\bf{}int} {\it{}i}=0; {\it{}i}\begin{math}<\end{math}6; {\it{}i}\protect\PP) {\nwlbrace}
 {\bf{}cycle2D}({\bf{}lst}(0,0),{\it{}e2}.{\it{}subs}({\it{}sign}\begin{math}\equiv\end{math}1),-{\it{}i}\begin{math}\ast\end{math}{\it{}i}\begin{math}\ast\end{math}.04)
         .{\it{}asy\_path}({\it{}asymptote}, -1.75, 1.75, -1.5, 2, 0, ({\it{}i}\begin{math}>\end{math}0));
 {\bf{}cycle2D}({\bf{}lst}(0,0),{\it{}e2}.{\it{}subs}({\it{}sign}\begin{math}\equiv\end{math}1),{\it{}i}\begin{math}\ast\end{math}{\it{}i}\begin{math}\ast\end{math}.12)
         .{\it{}asy\_path}({\it{}asymptote}, -1.75, 1.75, -1.5, 2, 0, {\bf{}true});
{\nwrbrace}
{\it{}asymptote} \begin{math}\ll\end{math} {\tt{}";"} \begin{math}\ll\end{math} {\it{}endl};
{\it{}asymptote} \begin{math}\ll\end{math} {\tt{}"path[] Ab="};
{\bf{}cycle2D}({\bf{}lst}(0,0),{\it{}e2}.{\it{}subs}({\it{}sign}\begin{math}\equiv\end{math}1),-1)
        .{\it{}asy\_path}({\it{}asymptote}, -1.75, 1.75, -1.5, 2, 0, {\bf{}false});
{\it{}asymptote} \begin{math}\ll\end{math} {\tt{}";"} \begin{math}\ll\end{math} {\it{}endl};

\nwused{\\{NWgqRer-3sxui-2}}\nwendcode{}\nwbegindocs{151}Hyperbolic orbits for reflected orbits. 
\nwenddocs{}\nwbegincode{152}\sublabel{NWgqRer-Ghke2-3}\nwmargintag{{\nwtagstyle{}\subpageref{NWgqRer-Ghke2-3}}}\moddef{Drawing pictures~{\nwtagstyle{}\subpageref{NWgqRer-Ghke2-1}}}\plusendmoddef\Rm{}\nwstartdeflinemarkup\nwusesondefline{\\{NWgqRer-3sxui-2}}\nwprevnextdefs{NWgqRer-Ghke2-2}{NWgqRer-Ghke2-4}\nwenddeflinemarkup
{\it{}asymptote} \begin{math}\ll\end{math} {\tt{}"path[] At="};
{\bf{}for}({\bf{}int} {\it{}i}=0; {\it{}i}\begin{math}<\end{math}6; {\it{}i}\protect\PP) {\nwlbrace}
        {\bf{}cycle2D}({\bf{}lst}(0,0),{\it{}e2}.{\it{}subs}({\it{}sign}\begin{math}\equiv\end{math}1),{\it{}i}\begin{math}\ast\end{math}{\it{}i}\begin{math}\ast\end{math}.04)
                .{\it{}asy\_path}({\it{}asymptote}, -1.75, 1.75, -1.5, 2, 0, ({\it{}i}\begin{math}>\end{math}0));
        {\bf{}cycle2D}({\bf{}lst}(0,0),{\it{}e2}.{\it{}subs}({\it{}sign}\begin{math}\equiv\end{math}1),-{\it{}i}\begin{math}\ast\end{math}{\it{}i}\begin{math}\ast\end{math}.12)
                .{\it{}asy\_path}({\it{}asymptote}, -1.75, 1.75, -1.5, 2, 0, {\bf{}true});
{\nwrbrace}
{\it{}asymptote} \begin{math}\ll\end{math} {\tt{}";"} \begin{math}\ll\end{math} {\it{}endl};
{\it{}asymptote} \begin{math}\ll\end{math} {\tt{}"path[] Abt="};
{\bf{}cycle2D}({\bf{}lst}(0,0),{\it{}e2}.{\it{}subs}({\it{}sign}\begin{math}\equiv\end{math}1),1)
                .{\it{}asy\_path}({\it{}asymptote}, -1.75, 1.75, -1.5, 2, 0, {\bf{}false});
{\it{}asymptote} \begin{math}\ll\end{math} {\tt{}";"} \begin{math}\ll\end{math} {\it{}endl};

\nwused{\\{NWgqRer-3sxui-2}}\nwendcode{}\nwbegindocs{153}Parabolic orbits, subgroup \(N\).
\nwenddocs{}\nwbegincode{154}\sublabel{NWgqRer-Ghke2-4}\nwmargintag{{\nwtagstyle{}\subpageref{NWgqRer-Ghke2-4}}}\moddef{Drawing pictures~{\nwtagstyle{}\subpageref{NWgqRer-Ghke2-1}}}\plusendmoddef\Rm{}\nwstartdeflinemarkup\nwusesondefline{\\{NWgqRer-3sxui-2}}\nwprevnextdefs{NWgqRer-Ghke2-3}{NWgqRer-Ghke2-5}\nwenddeflinemarkup
{\it{}asymptote} \begin{math}\ll\end{math} {\tt{}"path[] N="};
{\bf{}for}({\bf{}int} {\it{}i}=0; {\it{}i}\begin{math}<\end{math}6; {\it{}i}\protect\PP)
        {\bf{}cycle2D}(1,{\bf{}lst}(0,{\bf{}numeric}(1,2)),{\bf{}numeric}({\it{}i},2)-1,{\it{}e2}.{\it{}subs}({\it{}sign}\begin{math}\equiv\end{math}0))
                .{\it{}asy\_path}({\it{}asymptote}, -1.5, 1.5, -2, 2, 0, ({\it{}i}\begin{math}>\end{math}0));
{\it{}asymptote} \begin{math}\ll\end{math} {\tt{}";"} \begin{math}\ll\end{math} {\it{}endl};
{\it{}asymptote} \begin{math}\ll\end{math} {\tt{}"path[] Nb="};
{\bf{}cycle2D}(1,{\bf{}lst}(0,{\bf{}numeric}(1,2)),-1,{\it{}e2}.{\it{}subs}({\it{}sign}\begin{math}\equiv\end{math}0))
        .{\it{}asy\_path}({\it{}asymptote}, -1.5, 1.5, -2, 2, 0, {\bf{}false});
{\it{}asymptote} \begin{math}\ll\end{math} {\tt{}";"} \begin{math}\ll\end{math} {\it{}endl};

\nwused{\\{NWgqRer-3sxui-2}}\nwidentuses{\\{{\nwixident{numeric}}{numeric}}}\nwindexuse{\nwixident{numeric}}{numeric}{NWgqRer-Ghke2-4}\nwendcode{}\nwbegindocs{155} Parabolic orbits, subgroup \(N^\prime\).
\nwenddocs{}\nwbegincode{156}\sublabel{NWgqRer-Ghke2-5}\nwmargintag{{\nwtagstyle{}\subpageref{NWgqRer-Ghke2-5}}}\moddef{Drawing pictures~{\nwtagstyle{}\subpageref{NWgqRer-Ghke2-1}}}\plusendmoddef\Rm{}\nwstartdeflinemarkup\nwusesondefline{\\{NWgqRer-3sxui-2}}\nwprevnextdefs{NWgqRer-Ghke2-4}{NWgqRer-Ghke2-6}\nwenddeflinemarkup
{\it{}asymptote} \begin{math}\ll\end{math} {\tt{}"path[] N1="};
{\bf{}for}({\bf{}int} {\it{}i}=0; {\it{}i}\begin{math}<\end{math}5; {\it{}i}\protect\PP)
        {\bf{}cycle2D}(.5\begin{math}\ast\end{math}{\it{}i}\begin{math}\ast\end{math}{\it{}i}\begin{math}\ast\end{math}{\it{}i}+1,{\bf{}lst}(0,{\bf{}numeric}(1,2)),-1,{\it{}e2}.{\it{}subs}({\it{}sign}\begin{math}\equiv\end{math}0))
                .{\it{}asy\_path}({\it{}asymptote}, -1.5, 1.5, -1.5, 2, 0, ({\it{}i}\begin{math}>\end{math}0));
{\it{}asymptote} \begin{math}\ll\end{math} {\tt{}";"} \begin{math}\ll\end{math} {\it{}endl};
{\it{}asymptote} \begin{math}\ll\end{math} {\tt{}"path[] N1b="};
{\bf{}cycle2D}(1,{\bf{}lst}(0,{\bf{}numeric}(1,2)),-1,{\it{}e2}.{\it{}subs}({\it{}sign}\begin{math}\equiv\end{math}0))
        .{\it{}asy\_path}({\it{}asymptote}, -1.5, 1.5, -1.5, 2, 0, {\bf{}false});
{\it{}asymptote} \begin{math}\ll\end{math} {\tt{}";"} \begin{math}\ll\end{math} {\it{}endl};

\nwused{\\{NWgqRer-3sxui-2}}\nwidentuses{\\{{\nwixident{numeric}}{numeric}}}\nwindexuse{\nwixident{numeric}}{numeric}{NWgqRer-Ghke2-5}\nwendcode{}\nwbegindocs{157}Parabolic orbits, subgroup \(N\) (geodesic version?).
\nwenddocs{}\nwbegincode{158}\sublabel{NWgqRer-Ghke2-6}\nwmargintag{{\nwtagstyle{}\subpageref{NWgqRer-Ghke2-6}}}\moddef{Drawing pictures~{\nwtagstyle{}\subpageref{NWgqRer-Ghke2-1}}}\plusendmoddef\Rm{}\nwstartdeflinemarkup\nwusesondefline{\\{NWgqRer-3sxui-2}}\nwprevnextdefs{NWgqRer-Ghke2-5}{NWgqRer-Ghke2-7}\nwenddeflinemarkup
{\it{}asymptote} \begin{math}\ll\end{math} {\tt{}"path[] Ng="};
{\bf{}for}({\bf{}int} {\it{}i}=0; {\it{}i}\begin{math}<\end{math}6; {\it{}i}\protect\PP)
        {\bf{}cycle2D}(0.5,{\bf{}lst}(0,{\bf{}numeric}(1,2)),{\bf{}numeric}({\it{}i},2)-.5,{\it{}e2}.{\it{}subs}({\it{}sign}\begin{math}\equiv\end{math}0))
                .{\it{}asy\_path}({\it{}asymptote}, -1.5, 1.5, -2, 2, 0, ({\it{}i}\begin{math}>\end{math}0));
{\it{}asymptote} \begin{math}\ll\end{math} {\tt{}";"} \begin{math}\ll\end{math} {\it{}endl};
{\it{}asymptote} \begin{math}\ll\end{math} {\tt{}"path[] Ngb="};
{\bf{}cycle2D}(0.5,{\bf{}lst}(0,{\bf{}numeric}(1,2)),-.5,{\it{}e2}.{\it{}subs}({\it{}sign}\begin{math}\equiv\end{math}0))
        .{\it{}asy\_path}({\it{}asymptote}, -1.5, 1.5, -2, 2, 0, {\bf{}false});
{\it{}asymptote} \begin{math}\ll\end{math} {\tt{}";"} \begin{math}\ll\end{math} {\it{}endl};

\nwused{\\{NWgqRer-3sxui-2}}\nwidentuses{\\{{\nwixident{numeric}}{numeric}}}\nwindexuse{\nwixident{numeric}}{numeric}{NWgqRer-Ghke2-6}\nwendcode{}\nwbegindocs{159} Parabolic orbits, subgroup \(N^\prime\) (geodesic version?).
\nwenddocs{}\nwbegincode{160}\sublabel{NWgqRer-Ghke2-7}\nwmargintag{{\nwtagstyle{}\subpageref{NWgqRer-Ghke2-7}}}\moddef{Drawing pictures~{\nwtagstyle{}\subpageref{NWgqRer-Ghke2-1}}}\plusendmoddef\Rm{}\nwstartdeflinemarkup\nwusesondefline{\\{NWgqRer-3sxui-2}}\nwprevnextdefs{NWgqRer-Ghke2-6}{NWgqRer-Ghke2-8}\nwenddeflinemarkup
{\it{}asymptote} \begin{math}\ll\end{math} {\tt{}"path[] N1g="};
{\bf{}for}({\bf{}int} {\it{}i}=0; {\it{}i}\begin{math}<\end{math}5; {\it{}i}\protect\PP)
        {\bf{}cycle2D}(.25\begin{math}\ast\end{math}{\it{}i}\begin{math}\ast\end{math}{\it{}i}\begin{math}\ast\end{math}{\it{}i}+.5,{\bf{}lst}(0,{\bf{}numeric}(1,2)),-.5,{\it{}e2}.{\it{}subs}({\it{}sign}\begin{math}\equiv\end{math}0))
                .{\it{}asy\_path}({\it{}asymptote}, -1.5, 1.5, -1.5, 2, 0, ({\it{}i}\begin{math}>\end{math}0));
{\it{}asymptote} \begin{math}\ll\end{math} {\tt{}";"} \begin{math}\ll\end{math} {\it{}endl};
{\it{}asymptote} \begin{math}\ll\end{math} {\tt{}"path[] N1gb="};
{\bf{}cycle2D}(.5,{\bf{}lst}(0,{\bf{}numeric}(1,2)),-.5,{\it{}e2}.{\it{}subs}({\it{}sign}\begin{math}\equiv\end{math}0))
        .{\it{}asy\_path}({\it{}asymptote}, -1.5, 1.5, -1.5, 2, 0, {\bf{}false});
{\it{}asymptote} \begin{math}\ll\end{math} {\tt{}";"} \begin{math}\ll\end{math} {\it{}endl};

\nwused{\\{NWgqRer-3sxui-2}}\nwidentuses{\\{{\nwixident{numeric}}{numeric}}}\nwindexuse{\nwixident{numeric}}{numeric}{NWgqRer-Ghke2-7}\nwendcode{}\nwbegindocs{161} Parabolic orbits, subgroup \(N^\prime\) elliptic geodesics.
\nwenddocs{}\nwbegincode{162}\sublabel{NWgqRer-Ghke2-8}\nwmargintag{{\nwtagstyle{}\subpageref{NWgqRer-Ghke2-8}}}\moddef{Drawing pictures~{\nwtagstyle{}\subpageref{NWgqRer-Ghke2-1}}}\plusendmoddef\Rm{}\nwstartdeflinemarkup\nwusesondefline{\\{NWgqRer-3sxui-2}}\nwprevnextdefs{NWgqRer-Ghke2-7}{\relax}\nwenddeflinemarkup
{\bf{}numeric} {\it{}L}[9]={\nwlbrace}-2,-1.4,-.8,-.4,0,.4,.8,1.2,2{\nwrbrace}, {\it{}low}, {\it{}high};

{\it{}asymptote} \begin{math}\ll\end{math} {\tt{}"path[] N1geoE="};
{\bf{}for}({\bf{}int} {\it{}i}=0; {\it{}i}\begin{math}<\end{math}9; {\it{}i}\protect\PP){\nwlbrace}
        {\bf{}if} ({\it{}abs}({\it{}L}[{\it{}i}])\begin{math}<\end{math}.5) {\nwlbrace}
                        {\it{}low}={\it{}min}((4\begin{math}\ast\end{math}{\it{}L}[{\it{}i}]-2)\begin{math}\div\end{math}(4\begin{math}\ast\end{math}{\it{}L}[{\it{}i}]\begin{math}\ast\end{math}{\it{}L}[{\it{}i}]-1),(4\begin{math}\ast\end{math}{\it{}L}[{\it{}i}]+2)\begin{math}\div\end{math}(4\begin{math}\ast\end{math}{\it{}L}[{\it{}i}]\begin{math}\ast\end{math}{\it{}L}[{\it{}i}]-1));
                        {\it{}high}={\it{}max}((4\begin{math}\ast\end{math}{\it{}L}[{\it{}i}]-2)\begin{math}\div\end{math}(4\begin{math}\ast\end{math}{\it{}L}[{\it{}i}]\begin{math}\ast\end{math}{\it{}L}[{\it{}i}]-1),(4\begin{math}\ast\end{math}{\it{}L}[{\it{}i}]+2)\begin{math}\div\end{math}(4\begin{math}\ast\end{math}{\it{}L}[{\it{}i}]\begin{math}\ast\end{math}{\it{}L}[{\it{}i}]-1));
                {\nwrbrace} {\bf{}else} {\bf{}if} ({\it{}L}[{\it{}i}]\begin{math}>\end{math}.5) {\nwlbrace}
                        {\it{}low}={\bf{}numeric}(-3,2);
                        {\it{}high}={\it{}min}((4\begin{math}\ast\end{math}{\it{}L}[{\it{}i}]-2)\begin{math}\div\end{math}(4\begin{math}\ast\end{math}{\it{}L}[{\it{}i}]\begin{math}\ast\end{math}{\it{}L}[{\it{}i}]-1),(4\begin{math}\ast\end{math}{\it{}L}[{\it{}i}]+2)\begin{math}\div\end{math}(4\begin{math}\ast\end{math}{\it{}L}[{\it{}i}]\begin{math}\ast\end{math}{\it{}L}[{\it{}i}]-1));
                {\nwrbrace} {\bf{}else} {\nwlbrace}
                        {\it{}low}={\it{}max}((4\begin{math}\ast\end{math}{\it{}L}[{\it{}i}]-2)\begin{math}\div\end{math}(4\begin{math}\ast\end{math}{\it{}L}[{\it{}i}]\begin{math}\ast\end{math}{\it{}L}[{\it{}i}]-1),(4\begin{math}\ast\end{math}{\it{}L}[{\it{}i}]+2)\begin{math}\div\end{math}(4\begin{math}\ast\end{math}{\it{}L}[{\it{}i}]\begin{math}\ast\end{math}{\it{}L}[{\it{}i}]-1));
                        {\it{}high}={\bf{}numeric}(3,2);
                {\nwrbrace}
        {\bf{}cycle2D}(-1+4+4\begin{math}\ast\end{math}{\it{}L}[{\it{}i}]\begin{math}\ast\end{math}{\it{}L}[{\it{}i}],{\bf{}lst}(4\begin{math}\ast\end{math}{\it{}L}[{\it{}i}],2),0,{\it{}e2}.{\it{}subs}({\it{}sign}\begin{math}\equiv\end{math}0))
                .{\it{}asy\_path}({\it{}asymptote}, {\it{}max}({\bf{}numeric}(-3,2),{\it{}low}), {\it{}min}({\bf{}numeric}(3,2),{\it{}high}), -1.5, 2, 0,({\it{}i}\begin{math}>\end{math}0));
{\nwrbrace}
{\it{}asymptote} \begin{math}\ll\end{math} {\tt{}";"} \begin{math}\ll\end{math} {\it{}endl};

{\it{}asymptote} \begin{math}\ll\end{math} {\tt{}"path[] N1geoEf="};
{\bf{}for}({\bf{}int} {\it{}i}=0; {\it{}i}\begin{math}<\end{math}9; {\it{}i}\protect\PP){\nwlbrace}
        {\bf{}cycle2D}(-1+4+4\begin{math}\ast\end{math}{\it{}L}[{\it{}i}]\begin{math}\ast\end{math}{\it{}L}[{\it{}i}],{\bf{}lst}(4\begin{math}\ast\end{math}{\it{}L}[{\it{}i}],2),0,{\it{}e2}.{\it{}subs}({\it{}sign}\begin{math}\equiv\end{math}0))
                .{\it{}asy\_path}({\it{}asymptote}, -1.5, 1.5, -1.5, 2, 0,({\it{}i}\begin{math}>\end{math}0));
{\nwrbrace}
{\it{}asymptote} \begin{math}\ll\end{math} {\tt{}";"} \begin{math}\ll\end{math} {\it{}endl};

{\it{}asymptote}.{\it{}close}();

\nwused{\\{NWgqRer-3sxui-2}}\nwidentuses{\\{{\nwixident{numeric}}{numeric}}}\nwindexuse{\nwixident{numeric}}{numeric}{NWgqRer-Ghke2-8}\nwendcode{}\nwbegindocs{163}\nwdocspar
\subsubsection{Output routines}
\label{sec:output-routines}

We use standardised routines to output results of calculations.
\nwenddocs{}\nwbegincode{164}\sublabel{NWgqRer-4MrGQl-1}\nwmargintag{{\nwtagstyle{}\subpageref{NWgqRer-4MrGQl-1}}}\moddef{Output routines~{\nwtagstyle{}\subpageref{NWgqRer-4MrGQl-1}}}\endmoddef\Rm{}\nwstartdeflinemarkup\nwprevnextdefs{\relax}{NWgqRer-4MrGQl-2}\nwenddeflinemarkup
{\bf{}void} {\it{}formula\_out}({\it{}string} {\it{}S}, {\bf{}const} {\bf{}ex} & {\it{}F}, {\bf{}bool} {\it{}lineend}) {\nwlbrace}\nwindexdefn{\nwixident{formula{\_}out}}{formula:unout}{NWgqRer-4MrGQl-1}
        {\it{}cout} \begin{math}\ll\end{math} {\it{}S} \begin{math}\ll\end{math} ({\it{}latexout} ? {\tt{}"{\char92}{\char92}("} : {\tt{}""} ) \begin{math}\ll\end{math} {\it{}F} \begin{math}\ll\end{math} ({\it{}latexout} ? {\tt{}"{\char92}{\char92})"} : {\tt{}""} );
        {\bf{}if} ({\it{}lineend})
                {\it{}cout} \begin{math}\ll\end{math} ({\it{}latexout} ? {\tt{}"{\char92}{\char92}{\char92}{\char92}"} : {\tt{}""} ) \begin{math}\ll\end{math} {\it{}endl};
        {\bf{}else}
                {\it{}cout} \begin{math}\ll\end{math} {\tt{}"; "};
{\nwrbrace}

\nwalsodefined{\\{NWgqRer-4MrGQl-2}}\nwnotused{Output routines}\nwidentdefs{\\{{\nwixident{formula{\_}out}}{formula:unout}}}\nwendcode{}\nwbegindocs{165}This routine is used to check identities.
\nwenddocs{}\nwbegincode{166}\sublabel{NWgqRer-4MrGQl-2}\nwmargintag{{\nwtagstyle{}\subpageref{NWgqRer-4MrGQl-2}}}\moddef{Output routines~{\nwtagstyle{}\subpageref{NWgqRer-4MrGQl-1}}}\plusendmoddef\Rm{}\nwstartdeflinemarkup\nwprevnextdefs{NWgqRer-4MrGQl-1}{\relax}\nwenddeflinemarkup
{\bf{}void} {\it{}test\_out}({\it{}string} {\it{}S}, {\bf{}const} {\bf{}ex} & {\it{}T}) {\nwlbrace}\nwindexdefn{\nwixident{test{\_}out}}{test:unout}{NWgqRer-4MrGQl-2}
        {\it{}cout} \begin{math}\ll\end{math} {\it{}S} \begin{math}\ll\end{math} ({\it{}latexout} ? {\tt{}"{\char92}{\char92}textbf{\char123}"} : {\tt{}"*"} )
                 \begin{math}\ll\end{math} ({\it{}is\_a}\begin{math}<\end{math}{\bf{}dual\_number}\begin{math}>\end{math}({\it{}T}) ? {\it{}ex\_to}\begin{math}<\end{math}{\bf{}dual\_number}\begin{math}>\end{math}({\it{}T}).{\it{}normal}().{\it{}is\_zero}() :
                         {\it{}T}.{\it{}evalm}().{\it{}normal}().{\it{}is\_zero\_matrix}()) \begin{math}\ll\end{math} ({\it{}latexout} ? {\tt{}"{\char125}{\char92}{\char92}{\char92}{\char92}"} : {\tt{}"*"} )
                 \begin{math}\ll\end{math} {\it{}endl};
{\nwrbrace}

\nwidentdefs{\\{{\nwixident{test{\_}out}}{test:unout}}}\nwendcode{}

\nwixlogsorted{c}{{*}{NWgqRer-1p0Y9w-1}{\nwixd{NWgqRer-1p0Y9w-1}\nwixd{NWgqRer-1p0Y9w-2}}}%
\nwixlogsorted{c}{{Check identities}{NWgqRer-22TNkn-1}{\nwixd{NWgqRer-22TNkn-1}\nwixd{NWgqRer-22TNkn-2}\nwixd{NWgqRer-22TNkn-3}\nwixd{NWgqRer-22TNkn-4}\nwixd{NWgqRer-22TNkn-5}\nwixd{NWgqRer-22TNkn-6}\nwixd{NWgqRer-22TNkn-7}\nwixd{NWgqRer-22TNkn-8}\nwixd{NWgqRer-22TNkn-9}\nwixd{NWgqRer-22TNkn-A}\nwixd{NWgqRer-22TNkn-B}\nwixd{NWgqRer-22TNkn-C}\nwixd{NWgqRer-22TNkn-D}\nwixd{NWgqRer-22TNkn-E}\nwixd{NWgqRer-22TNkn-F}\nwixu{NWgqRer-nXe8t-5}}}%
\nwixlogsorted{c}{{Definition of variables}{NWgqRer-2iImO3-1}{\nwixu{NWgqRer-1p0Y9w-2}\nwixd{NWgqRer-2iImO3-1}\nwixd{NWgqRer-2iImO3-2}\nwixd{NWgqRer-2iImO3-3}\nwixd{NWgqRer-2iImO3-4}\nwixd{NWgqRer-2iImO3-5}\nwixd{NWgqRer-2iImO3-6}\nwixd{NWgqRer-2iImO3-7}}}%
\nwixlogsorted{c}{{Drawing pictures}{NWgqRer-Ghke2-1}{\nwixu{NWgqRer-3sxui-2}\nwixd{NWgqRer-Ghke2-1}\nwixd{NWgqRer-Ghke2-2}\nwixd{NWgqRer-Ghke2-3}\nwixd{NWgqRer-Ghke2-4}\nwixd{NWgqRer-Ghke2-5}\nwixd{NWgqRer-Ghke2-6}\nwixd{NWgqRer-Ghke2-7}\nwixd{NWgqRer-Ghke2-8}}}%
\nwixlogsorted{c}{{Induced representations}{NWgqRer-1h89LW-1}{\nwixd{NWgqRer-1h89LW-1}\nwixd{NWgqRer-1h89LW-2}\nwixd{NWgqRer-1h89LW-3}\nwixd{NWgqRer-1h89LW-4}\nwixd{NWgqRer-1h89LW-5}\nwixd{NWgqRer-1h89LW-6}\nwixd{NWgqRer-1h89LW-7}\nwixd{NWgqRer-1h89LW-8}\nwixu{NWgqRer-3sxui-4}}}%
\nwixlogsorted{c}{{Induced representations routines}{NWgqRer-11aTCz-1}{\nwixd{NWgqRer-11aTCz-1}\nwixd{NWgqRer-11aTCz-2}\nwixd{NWgqRer-11aTCz-3}\nwixd{NWgqRer-11aTCz-4}\nwixu{NWgqRer-1p0Y9w-2}}}%
\nwixlogsorted{c}{{Main procedure}{NWgqRer-3sxui-1}{\nwixu{NWgqRer-1p0Y9w-2}\nwixd{NWgqRer-3sxui-1}\nwixd{NWgqRer-3sxui-2}\nwixd{NWgqRer-3sxui-3}\nwixd{NWgqRer-3sxui-4}\nwixd{NWgqRer-3sxui-5}}}%
\nwixlogsorted{c}{{Output routines}{NWgqRer-4MrGQl-1}{\nwixd{NWgqRer-4MrGQl-1}\nwixd{NWgqRer-4MrGQl-2}}}%
\nwixlogsorted{c}{{Show expressions}{NWgqRer-1Pc9Jw-1}{\nwixd{NWgqRer-1Pc9Jw-1}\nwixd{NWgqRer-1Pc9Jw-2}\nwixd{NWgqRer-1Pc9Jw-3}\nwixd{NWgqRer-1Pc9Jw-4}\nwixd{NWgqRer-1Pc9Jw-5}\nwixd{NWgqRer-1Pc9Jw-6}\nwixd{NWgqRer-1Pc9Jw-7}\nwixd{NWgqRer-1Pc9Jw-8}\nwixd{NWgqRer-1Pc9Jw-9}\nwixd{NWgqRer-1Pc9Jw-A}\nwixd{NWgqRer-1Pc9Jw-B}\nwixd{NWgqRer-1Pc9Jw-C}\nwixd{NWgqRer-1Pc9Jw-D}\nwixu{NWgqRer-nXe8t-5}}}%
\nwixlogsorted{c}{{Test routine}{NWgqRer-nXe8t-1}{\nwixu{NWgqRer-1p0Y9w-2}\nwixd{NWgqRer-nXe8t-1}\nwixd{NWgqRer-nXe8t-2}\nwixd{NWgqRer-nXe8t-3}\nwixd{NWgqRer-nXe8t-4}\nwixd{NWgqRer-nXe8t-5}}}%
\nwixlogsorted{i}{{\nwixident{formula{\_}out}}{formula:unout}}%
\nwixlogsorted{i}{{\nwixident{main}}{main}}%
\nwixlogsorted{i}{{\nwixident{numeric}}{numeric}}%
\nwixlogsorted{i}{{\nwixident{parab{\_}rot{\_}sub}}{parab:unrot:unsub}}%
\nwixlogsorted{i}{{\nwixident{test{\_}out}}{test:unout}}%
\nwbegindocs{167}\nwdocspar

\nwenddocs{}

}{}

\end{document}